\documentclass[12pt,oneside]{amsart}
\usepackage[T1]{fontenc}
\usepackage{lmodern}
\input{glyphtounicode}
\usepackage{amsmath,amssymb,amsthm,mathrsfs}
\usepackage[margin=1.05in]{geometry}
\usepackage{enumitem}
\usepackage{needspace}
\usepackage[hidelinks]{hyperref}
\newtheorem{thm}{Theorem}[section]
\newtheorem{prop}[thm]{Proposition}
\newtheorem{lem}[thm]{Lemma}
\newtheorem{cor}[thm]{Corollary}
\newtheorem*{introthm}{Theorem A}
\newtheorem*{introclass}{Theorem B}
\newtheorem*{introgroup}{Theorem C}
\newtheorem*{introdye}{Theorem D}
\newtheorem*{introintrinsic}{Theorem E}
\theoremstyle{definition}
\newtheorem{defn}[thm]{Definition}
\newtheorem{example}[thm]{Example}
\theoremstyle{remark}
\newtheorem{rem}[thm]{Remark}
\newtheorem*{question*}{Question}
\numberwithin{equation}{section}
\newcommand{\norm}[1]{\left\|#1\right\|}
\newcommand{\abs}[1]{\left|#1\right|}
\newcommand{\eps}{\varepsilon}

\newcommand{\N}{\mathbb N}
\newcommand{\T}{\mathrm{T}}

\newcommand{\dist}{\operatorname{dist}}
\newcommand{\Ad}{\operatorname{Ad}}

\newcommand{\supp}{\operatorname{supp}}

\newcommand{\Prob}{\operatorname{Prob}}

\setlist[enumerate]{label=\textup{(\roman*)},leftmargin=*}
\setlist[itemize]{leftmargin=*}

\hypersetup{
  pdftitle={Uniform measure equivalence relations, hyperfiniteness, and Cartan pairs},
  pdfauthor={George A. Elliott and Zhuang Niu},
  pdfsubject={Uniform measure equivalence relations and completed Cartan pairs},
  pdfkeywords={equivalence relation, Cartan subalgebra, uniform tracial completion, small boundary property, amenability}
}

\title[Uniform measure relations]{Uniform measure equivalence relations, 
hyperfiniteness, and Cartan pairs}
\author{George A. Elliott}
\address{Department of Mathematics, University of Toronto,
40 St. George Street, Toronto, Ontario M5S 2E4, Canada}
\email{elliott@math.toronto.edu}
\author{Zhuang Niu}
\address{Department of Mathematics and Statistics, University of Wyoming,
Dept. 3036, 1000 E. University Avenue, Laramie, Wyoming 82071, USA}
\email{zniu@uwyo.edu}
\date{}
\begin{document}
\begin{abstract}
We introduce uniform measure equivalence relations and show that, under
the relative small boundary property, measurewise amenability of a regular
relation is equivalent to uniform AF hyperfiniteness.
For aperiodic amenable relations in this class equipped with their full
compact families of invariant probability measures, we classify the
uniform tracial completions of the associated Cartan pairs by their
trace simplexes. Every affine homeomorphism between these simplexes is
induced by a uniform orbit equivalence on Borel sets conull for every
invariant probability measure. Every topological isomorphism between
the corresponding uniform full groups is also induced by a uniform
orbit equivalence and is automatically isometric.
\end{abstract}
\maketitle
\raggedbottom
{\small\tableofcontents}

\section{Introduction}

Feldman and Moore showed that every countable Borel equivalence relation
on a standard Borel space is the orbit relation of a countable group of
Borel automorphisms~\cite[Theorem~1]{FMI}. In the presence of a
quasi-invariant measure, their correspondence describes a von Neumann
algebra with separable predual and a Cartan subalgebra by a countable
measured relation and the cohomology class of a circle-valued
$2$-cocycle~\cite[Theorem~1]{FMII}. Connes, Feldman, and Weiss proved
that such a relation is amenable if and only if, on a conull set, it is
an increasing union of subrelations with finite classes~\cite[Theorem~10]{CFW}.

The classification in~\cite[Theorem~1.1 and Corollary~1.3]{TC} treats
free minimal actions of countably infinite amenable groups on compact
metrizable spaces with the small boundary property. For two such actions,
every affine homeomorphism between their invariant-measure simplexes is
realized both by an isomorphism of the completed crossed-product Cartan
pairs and by a uniform measure orbit equivalence.
In particular, each such action admits a Cantor minimal $\mathbb Z$-model
up to uniform measure orbit equivalence~\cite[Corollary~6.23]{TC}.
The group elements give a regular atlas for the orbit relation of the
action. Its completed Cartan pair is canonically the completed
crossed-product pair
(Lemma~\ref{lem:uniform-relation-algebra-action-case}).

Starting from the relation itself, we ask whether amenability for each
measure in a designated compact convex family gives a single sequence
of compatible finite matrix charts approximating the relation uniformly
over the family. For regular relations with the relative small boundary
property, Theorem~A gives such charts even for a proper subfamily of
the invariant probability measures. When these relations are aperiodic
and carry their full compact families of invariant probabilities,
Theorem~B extends the simplex classification and the realization of
prescribed affine homeomorphisms to this setting, including nonminimal
relations.

Our data consist of a countable Borel equivalence relation $R$ on a
compact metrizable space $X$, a nonempty compact convex family $K$ of
$R$-invariant probability measures, and the bounded uniform completion
$D_K=\mathrm{C}(X)^K$. A regular countable atlas has domain, range, and
coincidence projections in $D_K$, and its partial bijections preserve
the corresponding diagonal corners. These conditions specify which
graph coefficients are uniformly observable. They hold for the global
atlas of a free continuous action, but are additional requirements for
general relations. We call a set universally null when it is null for
every measure in $K$. A regular atlas need not come from an \'etale
groupoid on the original topological unit space; an aperiodic example,
including a version with a nontrivial full invariant family, is given in
Example~\ref{ex:regular-atlas-no-original-etale-model}.
Preservation of the completed observable algebra does not by itself
ensure regular graph coefficients, even for a continuous
orbit-preserving homeomorphism with the full invariant family
(Example~\ref{ex:observable-preservation-not-admissibility}).

The relative small boundary property ($K$-SBP), introduced in
\cite[Definition~2.2]{EN}, asks for a basis of open sets with
universally null boundaries.
Its diagonal characterization is
bounded uniform $2$-norm density of finite projection diagonals.
By \cite[Theorem~8.4]{TC}, this is also equivalent to real rank zero of
$D_K$ in its $\mathrm{C}^*$-norm. Our
finite-dimensional criterion, Theorem~\ref{thm:finite-dimensional-relation-SBP},
shows that $K$-SBP is automatic when $\dim X<\infty$ and the classes
are infinite on a common invariant $K$-conull set. It uses the regular
atlas and invariance, with no amenability or minimality assumption.

The basic examples in Section~\ref{subsec:basic-uniform-relation-examples}
include orbit relations of free continuous actions of countable groups,
tail relations of Bratteli diagrams, and orbit relations of principal
\'etale groupoids with relative SBP. An irrational circle rotation and
the dyadic odometer have SBP for their unique invariant measures. Standard AF
diagonals and ample principal groupoids have SBP because their unit
spaces have clopen bases. The framework also includes $\mathrm{C}^*$-diagonals in
separable unital nuclear $\mathrm{C}^*$-algebras with a nonempty compact convex
family of traces whose restrictions to the diagonal have SBP. The
unique extension property gives a principal Weyl groupoid, and the
construction allows its Weyl twist
(Example~\ref{ex:SBP-Cartan-uniform-relation}).

Our main theorem combines small-boundary diagonal approximation with
amenability.

\begin{introthm}[Uniform Connes--Feldman--Weiss theorem]
Let $(X,R,K,D_K,\mathcal S)$ be a regular uniform measure relation
with relative $K$-SBP. The following conditions are equivalent:
\begin{enumerate}
\item $(X,R,\mu)$ is amenable for every $\mu\in K$;
\item $R$ has regular uniform Reiter amenability;
\item $R$ has local regular finite approximation;
\item $R$ is strongly uniformly AF hyperfinite in the compatible
finite matrix-chart sense.
\end{enumerate}
The AF model retains the designated compact measure family. No
aperiodicity, simplicity, or Bauer hypothesis is required.
\end{introthm}

The conditions are defined in
Definitions~\ref{defn:uniform-Reiter-amenability},
\ref{defn:local-regular-cuts}, and
\ref{defn:uniform-AF-hyperfinite-relation}; the theorem is proved as
Theorem~\ref{thm:uniform-CFW-SBP}. There are two uniformization steps.
First, convex separation converts measurewise amenability into regular
Reiter kernels with uniform error bounds. Second, relative SBP permits
a deterministic regular rounding of these kernels. A repair of the
resulting finite cuts supplies nested compatible matrix charts.
Thus no additional uniform amenability assumption is needed.
Example~\ref{ex:Reiter-without-AF-coordinates} shows that SBP cannot
be omitted from the strong AF conclusion.

For a principal \'etale groupoid, the finite cuts have a geometric
realization in the original topology. Theorem~\ref{thm:SBP-groupoid-castles}
constructs small subgroupoids, meaning relatively compact open
subgroupoids with finite multisection coordinates. Their cores capture
every arrow in a prescribed compact set inside an enlarged small
subgroupoid, their source fibres are uniformly F\o lner, and their unit
spaces cover uniformly almost all of $X$. The levels may refine any
prescribed open cover. The proof uses the regular-to-open lifting of
Lemma~\ref{lem:regular-open-castle-lifting}. For minimal groupoids carrying
their full invariant-probability family, this gives almost elementariness
in measure (Corollary~\ref{cor:SBP-groupoid-AE-measure}). The comparison
with action URP and intrinsic relation URP appears in
Remarks~\ref{rem:SBP-small-subgroupoids-URP} and
\ref{rem:SBP-small-subgroupoids-actions}.
For minimal topologically amenable principal groupoids with infinite
compact unit space and relative SBP for the full nonempty family of
invariant probability measures, Corollary~\ref{cor:SBP-COS-groupoid-Z} combines these small subgroupoids
with Cuntz comparison of open sets to prove $\mathcal Z$-absorption of
the untwisted reduced groupoid algebra. The comparison hypothesis retains
the fixed amplification allowed in action COS.

There is also a comparison consequence which requires neither
amenability nor minimality. For a principal \'etale groupoid
carrying its full nonempty invariant-probability family,
Theorem~\ref{thm:SBP-groupoid-property-C} shows that relative SBP
implies Property \textup{(C)} for the standard Cartan inclusion
$\mathrm{C}(X)\subseteq \mathrm{C}_r^*(\mathcal G)$, with actual unitary witnesses.
The proof uses uniform finite matching along bisections rather
than the amenability-dependent finite exhaustions.

The algebraic counterpart is an untwisted forward Feldman--Moore
construction. Its completed Cartan pair realizes admissible partial
bijections as partial normalizers modulo diagonal phases and has a
separable $\mathrm{C}^*$-Cartan core. The admissible full group is Polish, and
each admissible partial bijection extends to an element of this group.
Theorem~\ref{thm:uniform-relation-pair-isomorphism} identifies uniform
spatial relation isomorphisms with marked completed-pair isomorphisms
for the constructed models. The construction specializes to the
completed crossed product of a free action. Section~\ref{sec:intrinsic-reconstruction}
proves the inverse reconstruction and twisted uniqueness theorems for
uniformly separable completed pairs whose partial-isometry normalizers
generate the bounded completion. The additional fibrewise maximal
abelianness condition is intrinsic and necessary for a principal relation
model; it does not follow merely from global maximal abelianness and SBP.

Compatible finite charts reconstruct both an AF Cartan model and its
spatial tail relation
(Theorem~\ref{thm:uniform-AF-reconstruction}). If the classes are
infinite on a common invariant universally conull set, the AF model
has a single uniformly admissible generator. This gives a Cantor
dynamical Cartan model with designated measures
(Theorem~\ref{thm:local-cuts-dynamical-model}). For a full invariant
family, the AF reconstruction has a stronger consequence.

\begin{introclass}[Classification of completed Cartan pairs]
Let two regular uniform measure relations have relative SBP, amenability
for every invariant probability measure, and infinite classes on common
invariant universally conull sets. Suppose that their full invariant
probability families are nonempty and compact. Every affine homeomorphism
between these families is implemented by an isomorphism of their
uniform completed Cartan pairs and by a uniform measure relation
isomorphism. The full invariant-measure simplex is a complete invariant
of the marked completed pairs in this class. Each pair has a standard
simple infinite-dimensional AF model carrying the full trace simplex,
and a Cantor minimal
$\mathbb Z$-realization with that invariant-measure simplex.
\end{introclass}

The proof in Section~\ref{subsec:uniform-CFW} does not require minimality
or a Bauer simplex. The reconstructed family is a closed trace face:
a decomposition of one of its central measures remains on the common
invariant coding image and pulls back to invariant measures in the full
family. Aperiodicity makes its extreme traces type $\mathrm{II}_1$.
We prove a relative classification theorem for standard AF Cartan
completions over such faces, allowing nonsimple AF algebras and
nonunital finite-stage maps. This supplies the classification even when
the first reconstructed AF diagram is not simple or has additional
traces. An infinite minimal principal \'etale groupoid with SBP and a
nonempty full invariant family satisfies the aperiodicity hypothesis;
topological amenability supplies the measured amenability required here
(Corollary~\ref{cor:minimal-SBP-groupoid-classification}).

When two pairs covered by Theorem B have the same ambient completion
and the same designated traces, the identity on the trace simplex
lifts to an automorphism carrying one completed diagonal onto the other
and fixing every designated trace. A common ambient completion can also
contain nonconjugate uniform Cartan subalgebras, as shown in
Remark~\ref{rem:uniform-classification-scope}.

For $\mathrm{C}^*$-Cartan inclusions, the designated trace simplex classifies the
uniform tracial completions of all pairs satisfying the intrinsic
hypotheses of Theorem~\ref{thm:intrinsic-Cartan-classification}.
In particular, let $A$ be a simple, separable, unital, infinite-dimensional
nuclear $\mathrm{C}^*$-algebra with $\T(A)\ne\varnothing$. Consider all Cartan
subalgebras $D\subseteq A$ with SBP relative to $\T(A)|_D$ whose Weyl
groupoids have trivial isotropy almost everywhere for every restricted
trace. Their completions $D^u\subseteq A^u$ form a single conjugacy
class under automorphisms of $A^u$ fixing every trace in $\T(A)$.
Their Weyl relations are uniformly measure orbit equivalent with the
identity marking on $\T(A)$
(Example~\ref{ex:SBP-Cartan-uniform-relation} and
Corollary~\ref{cor:amenable-twisted-dynamical-model}).

This includes every $\mathrm{C}^*$-diagonal in $A$ with SBP, since its Weyl
groupoid is principal. In an irrational rotation $\mathrm{C}^*$-algebra the trace
is unique, so SBP is automatic. All $\mathrm{C}^*$-Cartan subalgebras whose
closures are Cartan in the tracial factor have conjugate completions.
These include all $\mathrm{C}^*$-diagonals, among them Kumjian's pairwise
nonconjugate examples whose spectra are disjoint unions of $n$
circles, one for each $n\geq 1$ (see~\cite[Remark~5.3]{LiRenault}).

For arbitrary designated subfamilies, the simple full-trace model
condition remains an explicit sufficient hypothesis, and the measure
family alone need not classify the pair.

The simple full-trace AF model also connects Theorem B with the
metric classification of the projection Boolean algebra
\cite[Theorems~8.13 and~8.14]{TC}.
The projection Boolean algebra of the completed diagonal, equipped
with the metric $d_K(p,q)=\|p-q\|_{2,K}^2$, recovers the trace simplex
up to affine homeomorphism. Thus distinct affine types of trace
simplex give completed diagonals that admit no $*$-isomorphism
isometric for the uniform $2$-norm. The underlying norm $\mathrm{C}^*$-algebras
can also differ, as in~\cite[Proposition~8.23]{TC}. In the case
of finitely many extreme traces, the underlying norm $\mathrm{C}^*$-algebra
is always the standard atomless $L^\infty([0,1])$ algebra
\cite[Theorem~8.18]{TC}.

The classification also permits AF models with arbitrarily large
matrix multiplicities. This gives a uniform version of the
Giordano--Pestov full-group theorem~\cite{GP}, including the
nonergodic measured case treated by Le Ma\^itre~\cite{LM}.

\begin{introgroup}[Uniform full-group amenability]
Let a regular uniform measure relation have relative SBP and infinite
classes on a common invariant universally conull set. Assume that its
full invariant-probability family is nonempty and compact. Then
amenability for every invariant probability is equivalent to
topological amenability and to extreme amenability of its admissible
uniform full group. Under these conditions the group is a L\'evy group,
with an increasing dense union of finite permutation subgroups whose
normalized counting measures concentrate in the uniform disagreement
metric.
\end{introgroup}

Theorem~\ref{thm:uniform-full-group-extreme-amenability} proves this
assertion. The concentration estimate is uniform over the entire
measure simplex. Its reverse implication uses the continuous dense
restriction maps to the single-measure full groups established in
Lemma~\ref{lem:uniform-full-group-measured-density}.

\begin{introdye}[Reconstruction from the uniform full group]
Let two regular uniform measure relations have relative SBP,
measurewise amenability, and infinite classes on common invariant
universally conull sets. Assume that each carries its full nonempty
compact invariant-probability family. Every topological isomorphism
between their admissible uniform full groups is implemented by a
uniform measure relation isomorphism, uniquely modulo universal
almost-everywhere equality, and is automatically isometric.
Consequently, uniform measure orbit equivalence, topological
full-group isomorphism, isometric full-group isomorphism, and affine
homeomorphism of the invariant-measure simplexes are equivalent.
\end{introdye}

Theorem~\ref{thm:uniform-Dye-reconstruction} proves this uniform Dye
reconstruction. Involutions recover the regular Boolean algebra,
using Fremlin's Boolean reconstruction theorem~\cite[\S384D]{Fremlin}.
The group topology recovers the bounded uniform diagonal topology,
and the full invariant family is characterized by its continuous
full-group-invariant states. Fullness matters:
Example~\ref{ex:full-group-topology-loses-designated-weights} gives
topologically isomorphic groups with different designated weights.
No automatic continuity assertion for abstract group isomorphisms
is needed.

The inverse construction makes the operator-algebraic content intrinsic.
Write $M_\tau=\pi_\tau(M)''$ and $D_\tau=\pi_\tau(D)''$.

\begin{introintrinsic}[Intrinsic reconstruction and uniform CFW]
Let $(D\subseteq M,K)$ be a marked uniformly separable completed Cartan
pair in the sense of Definition~\ref{defn:intrinsic-completed-Cartan}.
It admits a regular scalar-twisted uniform measure relation model if
and only if $D_\tau$ is maximal abelian in $M_\tau$ for every $\tau\in K$.
The model is unique up to marked uniform relation isomorphism and regular
cohomology of its twist, independently of the chosen separable core.

Assume this fibrewise condition and intrinsic SBP. Then every $M_\tau$
is injective if and only if the marked pair has a standard AF diagonal
model $(\mathrm{C}(Y)^L\subseteq B^L,L)$, for a unital AF algebra $B$ and a compact
designated family $L\subseteq \T(B)$. In this case one regular gauge removes
the reconstructed twist, and there are increasing finite-dimensional
subalgebras of $M$, with matrix units normalizing $D$, whose algebras and
diagonals are boundedly dense in the uniform $2$-norm.
\end{introintrinsic}

Theorems~\ref{thm:intrinsic-uniform-FM} and
\ref{thm:twisted-uniform-pair-isomorphism} give reconstruction and
uniqueness. Proposition~\ref{prop:intrinsic-fibre-Cartan-tests} characterizes
the fibrewise condition by finite averages of diagonal conjugations;
under SBP it is also equivalent to uniform convergence of finite diagonal
compressions to the expectation. Theorem~\ref{thm:intrinsic-uniform-CFW}
gives the AF conclusion without any extra uniform amenability, simplicity,
aperiodicity, or Bauer hypothesis. If the extreme designated fibres are
type $\mathrm{II}_1$ factors, the designated AF family is automatically a
diffuse closed trace face. Theorem~\ref{thm:intrinsic-Cartan-classification}
then gives simplex classification, trace-preserving automorphism
conjugacy in a common ambient algebra, and a simple full-trace AF model.
These intrinsic hypotheses do not require that an initial spatial
presentation carry its entire invariant-measure family.

The finite charts also remove regular scalar twists: one regular
gauge trivializes the cocycle simultaneously for all designated
measures (Theorem~\ref{thm:uniform-regular-untwisting}). Relative SBP
regularizes locally trivial circle twists on principal \'etale groupoids
(Proposition~\ref{prop:SBP-regularizes-twist}), so the same conclusion
applies to their completed pairs when the orbit relations are
measurewise amenable. Section~\ref{sec:etale-presentations} distinguishes
exact continuous countable-group presentations from uniform
$\mathbb Z$-models, and separates freeness, minimality, and full
invariant-measure control. Section~\ref{sec:controlled-RN} extends the
finite approximation and reconstruction results to quasi-invariant
families with a controlled common Radon--Nikodym cocycle. In that
setting the state completion retains the modular cocycle;
classification by the designated simplex alone is not asserted.

Finally, Section~\ref{sec:AH-without-SBP} studies approximation without
small boundaries. Regular bounded finite cuts give an ASH Cartan
core; compatible split matrix charts give an AH core. We define an
intrinsic relation uniform Rokhlin property (URP) by embedded
finite-cover kernel models, without connecting maps, and prove
\[
 \text{relation URP}
 \quad\Longleftrightarrow\quad
 \text{uniform local AH Cartan approximation}
\]
without SBP (Theorem~\ref{thm:relation-URP-local-AH}). This property is
invariant under uniform measure orbit equivalence and implies
measurewise amenability. Within this class, SBP is equivalent to strong
uniform AF hyperfiniteness
(Corollary~\ref{cor:local-AH-Reiter-and-SBP}). For free actions of countably infinite amenable groups, URP gives
split local homogeneous approximation with locally bounded dimension
ratio (Proposition~\ref{prop:URP-continuous-local-AH};
see also~\cite[Lemma~4.1]{NiuZ}). Amenability does not yet
supply these local homogeneous models without SBP in the general
regular relation setting. Passage from local blocks to a compatible
ASH Cartan core remains open in the stated generality. The intrinsic
reconstruction of Section~\ref{sec:intrinsic-reconstruction} concerns
tracial completed pairs with partial-isometry regularity and fibrewise
maximal abelianness; it does not supply a nonsingular inverse theorem or
remove either regularity requirement.

Throughout, completions are taken on operator-norm bounded sets.
Isomorphisms of completed pairs preserve their designated traces or,
in the quasi-invariant section, their designated states. Sections~2--4
establish the uniform measure framework; Sections~5--8 develop the
algebraic and AF constructions;
Sections~\ref{sec:local-regular-cuts}--\ref{sec:regular-twists}
prove the finite-cut, groupoid, uniform amenability, and untwisting results. Section~\ref{sec:controlled-RN} treats controlled measure classes.
Section~\ref{sec:intrinsic-reconstruction} proves intrinsic reconstruction
and its CFW consequences, followed by the discussion without SBP in
Section~\ref{sec:AH-without-SBP}. The final section records further
questions and related work.

\medskip
\noindent\textbf{Acknowledgements.}
The work was carried out using ChatGPT.
The first author is supported by an NSERC Discovery Grant, and the
second author is supported by a Simons Foundation Grant (MP-TSM-00002606).

\section{Bounded completions and Borel realization}

\begin{defn}\label{def:completion}
Let $C$ be a unital $\mathrm{C}^*$-algebra and let $K\subseteq\T(C)$ be nonempty,
compact, and convex. Write
\[
 \norm{c}_{2,K}=\sup_{\tau\in K}\tau(c^*c)^{1/2}.
\]
The bounded tracial completion $C^K$ is the quotient of the
$\mathrm{C}^*$-algebra of bounded $\norm{\cdot}_{2,K}$-Cauchy sequences by the
ideal of bounded $\norm{\cdot}_{2,K}$-null sequences. Every trace in $K$
extends to this quotient. Its unit ball is the $\norm{\cdot}_{2,K}$
completion of the unit ball of $C$. We write $C^u=C^{\T(C)}$.
For a subalgebra $F\subseteq C$, $F^K$ means its bounded closure inside
$C^K$, using the restrictions of the designated traces.
For a displayed pair $F\subset C$, the notation $F^u$ always uses
the restrictions of $\T(C)$ and means the bounded closure in $C^u$.
\end{defn}

A nonempty compact trace family may be replaced by its closed convex
hull without changing the uniform $2$-seminorm or the completion.
This convention also permits compact designated subfamilies below.

Multiplication and involution are continuous in the uniform $2$-norm
on bounded sets, since
\begin{equation}\label{eq:two-mult}
 \norm{xay}_{2,K}\leq\norm{x}\norm{y}\norm{a}_{2,K}.
\end{equation}
The unit-ball description in Definition~\ref{def:completion} follows
by cutting bounded representatives down to the unit ball with
continuous functional calculus. All isomorphisms of completed pairs
in this paper preserve the designated trace sets. They are therefore
uniform $2$-norm isometries. An abstract $\mathrm{C}^*$-algebra isomorphism
without the trace marking is not the notion being classified.

We will also use bounded density in the following form. If a
$\mathrm{C}^*$-subalgebra $C_0\subset C$ is $2$-norm dense in $C$, its unit
ball is $2$-norm dense in the unit ball of $C$. Indeed, replace an
approximant $x\in C_0$ by its radial truncation
$r(x)=x\min\{1,|x|^{-1}\}$, defined continuously at zero. For a
contraction $a$, one has
$\norm{r(x)-a}_{2,\tau}\leq\norm{x-a}_{2,\tau}$ for each trace.
To see this in the finite von Neumann algebra of the trace, write
$x=v|x|$ and $t=(|x|-1)_+$. Then
$\operatorname{Re}\tau((x-r(x))^*(r(x)-a))
=\tau(t)-\operatorname{Re}\tau(tv^*a)\geq0$.
Expanding $\norm{x-a}_{2,\tau}^2$ proves the inequality. Taking the
supremum proves the assertion, also for completed ambient algebras.

\subsection{Spatial realization of completed diagonals}

Let \(X\) be compact metrizable and let \(\Delta\) be a nonempty compact set of Borel probability measures on \(X\).
Put
\[
 \norm{f}_{2,\Delta}=\sup_{\mu\in\Delta}\left(\int_X\abs{f}^2\,d\mu\right)^{1/2},
 \qquad \mathcal D_\Delta=\overline{\mathrm{C}(X)}^{\Delta}.
\]
The completion is taken on norm-bounded sets, after quotienting by the seminorm kernel if necessary.
A Borel set is \emph{\(\Delta\)-null} if it has measure zero for every \(\mu\in\Delta\).
In this section, ``universally conull'' always refers to the designated measure family.
Write \(\mathcal B_\Delta(X)\) for the algebra of bounded Borel functions modulo functions vanishing outside a \(\Delta\)-null set.
Its norm is \(\sup_{\mu\in\Delta}\operatorname*{ess\,sup}_{\mu}\abs{b}\).

\begin{lem}\label{lem:orbit-borel-realization}
There is a canonical faithful unital homomorphism
\[
 \mathcal D_\Delta\longrightarrow\mathcal B_\Delta(X).
\]
Its image consists of the bounded Borel classes uniformly approximable in \(\norm{\cdot}_{2,\Delta}\) by a norm-bounded sequence of continuous functions.
For every \(b\in\mathcal D_\Delta\), the functions
\(\mu\mapsto\int b\,d\mu\) and \(\mu\mapsto\int\abs{b}^2\,d\mu\) are continuous on \(\Delta\).
Moreover, \(\mathcal D_\Delta\) embeds faithfully into \(L^\infty(X,\mu_0)\), where
\(\mu_0=\sum_{n\geq1}2^{-n}\mu_n\) and \((\mu_n)\) is any weak*-dense sequence in \(\Delta\).
\end{lem}

\begin{proof}
Let \((f_n)\) be a uniform \(2\)-Cauchy sequence with \(\norm{f_n}\leq M\).
Choose a subsequence satisfying
\(\sum_j\norm{f_{n_{j+1}}-f_{n_j}}_{2,\Delta}<\infty\).
For every \(\mu\in\Delta\),
\[
 \int_X\sum_j\abs{f_{n_{j+1}}-f_{n_j}}\,d\mu
 \leq\sum_j\norm{f_{n_{j+1}}-f_{n_j}}_{2,\Delta}<\infty.
\]
Thus the subsequence converges pointwise outside one Borel \(\Delta\)-null set.
Its limit \(b\), extended by zero on that set, is Borel and bounded by \(M\).
Minkowski's inequality and Fatou's lemma give
\[
 \norm{b-f_{n_j}}_{2,\Delta}
 \leq\sum_{k\geq j}\norm{f_{n_{k+1}}-f_{n_k}}_{2,\Delta}\longrightarrow0.
\]
Equivalent sequences give limits agreeing almost everywhere for every measure in \(\Delta\).
Addition and involution pass to these limits; multiplication does also, since on bounded sets
\[
 \norm{f_ng_n-fg}_{2,\Delta}
 \leq\norm{f_n}\norm{g_n-g}_{2,\Delta}+\norm{g}\norm{f_n-f}_{2,\Delta}.
\]
This proves the homomorphism assertion, faithfulness, and the image description.
The integrals of \(f_n\) and \(\abs{f_n}^2\) converge uniformly on \(\Delta\) to the stated integrals of \(b\), proving continuity.
For the last assertion, passing to \(L^\infty(X,\mu_0)\) is well-defined because every \(\Delta\)-null set is \(\mu_0\)-null.
If \(\int\abs{b}^2\,d\mu_0=0\), then this integral vanishes for every \(\mu_n\), and continuity makes it vanish on all of \(\Delta\).
Hence \(b=0\) in \(\mathcal D_\Delta\).
\end{proof}

\begin{rem}
The measure \(\mu_0\) need not dominate every measure in \(\Delta\), and a \(\mu_0\)-null Borel set need not be \(\Delta\)-null.
Faithfulness in the lemma concerns the completed diagonal.
We use \(\mathcal B_\Delta(X)\) for statements involving all measures simultaneously.
\end{rem}

\begin{thm}\label{thm:orbit-spatial-diagonal}
Let \(X,Y\) be compact metrizable, and let \(\Delta_X,\Delta_Y\) be nonempty compact sets of probability measures.
Suppose \(H:\Delta_X\to\Delta_Y\) is a homeomorphism and
\(\Phi:\mathcal D_{\Delta_X}\to\mathcal D_{\Delta_Y}\) is a unital isomorphism such that
\[
 \int_Y\Phi(b)\,d(H\mu)=\int_Xb\,d\mu
 \quad(\mu\in\Delta_X,\ b\in\mathcal D_{\Delta_X}).
\]
Then there are Borel universally conull sets \(X^0\subseteq X\), \(Y^0\subseteq Y\), and a Borel isomorphism \(\theta:X^0\to Y^0\), such that
\(\theta_*\mu=H\mu\) for every \(\mu\in\Delta_X\) and
\[
 \Phi([b])=[b\circ\theta^{-1}]\quad\text{in }\mathcal B_{\Delta_Y}(Y).
\]
Functions on the conull domains may be extended by zero to their complements.
\end{thm}

\begin{proof}
The trace correspondence makes \(\Phi\) an isometry for the uniform \(2\)-norms.
Choose a countable norm-dense unital *-subalgebra \(\mathcal A\subseteq\mathrm{C}(X)\), over \(\mathbb Q+i\mathbb Q\), containing a separating family.
Choose a Borel representative \(b_a\) of \(\Phi(a)\), bounded by \(\norm{a}\), for each \(a\in\mathcal A\).
All rational linear, multiplicative, involutive, and unital relations between these representatives hold outside a universally null set.
There are countably many relations, so they hold simultaneously on one Borel universally conull set \(S_Y\).
At each \(y\in S_Y\), evaluation \(a\mapsto b_a(y)\) extends to a character of \(\mathrm{C}(X)\).
It is therefore evaluation at a unique point \(\sigma(y)\in X\).
A countable continuous separating family shows that \(\sigma\) is Borel.
Extend it arbitrarily to the omitted set.
The trace identity and norm density imply
\(\sigma_*(H\mu)=\mu\), by equality of integrals of continuous functions.

If bounded continuous \(a_n\) approximate a Borel representative \(b\) of an arbitrary completed element, then
\[
 \sup_{\nu\in\Delta_Y}\int_Y\abs{a_n\circ\sigma-b\circ\sigma}^2\,d\nu
 =\sup_{\mu\in\Delta_X}\int_X\abs{a_n-b}^2\,d\mu\longrightarrow0.
\]
Consequently \(\Phi([b])=[b\circ\sigma]\) for all completed \(b\).
Apply the same argument to \(\Phi^{-1}\), obtaining a Borel map \(\rho:X\to Y\), a good set \(S_X\), and \(\rho_*\mu=H\mu\).
Applying the two compositions to countably many separating functions gives
\(\sigma\rho(x)=x\) universally almost everywhere and \(\rho\sigma(y)=y\) universally almost everywhere.
Preimages of the respective universally null sets are universally null, by the pushforward identities.
Set
\[
 \begin{split}
 X^0&=S_X\cap\rho^{-1}(S_Y)\cap\{x:\sigma\rho(x)=x\},\\
 Y^0&=S_Y\cap\sigma^{-1}(S_X)\cap\{y:\rho\sigma(y)=y\}.
 \end{split}
\]
These are Borel universally conull sets, and \(\theta=\rho|_{X^0}\) is a Borel bijection onto \(Y^0\), with inverse \(\sigma|_{Y^0}\).
\end{proof}

\begin{rem}
The pullback identity is an equality of Borel classes.
For independently chosen representatives its pointwise exceptional set may depend on the completed element.
Exact pointwise equality on the fixed domain is obtained by defining image representatives by pullback, without intersecting uncountably many exceptional sets.
\end{rem}

\section{Uniform measure equivalence relations}
\label{sec:uniform-measured-relations}

A countable equivalence relation replaces the action of a particular group.
Besides the relation and its invariant measures, we retain the uniform
regularity of the diagonal and the graph coefficients.  The definition
specializes to the usual probability-measure-preserving relation for one
measure, and to uniform measure orbit equivalence for free continuous actions.

Let $X$ be a compact metrizable space, let $R\subseteq X\times X$ be a
countable Borel equivalence relation, and let $K$ be a nonempty compact convex
set of $R$-invariant Borel probability measures, with its weak star topology.
Here invariance means that every Borel partial bijection with graph contained
in $R$ preserves each measure in $K$ on its domain.  A Borel set is
\emph{$K$-null} if it is null for every measure in $K$.  We identify bounded
Borel functions that agree outside a $K$-null set and put
\[
 D_K=\mathrm{C}(X)^K,
 \qquad
 \|f\|_{2,K}=\sup_{\mu\in K}
       \left(\int_X |f|^2\,d\mu\right)^{1/2}.
\]
As before, $D_K$ denotes the bounded uniform tracial completion of $\mathrm{C}(X)$:
its elements are the uniform $2$-norm limits of uniformly bounded sequences
from $\mathrm{C}(X)$.  They admit bounded Borel representatives, modulo $K$-null sets.
In particular, $\mu\mapsto\int f\,d\mu$ is continuous on $K$ for $f\in D_K$.
If $p\in D_K$ is a projection, it is represented by a Borel characteristic
function and $\mu\mapsto\mu(p)$ is continuous.  Since $R$ is countable and
the measures are $R$-invariant, the $R$-saturation of a $K$-null Borel set is
again a $K$-null Borel set.  Thus countably many representative choices can
always be made on one $R$-invariant $K$-conull Borel set.

\subsection{Regular countable atlases}

\begin{defn}[Regular countable atlas]
\label{defn:regular-relation-atlas}
A \emph{regular countable atlas} for $(X,R,K)$ is a countable collection
$\mathcal S$ of Borel partial bijections $s:A_s\to B_s$ satisfying the
following conditions, after restriction to one $R$-invariant $K$-conull set:
\begin{enumerate}
\item $\mathcal S$ contains the identity and the empty map, is closed under
inversion and composition, and its graphs cover $R$.
\item The domain and range projections
$p_s=1_{A_s}$ and $q_s=1_{B_s}$ belong to $D_K$.
\item Pullback defines an isomorphism of corners
\[
 \alpha_s:q_sD_K\longrightarrow p_sD_K,
 \qquad
 \alpha_s(f)(x)=1_{A_s}(x)f(sx).
\]
\item For $s,t\in\mathcal S$, the coincidence projection
\[
 c(s,t)=1_{\{x\in A_s\cap A_t:sx=tx\}}
\]
belongs to $D_K$.
\end{enumerate}
The data $(X,R,K,D_K,\mathcal S)$ will be called a
\emph{uniform measure countable equivalence relation}.
Atlases giving the same admissible partial maps below are regarded as
equivalent.
\end{defn}

Condition (4) records the regularity of graph intersections.  It is needed
in addition to invariance of the diagonal under pullback.  In a free group
action, the global maps $x\mapsto gx$, together with the empty map, form such an atlas: distinct maps have
empty coincidence sets.  More generally, for an action with stabilizers,
the same atlas is regular provided all coincidence projections belong to
$D_K$.

\subsection{Basic examples}
\label{subsec:basic-uniform-relation-examples}

The following constructions give the main sources of regular uniform
measure relations used in this paper. The designated measure family is
part of each example. When relative SBP is invoked below, it means that
the unit space has a basis of open sets whose boundaries are null for
every designated measure; see Definition~\ref{defn:relative-K-SBP}.

\begin{example}[Free continuous actions]
\label{ex:free-action-uniform-relation}
Let a countable group $\Gamma$ act freely and continuously on a compact
metrizable space $X$, and let $K$ be a nonempty compact convex family of
invariant probability measures. Its orbit relation
\[
 R_\Gamma=\{(gx,x):g\in\Gamma,\ x\in X\}
\]
has the regular atlas consisting of the global group maps and the empty
map. Domain and range projections are $1$ or $0$, pullbacks preserve
$\mathrm{C}(X)$ and its bounded uniform completion, and freeness makes distinct
group maps have empty coincidence sets. Neither amenability, minimality,
nor SBP is needed for this regularity assertion. The completed relation
pair is the completed crossed-product pair, by
Lemma~\ref{lem:uniform-relation-algebra-action-case}.

For concrete aperiodic amenable examples, take an irrational rotation
of the circle or the adding-one homeomorphism of the dyadic integers,
with its unique invariant probability. Both have relative SBP and
therefore satisfy the uniform AF conclusion of
Theorem~\ref{thm:uniform-CFW-SBP}. Their singleton invariant-measure
simplexes are affinely homeomorphic, so
Theorem~\ref{thm:uniform-Cartan-classification} makes these two relations
uniformly measure orbit equivalent. This is a statement about the
completed observable algebras, not a topological conjugacy.
\end{example}

\begin{example}[AF tail relations and families with many measures]
\label{ex:AF-tail-uniform-relation}
Let $B$ be a separable unital AF algebra with its standard diagonal
$\mathrm{C}(Y)$, where $Y$ is the path space of a coordinate Bratteli diagram.
Let $R_B$ be tail equivalence, and let $K$ be any nonempty compact
convex family of central probability measures on $Y$. Finite-prefix
substitutions form a regular countable atlas: their domains, ranges,
and coincidence sets are clopen, and their pullbacks preserve the
completed corners. The level relations give compatible finite matrix
charts, while the finite cylinder diagonals are norm dense in $\mathrm{C}(Y)$.
Thus this is a strongly uniformly AF hyperfinite relation, with the
designated completed standard AF pair; see
Corollary~\ref{cor:uniform-AF-spatial-characterization}.
The family $K$ need not be the full central-measure family, and the
diagram need not be simple.

An explicit example with a nontrivial full invariant family is obtained
by taking a Cantor space $Z$, putting $\Omega=\{0,1\}^{\N}$ and
$Y=Z\times\Omega$, and declaring
\[
 (z,\omega)R(z',\omega')
 \quad\Longleftrightarrow\quad
 z=z'\ \text{and}\ \omega_n=\omega'_n\ \text{for all sufficiently large }n.
\]
Let $m$ be the fair Bernoulli product measure on $\Omega$. The full
invariant family is
\[
 K=\{\nu\otimes m:\nu\in\Prob(Z)\}.
\]
Indeed, for an invariant measure $\mu$, every Borel $E\subseteq Z$,
and every length-$n$ cylinder $[w]\subseteq\Omega$, prefix invariance gives
$\mu(E\times[w])=2^{-n}\mu(E\times\Omega)$. These equalities determine
the displayed measures. For $f\in \mathrm{C}(Y)$,
\[
 \|f\|_{2,K}^2
 =\sup_{z\in Z}\int_\Omega |f(z,\omega)|^2\,dm(\omega).
\]
This explicitly exhibits the uniformity over the parameter space.
The corresponding standard AF pair is
\[
 \mathrm{C}(Z\times\Omega)\subseteq \mathrm{C}(Z)\otimes \mathrm{M}_{2^\infty},
\]
where the diagonal of the dyadic UHF algebra $\mathrm{M}_{2^\infty}$ is
$\mathrm{C}(\Omega)$. Every class is infinite, but the relation is not minimal:
the parameter $z$ is unchanged. Thus the aperiodic classification
theorem applies without a minimality assumption.
\end{example}

\begin{example}[Principal \'etale groupoids]
\label{ex:principal-groupoid-uniform-relation}
Let $\mathcal G$ be a second-countable locally compact Hausdorff
principal \'etale groupoid with compact metrizable unit space $X$,
and let $K$ be a nonempty compact convex family of invariant
probability measures. If $(X,K)$ has relative SBP, its orbit relation
\[
 R_{\mathcal G}=\{(r\gamma,s\gamma):\gamma\in\mathcal G\}
\]
has a regular countable atlas obtained from small-boundary restrictions
of open bisections. Principality identifies graph coincidences with
bisection intersections; relative SBP makes their projections regular.
This is proved in Proposition~\ref{prop:principal-groupoid-SBP-atlas}.
For an ample groupoid the clopen basis already supplies relative SBP.
Topological amenability is not needed to obtain the atlas, but it does
give the measured amenability needed for the uniform AF conclusion.
Locally trivial circle twists are retained as regular twists of the
relation and, in the amenable case, removed in the completion by
Proposition~\ref{prop:SBP-regularizes-twist} and
Theorem~\ref{thm:uniform-regular-untwisting}.
\end{example}

\begin{example}[A regular atlas without an \'etale presentation on the unit space]
\label{ex:regular-atlas-no-original-etale-model}
Let $C=\{0,1\}^{\N}$ have its fair Bernoulli probability $\nu$, let
$m$ be Haar probability on the circle $\mathbb T$, and put
\[
 X=C\sqcup\mathbb T,\qquad \mu=\tfrac12\nu+\tfrac12m.
\]
Choose a measure-preserving Borel bijection
$\phi:(C,\nu)\to(\mathbb T,m)$. Binary expansion gives such a
bijection away from countable exceptional sets; modifying it on a
countable set gives a bijection everywhere. Fix an irrational rotation
$a$ of $\mathbb T$. Define Borel automorphisms $S,T$ of $X$ by
\[
 \begin{array}{lll}
 S(c)=\phi(c),&T(c)=\phi^{-1}a\phi(c)&(c\in C),\\
 S(z)=\phi^{-1}(z),&T(z)=a(z)&(z\in\mathbb T).
 \end{array}
\]
They preserve $\mu$, commute, and generate a free Borel action of
$\Gamma=\mathbb Z\times\mathbb Z/2$. Let $R$ be its orbit relation.
Every class is infinite and meets both components. Moreover,
$M(R)=\{\mu\}$: invariance under $S$ forces equal component masses,
and the normalized restriction to $\mathbb T$ must be the unique
invariant probability of the irrational rotation. Its restriction to
$C$ is then determined by $\phi$.

For $K=\{\mu\}$, bounded $L^2(\mu)$ approximation by continuous
functions gives $D_K=L^\infty(X,\mu)$. Thus the global group maps,
together with the empty map, form a regular countable atlas: their
pullbacks preserve $D_K$, their support projections are $0$ or $1$,
and distinct group maps have empty coincidence sets. In particular,
this is an aperiodic regular uniform measure relation with an atomless
full invariant family and $\dim X=1$.

Nevertheless, no \'etale groupoid with unit space $X$ carrying its
given topology can have orbit relation $R$. An arrow from $c\in C$
to $z\in\mathbb T$ would, by restriction of an open bisection, give
a homeomorphism between open neighborhoods contained in the two
components. A nonempty open subset of $C$ is zero dimensional,
whereas every nonempty open subset of $\mathbb T$ contains an arc.
Such a homeomorphism is impossible. The obstruction persists even
if agreement of the orbit relations is required only on a common
$\mu$-conull set while the groupoid unit space remains $X$:
$S$ preserves $\mu$, so a conull set contains both $c$ and $S(c)$
for almost every $c\in C$.

There is no obstruction to another topological presentation. The
Borel bijection
\[
 c\longmapsto(\phi(c),0),\qquad
 z\longmapsto(z,1)
\]
identifies this measured action with the continuous action on
$\mathbb T\times\{0,1\}$ generated by rotation and the flip.
Thus $R$ is a principal Borel groupoid, as is every Borel equivalence
relation, and has another \'etale model. It has no \'etale presentation
on the prescribed topological unit space.

For a nontrivial compact full invariant family, take a Cantor space
$P$, let $\Gamma$ act trivially on $P$ and as above on $X$, and set
\[
 \widetilde X=P\times X,\qquad
 \widetilde K=\{\lambda\otimes\mu:\lambda\in\Prob(P)\}.
\]
This is the full invariant family. Indeed, disintegration over $P$
has invariant conditional probabilities for almost every parameter,
and the invariant probability on each fiber is uniquely $\mu$.
For a bounded Borel $f$,
\[
 \|f\|_{2,\widetilde K}^2
 =\sup_{p\in P}\int_X|f(p,x)|^2\,d\mu(x).
\]
The same global atlas is regular. To check pullbacks of continuous
functions, approximate in the parameter by a finite clopen step
function $\sum_i1_{P_i}(p)f_i(x)$. After pullback by a fixed group
map, approximate each of the finitely many bounded Borel fiber
functions in $L^2(\mu)$ by bounded continuous functions. This gives
uniform $2$-norm approximation over all parameters. Inverses give
the corner isomorphisms, and freeness gives the coincidence
projections. The topological obstruction remains: $P\times C$ is
zero dimensional, while every nonempty open subset of
$P\times\mathbb T$ contains an arc. All designated measures are
atomless, and every class is infinite. Theorem~\ref{thm:finite-dimensional-relation-SBP}
also applies to both presentations on $X$ and $\widetilde X$.
\end{example}

\begin{example}[$\mathrm{C}^*$-Cartan pairs with SBP and universally null isotropy]
\label{ex:SBP-Cartan-uniform-relation}
Let $D=\mathrm{C}(X)\subseteq A$ be a separable unital $\mathrm{C}^*$-Cartan pair, let
$L\subseteq \T(A)$ be nonempty, compact, and convex, and set
$K=\{\tau|_D:\tau\in L\}$. Assume that $(X,K)$ has relative SBP.
Renault's reconstruction~\cite{RenaultC} identifies the pair with
$\mathrm{C}(X)\subseteq \mathrm{C}_r^*(\mathcal G,\Sigma)$ for its Hausdorff
\'etale Weyl groupoid and Weyl twist. In addition, assume
\[
 \mu\{x\in X:\mathcal G_x^x\ne\{x\}\}=0
 \qquad(\mu\in K).
\]
The isotropy-free units form an invariant Borel set $X_0$ conull for
every designated measure. On the compact unit space $X$, the Weyl
orbit relation admits a regular atlas after restriction to $X_0$.
Its regular twisted completion is the designated completion
$(D_K,A^L)$ of the original pair.

The restriction of a trace is invariant. Indeed, on a trivializing
bisection the trace identity $\tau(v^*v)=\tau(vv^*)$ transports the
squared absolute value of each compactly supported section $v$ from
source to range. Conversely every invariant probability gives a trace
by composition with $\mathbb E$: on compactly supported sections the
diagonal coefficients of $a^*a$ and $aa^*$ integrate to the same value
by invariance and inversion.

To check the regular atlas assertion, choose countably many bisections carrying
continuous local twist sections and restrict their sources to
small-boundary open sets with compact closure in the original sources.
Invariance gives small boundaries for the ranges, and continuous
cutoffs give the domain projections and corner pullbacks in $D_K$.
For two such charts, equality of endpoints on $X_0$ is equality of
arrows. Hence their coincidence set agrees, modulo a $K$-null set,
with the source of their bisection intersection. Equality of the two
arrow-valued sections is closed in their common domain by Hausdorffness
and open by the bisection property. The same argument applies to finite
words in the restricted bisections and their inverses, whose source and
range boundaries remain null as in
Proposition~\ref{prop:principal-groupoid-SBP-atlas}. The word charts
thus form a regular atlas on $X_0$, using the completed diagonal of the
original compact space $X$. The local sections give a regular scalar
twist by the argument of Proposition~\ref{prop:SBP-regularizes-twist}.

To identify the trace marking, fix a designated trace $\tau$ and write
$\mu=\tau|_D$. For a compactly supported twisted section $f$ on an
off-unit bisection $U$, the functional $g\mapsto\tau(fg)$ on $\mathrm{C}(X)$
is bounded for $L^2(\mu)$ by Cauchy--Schwarz. Its representing density
vanishes off $s(U)$ by localization. Traciality and the bisection
covariance imply that on $s(U)$ the density is supported where
$h(x)=h(\beta_Ux)$ for all $h$ in a countable point-separating family
in $\mathrm{C}(X)$. This fixed set is contained in the isotropy locus and is
$\mu$-null. Hence $\tau(f)=0$. Decomposing compactly supported sections
over the units and off-unit bisections gives
$\tau=\mu\circ\mathbb E$. Continuous cutoffs and finite sums of local
bisection sections now identify the two bounded completions, including
the expectation and the designated traces.

In particular, this construction applies to every $\mathrm{C}^*$-diagonal with
relative SBP: the unique extension property for pure states is
equivalent to principality of the Weyl groupoid~\cite{RenaultC}.
If $A$ is nuclear, the reduced algebra of the Fell line bundle of its
Weyl twist is nuclear. Takeishi's theorem~\cite[Theorem~5.4]{Takeishi}
then gives amenability of the Weyl groupoid and hence amenability of
the orbit relation for each designated measure. Theorem~\ref{thm:uniform-CFW-SBP} and
Theorem~\ref{thm:uniform-regular-untwisting} then give a standard AF
Cartan model for $(D_K,A^L)$. No aperiodicity or simplicity is needed
for this AF-model conclusion.
\end{example}

\begin{rem}[Why the isotropy qualification is necessary]
\label{rem:SBP-Cartan-isotropy-obstruction}
SBP alone does not justify the preceding completed-pair assertion.
Let $X=\mathbb Z\cup\{\infty\}$ be the one-point compactification,
let $T(n)=n+1$ and $T(\infty)=\infty$, and put
$D=\mathrm{C}(X)\subseteq A=\mathrm{C}(X)\rtimes_T\mathbb Z$.
The action is topologically free, so this is a $\mathrm{C}^*$-Cartan pair, and
the zero-dimensional space $X$ has SBP. However,
\[
 0\longrightarrow\mathcal K(\ell^2(\mathbb Z))
 \longrightarrow A\longrightarrow \mathrm{C}(\mathbb T)\longrightarrow0.
\]
Every trace vanishes on the compact-operator ideal. Using the full
trace family, the uniform seminorm on the quotient is the supremum
norm, and consequently
\[
 (D^u,A^u)\cong\bigl(\mathbb C,\mathrm{C}(\mathbb T)\bigr).
\]
The completed diagonal is not maximal abelian. All invariant
probabilities on $X$ are concentrated at $\infty$; the equality
relation there forgets the isotropy retained by the algebra.
\end{rem}

Regular uniform measure relations need not have SBP: the split finite-fiber relation
and the double-cover relation in
Examples~\ref{ex:split-AH-without-SBP} and
\ref{ex:cover-AH-without-transversal} have finite regular atlases and
exact homogeneous models but fail relative SBP. Their details are
kept in Section~\ref{sec:AH-without-SBP}, where they explain the need
for homogeneous rather than AF approximation.

\subsection{Admissible partial bijections}

\begin{defn}[Uniform measure full pseudogroup]
\label{defn:uniform-relation-pseudogroup}
For a regular atlas $\mathcal S$, let
$\mathopen{[\![}R\mathclose{]\!]}^u_{K,\mathcal S}$ consist of Borel partial bijections
$\varphi:A\to B$ with graph contained in $R$ such that
$1_A,1_B\in D_K$ and
\[
 c_s(\varphi)
 =1_{\{x\in A\cap A_s:\varphi(x)=sx\}}\in D_K
 \qquad(s\in\mathcal S).
\]
Partial maps are identified when their domains and their values agree
outside a $K$-null set.  This is the
\emph{uniform measure full pseudogroup} of the regular relation.
Its unit group, denoted $[R]^u_{K,\mathcal S}$, consists of the maps with
domain and range equal to $X$ modulo $K$-null sets.  We call it the
\emph{uniform measure full group}.
\end{defn}

When the atlas is fixed, we also write
\[
 \mathcal D=D_K,\qquad G=[R]^u_{K,\mathcal S},\qquad
 \mathscr P=\mathopen{[\![}R\mathclose{]\!]}^u_{K,\mathcal S},\qquad
 [R]^u_{K,\mathcal D}=[R]^u_{K,\mathcal S}.
\]
The latter notation always retains the specified equivalence class of
atlases.  It does not assert that $R$, $K$, and $D_K$ alone determine
that class.  For a free action we abbreviate the full group further to
$[R_\Gamma]^u_K$.

For statements involving total maps, representatives may be chosen as
bijections of an $R$-invariant $K$-conull Borel set and extended by the
identity on its complement.  This convention avoids separate domains for
different representatives.

\begin{prop}[Uniform tails and closure]
\label{prop:regular-relation-pseudogroup-closure}
Let $\mathcal S=(s_j)_{j\geq1}$ be a regular atlas.  Every
$\varphi\in\mathopen{[\![}R\mathclose{]\!]}^u_{K,\mathcal S}$ preserves each measure
in $K$ and induces a uniform $2$-norm isometric isomorphism
$1_BD_K\to1_AD_K$ by pullback.  The pseudogroup contains $\mathcal S$ and
is closed under inversion, composition, and restriction to a domain whose
characteristic function belongs to $D_K$.  It is also closed under
countable disjoint unions whose total domain projection belongs to $D_K$;
the total range projection then belongs to $D_K$ automatically.
\end{prop}

\begin{proof}
Put $a=1_A$ and disjointize the graph coefficients in the order of the atlas:
\[
 e_j=c_{s_j}(\varphi)
       \prod_{i<j}\bigl(1-c_{s_i}(\varphi)\bigr).
\]
These are orthogonal projections in $D_K$, their pointwise sum is $a$, and
$\varphi=s_j$ on the support of $e_j$.  The functions
\[
 \mu\longmapsto \mu\!\left(a-\sum_{j\leq n}e_j\right)
\]
are continuous on $K$ and decrease pointwise to zero.  Dini's theorem gives
\begin{equation}
\label{eq:regular-relation-uniform-tails}
 \sup_{\mu\in K}
    \mu\!\left(a-\sum_{j\leq n}e_j\right)\longrightarrow0.
\end{equation}
Measure preservation follows either from $R$-invariance or directly by
partitioning the domain into these pieces.  The range projections
$\alpha_{s_j^{-1}}(e_j)$ are orthogonal, and their omitted trace mass is
the same as that of the omitted $e_j$.  Thus their bounded uniform
$2$-norm sum belongs to $D_K$ and is $1_B$.

For $f\in 1_BD_K$, pullback is the bounded uniform $2$-norm sum
\[
 f\circ\varphi
   =\sum_j e_j\alpha_{s_j}(q_{s_j}f).
\]
Indeed, finite partial sums have norm at most $\|f\|$ and their squared
uniform $2$-norm error is at most $\|f\|^2$ times the tail in
\eqref{eq:regular-relation-uniform-tails}.  Invariance of every measure
gives the asserted isometry.

For inversion, fix $t\in\mathcal S$.  On the range of the $j$th piece,
$\varphi^{-1}$ agrees with $t$ precisely when
$t s_j(x)=x$.  Consequently
\[
 c_t(\varphi^{-1})
  =\sum_j\alpha_{s_j^{-1}}
           \bigl(e_j c(t s_j,\mathrm{id})\bigr).
\]
Every term belongs to $D_K$, and the sum has the same uniform tail bound
as above.  Thus the inverse is admissible and gives the inverse corner
isomorphism.

For composition, let $\psi$ have a disjoint chart decomposition
$f_j$ with $\psi=t_j$ on its $j$th piece, and retain $e_i$ for $\varphi$.
The source pieces of $\varphi\psi$ are
\[
 h_{ij}=f_j\alpha_{t_j}(q_{t_j}e_i)\in D_K.
\]
They are orthogonal and their sum is the domain projection of the
composition, which belongs to $D_K$ by the corner pullback already proved.
Keeping only terms with $i,j\leq n$ loses, for each measure, at most the omitted
source mass for $\psi$ plus the omitted source mass for $\varphi$.
This tends uniformly to zero.  On the $ij$th piece the composition is
$s_i t_j$, so for $t\in\mathcal S$ its coincidence projection is the
bounded uniform $2$-norm sum
\[
 c_t(\varphi\psi)=\sum_{i,j}h_{ij}c(s_i t_j,t)\in D_K.
\]
The range projection is explicitly
$\alpha_{\varphi^{-1}}(1_{A_\varphi}1_{B_\psi})$, so it too belongs
 to $D_K$.  This proves closure under composition.  Restriction
multiplies the source coefficients by the restricting projection; its
range projection is the pullback by the inverse.

Finally, consider partial maps with disjoint domains and disjoint ranges.
If their total domain projection is in $D_K$, Dini's theorem applied to
the finite sums of their domain projections gives uniform tails.  The
range tails have equal trace mass.  The union's graph coefficients are
bounded uniform $2$-norm sums of the coefficients of its pieces and hence
belong to $D_K$; the same holds for its range projection.  The atlas itself
is contained in the pseudogroup by its coincidence condition.
\end{proof}

\begin{cor}\label{cor:regular-relation-equalizers}
If $\varphi,\psi\in\mathscr P$, their coincidence projection belongs
to $D_K$.
\end{cor}
\begin{proof}
On their common domain the coincidence projection is
$c_{\mathrm{id}}(\psi^{-1}\varphi)$, which belongs to $D_K$ by the
closure proposition and the identity chart.
\end{proof}

Proposition~\ref{prop:regular-relation-pseudogroup-closure} shows that
regular graph coefficients imply preservation of the completed
observable algebra. The converse would give a simpler definition of
the full pseudogroup, using only the diagonal pullback. The following
example explains why the graph-coefficient condition in
Definition~\ref{defn:uniform-relation-pseudogroup} must be retained.

The two conditions record different information. Pullback tests the
values of observables after transport, whereas graph coefficients
record which atlas chart performs that transport. Even a continuous
orbit-preserving homeomorphism can have topologically small
displacements whose integer orbit labels escape every finite set.
This occurs with the full invariant-measure family, as follows.

\begin{example}[Observable preservation does not imply admissibility]
\label{ex:observable-preservation-not-admissibility}
Let
\[
 B=\{0\}\cup\{1/n:n\geq1\},\qquad X=B\times\mathbb T,
\]
let $\alpha\in\mathbb R\setminus\mathbb Q$, and let $R$ be the orbit
relation of the fiberwise rotation
$S(t,z)=(t,z+\alpha)$. Use the regular atlas consisting of the
maps $S^k$, $k\in\mathbb Z$, and the empty map. Put
\[
 K=\{\nu\otimes m:\nu\in\operatorname{Prob}(B)\},
 \qquad \mathcal D=D_K=\mathrm{C}(X)^K,
\]
where $m$ is normalized Lebesgue measure on $\mathbb T$.
This is the full invariant-probability family: $S$ fixes the first
coordinate and is uniquely ergodic on each fiber. The supremum in
the uniform $2$-norm is already attained as a supremum over the
individual fiber measures $\delta_t\otimes m$.

Choose distinct nonzero integers $q_n$ with
$q_n\alpha\longrightarrow0$ in $\mathbb T$. Define
\[
 T(0,z)=(0,z),\qquad
 T(1/n,z)=(1/n,z+q_n\alpha).
\]
The convergence of the rotations to the identity makes both $T$ and
$T^{-1}$ continuous, including on the limiting fiber. The map $T$
preserves every measure in $K$, has graph in $R$, and preserves
$\mathrm{C}(X)$ and $\mathcal D$ under pullback. Thus it satisfies the
observable condition in both directions.

Nevertheless, its graph coefficient for the identity chart is
\[
 c_{\mathrm{id}}(T)
   =1_{\{0\}\times\mathbb T}=:p\notin\mathcal D.
\]
Indeed, if $f\in \mathrm{C}(X)$, then
$f(1/n,\cdot)\longrightarrow f(0,\cdot)$ in $L^2(\mathbb T,m)$.
Writing $\delta=\|f-p\|_{2,K}$ gives
\[
 \|f(0,\cdot)-1\|_2\leq\delta,
 \qquad \|f(0,\cdot)\|_2
      =\lim_n\|f(1/n,\cdot)\|_2\leq\delta.
\]
The triangle inequality yields $1\leq2\delta$, so $p$ cannot be
approximated in uniform $2$-norm by continuous functions. For
$k\ne0$, by contrast, $c_{S^k}(T)$ is either zero or the
characteristic function of one isolated clopen fiber, because the
$q_n$ are distinct. The identity coefficient alone already prevents
$T$ from being admissible, or from being added to this regular atlas.

The same construction displays the failure of uniform chart tails.
For every finite $F\subseteq\mathbb Z$, choose $n$ with $q_n\notin F$.
On the fiber $\{1/n\}\times\mathbb T$, the map $T$ agrees with none
of the charts $S^k$, $k\in F$. Consequently,
\[
 \sup_{\mu\in K}
 \mu\{x:T(x)\notin\{S^k(x):k\in F\}\}=1.
\]
This contradicts the uniform-tail conclusion
\eqref{eq:regular-relation-uniform-tails} that admissibility would
give. Although the displacements tend to zero in the topology of
$X$, their chart labels are not uniformly captured by finite sets.
Regular graph coefficients are what permit the passage from
chartwise transport to uniform finite approximation.
\end{example}

\begin{prop}[Independence of equivalent atlases]
\label{prop:uniform-relation-atlas-independence}
Suppose $\mathcal S$ and $\mathcal T$ are regular countable atlases for the
same $(X,R,K,D_K)$.  If every chart of either atlas is admissible for the
other, then
\[
 \mathopen{[\![}R\mathclose{]\!]}^u_{K,\mathcal S}
       =\mathopen{[\![}R\mathclose{]\!]}^u_{K,\mathcal T}.
\]
In particular, the uniform measure full group does not depend on the
choice of an equivalent atlas.
\end{prop}

\begin{proof}
Write an $\mathcal S$-admissible map $\varphi$ as the disjoint union of the
restrictions $s_j|_{e_j}$ used above.  For $t\in\mathcal T$,
\[
 c_t(\varphi)=\sum_j e_j c_t(s_j).
\]
The summands belong to $D_K$ by mutual admissibility and have uniformly
vanishing omitted source mass.  Thus the sum belongs to $D_K$, so
$\varphi$ is $\mathcal T$-admissible.  Interchanging the atlases proves
the reverse inclusion.
\end{proof}

\subsection{The full-group metric and spatial isomorphisms}

\begin{prop}[The complete uniform full-group metric]
\label{prop:uniform-relation-full-group-complete}
The formula
\[
 d_K(T,S)=\sup_{\mu\in K}\mu\{x:T(x)\ne S(x)\}
\]
defines a complete bi-invariant metric on $[R]^u_{K,\mathcal S}$.
This group is separable and hence Polish.
\end{prop}

\begin{proof}
All its elements preserve every measure in $K$.  Therefore right
translation preserves the metric, while injectivity of a left translating
map preserves the disagreement set.  Inversion is also an isometry: the
images under $T$ and $S$ of their agreement set coincide and have the
same measure as that set.

To prove completeness, pass from a Cauchy sequence to a subsequence,
still denoted $T_n$, for which
\[
 \sum_n d_K(T_n,T_{n+1})<\infty.
\]
The limsup of the consecutive disagreement sets is $K$-null by the
Borel--Cantelli lemma; the same holds for the inverses, by inversion
invariance of the metric.  Remove their $R$-saturations.  On the remaining
$R$-invariant $K$-conull Borel set, $T_n(x)$ and $T_n^{-1}(x)$ eventually
stabilize at every point.  Their eventual values define inverse Borel
bijections $T$ and $T^{-1}$ with graph in $R$.  Moreover,
\[
 d_K(T_n,T)\leq\sum_{j\geq n}d_K(T_j,T_{j+1})\longrightarrow0.
\]
For any fixed chart $s$,
\[
 \|c_s(T_n)-c_s(T)\|_{2,K}^2\leq d_K(T_n,T).
\]
Closedness of bounded sets of $D_K$ under uniform $2$-norm limits shows
that all graph coefficients of $T$ belong to $D_K$.  Thus $T$ is
admissible.  The original Cauchy sequence has the same limit.

For separability, use the first-chart disjointization $e_j(T)$ for a
total map $T$.  The first chart agreeing with an arrow is determined by
the arrow itself.  Hence, pointwise,
\[
 \sum_j|e_j(T)-e_j(S)|^2=2\,1_{\{T\ne S\}}.
\]
Each coefficient sequence has uniformly vanishing tails.  It therefore
belongs to the completion of finite sequences in $D_K$ for the norm
\[
 \|(f_j)\|=
    \sup_{\mu\in K}\left(\int_X\sum_j|f_j|^2\,d\mu\right)^{1/2}.
\]
This metric space is separable: $D_K$ is separable in the uniform
$2$-norm, since $\mathrm{C}(X)$ is separable in norm, and finite sequences from a
countable dense subset are dense in this completion.  The displayed
identity embeds the full group with metric $(2d_K)^{1/2}$ into that
separable metric space.  The two metrics give the same topology, so the
full group is separable.
\end{proof}

\begin{defn}[Uniform measure relation isomorphism]
\label{defn:uniform-measured-relation-isomorphism}
Consider two uniform measure countable equivalence relations,
\[
 (X,R,K,D_K,\mathcal S)
 \quad\text{and}\quad (Y,Q,L,D_L,\mathcal T).
\]
A \emph{uniform measure relation
isomorphism} implementing an affine homeomorphism $H:K\to L$ is a Borel
bijection $\theta:X_0\to Y_0$ between invariant universally conull Borel
sets such that:
\begin{enumerate}
\item $\theta$ carries $R|_{X_0}$ onto $Q|_{Y_0}$ and
$\theta_*\mu=H(\mu)$ for every $\mu\in K$;
\item pullback by $\theta$ carries $D_L$ onto $D_K$;
\item every $\theta s\theta^{-1}$, $s\in\mathcal S$, is
$\mathcal T$-admissible, and every $\theta^{-1}t\theta$,
$t\in\mathcal T$, is $\mathcal S$-admissible.
\end{enumerate}
Here ``universally conull'' refers to all measures in the indicated family.
\end{defn}

\begin{defn}\label{defn:uniform-measure-oe}
Let countable groups \(\Gamma_i\) act freely and continuously on compact metrizable \(X_i\), and let \(K_i=M_{\Gamma_i}(X_i)\neq\varnothing\).
A \emph{uniform measure orbit equivalence} implementing an affine homeomorphism \(H:K_1\to K_2\) is a Borel bijection \(\theta:X_1^0\to X_2^0\) between invariant universally conull Borel sets such that:
\begin{enumerate}
\item \(\theta(\Gamma_1x)=\Gamma_2\theta(x)\) and \(\theta_*\mu=H\mu\);
\item \(f\circ\theta^{-1}\in\mathcal D_2\) for \(f\in\mathrm{C}(X_1)\), and \(k\circ\theta\in\mathcal D_1\) for \(k\in\mathrm{C}(X_2)\);
\item for every \(g\in\Gamma_1\) and \(h\in\Gamma_2\), the cocycle fiber \(\{x:\theta(gx)=h\theta(x)\}\) has indicator in \(\mathcal D_1\), and the corresponding fibers for \(\theta^{-1}\) have indicators in \(\mathcal D_2\).
\end{enumerate}
Here \(\mathcal D_i=\overline{\mathrm{C}(X_i)}^{K_i}\). Membership means uniform $L^2(K_i)$ approximability by a bounded sequence of continuous functions, so all conditions refer to the spaces and measure families.
\end{defn}

\begin{prop}[Spatial transport and the action case]
\label{prop:uniform-relation-isomorphism-transport}
A uniform measure relation isomorphism induces, by conjugation,
an isomorphism of the full pseudogroups and an isometric isomorphism
of their full groups.  For free continuous actions with the atlases consisting of
the global group maps and the empty map, this notion is precisely uniform measure orbit
equivalence from Definition~\ref{defn:uniform-measure-oe}.
\end{prop}

\begin{proof}
Pullback by $\theta$ and its inverse preserve the uniform $2$-norm because
$\theta_*\mu=H(\mu)$.  Decompose an admissible partial map into chart
restrictions as in Proposition~\ref{prop:regular-relation-pseudogroup-closure}.
Their conjugates are admissible restrictions of the conjugated charts.
Their total domain projection belongs to the target diagonal, so their
disjoint union is admissible by that proposition.  Applying the same
argument to $\theta^{-1}$ gives an isomorphism of the pseudogroups.
For total maps, measure transport gives
\[
 d_L(\theta T\theta^{-1},\theta S\theta^{-1})=d_K(T,S).
\]
In the action case, diagonal pullback is equivalent to the two observable
conditions in Definition~\ref{defn:uniform-measure-oe}: bounded uniform
$2$-norm closure extends pullback from continuous functions to the
completed diagonals.  The chart coefficients of $\theta g\theta^{-1}$
are exactly the cocycle-fiber projections in that definition, transported
by $\theta$; the inverse condition gives its reverse cocycle-fiber
projections.  All remaining conditions are the common spatial orbit and
measure conditions.
\end{proof}

\begin{rem}[The classical single-measure case]
\label{rem:uniform-relation-single-measure}
If $K=\{\mu\}$, then $D_K=L^\infty(X,\mu)$, since bounded continuous
functions are dense in bounded measurable functions for the $2$-norm.
Every countable Borel atlas of the probability-measure-preserving relation
is regular: all its support and coincidence projections are measurable,
and measure-preserving pullbacks preserve $L^\infty$.  Every Borel partial
bijection with graph in $R$ is admissible.  Thus the full pseudogroup and
full group are the usual measured ones, and the isomorphism notion is
ordinary probability-measure-preserving relation isomorphism.  No
continuity of the Borel charts is required.  Compactness of the chosen
topological presentation introduces no additional restriction in this case.
\end{rem}

\begin{rem}
The bare data $(R,K)$ do not specify which Borel graph coefficients
are uniformly approximable by the continuous diagonal.  That
regularity is retained by $D_K$ and the equivalence class of the
regular atlas. Section~\ref{sec:intrinsic-reconstruction} reconstructs
these regular data and their twist from uniformly separable completed
pairs under the partial-normalizer and fibrewise Cartan hypotheses.
For a given regular twist,
Theorem~\ref{thm:uniform-regular-untwisting} removes it under the
finite-chart hypotheses.
\end{rem}
\subsection{Extensions and countable group presentations}

\begin{lem}[Extension of an admissible partial map]
\label{lem:uniform-partial-map-extension}
Every $\varphi:A\to B$ in $\mathscr P$ agrees on its domain, modulo a
$K$-null set, with an element of $G$.
\end{lem}

\begin{proof}
Put $U=A\setminus B$ and $V=B\setminus A$.  The directed graph of
$\varphi$ consists of chains and cycles, since every vertex has at most
one incoming and one outgoing arrow.  Its finite nontrivial chains start
in $U$ and end in $V$.  For $m\geq1$, let
\[
 V_m=\{x\in V:x\in\operatorname{dom}(\varphi^{-m}),\quad
                         \varphi^{-m}x\in U\}.
\]
The domain projection of $\varphi^{-m}$ belongs to $\mathcal D$, and
corner pullback by that partial map shows that $1_{V_m}\in\mathcal D$.
The maps $\varphi^{-m}|_{V_m}$ close the chains of length $m$.  Their
domains are pairwise disjoint, as are their ranges.

The set $N_-=V\setminus\bigcup_mV_m$ consists of terminal vertices with
an infinite past.  All inverse iterates $\varphi^{-j}N_-$ are defined,
are pairwise disjoint, and have the same measure as $N_-$ for every
$\mu\in K$.  Hence $\mu(N_-)=0$ for every such measure.  Similarly,
the initial vertices in $U$ with an infinite future form a Borel set
$N_+$ which is null for every measure in $K$.  In particular, the total
domain projection of the disjoint union
\[
 \psi=\bigsqcup_{m\geq1}\varphi^{-m}|_{V_m}
\]
is $1_V$ modulo $K$-null sets.  Its total range projection is $1_U$.
By Proposition~\ref{prop:regular-relation-pseudogroup-closure},
$\psi\in\mathscr P$.  The union
\[
 \widetilde\varphi
   =\varphi\sqcup\psi\sqcup
                   \operatorname{id}_{X\setminus(A\cup B)}
\]
has disjoint domains and disjoint ranges, both of total projection $1$,
and is therefore an element of $G$.  It extends $\varphi$.  Removing
the $R$-saturation of $N_-\cup N_+$ gives an actual bijection on a
common invariant conull set.
\end{proof}

\begin{cor}[An admissible countable group presentation]
\label{cor:uniform-relation-countable-group-presentation}
Every uniform measure countable equivalence relation admits, on a
common invariant conull Borel set, a presentation as the orbit
relation of a countable group of admissible Borel automorphisms.
Its global word atlas, with the empty map adjoined, is equivalent to
the original regular atlas.
The action need not be free.
\end{cor}

\begin{proof}
Extend each chart of the regular countable atlas to a unit
$T_j\in G$ by Lemma~\ref{lem:uniform-partial-map-extension}.  The
subgroup generated by the $T_j$ is countable.  Choose representatives
for its elements and remove the $R$-saturations of the countably many
exceptional sets for their group identities and the extension
identities.  On the remaining invariant conull Borel set they form
an actual group action.  Its maps have graph in $R$, and their
graphs contain all the original atlas graphs, so their orbits are
exactly the $R$-classes.  Pseudogroup closure makes every word map
admissible; equalizers of admissible partial maps are regular since
they are fixed-point coefficients of inverse compositions.  Thus
the word atlas is regular, and every chart of either atlas is
admissible for the other.  Atlas independence proves the assertion.
\end{proof}

\section{Small boundaries and finite approximation}
\label{subsec:uniform-hyperfiniteness}

Let $(X,R,K,\mathcal D,\mathscr P)$ be a uniform measure relation as
above, and write $G=[R]^u_{K,\mathcal D}$.  All the statements below
are understood after restriction to an $R$-invariant Borel set which
has measure one for every $\mu\in K$.  The same set is used for the
whole sequence of relations or maps under consideration.

\subsection{Small boundaries and enough diagonal projections}
\label{subsec:relative-SBP-diagonal-projections}

The projection requirement in uniform AF
hyperfiniteness has a precise small-boundary interpretation.  This condition is needed for finite
diagonal approximation, rather than for the definition of a regular
uniform measure relation or its untwisted algebra.

\begin{defn}[Relative small boundary property {\cite[Definition~2.2]{EN}}]\label{defn:relative-K-SBP}
The compact metrizable pair $(X,K)$ has the \emph{relative small boundary
property}, or \emph{$K$-SBP}, if $X$ has a basis of open sets $U$ such that
\[
 \mu(\partial U)=0\qquad(\mu\in K).
\]
\end{defn}

For $K=\{\mu\}$ this condition is automatic.  At each point, the
metric spheres of distinct radii are disjoint, so only countably many
have positive $\mu$-measure.  Arbitrarily small balls of the remaining
radii give a small-boundary basis.  The same argument works for a
countable designated family.

\begin{prop}[Small boundaries and finite diagonals]
\label{prop:relative-SBP-finite-diagonals}
The following conditions are equivalent:
\begin{enumerate}
\item $(X,K)$ has $K$-SBP;
\item there are increasing finite-dimensional unital subalgebras
$F_n\subseteq D_K$, generated by finite partitions of unity of
projections, such that every $f\in D_K$ has approximants $f_n\in F_n$
with
\[
 \|f_n\|\leq\|f\|,\qquad \|f-f_n\|_{2,K}\longrightarrow0.
\]
\end{enumerate}
Under (1), the generating partitions may be represented by nested
finite Borel partitions of $X$ whose mesh tends to zero, and the norm
closure of their finite diagonals contains the image of $\mathrm{C}(X)$.
\end{prop}

\begin{proof}
Assume (1).  First, the characteristic function of every small-boundary
open set belongs to $D_K$.  For $\varnothing\ne U\ne X$, put
\[
 a_n(x)=\min\{1,n\,d(x,X\setminus U)\},\qquad
 b_n(x)=\max\{0,1-n\,d(x,\overline U)\}.
\]
These are continuous, satisfy $0\leq a_n\leq1_U\leq b_n\leq1$,
and $b_n-a_n$ decreases pointwise to $1_{\partial U}$.
The continuous functions
$\mu\mapsto\int(b_n-a_n)\,d\mu$ on the compact set $K$ consequently
decrease to zero.  Dini's theorem gives
\[
 \|1_U-a_n\|_{2,K}^2
       \leq\sup_{\mu\in K}\int(b_n-a_n)\,d\mu\longrightarrow0.
\]
The cases $U=\varnothing,X$ are immediate.  Hence $1_U\in D_K$.

For each $n$, take a finite cover by small-boundary basis sets of
diameter less than $2^{-n}$.  Ordered differences make this a finite
Borel partition: replace a cover $U_1,\ldots,U_m$ by
\[
 U_i\setminus\bigcup_{j<i}U_j\qquad(1\leq i\leq m)
\]
and omit empty pieces.  Its projections belong to $D_K$ by finite
Boolean operations.  Joining the first $n$ such partitions gives
nested partitions $\mathcal P_n$ of mesh at most $2^{-n}$.
Let $F_n$ be their finite diagonal algebras.

For $f\in \mathrm{C}(X)$, choose a sample in each non-null partition piece and
use its $f$-value on that piece.  These step functions approximate $f$
in norm, by uniform continuity.  They can be chosen with norm at most
$\|f\|$ in the common trace-faithful quotient. Indeed, let $\mu_0$
be a dense barycenter as in Lemma~\ref{lem:orbit-borel-realization}.
Each nonzero partition projection has positive $\mu_0$-measure, so
its piece admits a sample in $\operatorname{supp}(\mu_0)$.
Trace-null pieces contribute zero.
Thus $\overline{\bigcup_nF_n}^{\|\cdot\|}$ contains the continuous
diagonal.  Bounded uniform completion then gives (2).  Explicitly,
approximate $f\in D_K$ by bounded continuous functions, clip them to
the disk of radius $\|f\|$, and approximate the resulting functions
by these finite step algebras.  Nesting permits choices $f_n\in F_n$
with the required norm bound and uniform $2$-norm convergence.

Conversely, assume (2).  Every self-adjoint $f\in D_K$ is approximated
in uniform $2$-norm by self-adjoint elements of the finite algebras
$F_n$ (replace an approximant by its real part).  Such elements have
finite spectrum.  An arbitrarily small real scalar perturbation
avoiding the negatives of their finitely many spectral values makes
them invertible in $D_K$.  Therefore invertible self-adjoint elements
of $D_K$ are uniformly $2$-norm dense in its self-adjoint part.
The small-boundary characterization in
\cite[Theorem~2.16 and Corollary~2.17]{EN} gives $K$-SBP.
\end{proof}

The general norm real rank zero lemma in \cite[Lemma~8.3]{TC}
applies to the bounded uniform completion $D_K$: its closed unit ball
is uniformly $2$-norm complete, so the finite-spectrum approximation
above implies norm real rank zero. Conversely, norm real rank zero
implies $K$-SBP by \cite[Theorem~2.16 and Corollary~2.17]{EN}.
Thus $K$-SBP is equivalent to real rank zero of the entire completion,
answering \cite[Remark~2.18]{EN}. The proof and the resulting
Boolean-algebra consequences are given in \cite[Section~8]{TC}.
Neither this equivalence nor the finite-diagonal characterization
requires amenability, aperiodicity, or a Bauer hypothesis on $K$.

\begin{cor}[Invariance of relative SBP]
\label{cor:relative-SBP-uniform-invariance}
Relative SBP is preserved by a designated-trace-preserving isomorphism
of completed diagonals, and hence by a uniform measure relation
isomorphism.
\end{cor}
\begin{proof}
Such an isomorphism transports the finite projection diagonals,
their norm bounds, and their uniform $2$-norm density.
Apply Proposition~\ref{prop:relative-SBP-finite-diagonals} in both directions.
\end{proof}

\begin{cor}[An AF diagonal core adapted to the atlas]
\label{cor:relative-SBP-AF-diagonal-core}
If $(X,K)$ has $K$-SBP, the separable Cartan core in
Proposition~\ref{prop:uniform-relation-separable-Cartan-core} may be
chosen with a commutative AF diagonal $D_0\subseteq D_K$ containing
the continuous diagonal.  Its support and coincidence projections
belong to $D_0$, it is invariant under every atlas corner pullback,
and its bounded uniform completion is $D_K$.
\end{cor}

\begin{proof}
Start with the countably many partition projections just constructed
and all the countably many support and coincidence projections of
the regular atlas.  Repeatedly adjoin
\[
 \alpha_s(q_sp)\qquad(s\in\mathcal S)
\]
for every projection $p$ already in this countable collection.
These are again projections.  The norm-closed algebra generated by
the resulting countable set of commuting projections is separable
and AF, and contains $\mathrm{C}(X)$ by the preceding proof.  It is invariant
under the corner pullbacks by continuity and its defining closure
operation.  The Cartan-core construction therefore applies with this
$D_0$.  Its bounded uniform completion is $D_K$ because it contains
the continuous diagonal and is contained in $D_K$.
\end{proof}

\begin{rem}\label{rem:relative-SBP-scope}
In particular, uniform AF hyperfiniteness implies $K$-SBP, by its
finite-diagonal density requirement.  Conversely, $K$-SBP supplies
enough regular diagonal projections, and even an AF diagonal core;
it does not by itself supply finite orbit relations, regular
transversals, or compatible matrix charts in the ambient algebra.
Those retain the additional orbit information.
The equivalence with norm real rank zero of $D_K$ is proved in
\cite[Theorem~8.4]{TC}; the finite
orbit relations and compatible matrix charts require the separate
relation-theoretic arguments below.
\end{rem}

\subsection{Finite-dimensional unit spaces}
\label{subsec:finite-dimensional-relation-SBP}

Finite covering dimension supplies relative small boundaries when the
regular relation is aperiodic. The Borel charts first admit compact
homeomorphic representatives; general position is then applied on
those compact sets. Write $c_K(E)=\sup_{\mu\in K}\mu(E)$ for a Borel
set $E$.

\begin{lem}[Compact homeomorphic representatives]
\label{lem:regular-atlas-compact-homeomorphic-pieces}
Let $(X,R,K,D_K,\mathcal S)$ have a regular countable atlas. There
are countably many homeomorphisms $h_j:C_j\to D_j$ between compact
subsets of $X$, and one invariant $K$-conull Borel set $Y$, such that:
\begin{enumerate}
\item $x\in C_j\cap Y$ implies $h_j(x)\in Y$ and $xRh_j(x)$;
\item every pair $xRy$ in $Y$ is represented by some $h_j$:
$x\in C_j$ and $h_j(x)=y$.
\end{enumerate}
The compact representatives need not have graphs in $R$ outside
$Y$, and their compact domains need not be open.
\end{lem}

\begin{proof}
We first record the capacity quasi-continuity argument; compare
\cite[Proposition~24]{DHP}. Every $u\in D_K$ has a Borel representative
$\widehat u$ continuous outside open sets of arbitrarily small
$c_K$-capacity. Indeed, take bounded continuous approximants $u_n$
with $\|u_{n+1}-u_n\|_{2,K}^2\leq2^{-4n}$ and put
\[
 O_N=\bigcup_{n\geq N}\{|u_{n+1}-u_n|>2^{-n}\}.
\]
Chebyshev's inequality gives
\[
 c_K(O_N)\leq\sum_{n\geq N}2^{2n}\|u_{n+1}-u_n\|_{2,K}^2
 \leq\sum_{n\geq N}2^{-2n}\longrightarrow0.
\]
The sequence converges uniformly on each $X\setminus O_N$.
Its pointwise limit off the $K$-null set $\bigcap_NO_N$, defined
arbitrarily on that set, represents $u$ and has the asserted continuity.
It agrees with any given Borel representative outside a $K$-null set.
For a countable family, choose exceptional open sets whose capacities
sum to less than $\varepsilon$. Their complement is one compact set
on which all chosen representatives are continuous. Representatives
of real elements taking values in $[0,1]$ may be clipped to that interval.

Fix $s:A\to B$ and a continuous embedding
$e=(e_k):X\to[0,1]^{\N}$. Regularity gives
\[
 1_A\in D_K,\qquad
 u_k=\alpha_s(1_Be_k)=1_A(e_k\circ s)\in D_K.
\]
Choose quasi-continuous representatives and a compact $F^0$ of
complement capacity less than $\varepsilon$ on which they are
all continuous. Since $e(X)$ is closed, the set
\[
 F=\{x\in F^0:\widehat{1_A}(x)=1,\
                   (\widehat u_k(x))_{k\geq1}\in e(X)\}
\]
is compact, and $H(x)=e^{-1}((\widehat u_k(x))_{k\geq1})$ is
continuous on $F$. Outside one $K$-null set, $F\subseteq A$ and
$H=s$ on $F$. More precisely, $F\setminus A$ and the disagreement
set in $F\cap A$ are $K$-null, while
$c_K(A\setminus F)<\varepsilon$.

Apply the same construction to $s^{-1}$, obtaining a compact $L$
and a continuous $J:L\to X$ with $c_K(B\setminus L)<\varepsilon$
and $J=s^{-1}$ outside a $K$-null set. Put
\[
 C=\{x\in F:H(x)\in L,\ J(H(x))=x\}.
\]
This is compact, and $h=H|_C$ is a homeomorphism onto its compact
image, with inverse $J$. Outside the exceptional sets, a point of
$A\cap F$ fails to lie in $C$ only if its $s$-image fails to lie
in $L$. Invariance therefore gives
\[
 c_K(A\setminus C)
 \leq c_K(A\setminus F)+c_K(s^{-1}(B\setminus L))
 <2\varepsilon.
\]
The points of $C\setminus A$ and the disagreements between $h$ and
$s$ on $C\cap A$ are $K$-null.

Take $\varepsilon\to0$ along a sequence for every atlas chart.
For each $s$, the resulting compact domains cover $A_s$ modulo
a $K$-null set, since the missed domain has capacity bounded by
each of these errors. Remove all missed domains, all source
discrepancies just described, and the original atlas exceptional
set, together with their $R$-saturations. This removes one
$K$-null Borel set. Its complement $Y$ is invariant and satisfies
the two assertions.
\end{proof}

\begin{thm}[Finite dimensionality and small boundaries]
\label{thm:finite-dimensional-relation-SBP}
Let $(X,R,K,D_K,\mathcal S)$ be a regular uniform measure relation.
Suppose that $\dim X=d<\infty$ and the classes are infinite on
an invariant Borel set conull for every $\mu\in K$. Then $(X,K)$
has relative $K$-SBP.

More precisely, there are an invariant $K$-conull Borel set $Y$
and a dense $G_\delta$ subset $\mathcal H\subseteq \mathrm{C}(X,\mathbb R)$
such that
\begin{equation}\label{eq:finite-dimensional-orbit-zero-bound}
 \#(g^{-1}(0)\cap[x]_R)\leq d
 \qquad(g\in\mathcal H,\ x\in Y).
\end{equation}
Consequently $\mu(g^{-1}(0))=0$ for every $g\in\mathcal H$ and
$\mu\in K$. No amenability or minimality assumption is required.
\end{thm}

\begin{proof}
Use Lemma~\ref{lem:regular-atlas-compact-homeomorphic-pieces} and
intersect its invariant conull set with the infinite-class locus.
Enlarge its countable list by restricting every $h_j$ to $C_j\cap L$,
where $L$ ranges over a countable family of closed neighborhoods
that can separate finitely many distinct points. Closed balls
with centers in a countable dense subset and positive rational
radii suffice. The resulting compact homeomorphisms still cover
every pair in $R|_Y$.

Put $q=d+1$. For each $q$-tuple of listed homeomorphisms
$h_i:C_i\to D_i$ with pairwise disjoint source domains, let
$\mathcal O_{(h_i)}$ consist of the functions $g\in \mathrm{C}(X,\mathbb R)$
satisfying
\begin{equation}\label{eq:compact-chart-zero-avoidance}
 \bigcap_{i=1}^{q}h_i(C_i\cap g^{-1}(0))=\varnothing.
\end{equation}
This condition is open and dense for the uniform norm. Put
$B=\bigcap_iD_i$. If $B=\varnothing$ there is nothing to prove.
Otherwise the condition says that the continuous map
\[
 y\longmapsto
 \bigl(g(h_1^{-1}y),\ldots,g(h_q^{-1}y)\bigr),\qquad y\in B,
\]
avoids zero. Compactness of $B$ makes this an open condition.

For density, fix $f\in \mathrm{C}(X,\mathbb R)$ and $\delta>0$, and put
\[
 F:B\longrightarrow\mathbb R^q,\qquad
 F(y)=\bigl(f(h_1^{-1}y),\ldots,f(h_q^{-1}y)\bigr).
\]
Since $\dim B\leq d<q$, there is a continuous
$G:B\to\mathbb R^q\setminus\{0\}$ with
$\sup_{y\in B}\max_i|G_i(y)-F_i(y)|<\delta$.
Here is the dimension argument explicitly. Refine a finite cover
controlling the oscillation of $F$ to one of multiplicity at most
$d+1$, and take a subordinate partition of unity. Approximate $F$
by the corresponding convex combinations of finitely many vectors.
Perturb these vectors slightly so that zero belongs to the affine
span of no subset of at most $d+1$ of them. Every active convex
combination then avoids zero, giving the required $G$.

The scalar corrections $G_i(h_i(x))-f(x)$ are continuous on the
pairwise disjoint compact sets $h_i^{-1}(B)$. They define one
continuous correction on their union. Tietze extension gives
$b\in \mathrm{C}(X,\mathbb R)$ with $\|b\|<\delta$ agreeing with these
corrections. Then $g=f+b$ satisfies
\eqref{eq:compact-chart-zero-avoidance}, proving density.

There are countably many such tuples. Baire's theorem gives a dense
$G_\delta$ intersection $\mathcal H$ of their open dense sets.
If $q$ distinct zeros $x_1,\ldots,x_q$ of some $g\in\mathcal H$
belonged to one class in $Y$, choose a common range point $y$
and compact representatives sending $x_i$ to $y$. Restrict their
domains to pairwise disjoint closed neighborhoods of the $x_i$.
These maps are in our list, and $y$ belongs to every set on the
left of \eqref{eq:compact-chart-zero-avoidance}, a contradiction.
This proves \eqref{eq:finite-dimensional-orbit-zero-bound}.

Set $Z=g^{-1}(0)$ and fix $\mu\in K$. Counting-measure symmetry
for the probability-measure-preserving relation $R|_Y$ gives
\[
 \int_Y\#(Z\cap[x]_R)\,d\mu(x)
 =\int_{Z\cap Y}\#[x]_R\,d\mu(x).
\]
This follows by decomposing the countable Borel relation into
Borel partial bijections and using invariance; the integrals are
extended nonnegative integrals. The left side is at most $d$.
The right side would be infinite if $\mu(Z)>0$, since the classes
in $Y$ are infinite. Thus $\mu(Z)=0$ for every $\mu\in K$.

Finally, let $A\subseteq O$ with $A$ closed and $O$ open. Choose
$f\in \mathrm{C}(X,\mathbb R)$ equal to $-1$ on $A$ and to $1$ on
$X\setminus O$, and choose $g\in\mathcal H$ with $\|g-f\|<1/2$.
Then $W=\{g<0\}$ satisfies
\[
 A\subseteq W\subseteq\overline W\subseteq O,\qquad
 \partial W\subseteq g^{-1}(0).
\]
Its boundary is null for every designated measure. Taking $A$
to be a singleton gives a small-boundary basis.
\end{proof}

\begin{cor}[Atomless full invariant families]
\label{cor:finite-dimensional-atomless-full-family-SBP}
Suppose $\dim X<\infty$, $K=M(R)$ is nonempty and compact, and
$(X,R,K)$ has a regular countable atlas. If every measure in $K$
is atomless, then $(X,K)$ has relative SBP.
\end{cor}
\begin{proof}
A finite $R$-class would support an atomic $R$-invariant probability,
namely its normalized counting measure. Thus every class is infinite,
and Theorem~\ref{thm:finite-dimensional-relation-SBP} applies.
\end{proof}

\begin{rem}[Scope of the finite-dimensional criterion]
\label{rem:finite-dimensional-SBP-scope}
For a designated subfamily $K$, atomlessness alone does not rule
out finite classes: the equality relation preserves every probability
measure. The theorem retains the hypothesis of infinite classes
on a common $K$-conull set. Conversely, under that hypothesis
invariance already forces every designated measure to be atomless,
since an atom would give infinitely many atoms of the same positive
mass in its class. The compact pieces in the proof need no open
extensions. Example~\ref{ex:regular-atlas-no-original-etale-model}
shows why regularity does not require an \'etale presentation on
the original unit space.
\end{rem}

\subsection{Uniform hyperfiniteness and finite matrix charts}

\begin{defn}\label{defn:uniform-hyperfinite-relation}
The uniform measure relation is \emph{uniformly hyperfinite} if there
are increasing finite Borel equivalence relations
\[
 R_1\subseteq R_2\subseteq\cdots\subseteq R
\]
such that the following conditions hold:
\begin{enumerate}
\item each $R_n$ is generated by finitely many elements of the
admissible pseudogroup $\mathscr P$;
\item $R=\bigcup_nR_n$ on one $R$-invariant, $K$-universally conull
Borel set;
\item for every $T\in G$,
\[
 \lim_{n\to\infty}\sup_{\mu\in K}
 \mu\{x:(x,Tx)\notin R_n\}=0.
 \tag{UH}\label{eq:uniform-hyperfinite-escape}
\]
\stepcounter{equation}
\end{enumerate}
Here ``finite'' means that each equivalence class is finite.  When
there are integers $N_n$ with $|[x]_{R_n}|\leq N_n$ for every $x$, we
say that the exhaustion has \emph{bounded finite levels}.
\end{defn}

For a single probability measure, the increasing exhaustion already
implies \eqref{eq:uniform-hyperfinite-escape}.  For a family of
measures, condition \eqref{eq:uniform-hyperfinite-escape} requires
the convergence to be uniform over the family.  It also refers to
the admissible full group, and hence to the prescribed uniform
observable algebra.

By Lemma~\ref{lem:uniform-partial-map-extension}, condition
\eqref{eq:uniform-hyperfinite-escape} is equivalent to the same
escape condition for every admissible partial map, restricted to its
domain: extend it to a total admissible map and compare escape sets.

The following stronger spatial condition records the finite matrix
charts used by an AF model.  It is useful to distinguish this
condition from Definition~\ref{defn:uniform-hyperfinite-relation}.

\begin{defn}\label{defn:uniform-AF-hyperfinite-relation}
A \emph{finite $\mathcal D$-regular matrix chart system} consists of
a finite Borel partition
\[
 X=\bigsqcup_{a=1}^{r}\bigsqcup_{i=1}^{k_a}Y_{a,i},
 \qquad 1_{Y_{a,i}}\in\mathcal D,
\]
and admissible partial bijections
$s_{a,i}:Y_{a,1}\longrightarrow Y_{a,i}$, with
$s_{a,1}=\operatorname{id}_{Y_{a,1}}$.  It defines the finite relation
\[
 S=\bigcup_{a=1}^{r}\bigcup_{i,j=1}^{k_a}
       \operatorname{graph}(s_{a,i}s_{a,j}^{-1}).
\]
Thus $Y_{a,1}$ is a Borel transversal for the corresponding block,
and its partial bijections give $k_a$ by $k_a$ matrix charts.  The
associated finite diagonal is
\[
 E_S=\operatorname{span}_{\mathbb C}
       \{1_{Y_{a,i}}:1\leq a\leq r,\ 1\leq i\leq k_a\}.
\]

The relation is \emph{uniformly AF hyperfinite} if it has a uniformly
hyperfinite exhaustion $(R_n)$ coming from finite
$\mathcal D$-regular matrix chart systems, with the following two
additional properties:
\begin{enumerate}
\item the chart systems have coordinate refinements: the level
$n+1$ partition refines the level $n$ partition, and every level
$n$ matrix-chart map is a finite disjoint union of matching level
$n+1$ matrix-chart maps;
\item the finite diagonals $E_n=E_{R_n}$ exhaust $\mathcal D$ in the
bounded uniform $2$-norm sense: for every $f\in\mathcal D$ there are
$f_n\in E_n$ with $\|f_n\|\leq\|f\|$ and
$\|f-f_n\|_{2,K}\longrightarrow0$.
\end{enumerate}
The coordinate refinements are the spatial form of the usual
diagonal-preserving, phase-free embeddings of finite-dimensional
matrix algebras.
\end{defn}

The finite matrix charts provide more information than the sizes of
the finite equivalence classes.  In particular, they supply regular
transversals, compatible matrix coordinates, and enough finite
diagonal functions to approximate all uniform observables.  These
are precisely the data to be retained when one uses an AF
presentation, rather than only a finite-relation exhaustion.

\begin{prop}\label{prop:dense-locally-finite-uniform-hyperfinite}
If
\[
 L_1\subseteq L_2\subseteq\cdots\subseteq G
\]
are finite subgroups and $\bigcup_nL_n$ is dense in $G$ for $d_K$,
then $R$ is uniformly hyperfinite with bounded finite levels.  The
level $n$ relation may be taken to be the orbit relation of $L_n$;
its classes have size at most $|L_n|$.
\end{prop}

\begin{proof}
By Corollary~\ref{cor:uniform-relation-countable-group-presentation},
choose a countable admissible global graphing $T_1,T_2,\ldots$ of $R$
on a common invariant conull set.
Put $L=\bigcup_nL_n$.  This is a countable group of equivalence
classes of admissible Borel maps.  Choose one Borel representative
$S_\ell$ for every $\ell\in L$, and representatives for all the
graphing maps $T_i$.  First restrict to the intersection of their
domains and of the common invariant conull set on which the
graphing generates $R$.  Each of these domains is $R$-invariant and
has measure one for every $\mu\in K$.

For $\ell,m\in L$, the set
\[
 B_{\ell,m}=
 \{x:S_\ell S_mx\ne S_{\ell m}x\}
\]
is null for every $\mu\in K$.  The same is true of
$B_e=\{x:S_ex\ne x\}$.  There are countably many such sets.  Their
union therefore remains $K$-universally null.  Saturating that
union under the countable graphing $T_i,T_i^{-1}$ still gives a
Borel, $K$-universally null set, since every graphing map preserves
every measure in $K$.  Its complement is $R$-invariant.  On this
complement the maps $S_\ell$ satisfy all the group identities
pointwise and give an actual Borel action of $L$.  In particular,
each $L_n$ acts as a finite group.  The representatives retain their
admissibility after this restriction.

For each $i,m\geq1$, use density to choose
$\ell_{i,m}\in L$ such that
\[
 d_K(T_i,S_{\ell_{i,m}})<2^{-m}.
\]
Let $A_{i,m}=\{x:T_ix\ne S_{\ell_{i,m}}x\}$.  For every $\mu\in K$,
\[
 \sum_{m=1}^{\infty}\mu(A_{i,m})
 \leq\sum_{m=1}^{\infty}2^{-m}<\infty.
\]
Thus the Borel set $\limsup_mA_{i,m}$ is null for every $\mu\in K$.
Take the union of these limsup sets over all $i$ and remove its
$R$-saturation as well.  As before, this removes one Borel set which
is null for every measure in $K$.  Denote the remaining invariant
set by $X_0$.

On $X_0$, every $T_ix$ is equal to $S_{\ell_{i,m}}x$ for all
sufficiently large $m$.  Define
\[
 R_n=\{(x,S_\ell x):x\in X_0,\ \ell\in L_n\}.
\]
These are increasing Borel equivalence relations, generated by
finitely many admissible maps, and their classes have size at most
$|L_n|$.  Their union contains the graph of every $T_i$ on $X_0$.
It therefore contains $R|_{X_0}$.  Conversely, each $S_\ell$ has
graph in $R$, so
\[
 R|_{X_0}=\bigcup_nR_n.
\]

Finally, let $T\in G$.  For every $\ell\in L_n$, up to a common null
set,
\[
 \{x:(x,Tx)\notin R_n\}
 \subseteq\{x:Tx\ne S_\ell x\}.
\]
Consequently,
\[
 \sup_{\mu\in K}\mu\{x:(x,Tx)\notin R_n\}
 \leq\inf_{\ell\in L_n}d_K(T,S_\ell)
 \longrightarrow0,
\]
because the subgroups are increasing and their union is dense.
This proves all the required conditions.
\end{proof}

\subsection{Regular uniform Reiter amenability}
\label{subsec:uniform-Reiter-amenability}

Let $(X,R,K,\mathcal D,\mathscr P)$ be a uniform measure relation
with regular countable atlas $\mathcal S$.  For two probability
measures on the countable class $[x]_R$, use the convention
\[
 \|\eta-\zeta\|_1
     =\sum_{y\in[x]_R}|\eta(\{y\})-\zeta(\{y\})|.
\]
This distance is at most two.  All kernels and partial maps below
are represented on one invariant Borel set conull for every
measure in $K$.

\begin{defn}\label{defn:uniform-Reiter-amenability}
The relation has \emph{regular uniform Reiter amenability} if there
are probability kernels $x\mapsto m_n^x$ on
$[x]_R$ of the form
\[
 m_n^x=\sum_{j=1}^{J_n}f_{n,j}(x)\delta_{\varphi_{n,j}x},
 \qquad \varphi_{n,j}\in\mathscr P,
\]
where $J_n<\infty$, the coefficients $f_{n,j}\in\mathcal D$ are
nonnegative, are supported on $\operatorname{dom}(\varphi_{n,j})$,
and satisfy $\sum_jf_{n,j}=1$, such that, for every $s\in\mathcal S$,
\begin{equation}\label{eq:uniform-Reiter-condition}
 \sup_{\mu\in K}\int_{\operatorname{dom}(s)}
          \|m_n^{sx}-m_n^x\|_1\,d\mu(x)\longrightarrow0.
\end{equation}
An atom with zero coefficient is omitted, so its partial map need
not be defined at $x$.  Repeated atoms are combined in the kernel.
\end{defn}

\begin{prop}\label{prop:uniform-Reiter-atlas-independence}
For kernels as in Definition~\ref{defn:uniform-Reiter-amenability},
condition \eqref{eq:uniform-Reiter-condition} is equivalent to its
validity for every $s\in\mathscr P$.  It is consequently independent
of the choice of an equivalent regular atlas.
\end{prop}

\begin{proof}
First observe that, for a fixed $n$ and $s\in\mathscr P$, the function
\[
 d_{n,s}(x)=1_{\operatorname{dom}(s)}(x)
                          \|m_n^{sx}-m_n^x\|_1
\]
belongs to $\mathcal D$.  The coefficients of the first kernel,
pulled back to the source of $s$, are corner pullbacks of the
$f_{n,j}$; its support maps are $\varphi_{n,j}s$.  The second kernel
has support maps $\varphi_{n,j}$.  Partition the source into the
finitely many Boolean pieces determined by the domains and the
pairwise equalizers of these partial maps.  All these projections
belong to $\mathcal D$.  On each piece the coinciding atoms form
fixed groups, and the distance is the sum of the absolute values
of the differences of their total coefficients.  Finite sums and
continuous functional calculus show that $d_{n,s}\in\mathcal D$.
In particular $0\leq d_{n,s}\leq2$.

Enumerate the original atlas as $(s_j)$, and let
$\varphi\in\mathscr P$.  Its first-chart decomposition gives
orthogonal projections $e_j\in\mathcal D$ on which
$\varphi=s_j$, with
\[
 \sup_{\mu\in K}\mu\left(1_{\operatorname{dom}(\varphi)}
                              -\sum_{j\leq N}e_j\right)
       \longrightarrow0\quad(N\to\infty).
\]
Since $d_{n,\varphi}=d_{n,s_j}$ on the $j$th piece, for every $N$,
\begin{align*}
 \sup_{\mu\in K}\int d_{n,\varphi}\,d\mu
 &\leq\sum_{j\leq N}\sup_{\mu\in K}
                     \int d_{n,s_j}\,d\mu\\
 &\quad+2\sup_{\mu\in K}
       \mu\left(1_{\operatorname{dom}(\varphi)}
                              -\sum_{j\leq N}e_j\right).
\end{align*}
For fixed $N$ the first term tends to zero as $n\to\infty$.
Letting $N\to\infty$ proves the condition for $\varphi$.
The converse follows because atlas charts are admissible.
Equivalent atlases have the same admissible pseudogroup, so they
have exactly the same finite regular kernels and tests.
\end{proof}

\begin{prop}\label{prop:uniform-AF-implies-Reiter}
Every uniformly AF hyperfinite relation has regular uniform Reiter
amenability.
\end{prop}

\begin{proof}
At level $n$, write its finite matrix charts as
$s_{a,i}:Y_{a,1}\to Y_{a,i}$, with $1\leq i\leq k_a$.
For $x\in Y_{a,j}$ define
\[
 m_n^x=\frac1{k_a}\sum_{i=1}^{k_a}
                \delta_{s_{a,i}s_{a,j}^{-1}x}.
\]
This is a finite regular kernel: its coefficients are
$k_a^{-1}1_{Y_{a,j}}\in\mathcal D$, and its support maps are the
admissible matrix-chart maps.  It is exactly the uniform probability
measure on the finite class $[x]_{R_n}$.

For $\varphi\in\mathscr P$, the two $R_n$-classes at $x$ and
$\varphi x$ are either equal or disjoint.  Therefore, pointwise on
its domain,
\[
 \|m_n^{\varphi x}-m_n^x\|_1
       =2\,1_{\{(x,\varphi x)\notin R_n\}}.
\]
The uniform escape condition applies to partial maps by
Lemma~\ref{lem:uniform-partial-map-extension}.  Integrating the last
identity proves all the Reiter tests.  This estimate uses only the
regular finite charts and uniform escape, rather than the further
coordinate-refinement and diagonal-density conditions.
\end{proof}

\begin{rem}\label{rem:uniform-Reiter-CFW-program}
The same proof applies to an exhaustion by finite admissible group
orbits: averaging $\delta_{\ell x}$ over the finite group gives the
uniform measure on each orbit, since every point in an orbit has
the same number of preimages under $\ell\mapsto\ell x$.

Under relative SBP, the converse finite-approximation direction is
Theorem~\ref{thm:uniform-CFW-SBP}: regular Reiter kernels yield
deterministic regular finite cuts, and
Theorem~\ref{thm:local-cuts-AF} supplies compatible refinements.
The unconditional strong AF conclusion fails without relative SBP,
as Example~\ref{ex:Reiter-without-AF-coordinates} shows.
\end{rem}

\subsection{From measurewise hyperfiniteness to uniform kernels}
\label{sec:uniform-Reiter-minimax}

Let $(X,R,K,D_K,\mathcal S)$ be a uniform measure relation.  Enumerate
its regular atlas as $(s_j)$ with $s_1=\operatorname{id}$, and put
\[
 r_j=p_{s_j}\prod_{i<j}(1-c(s_j,s_i)).
\]
The restricted graphs $s_j|_{\operatorname{supp}(r_j)}$ are disjoint
and cover $R$, and $r_1=1$.  These are the disjoint arrow coordinates
used in the untwisted relation algebra.

We use the following sufficient subclass of the kernels in
Definition~\ref{defn:uniform-Reiter-amenability}.  A disjoint finite
regular probability kernel is a finitely supported sequence $k=(k_j)$
with $k_j\in r_jD_K$, $k_j\geq0$, and $\sum_jk_j=1$.
Its row is $m_k(x)=\sum_jk_j(x)\delta_{s_jx}$.
The restricted source maps belong to $\mathscr P$.
Write $e_s(k)(x)=1_{A_s}(x)\|m_k(sx)-m_k(x)\|_1$.

For a finite kernel, $e_s(k)\in D_K$.  Indeed, the two rows use only
finitely many output charts, namely restrictions of $s_j$ and $s_j s$.
Disjointize this finite family by its regular coincidence projections.
On each resulting arrow coordinate the difference coefficient is a
finite sum of diagonal coefficients and their corner pullbacks.  The
row $\ell^1$ norm is the sum of their absolute values and therefore
belongs to $D_K$.  In particular,
$\mu\mapsto\int e_s(k)\,d\mu$ is continuous and affine on $K$.

\begin{thm}\label{thm:measurewise-hyperfinite-uniform-Reiter}
If $(X,R,\mu)$ is hyperfinite as an ordinary measured relation for
every $\mu\in K$, then the uniform measure relation has regular
uniform Reiter amenability.
\end{thm}

\begin{proof}
Fix first one $\mu\in K$ and finitely many test charts
$s_{i_1},\ldots,s_{i_t}$.  Choose a $\mu$-conull invariant set on which
$R=\bigcup_nR_n$ with increasing finite Borel subrelations $R_n$.
The probability kernel which is uniform on $[x]_{R_n}$ is measurable.
Its rows at $x$ and $sx$ agree whenever $(x,sx)\in R_n$; elsewhere their
$\ell^1$ distance is at most two.  Dominated convergence consequently
gives kernels whose integrated errors tend to zero for all the chosen
tests.  Extend their rows by the identity probability on the omitted
$\mu$-null set.

We show that these measurable kernels can be replaced, for this $\mu$,
by finite regular ones.  For any measurable probability kernel $m$,
retain its weights on the first $J$ disjoint arrow coordinates and
place the omitted mass on the identity coordinate.  This gives a
finite measurable probability kernel $b=(b_1,\ldots,b_J)$ with
$b_j=r_jb_j$, and
\[
 \int_X\|m_b(x)-m(x)\|_1\,d\mu(x)\longrightarrow0
 \qquad(J\to\infty).
\]
The estimate follows from twice the omitted row mass and dominated
convergence.  By $L^1(\mu)$ density of continuous functions, choose
$a_j=r_jh_j\in r_jD_K$, with $h_j\in \mathrm{C}(X)$ and $0\leq h_j\leq1$,
so that $\sum_j\int|a_j-b_j|\,d\mu$ is arbitrarily small.  Put
\[
 A=\sum_{j=1}^Ja_j,\qquad c=\max\{1,A\},\qquad
 k_j=\frac{a_j}{c}
       +\left(1-\frac{A}{c}\right)1_{\{j=1\}}.
\]
Continuous functional calculus in $D_K$ makes these regular positive
coefficients, and their row sum is exactly one for all designated
measures.  The fallback term is allowed because $r_1=1$.  Pointwise,
\[
 \sum_j|k_j-b_j|\leq2\sum_j|a_j-b_j|.
\]
Indeed, normalization changes the vector $(a_j)$ by $|A-1|$ in
$\ell^1$, and $|A-1|\leq\sum_j|a_j-b_j|$.

For two probability kernels $m,m'$, invariance of $\mu$ under $s$ gives
\[
 \left|\int e_s(m)\,d\mu-\int e_s(m')\,d\mu\right|
 \leq2\int_X\|m(x)-m'(x)\|_1\,d\mu(x).
\]
This is the triangle inequality integrated over $A_s$, followed by the
change of variables $x\mapsto sx$ for the first row.  The preceding
truncation and normalization therefore show, for every $\mu\in K$,
\begin{equation}\label{eq:single-measure-regular-Reiter-infimum}
 \inf_k F_k(\mu)=0,\qquad
 F_k(\mu)=\sum_{\ell=1}^t\int e_{s_{i_\ell}}(k)\,d\mu,
\end{equation}
where the infimum runs over finite regular probability kernels.

We now make the choice uniform.  Each $F_k$ is continuous and affine
on the compact convex set $K$.  Kernels form a convex set, and the
triangle inequality gives
$F_{\lambda k+(1-\lambda)k'}\leq\lambda F_k+(1-\lambda)F_{k'}$.
Suppose $\inf_k\sup_KF_k>0$, and choose $\delta>0$ smaller than this
infimum.  In the real Banach space $\mathrm{C}(K)$ the convex upward set
\[
 \mathcal C=\{f:\ F_k\leq f\text{ for some }k\}
\]
is disjoint from the open convex set
$\mathcal U=\{f:\sup_Kf<\delta\}$.
Hahn--Banach separation gives a nonzero continuous functional
$\Lambda$ with $\sup_{\mathcal U}\Lambda\leq\inf_{\mathcal C}\Lambda$.
Upward closure of $\mathcal C$ forces $\Lambda$ to be positive.
Normalize $\Lambda(1)=1$, and represent it by a probability measure
$\rho$ on $K$.  Since constants below $\delta$ belong to $\mathcal U$,
$\Lambda(F_k)\geq\delta$ for every $k$.  The barycenter
$\bar\mu=\int_K\mu\,d\rho(\mu)$ belongs to $K$, and affinity gives
$F_k(\bar\mu)=\Lambda(F_k)\geq\delta$, contradicting
\eqref{eq:single-measure-regular-Reiter-infimum}.
Thus $\inf_k\sup_KF_k=0$.  Choose a kernel for the first $n$ atlas
tests with supremum of their summed errors below $1/n$.
These choices give the required sequence.
\end{proof}

\begin{cor}\label{cor:measurewise-amenable-uniform-Reiter}
If $(X,R,\mu)$ is amenable as an ordinary measured relation for every
$\mu\in K$, then the uniform measure relation has regular uniform Reiter amenability.
\end{cor}

\begin{proof}
Classical Connes--Feldman--Weiss~\cite{CFW} gives measured hyperfiniteness
for each individual $\mu$.  Apply
Theorem~\ref{thm:measurewise-hyperfinite-uniform-Reiter}.
\end{proof}

\begin{rem}\label{rem:Reiter-versus-regular-cuts}
The minimax argument itself supplies probability kernels rather than
finite subrelations.  Theorem~\ref{thm:regular-Reiter-rounding}
supplies the additional deterministic construction under relative
SBP, using finite rational kernels and local conditional expectations.
The resulting cuts are regular in the prescribed diagonal, a
requirement not contained in ordinary measured hyperfiniteness.
\end{rem}

\begin{example}[Relative SBP is necessary for the AF conclusion]
\label{ex:Reiter-without-AF-coordinates}
Let $X=[0,1]$, let $R$ be equality, and let $K=\operatorname{Prob}(X)$.
The atlas consisting of the identity and empty map is regular.
Since $K$ contains every point mass,
\[
 \|f\|_{2,K}=\|f\|_\infty\qquad(f\in \mathrm{C}(X)),
\]
so $D_K=\mathrm{C}([0,1])$.  The exact kernel $m^x=\delta_x$ has zero
Reiter error, and the constant finite relation exhaustion $R_n=R$
is uniformly hyperfinite.  Nevertheless the relation is not uniformly
AF hyperfinite: connectedness implies that the only projections in
$D_K$ are $0$ and $1$, so every finite regular diagonal $E_n$ consists
of constants and cannot be uniformly $2$-norm dense in $D_K$.
Thus Reiter amenability alone does not imply the strong AF condition.
The missing hypothesis here is relative SBP: the pair
$([0,1],\operatorname{Prob}([0,1]))$ does not have it.  Equality
itself is generated by the identity.
\end{example}

\section{The untwisted uniform relation algebra}
\label{sec:uniform-relation-algebra}

We now construct an operator algebra from a regular countable relation.
This is a forward construction: no reconstruction assertion for an
arbitrary uniformly completed Cartan pair is needed.  Throughout this
section, $(X,R,K,D_K,\mathcal S)$ is as in
Definition~\ref{defn:regular-relation-atlas}.  We first restrict to the
common invariant conull set on which its countable atlas covers $R$.
All Borel equalities below are understood modulo sets null for every
measure in $K$ when representatives are changed.

\subsection{The regular representation and the designated traces}

Choose a dense sequence $(\mu_n)$ in $K$ and set
\[
 \mu_0=\sum_{n=1}^{\infty}2^{-n}\mu_n\in K.
\]
The state $f\mapsto\int f\,d\mu_0$ is faithful on $D_K$.  Indeed, if
$\int |f|^2\,d\mu_0=0$, then the continuous function
$\mu\mapsto\int |f|^2\,d\mu$ vanishes at every $\mu_n$, hence on $K$.
This observation concerns $D_K$; it does not say that every
$\mu_0$-null Borel set is null for all measures in $K$.

Put
\[
 \nu_0(F)=\int_X\sum_{y\in[x]_R}1_F(y,x)\,d\mu_0(x),
 \qquad H_0=L^2(R,\nu_0).
\]
The counting measure $\nu_0$ is sigma finite, since the graphs of the
countable atlas cover $R$.  A bounded Borel function $f$ acts by
$(f\xi)(y,x)=f(y)\xi(y,x)$.  For a chart $s:A_s\to B_s$, define
\[
 (V_s\xi)(y,x)=1_{B_s}(y)\xi(s^{-1}y,x).
\]
Then $V_s^*=V_{s^{-1}}$, $V_s^*V_s=1_{A_s}$,
$V_sV_s^*=1_{B_s}$, and $V_s$ implements its corner pullback on
the diagonal.  Let
\[
 \mathcal N_0=\mathrm{W}^*(L^\infty(X,\mu_0),V_s:s\in\mathcal S).
\]
Its diagonal expectation is denoted $\mathbb{E}_0$.  It is characterized on
finite chart kernels by restriction to the diagonal, and in particular
\begin{equation}
\label{eq:uniform-relation-regular-expectation}
 \mathbb{E}_0(V_sf)=c(s,\mathrm{id})f.
\end{equation}
It is normal, faithful, and completely positive, and
$\tau_0=\mu_0\circ \mathbb{E}_0$ is a faithful normal trace.

These elementary properties can be seen directly in the regular
representation.  Decomposing $H_0$ over its second coordinate identifies
$\mathbb{E}_0(a)(x)$ with the $(x,x)$ matrix entry of $a$ on
$\ell^2([x]_R)$, which proves positivity and diagonal bimodularity.
The vector $\delta=1_{\{(x,x):x\in X\}}$ is separating: the analogous
right chart operators, together with multiplication by bounded
functions of the second coordinate, commute with the left operators
and generate a dense subspace from $\delta$.  Thus $\mathbb{E}_0(a^*a)=0$ implies $a=0$.
Finally, $R$-invariance makes the counting measure symmetric under
$(y,x)\mapsto(x,y)$, which proves the trace identity on finite chart
kernels and then on their von Neumann closure.  This construction uses
the given relation and its untwisted graph operators only.

\begin{prop}[Designated traces]
\label{prop:uniform-relation-designated-traces}
Let
\[
 A_R=\mathrm{C}^*(D_K,V_s:s\in\mathcal S)\subseteq\mathcal N_0.
\]
The restriction of $\mathbb{E}_0$ is a faithful conditional expectation
$\mathbb{E}:A_R\to D_K$.  For every $\mu\in K$,
\[
 \tau_\mu(a)=\int_X \mathbb{E}(a)\,d\mu
\]
is a tracial state, and $\mu\mapsto\tau_\mu(a)$ is continuous on $K$
for every $a\in A_R$.
\end{prop}

\begin{proof}
The algebraic span of $V_sf$, $s\in\mathcal S$, $f\in D_K$, is a
unital star algebra.  This follows from closure of the atlas under
composition and inversion and from its corner pullback condition.
Formula~\eqref{eq:uniform-relation-regular-expectation} and regularity
of coincidence projections show that its expectation belongs to $D_K$.
Norm continuity of $\mathbb{E}_0$ and norm closedness of $D_K$ give the same
conclusion on $A_R$.  Faithfulness follows from faithfulness of $\mathbb{E}_0$.

Integration against $\mu$ is a state on $D_K$, so $\tau_\mu$ is a state
on $A_R$.  To check the trace identity, write $a,b$ as finite chart
kernels.  The diagonal of their product is
\[
 \mathbb{E}(ab)(x)=\sum_{y\in[x]_R}a(x,y)b(y,x).
\]
For these kernels the relevant sums have finite row support.  Symmetry
of the counting measure associated with the $R$-invariant measure
$\mu$ gives $\tau_\mu(ab)=\tau_\mu(ba)$.  This identity extends by
norm continuity.  Continuity on $K$ follows from $\mathbb{E}(a)\in D_K$.
\end{proof}

\begin{lem}[Faithful bounded uniform completion]
\label{lem:uniform-relation-faithful-completion}
The bounded uniform completion of $A_R$ for the designated traces
embeds faithfully into $\mathcal N_0$.  Denote its image by $M_R$.
It is a unital $\mathrm{C}^*$-algebra, is complete on each norm bounded set for
\[
 \|a\|_{2,K}=\sup_{\mu\in K}\tau_\mu(a^*a)^{1/2},
\]
and $\mathbb{E}$ extends to a faithful conditional expectation
$\mathbb{E}:M_R\to D_K$.  The designated traces extend to $M_R$ and remain
continuous as functions of $\mu\in K$.  The uniform $2$-norm topology
on its unit ball is separable.
\end{lem}

\begin{proof}
We give details because the measures in $K$ need not be absolutely
continuous with respect to $\mu_0$.  A norm bounded uniform
$2$-norm Cauchy sequence $(a_j)$ is Cauchy for the trace $\tau_0$,
and hence has a bounded $2$-norm limit in the finite von Neumann algebra
$\mathcal N_0$.  Also,
\[
 \tau_{\mu_n}\leq 2^n\tau_0
\]
on positive elements of $A_R$.  These dominated traces extend normally
to $\mathcal N_0$.  If the $\tau_0$-limit of $(a_j)$ is zero, then
$\lim_j\tau_{\mu_n}(a_j^*a_j)=0$ for every $n$.  On the other hand,
\[
 \mu\longmapsto\lim_j\tau_\mu(a_j^*a_j)
\]
is continuous on $K$: the functions being limited are uniformly
Cauchy by the Cauchy--Schwarz inequality and the common norm bound.
It vanishes on the dense sequence $(\mu_n)$, so it vanishes on $K$.
This proves injectivity of the map from the bounded uniform completion.
The same argument shows that the designated trace values are independent
of the approximating sequence and continuous on $K$.

For clarity, the $\mathrm{C}^*$-norm is the inherited operator norm in
$\mathcal N_0$.  If $x$ belongs to the image, its normally extended
$\mu_n$-traces satisfy
$\tau_{\mu_n}(x^*x)\leq\|x\|^2$.  Continuity on $K$ gives this
inequality for every $\mu\in K$.  Thus
$\|x\|_{2,K}\leq\|x\|$.  Representatives of $x$ may be chosen
with norm at most $\|x\|+\varepsilon$: truncate a bounded representing
sequence by
\[
 z\longmapsto z\,h_c(z^*z),\qquad
 h_c(t)=\min\{1,c/\sqrt t\},\quad h_c(0)=1,
 \qquad c=\|x\|+\varepsilon.
\]
Continuous functional calculus on uniformly bounded sets preserves
uniform $2$-norm Cauchy sequences, by polynomial approximation and
the two multiplication estimates.  The truncated sequence has norm
at most $c$ and the same limit in $\mathcal N_0$, hence the same
uniform limit by the injectivity already proved.

Products and adjoints are therefore represented by bounded products
and adjoints of representing sequences.  The image is closed in
operator norm: an operator norm Cauchy sequence in the image is
uniformly $2$-norm Cauchy, and a diagonal choice of these truncated
representatives produces its bounded uniform limit.  Thus it is a
$\mathrm{C}^*$-algebra.  The same diagonal argument proves uniform $2$-norm
completeness of every norm ball; truncation with $c$ equal to its
radius shows that the limit remains in that ball.

Schwarz's inequality for $\mathbb{E}$ gives
$\|\mathbb{E}(a)\|_{2,K}\leq\|a\|_{2,K}$ on $A_R$.  Its norm bounded extension
therefore takes values in $D_K$, agrees with $\mathbb{E}_0$, and is a faithful
conditional expectation.  Positivity and traciality of the extended
states follow by bounded uniform approximation.

Finally, norm bounded sets of $D_K$ are separable in the uniform
$2$-norm.  Finite chart polynomials with coefficients from countable
dense subsets of these sets are uniformly $2$-norm dense in $A_R$.
Truncation as above supplies a countable dense subset of its unit
ball, whose bounded uniform completion is the unit ball of $M_R$.
\end{proof}

\subsection{Disjoint graph coordinates and the diagonal}

Enumerate the atlas as $(s_j)_{j\geq1}$ with $s_1=\mathrm{id}$, and put
\[
 r_j=p_{s_j}\prod_{i<j}\bigl(1-c(s_j,s_i)\bigr),
 \qquad w_j=V_{s_j}r_j.
\]
The graphs of the restricted maps represented by $w_j$ are disjoint
and cover $R$.  Notice that the projections $r_j$ need not be orthogonal
on $X$: they select different arrows, not different source points.
Nevertheless,
\begin{equation}
\label{eq:uniform-relation-orthogonal-Fourier-coordinates}
 \mathbb{E}(w_i^*w_j)=\begin{cases}r_j,&i=j,\\0,&i\ne j.\end{cases}
\end{equation}
Every original chart $V_{s_j}$ is a finite sum of restrictions of
$w_1,\ldots,w_j$, because its first matching chart has index at most $j$.

\begin{lem}[Uniform Fourier approximation]
\label{lem:uniform-relation-Fourier-approximation}
For $x\in M_R$, put
\[
 b_j(x)=\mathbb{E}(w_j^*x)\in r_jD_K,
 \qquad Q_nx=\sum_{j\leq n}w_j b_j(x).
\]
Then $\|x-Q_nx\|_{2,K}\to0$, and for each $\mu\in K$,
\begin{equation}
\label{eq:uniform-relation-Fourier-tail}
 \|x-Q_nx\|_{2,\mu}^2
  =\tau_\mu(x^*x)-\sum_{j\leq n}\int_X|b_j(x)|^2\,d\mu.
\end{equation}
The diagonal $D_K$ is a maximal abelian subalgebra of $M_R$.
\end{lem}

\begin{proof}
Equation~\eqref{eq:uniform-relation-orthogonal-Fourier-coordinates}
and diagonal bimodularity of $\mathbb{E}$ give
\[
 \mathbb{E}((x-Q_nx)^*(x-Q_nx))
       =\mathbb{E}(x^*x)-\sum_{j\leq n}|b_j(x)|^2\geq0.
\]
They also show that $Q_n$ is contractive for every designated $2$-norm.
Finite chart polynomials have finite support in these coordinates.
Since they are uniformly $2$-norm dense in $M_R$, this proves the
convergence and the stated identity.  No uniform operator norm bound
on $Q_nx$ is asserted for a general $x$.

For maximal abelianness, suppose $x\in M_R$ commutes with $D_K$.
Its regular kernel $x\delta$ satisfies
\[
 (f(y)-f(z))(x\delta)(y,z)=0\qquad(f\in \mathrm{C}(X)).
\]
A countable family of continuous functions separating the points of
$X$ shows that this kernel is supported on the diagonal.  Since
$\delta$ is separating, $x=\mathbb{E}_0(x)$.  The expectation extension takes
values in $D_K$, so $x\in D_K$.
\end{proof}

\begin{prop}[A separable Cartan core]
\label{prop:uniform-relation-separable-Cartan-core}
There is a separable unital $\mathrm{C}^*$-subalgebra $D_0\subseteq D_K$ containing
the image of $\mathrm{C}(X)$ such that
\[
 A_0=\mathrm{C}^*(D_0,V_s:s\in\mathcal S)
\]
is separable, $D_0\subseteq A_0$ is a Cartan pair, and its bounded
completion for the designated traces is canonically $M_R$.
In particular, the designated trace family is affinely homeomorphic
to $K$.
\end{prop}

\begin{proof}
Start with the image of $\mathrm{C}(X)$ and the countably many support and
coincidence projections of the atlas.  From a separable diagonal
algebra already chosen, adjoin the images of a countable dense subset
under every corner map $f\mapsto\alpha_s(q_sf)$, $s\in\mathcal S$.
Iterate this operation and take the norm closure of the increasing
union.  The resulting $D_0$ is separable, contains all the indicated
projections, and is invariant under every corner pullback, since
these maps are norm continuous.  Its bounded uniform completion is
$D_K$: it contains $\mathrm{C}(X)$ and is contained in $D_K$.

The finite chart calculation shows that $\mathbb{E}(A_0)\subseteq D_0$.
Its restricted expectation is faithful, and the chart operators
normalize $D_0$ and generate $A_0$ together with the diagonal.
If an element of $A_0$ commutes with $D_0$, it commutes with the
continuous functions.  The diagonal-kernel argument in
Lemma~\ref{lem:uniform-relation-Fourier-approximation} shows that it
equals its expectation and hence belongs to $D_0$.  Thus $D_0$
is maximal abelian.  It contains the unit of $A_0$, so these properties
are precisely the Cartan conditions.

The same faithful-completion argument embeds the designated completion
of $A_0$ into $\mathcal N_0$.  It contains $D_K$ as the completion of
$D_0$, and therefore contains $A_R$ and its bounded completion.
The reverse inclusion follows from $A_0\subseteq A_R$.  Finally,
$\mu\mapsto\tau_\mu|_{A_0}$ is affine and continuous and is injective
because $D_0$ contains the continuous functions on $X$.  Compactness
of $K$ makes it a homeomorphism onto the designated trace family;
no assertion that this family exhausts all traces of $A_0$ is made.
\end{proof}

\subsection{The full pseudogroup and normalizers}

For an admissible partial map $\varphi:A\to B$, use its first-chart
disjointization $e_j(\varphi)$ from
Proposition~\ref{prop:regular-relation-pseudogroup-closure}.  Define its
\emph{canonical graph partial isometry} by
\begin{equation}
\label{eq:uniform-relation-canonical-graph-lift}
 v_\varphi=\sum_j V_{s_j}e_j(\varphi).
\end{equation}
Its finite partial sums have norm at most one, because both their
initial projections and their range projections are orthogonal.
The uniformly vanishing tails prove convergence in $M_R$.

\begin{prop}[Untwisted uniform normalizer correspondence]
\label{prop:uniform-full-group-normalizers}
The map $\varphi\mapsto v_\varphi$ identifies the uniform measure full
pseudogroup with the partial-isometry normalizers of $D_K$ in $M_R$,
modulo their diagonal source phases.  More precisely, a partial isometry
$v\in M_R$ satisfying
\[
 p=v^*v\in D_K,\quad q=vv^*\in D_K,\quad
 vD_Kv^*\subseteq D_K,\quad v^*D_Kv\subseteq D_K
\]
has a unique admissible partial map $\varphi$ modulo $K$-null sets,
and a unique $a\in pD_K$ with $a^*a=aa^*=p$, such that
\[
 v=v_\varphi a.
\]
Canonical lifts satisfy
\[
 v_\varphi^*=v_{\varphi^{-1}},\qquad
 v_\varphi v_\psi=v_{\varphi\psi}.
\]
In particular, there is a split exact sequence
\begin{equation}
\label{eq:uniform-relation-split-normalizer-sequence}
 1\longrightarrow\mathcal U(D_K)
 \longrightarrow\mathcal N_{\mathcal U(M_R)}(D_K)
 \longrightarrow[R]^u_{K,\mathcal S}
 \longrightarrow1,
\end{equation}
split by $T\mapsto v_T$.  Its canonical and quotient metrics satisfy
\begin{align}
 \|v_T-v_S\|_{2,K}^2&=2d_K(T,S),\label{eq:uniform-relation-canonical-metric}\\
 \inf_{a\in\mathcal U(D_K)}\|v-av'\|_{2,K}^2
     &=2d_K(T,S),\label{eq:uniform-relation-quotient-metric}
\end{align}
where $v,v'$ induce $T,S$, respectively.  Thus the normalizer
quotient with its uniform $2$-norm quotient distance is isometric to
the full group with distance $(2d_K)^{1/2}$, and the canonical section
is continuous.
\end{prop}

\begin{proof}
The canonical graph operators implement corner pullbacks on $D_K$.
Their initial and range projections are $1_A$ and $1_B$.
Composition and adjoint of the untwisted graph kernels give the two
displayed identities first in the regular representation and hence
in $M_R$.  Thus they are partial normalizers.

Conversely, let $v$ be a partial normalizer with initial projection
$p$ and range projection $q$.  Its corner isomorphism extends normally
to the corresponding corners of $L^\infty(X,\mu_0)$, since $D_K$
contains the continuous functions and is strongly dense in that
diagonal.  The countable-generator spatial realization argument
gives a Borel partial bijection $\varphi$ between the supports of
$p$ and $q$, modulo $\mu_0$-null sets, such that
\[
 v(pf)v^*=q(f\circ\varphi^{-1})
 \qquad(f\in \mathrm{C}(X)).
\]
For example, choose simultaneous representatives of the images of a
countable dense unital star subalgebra of $\mathrm{C}(X)$.  At almost every
point their values form a character of $\mathrm{C}(X)$ and hence determine a
point of $X$; the inverse corner map gives the inverse point map.
Trace preservation puts these point maps in the indicated supports.

The relation kernel $v\delta$ satisfies
\[
 (v\delta)(y,x)f(x)
  =f(\varphi^{-1}y)(v\delta)(y,x).
\]
The same countable point-separating family therefore forces its
support to lie on $y=\varphi x$.  Also,
$\mathbb{E}_0(v^*v)=p$ says that the squared column norm is one on that
support.  Consequently the point map has graph in $R$ for
$\mu_0$-almost every source point, and its graph coefficient has
modulus one.  For the disjoint graph coordinates above,
\[
 e_j=|b_j(v)|^2\in D_K
\]
are therefore orthogonal projections, dominated by $p$, for
$\mu_0$-almost every point.  Faithfulness of the representation of
$D_K$ upgrades these projection identities and orthogonality
identities to identities in $D_K$.

It is important to establish coverage for all measures, rather than
only for $\mu_0$.  Formula~\eqref{eq:uniform-relation-Fourier-tail}
and uniform Fourier approximation give
\begin{equation}
\label{eq:uniform-relation-normalizer-universal-tails}
 \sup_{\mu\in K}\mu\!\left(p-\sum_{j\leq n}e_j\right)
      \longrightarrow0.
\end{equation}
Define $\varphi=s_j$ on the support of $e_j$.  Its range pieces are
$\alpha_{s_j^{-1}}(e_j)\in D_K$.  They are orthogonal and dominated
by $q$, as may again be checked in the faithful $\mu_0$
representation.  Their traces equal the traces of $e_j$, while
$\tau_\mu(p)=\tau_\mu(q)$ for every $\mu\in K$.  Thus they exhaust
$q$ with uniformly vanishing tails as well.  Remove the countable
exceptional sets and their $R$-saturations.  These source and range
pieces give a Borel partial bijection on one $R$-invariant
$K$-conull set, with initial and final projections $p,q$.

For every original chart $s$, its coincidence coefficient is
\[
 c_s(\varphi)=\sum_j e_j c(s_j,s)\in D_K,
\]
where convergence follows from
\eqref{eq:uniform-relation-normalizer-universal-tails}.
Hence $\varphi$ is admissible.  Its canonical lift has the same
spatial corner map as $v$.  The element $a=v_\varphi^*v$ consequently
commutes with $pD_K$ and is a unitary in the corner $pM_Rp$.
Maximal abelianness gives $a\in pD_K$, proving the factorization.
The same argument and faithfulness give uniqueness.

For unitaries this gives the exact sequence and its homomorphic
section.  The kernel is precisely the diagonal unitary group.
Since graph lifts have phase one,
\[
 \mathbb{E}(v_T^*v_S)=1_{\{T=S\}}.
\]
This proves~\eqref{eq:uniform-relation-canonical-metric}.  For any
two normalizers inducing $T,S$, the real part of the mixed diagonal
coefficient is at most $1_{\{T=S\}}$, pointwise, regardless of the
diagonal phases.  Therefore the quotient distance is at least the
right side of~\eqref{eq:uniform-relation-quotient-metric}.
Writing $v=v_Tb$ and $v'=v_Sc$, the source phase $c^*b$ makes
$v'(c^*b)=v_Sb$.  It achieves equality because right multiplication
by $b$ preserves the uniform $2$-norm.  Right and left diagonal
phases give the same infimum, since $v'$ normalizes $D_K$.
The metric formulas prove continuity of the quotient and of its
canonical section.
\end{proof}

\subsection{Isomorphisms within the untwisted relation models}
\label{subsec:uniform-relation-pair-isomorphisms}

The preceding construction also retains the spatial isomorphism class of
the regular relation. We first prove the isomorphism theorem within
these untwisted models. Section~\ref{sec:intrinsic-reconstruction} extends
it to twisted models and supplies the inverse construction under its
intrinsic hypotheses.

\begin{lem}[Independence of the dense barycenter]
\label{lem:uniform-relation-base-measure-independence}
For a finite regular graph kernel $a$, its norm in $A_R$ is
\[
 \|a\|=\sup_{\mu\in K}\|\lambda_\mu(a)\|,
\]
where $\lambda_\mu$ is its regular representation on
$L^2(R,\nu_\mu)$.  Consequently $A_R$ and its designated bounded
completion $M_R$ do not depend on the dense sequence used to choose
$\mu_0$.
\end{lem}

\begin{proof}
For a designated trace $\tau_\mu$ from
Proposition~\ref{prop:uniform-relation-designated-traces}, the map
$a\mapsto\lambda_\mu(a)\delta_\mu$ has inner product
$\tau_\mu(b^*a)$.  Its range is dense in $L^2(R,\nu_\mu)$.
Indeed, the disjoint graph coordinates cover $R$, and continuous
coefficient functions are dense for the $2$-norm of the single measure
on each graph.  Thus this regular representation is the GNS
representation of $\tau_\mu$ and $\|\lambda_\mu(a)\|\leq\|a\|$
for every $\mu\in K$. Equality holds for $\mu=\mu_0\in K$ by the
definition of $A_R$.
The supremum is intrinsic to the designated family, and its
uniform $2$-norm completion is intrinsic as well.
\end{proof}

\begin{thm}[Spatial and completed-pair isomorphisms]
\label{thm:uniform-relation-pair-isomorphism}
Let
\[
 (X,R,K,D_K,\mathcal S),\qquad
 (Y,Q,L,D_L,\mathcal T)
\]
be two regular uniform measure countable equivalence relations,
with constructed untwisted completed pairs $(D_K,M_R)$ and
$(D_L,M_Q)$.  Fix an affine homeomorphism $H:K\to L$.
The following are equivalent:
\begin{enumerate}
\item There is a uniform measure relation isomorphism $\theta$
implementing $H$.
\item There is an isomorphism $\Phi:M_R\to M_Q$ with
$\Phi(D_K)=D_L$ and
\[
 \tau_{H(\mu)}(\Phi(x))=\tau_\mu(x)
       \qquad(x\in M_R,\ \mu\in K).
\]
\end{enumerate}
Given $\theta$ in (1), there is a canonical such isomorphism satisfying
\[
 \Phi_\theta(f)=f\circ\theta^{-1},\qquad
 \Phi_\theta(v_\varphi)=v_{\theta\varphi\theta^{-1}}
 \quad
 \left(\varphi\in
       \mathopen{[\![}R\mathclose{]\!]}^u_{K,\mathcal S}\right).
\]
For each isomorphism $\Phi$ in (2), its restriction to the diagonal
determines a spatial map $\theta$ up to universally almost-everywhere
equality. The image of each canonical lift agrees with its canonical
transport by this $\theta$ up to a diagonal phase on its support.
\end{thm}

\begin{proof}
Suppose first that $\theta$ is a uniform measure relation
isomorphism.  By
Lemma~\ref{lem:uniform-relation-base-measure-independence}, choose
compatible dense barycenters $\mu_0$ and
$\eta_0=H(\mu_0)=\theta_*\mu_0$.  The invariant conull sets on
which $\theta$ is defined give a counting-space unitary
\[
 U_\theta:L^2(R,\nu_{\mu_0})\longrightarrow
                   L^2(Q,\nu_{\eta_0}),
 \qquad
 (U_\theta\xi)(\theta y,\theta x)=\xi(y,x).
\]
It is unitary because $\theta$ bijects the equivalence classes and
transports the base measure.  It carries multiplication by $f$ to
multiplication by $f\circ\theta^{-1}$ and every graph operator
to the canonical graph operator of its conjugated partial map.
Those conjugated maps are admissible by
Proposition~\ref{prop:uniform-relation-isomorphism-transport},
so their canonical graph operators belong to $M_Q$.
Conjugation by $U_\theta$ consequently maps $A_R$ into $M_Q$.
It preserves all designated trace values, as follows by restricting
kernels to the diagonal and using measure transport.
Hence it preserves the uniform $2$-norm and extends by bounded
completion to a homomorphism $M_R\to M_Q$.
The same argument for $\theta^{-1}$ gives its inverse.  Conjugation
also preserves the diagonal expectations and the canonical graph
operators of all admissible partial maps, proving the asserted formula.

Conversely, let $\Phi$ satisfy (2).  Its diagonal restriction and
the spatial realization theorem
(Theorem~\ref{thm:orbit-spatial-diagonal}) give a Borel bijection
\[
 \theta:X_d\longrightarrow Y_d
\]
between sets conull for the respective designated families such that
$\theta_*\mu=H(\mu)$ and
$\Phi(f)=f\circ\theta^{-1}$ for $f\in D_K$.
At this point the sets need not be invariant for their relations.

For every source chart $s\in\mathcal S$,
$\Phi(V_s)$ is a partial normalizer of $D_L$.  By
Proposition~\ref{prop:uniform-full-group-normalizers}, it has an
admissible target partial map $\sigma_s$, with source and range
projections $\Phi(p_s)$ and $\Phi(q_s)$.  Its diagonal phase does
not affect its corner pullback.  Transporting the covariance identity
for $V_s$ therefore gives
\[
 \alpha_{\sigma_s}(\Phi(q_sf))
       =\Phi(\alpha_s(q_sf))
       \qquad(f\in D_K).
\]
Apply this identity to a countable family of continuous functions
separating the points of $X$, and use simultaneous Borel
representatives.  For almost every $x\in A_s$, for every measure
in $K$, it implies
\begin{equation}
\label{eq:uniform-relation-pair-chart-covariance}
 \sigma_s(\theta x)=\theta(sx).
\end{equation}
In making this assertion, omit the exceptional points where $x$ or
$sx$ is outside $X_d$, where the corresponding source or range
projections do not match their transported representatives, or
where $\sigma_s(\theta x)$ is outside $Y_d$.  These exceptions are
null for all measures in $K$: chart maps preserve the measures on
their domains, $\theta$ transports them, and $\sigma_s$ preserves
the target measures on its domain.  The displayed identities for
countably many functions add only a countable union of null sets.

Apply the same argument to $\Phi^{-1}$ and every target chart
$t\in\mathcal T$.  It supplies an admissible source partial map
$\rho_t$ and the reverse identities
\[
 \rho_t(\theta^{-1}y)=\theta^{-1}(ty)
\]
outside one universally null Borel set for each chart.
Because both atlases are countable, there are conull Borel subsets
$P\subseteq X_d$ and $Q_0=\theta(P)\subseteq Y_d$ such that
$\theta$ carries $R|_P$ onto $Q|_{Q_0}$: for a source-related pair
choose a chart whose graph contains it and apply
\eqref{eq:uniform-relation-pair-chart-covariance}; use the reverse
identities for a target-related pair.

We now retain invariant sets on both sides.  Put
\[
 P^{\mathrm{inv}}=X\setminus R(X\setminus P),\qquad
 Q^{\mathrm{inv}}=Y\setminus Q(Y\setminus Q_0),
\]
where $R(C)$ and $Q(C)$ denote relation saturations.  They are Borel,
invariant, and conull for every designated measure, by countability
and invariance of the relations.  Define
\[
 X_0=P^{\mathrm{inv}}\cap\theta^{-1}(Q^{\mathrm{inv}}),
 \qquad Y_0=\theta(X_0).
\]
These sets are conull and Borel.  They are invariant on both sides:
if $x\in X_0$ and $x'Rx$, then $x'\in P^{\mathrm{inv}}\subseteq P$,
and restricted relation transport puts $\theta x'$ in
$Q^{\mathrm{inv}}$.  Thus $x'\in X_0$.  Conversely, if
$y\in Y_0$ and $y'Qy$, then $y'\in Q^{\mathrm{inv}}\subseteq Q_0$;
the reverse transport puts $\theta^{-1}y'$ in
$P^{\mathrm{inv}}$, hence $y'\in Y_0$.

On these invariant sets, every conjugated source chart is the
restriction of $\sigma_s$, and every inversely conjugated target
chart is the restriction of $\rho_t$.  Restriction changes their
domains only by universally null sets.  They are therefore
admissible.  Together with the diagonal and measure transport already
proved, these are precisely the conditions for a uniform measure
relation isomorphism implementing $H$.
\end{proof}

\begin{rem}[The residual phases of a pair isomorphism]
\label{rem:uniform-relation-pair-residual-phases}
For the spatial map induced by a marked pair isomorphism $\Phi$,
the normalizer correspondence gives source-corner unitaries $a_s$
such that
\[
 \Phi(V_s)=v_{\theta s\theta^{-1}}a_s.
\]
Multiplicativity gives their cocycle coherence.  If
$\sigma_s=\theta s\theta^{-1}$ and the composite is defined at $y$,
then
\[
 a_{st}(y)=a_s(\sigma_t y)a_t(y).
\]
The phases agree on common graph restrictions because the
corresponding restricted source graph operators agree.  Thus the
additional phase information is retained rather than discarded.
The canonical forward isomorphism has all these phases equal to
one on their source corners.  The theorem identifies existence of
isomorphisms within the untwisted regular relation models; it does
not identify every algebra isomorphism with this canonical choice.
These residual $1$-cocycle phases are distinct from the $2$-cocycle
twist removed in Theorem~\ref{thm:uniform-regular-untwisting}.
\end{rem}

\subsection{Specialization to free actions}

\begin{lem}[The crossed-product specialization]
\label{lem:uniform-relation-algebra-action-case}
Let a countable group $\Gamma$ act freely and continuously on a compact
metrizable space $X$, and let $K$ be a nonempty compact convex family of
invariant probability measures.  Equip its orbit relation with the
atlas consisting of the global group maps and the empty map. Then the designated bounded uniform
completion of $\mathrm{C}(X)\rtimes_r\Gamma$ is canonically isomorphic to
$M_R$, carrying its completed diagonal onto $D_K$ and its canonical
expectation onto $\mathbb{E}$.
\end{lem}

\begin{proof}
Freeness identifies the relation coordinates $(gx,x)$ uniquely with
$(g,x)$.  Under this identification, the graph operators and the
continuous diagonal form the regular covariant representation of
$\mathrm{C}(X)\rtimes_r\Gamma$ associated with $\mu_0$.  Its expectation is
the identity Fourier coefficient, because every nonidentity group
map has empty coincidence set with the identity.  The designated
traces and their uniform $2$-norms therefore agree with those of the
crossed product.  If necessary, take the common tracial kernel
quotient first; its regular representation is faithful.  Indeed,
faithfulness of the canonical expectation and continuity of the
designated trace values show that the $\mu_0$-trace kernel is
exactly the common kernel.

Write $M$ for this crossed-product completion.  Its diagonal
completion is $D_K$, so it contains $D_K$ and every graph operator
$V_g$.  As a $\mathrm{C}^*$-algebra it therefore contains $A_R$, and bounded
uniform completeness gives $M_R\subseteq M$.  Conversely, the
crossed product is contained in $A_R$, so its bounded completion
is contained in $M_R$.  These inclusions identify the completions
and their expectations.  The normalizer correspondence above thus
specializes to the action full group without an amenability assumption.
\end{proof}

\begin{rem}[Scope of the forward construction]
\label{rem:uniform-relation-forward-construction-scope}
When $K=\{\mu\}$, one has $D_K=L^\infty(X,\mu)$, and
Kaplansky density shows that $M_R=\mathcal N_0$: the bounded
strong closure of the graph algebra is its bounded $2$-norm closure
for the single faithful trace.  Thus the construction specializes
to the usual untwisted von Neumann algebra of a measured relation.

This construction supplies an untwisted algebra and a maximal
abelian diagonal for every regular atlas in the stated sense.  It
also identifies the whole admissible full pseudogroup and full group
with its partial and unitary normalizers. The inverse theorem in
Section~\ref{sec:intrinsic-reconstruction} requires uniform separability,
partial-isometry normalizer regularity, and maximal abelianness in every
designated GNS closure. For a given regular relation
with the finite-chart hypotheses,
Theorem~\ref{thm:uniform-regular-untwisting} removes every regular
scalar twist uniformly over the designated trace family.
\end{rem}

\section{Reconstruction from finite matrix charts}
\label{subsec:uniform-AF-reconstruction}

We use the untwisted relation algebra $M_R$ constructed in
Section~\ref{sec:uniform-relation-algebra}.  In this section we also
write $\mathcal M_R=M_R$ and $\mathcal D=D_K$.
By a \emph{completed Cartan pair} we mean a maximal abelian completed
diagonal with a faithful designated-trace-preserving conditional
expectation, whose partial-isometry normalizers generate an algebra
that is dense in the bounded uniform $2$-norm topology.
This regularity convention differs from the norm regularity of the
separable $\mathrm{C}^*$-Cartan core in
Proposition~\ref{prop:uniform-relation-separable-Cartan-core}.

\begin{thm}[Forward AF reconstruction]
\label{thm:uniform-AF-reconstruction}
Suppose that $(X,R,K,\mathcal D,\mathscr P)$ is uniformly AF
hyperfinite.  Then there are a unital AF algebra $B$, its standard
diagonal $E=\mathrm{C}(Y)$, and a compact convex family $L\subseteq \T(B)$ with
the following properties:
\begin{enumerate}
\item restriction of traces gives an affine homeomorphism
$H:K\to L$;
\item there is an isomorphism of completed untwisted Cartan pairs
\[
 (\mathcal D,\mathcal M_R)\cong(E^L,B^L);
\]
\item $R$ is uniform measure isomorphic to the AF tail relation on
$Y$, with its standard prefix-substitution atlas and the central
measure family corresponding to $L$.
\end{enumerate}
The spatial isomorphism in (3) is a Borel bijection between invariant
sets conull for every measure in the designated families.  Neither
simplicity of $B$ nor the equality $L=\T(B)$ is asserted.
\end{thm}

\begin{proof}
For the level $n$ finite chart system, use its notation
$Y^{(n)}_{a,i}$ and $s^{(n)}_{a,i}$.  The graph operators
\[
 e^{(n,a)}_{ij}
    =v_{s^{(n)}_{a,i}(s^{(n)}_{a,j})^{-1}}
\]
satisfy the matrix-unit identities, and
\[
 e^{(n,a)}_{ii}=1_{Y^{(n)}_{a,i}},
 \qquad \sum_{a,i}e^{(n,a)}_{ii}=1.
\]
There is no phase factor because the graph realization is untwisted.
Discard any block whose source projection is $K$-null.  This does not
change the completed relation or diagonal.  All remaining matrix
blocks are nonzero, and their span is a unital finite-dimensional
algebra
\[
 B_n\cong\bigoplus_a \mathrm{M}_{k_{n,a}}.
\]
Coordinate refinement says precisely that each old matrix unit is a
finite sum of matching new matrix units.  Thus $B_n\subseteq B_{n+1}$
are the standard phase-free embeddings.  Put
\[
 B=\overline{\bigcup_nB_n}^{\|\cdot\|},
 \qquad E=\overline{\bigcup_nE_n}^{\|\cdot\|}.
\]
This is a standard unital AF pair.  Its diagonal spectrum $Y$ is a
compact metrizable zero-dimensional space, realized as the path space
of the coordinate Bratteli diagram.  It need not be a Cantor space.

For $\mu\in K$, set $H(\mu)=\tau_\mu|_B$.  The trace values on the
finite stages are finite linear combinations of
$\mu(Y^{(n)}_{a,i})$, so $H$ is continuous and affine.  It is
injective.  Indeed, if $H(\mu)=H(\mu')$, the two measures agree on
$E_n$ for every $n$.  Given $f\in \mathrm{C}(X)$, choose uniformly bounded
$f_n\in E_n$ converging to $f$ in uniform $2$-norm.  Then
\[
 \left|\int f\,d\mu-\int f\,d\mu'\right|
 \leq 2\|f-f_n\|_{2,K}\longrightarrow0.
\]
Thus $\mu=\mu'$.  Let $L=H(K)$; compactness gives the asserted
affine homeomorphism.  In what follows we use the same letter for a
trace in $L$ and its central measure on $Y$.

We next prove bounded uniform $2$-norm density of $B$ in
$\mathcal M_R$.  Diagonal density is part of the hypothesis.  Fix
$\varphi:A\to B'$ in $\mathscr P$.  Let $\varphi_n$ be its restriction
to those $x\in A$ with $(x,\varphi x)\in R_n$.  This restriction is
admissible: its domain is the finite union of the equalizer sets of
$\varphi$ with the level $n$ matrix-chart maps.  Its graph operator
belongs to
\[
 C_n=\mathrm{C}^*(B_n,\mathcal D)
     \cong\bigoplus_a
             \mathrm{M}_{k_{n,a}}(1_{Y^{(n)}_{a,1}}\mathcal D).
\]
More explicitly, it is the sum over $a,i,j$ of
$e^{(n,a)}_{ij}$ times the source projection on which $\varphi$
agrees with that chart.  It is a contraction.  By
Lemma~\ref{lem:uniform-partial-map-extension} and uniform
hyperfiniteness,
\[
 \|v_\varphi-v_{\varphi_n}\|_{2,K}^{\,2}
   =\sup_{\mu\in K}
       \mu\{x\in A:(x,\varphi x)\notin R_n\}
   \longrightarrow0.
\]

Every element of $C_n$ is a finite matrix with entries in diagonal
corners.  Approximate those entries by the finite diagonals $E_m$,
with $m\geq n$, to approximate it in uniform $2$-norm by elements of
$B_m$.  A uniform norm bound can be retained.  For completeness,
if $z\in C_n$ is a contraction and $b\in B_m$ approximates it, replace
$b$ by its radial truncation $b^{[2]}$ at norm two.  With
$p=1_{(2,\infty)}(|b|)$, for every trace in the designated family,
\[
 \|bp\|_2\geq2\|p\|_2,
 \qquad
 \|bp\|_2\leq\|(b-z)p\|_2+\|p\|_2.
\]
Hence $\|b-b^{[2]}\|_2\leq2\|b-z\|_2$, and
\[
 \|b^{[2]}\|\leq2,
 \qquad \|b^{[2]}-z\|_{2,K}\leq3\|b-z\|_{2,K}.
\]
Applying this to $z=v_{\varphi_n}$ gives uniformly norm-bounded
finite-stage approximations to $v_\varphi$.  Thus the complete
subalgebra $B^L\subseteq\mathcal M_R$ contains the diagonal and all
the original atlas graph operators.  It contains $A_R$, and bounded
uniform $2$-norm closure gives $B^L=\mathcal M_R$.  Similarly,
$E^L=\mathcal D$.

The usual finite diagonal-phase averages of a standard AF pair
converge to its diagonal expectation in uniform $2$-norm: on any
fixed finite stage, a sufficiently large diagonal average already
kills all its off-diagonal matrix units.  Bounded density extends
this to $B^L$.  An element commuting with $E^L$ is fixed by these
averages and therefore belongs to $E^L$.  This proves the MASA
property.  The expectation is faithful as above, and finite
coordinate permutation unitaries, together with the diagonal,
generate the finite stages.  Their normalizers therefore have linear
span dense in the uniform $2$-norm, which is the regularity meant for
the completed Cartan pair.  This proves (2).

It remains to retain the actual spatial relation.  The nested finite
partitions define a Borel coding map $\theta:X\to Y$, outside the
common null set already removed for the chart identities: evaluation
on all the $E_n$ gives the corresponding character of $E=\mathrm{C}(Y)$.
Its pushforward satisfies $\theta_*\mu=H(\mu)$, since this equality
holds on every finite diagonal.  Coordinate refinement also implies
that each level $n$ matrix-chart map is carried to the corresponding
prefix substitution in the path diagram.

The coding map is injective on a common conull Borel set.  To see
this, choose a countable family $f_j\in \mathrm{C}(X)$ separating points.
For each $j$, choose finite-diagonal approximations with summably
small squared uniform $2$-errors.  Markov's inequality and
Borel--Cantelli give a fixed Borel exceptional set, null for every
$\mu\in K$, outside which those approximations converge pointwise
to $f_j$, simultaneously for every $j$.  Their values are determined
by $\theta(x)$.  Equal codes consequently have equal values of all
$f_j$ and are equal.  Remove the $R$-saturation of the exceptional
set and denote the remaining invariant set by $X_0$.

The image $Y_0=\theta(X_0)$ is Borel, and the injective Borel map
$\theta:X_0\to Y_0$ has a Borel inverse.  It is conull for all central
measures in $L$.  Each finite prefix permutation is the image of a
permutation of a level $n$ finite chart system, so $Y_0$ is invariant
under tail equivalence.  The finite chart relations $R_n$ correspond
exactly to equality of tails beyond level $n$ on these sets.  Since
$R=\bigcup_nR_n$ there, $\theta$ carries $R|_{X_0}$ onto the AF tail
relation restricted to $Y_0$.

Pullback by $\theta$ identifies $E^L$ with $\mathcal D$, as already
proved by finite-diagonal density.  The inverse images of all
prefix-substitution charts are admissible.  Conversely, the graph
coefficients of a transported original atlas chart, relative to a
prefix-substitution chart, are the transported equalizer projections
of two original admissible partial maps.  They belong to $E^L$.
The spatial transport proposition therefore gives the uniform
measure relation isomorphism in (3).
\end{proof}

\begin{rem}\label{rem:uniform-AF-reconstruction-scope}
The AF algebra in Theorem~\ref{thm:uniform-AF-reconstruction} is the
algebra of the given coordinate diagram.  Its simplicity requires
the usual simplicity condition on that diagram; it does not follow
from uniform hyperfiniteness.  Likewise, $L=H(K)$ is the designated
compact trace family and may be a proper subset of $\T(B)$.  Atomic
central measures and finite tail classes are allowed in this
construction. With aperiodicity and the full compact invariant
family, Lemma~\ref{lem:full-family-AF-face} makes $L$ a diffuse closed
trace face, and Corollary~\ref{cor:full-family-simple-AF-model} replaces
this model by a simple infinite-dimensional full-trace model. This
does not assert simplicity or full trace control for the original
coordinate diagram.

The finite-permutation normalizer-density argument also applies to
a standard AF pair with a designated compact trace family: its
finite-stage conditional expectations, matrix-entry decompositions,
and approximation of diagonal projection partitions do not require
simplicity or faithfulness of each individual trace.  One works in
the quotient by the common trace-null ideal if necessary.  This
observation should not be confused with an assertion that every
uniformly hyperfinite relation has a compatible AF presentation.
The forward theorem uses all the matrix-chart and diagonal-density
conditions in Definition~\ref{defn:uniform-AF-hyperfinite-relation}.
\end{rem}

\begin{cor}[A spatial characterization of the AF condition]
\label{cor:uniform-AF-spatial-characterization}
A uniform measure relation is uniformly AF hyperfinite if and only if
it admits a uniform measure isomorphism to the tail relation of a
standard unital AF pair, equipped with a compact designated
central-measure family.
The family need not comprise all central measures.
\end{cor}

\begin{proof}
One implication is Theorem~\ref{thm:uniform-AF-reconstruction}.
For the converse, let $B$ be a standard unital AF algebra with diagonal
$\mathrm{C}(Y)$ and a nonempty compact convex family $L\subseteq \T(B)$.
Identify $L$ with its central probability measures $K_Y$ on $Y$.
The finite-prefix substitutions, with the identity and empty map,
form a regular countable atlas: their supports and coincidence sets
are clopen, their corner maps preserve continuous functions, and
bounded uniform $2$-norm closure extends those maps to $\mathrm{C}(Y)^{K_Y}$.
Compositions and inverses are again prefix substitutions or their
finite restrictions, which may be included in the atlas.

The length-$n$ prefix cylinders, grouped by their terminal vertex,
and the maps replacing one prefix by another give finite regular
matrix charts.  Their coordinate refinements are exactly those of
the Bratteli diagram.  The finite diagonals are norm dense in $\mathrm{C}(Y)$,
and hence bounded uniformly $2$-norm dense in $\mathrm{C}(Y)^{K_Y}$.
The level relations increase to tail equivalence pointwise.

For an admissible partial map $\varphi$, take its first-chart
partition $(e_j)$ with uniformly vanishing tails.  Given $\varepsilon>0$,
retain finitely many charts so that the omitted source mass is at most
$\varepsilon$ for every designated measure.  Those finitely many
prefix maps all lie in the level relation $R_N$ for some $N$.
Consequently, for $n\geq N$,
\[
 \sup_{\nu\in K_Y}
 \nu\{y\in\operatorname{dom}\varphi:(y,\varphi y)\notin R_n\}
 \leq\varepsilon.
\]
This proves uniform escape and all the AF hyperfiniteness conditions.
A uniform measure isomorphism transports the finite charts,
their refinements, bounded diagonal density, and these estimates.
Common trace-null cylinders may be discarded throughout; atomic
central measures and finite tail classes are allowed.
\end{proof}

\section{Finite permutations in completed AF pairs}

The standard diagonal expectation of a unital AF algebra extends to
any designated compact trace family.  Its completed diagonal is a MASA.
Indeed, averaging finite diagonal phase groups gives the expectation
exactly on a finite stage.  For a bounded completed element $x$, first
approximate it in uniform $2$-norm in a finite stage and apply that finite
average.  If $x$ commutes with the diagonal, the average fixes $x$, while
its value can be made arbitrarily close to the diagonal expectation of
$x$.  Thus $x$ belongs to the completed diagonal.  The expectation is
faithful, since it preserves all designated traces.  This argument also
shows that the completed pair is uniformly regular.

\begin{lem}\label{lem:AF-normalizer-density}
Let $(E,B)$ be the standard diagonal pair of a unital AF algebra,
with coordinate Bratteli inclusions $B_n\subseteq B_{n+1}$.
Let $L\subseteq\T(B)$ be nonempty and compact.
Put $\mathcal D=E^L$ and $\mathcal M=B^L$.
Let $P_n$ be the finite group of coordinate permutation unitaries
in $B_n$. For every
$v\in\mathcal N_{\mathcal U(\mathcal M)}(\mathcal D)$ and every
$\eta>0$, there are $m$, $p\in P_m$, and
$\lambda\in\mathcal U(\mathcal D)$ such that
\[
 \norm{v-\lambda p}_{2,L}<\eta.
\]
The groups $P_n$ are increasing, and their union is dense in the
completed normalizer quotient modulo diagonal unitaries.
\end{lem}

\begin{proof}
Write $B_n=\bigoplus_{a=1}^s\mathrm{M}_{k_a}$, with coordinate matrix units
$e^a_{ij}$ and central projections $z_a=\sum_i e^a_{ii}$.
In the completion we use the image of each stage, omitting blocks
annihilated by every designated trace. A coordinate permutation of
the surviving blocks lifts to the original stage by taking the
identity permutation on every omitted block.
Put $\mathcal D_a=e^a_{11}\mathcal D$ and
\[
 C_n=\mathrm{C}^*(B_n,\mathcal D)
 \cong\bigoplus_{a=1}^s\mathrm{M}_{k_a}(\mathcal D_a).
\]
The diagonal in this matrix description is $\mathcal D$.
If $\mathbb{E}:\mathcal M\to\mathcal D$ is the extended diagonal expectation,
define
\[
 \mathbb{E}_n(x)=
 \sum_{a=1}^s\sum_{i,j=1}^{k_a}
 e^a_{i1}\,\mathbb{E}(e^a_{1i}xe^a_{j1})\,e^a_{1j}.
\]
To see that this is a conditional expectation onto $C_n$, first
pinch $x$ into the corners $z_a\mathcal Mz_a$, identify each corner
with $\mathrm{M}_{k_a}(e^a_{11}\mathcal Me^a_{11})$, and apply the corner
expectation $\mathbb{E}$ entrywise. This is unital and completely positive,
fixes $C_n$, and is $C_n$-bimodular.
For every designated trace $\tau$, cyclicity shows
\[
 \tau(\mathbb{E}_n(x))
 =\sum_{a,i}\tau(\mathbb{E}(e^a_{1i}xe^a_{i1}))
 =\sum_{a,i}\tau(xe^a_{ii})=\tau(x).
\]
Thus $\mathbb{E}_n$ is the $L^2(\tau)$-orthogonal projection onto
the closure of $C_n$, and is uniformly $2$-contractive.
Since it fixes $B_n$ and $\bigcup_n B_n$ is uniformly $2$-dense,
\[
 \norm{v-\mathbb{E}_n(v)}_{2,L}\longrightarrow0.
\]

Put $b=\mathbb{E}_n(v)$ and $\beta=\operatorname{Ad}(v)|_{\mathcal D}$.
The $\mathcal D$-bimodule property gives $bd=\beta(d)b$.
Consequently $r=v^*b$ commutes with $\mathcal D$.
The completed diagonal is maximal abelian, so $r\in\mathcal D$.
Idempotence now gives
\[
 b=\mathbb{E}_n(b)=\mathbb{E}_n(vr)=\mathbb{E}_n(v)r=vr^2,
\]
and hence $r^2=r$. An idempotent in a commutative $\mathrm{C}^*$-algebra
is a projection. Thus
\[
 b=vr,\qquad
 b^*b=r,\qquad bb^*=\beta(r),\qquad
 \norm{v-b}_{2,\tau}^2=\tau(1-r).
\]
In particular $b$ is an exact partial normalizer of $\mathcal D$
inside $C_n$.

We next complete $b$ to a monomial unitary in $C_n$.
In a block $\mathrm{M}_{k_a}(\mathcal D_a)$ write its entries as $b_{ij}$.
The diagonal conditions for $b e_{jj}b^*$ and $b^*e_{ii}b$
give, respectively,
\[
 b_{ij}b_{\ell j}^*=0\quad(i\ne\ell),\qquad
 b_{ij}^*b_{i\ell}=0\quad(j\ne\ell).
\]
The positives $p_{ij}=|b_{ij}|^2$ therefore have orthogonal supports
in every row and column. Their column sums are the diagonal entries
of the projection $r$, so each $p_{ij}$ is a projection.
The finite Boolean algebra generated by these projections partitions
the corner unit $e^a_{11}$ into projections on which the nonzero
entries have one fixed partial-permutation pattern.
Extend every such pattern to a permutation by matching the unused
columns and rows in increasing order. Give each added entry phase
one, and retain the existing entries of $b$.
This produces a monomial unitary $w\in C_n$ with $wr=b$.
Put $c=w-b=w(1-r)$.
Then $c\in C_n$, $c^*c=1-r$, and
\[
 \tau((v-b)^*c)
 =\tau\bigl(\mathbb{E}_n(v-b)^*c\bigr)=0.
\]
It follows that
\begin{equation}\label{eq:AF-partial-normalizer-rounding}
 \norm{v-w}_{2,\tau}^2=2\tau(1-r).
\end{equation}
In each block its row phase is explicitly
\[
 \lambda_i=\sum_j b_{ij}+\left(1-\sum_j p_{ij}\right)
 \in\mathcal U(\mathcal D_a).
\]
Orthogonality of the row supports proves that these are unitaries.
Removing them writes $w=\lambda\sigma$, where
$\lambda\in\mathcal U(\mathcal D)$ and $\sigma\in C_n$ is a
permutation matrix over a finite projection partition in each
$\mathcal D_a$.

Here are explicit bounds for approximating $\sigma$ by a finite-stage
coordinate permutation. Let
$(q_{a,\ell})_{\ell=1}^{L_a}$ be the tail projection partition
selecting its permutations in block $a$.
At a common stage $m\geq n$, choose positive diagonal contractions
$f_{a,\ell}\in e^a_{11}E_m e^a_{11}$ approximating
$q_{a,\ell}$ uniformly in $2$-norm.
For $\ell<L_a$ put
\[
 p_{a,\ell}=\mathbf1_{[1/2,1]}(f_{a,\ell}).
\]
These are finite-stage diagonal projections and satisfy, pointwise,
\[
 |p_{a,\ell}-q_{a,\ell}|^2
 \leq4|f_{a,\ell}-q_{a,\ell}|^2.
\]
Disjointize the $p_{a,\ell}$ in order and put the remaining
projection into the last part. The resulting partition
$(\widehat q_{a,\ell})$ agrees with $(q_{a,\ell})$ outside
\[
 t_a=\bigvee_{\ell<L_a}|p_{a,\ell}-q_{a,\ell}|,
 \qquad
 \tau(t_a)\leq
 4\sum_{\ell<L_a}\norm{f_{a,\ell}-q_{a,\ell}}_{2,\tau}^2.
\]
Replace the selectors of $\sigma$ by these finite-stage selectors,
obtaining $p\in P_m$. This is a coordinate permutation because
the coordinate Bratteli maps copy each old matrix unit along identical
tail coordinates. The two unitaries agree off the central projection
$\widetilde t_a=\sum_i e^a_{i1}t_ae^a_{1i}$, and
$\tau(\widetilde t_a)=k_a\tau(t_a)$. Hence
\begin{equation}\label{eq:AF-selector-rounding}
 \norm{\sigma-p}_{2,L}^2
 \leq16\sum_{a=1}^s k_a\sum_{\ell<L_a}
 \norm{f_{a,\ell}-q_{a,\ell}}_{2,L}^2.
\end{equation}
Choose $n$ so late that
$\norm{v-w}_{2,L}<\eta/2$, using
\eqref{eq:AF-partial-normalizer-rounding}.
There are only finitely many selectors, so choose their approximants
to make \eqref{eq:AF-selector-rounding} less than $\eta^2/4$.
Then
\[
 \norm{v-\lambda p}_{2,L}
 \leq\norm{v-w}_{2,L}+\norm{\sigma-p}_{2,L}<\eta.
\]
Finally coordinate Bratteli inclusions send permutation matrices to
permutation matrices, so $P_n\subseteq P_{n+1}$.
Every finitely generated subgroup of $\bigcup_nP_n$ is contained in
one finite $P_n$, proving local finiteness as well as density in
the normalizer quotient.
\end{proof}

\section{Uniform Cartan classification and CFW models}
\label{subsec:uniform-CFW}

We first classify standard AF Cartan completions over diffuse closed
trace faces. Combined with uniform AF reconstruction, this proves the
classification theorem for aperiodic SBP relations carrying their full
compact invariant-measure families. We then give the single-generator
and minimal dynamical realizations. All conull sets below are conull
for every measure in the indicated family.

\subsection{AF completions over diffuse closed trace faces}
\label{subsec:diffuse-face-AF-classification}

In this subsection $(E_C,C)$ and $(E_B,B)$ are standard diagonal pairs
of separable unital AF algebras. Their coordinate presentations are
$C=\overline{\bigcup_n C_n}$ and $B=\overline{\bigcup_n B_n}$,
with unital phase-free coordinate inclusions. The algebras need not be
simple. A regular map between finite stages below means a possibly
nonunital homomorphism obtained by placing coordinate matrix-unit
copies of each source block in target blocks. It preserves the standard
diagonals and intertwines their diagonal expectations. Its multiplicities
may be zero. A coordinate permutation conjugacy preserves these properties.

Let $L_C\subseteq \T(C)$ and $L_B\subseteq \T(B)$ be nonempty closed faces.
Assume that every extremal trace in either face has a diffuse GNS factor.
The latter condition means that the finite von Neumann algebra generated
in that trace representation has no minimal projections. Completions use
the designated faces, and consequently include the quotient by their
common trace-null ideal.

\begin{lem}[Uniform disappearance of small blocks]
\label{lem:diffuse-face-small-blocks}
Let $A=\overline{\bigcup_m A_m}$ be a unital AF algebra with coordinate
stages $A_m=\bigoplus_j \mathrm{M}_{q_{m,j}}$, and let $L\subseteq \T(A)$ be a
nonempty compact face all of whose extremal GNS factors are diffuse.
Write $z_{m,j}$ for the unit of the $j$th block and put
\[
 z_{m,N}=\sum_{q_{m,j}\leq N}z_{m,j}.
\]
For every fixed $N$,
\[
 \sup_{\tau\in L}\tau(z_{m,N})\longrightarrow0.
\]
In particular,
\begin{equation}\label{eq:diffuse-face-rounding-mass}
 \sup_{\tau\in L}\sum_j\frac{\tau(z_{m,j})}{q_{m,j}}
 \longrightarrow0.
\end{equation}
\end{lem}

\begin{proof}
An old block of size greater than $N$ cannot map nontrivially into a
new block of size at most $N$. Thus $z_{m+1,N}\leq z_{m,N}$, where the
projections are compared using the coordinate inclusions.
Fix an extremal $\sigma\in L$ and work in its GNS von Neumann algebra
$\mathcal V_\sigma$. Since $L$ is a face, $\sigma$ is an extremal trace
of $A$, so $\mathcal V_\sigma$ is a factor. The strong limit $z$ of this
decreasing sequence of projections commutes with every fixed $A_k$,
because $z_{m,N}$ is central in $A_m$ for $m\geq k$. Thus $z$ is central
in $\mathcal V_\sigma$, and is either zero or one. If $z=1$, then every
image of $z_{m,N}$ is one. All the represented finite stages would then
be $N$-subhomogeneous. They satisfy the standard polynomial identity of
$\mathrm{M}_N$, which passes to their strong closure by bounded strong density
and continuity of multiplication. A diffuse finite factor does not
satisfy this identity: it contains a unital matrix algebra $\mathrm{M}_{N+1}$.
This contradiction shows $z=0$ and $\sigma(z_{m,N})\downarrow0$.

Every trace in the metrizable simplex $L$ is the barycenter of a
probability giving full measure to $\partial_eL\subseteq\partial_e\T(A)$.
Dominated convergence therefore
gives $\tau(z_{m,N})\downarrow0$ for every $\tau\in L$. These are decreasing
continuous affine functions on the compact set $L$. For any $\eta>0$,
the increasing open sets $\{\tau:\tau(z_{m,N})<\eta\}$ cover $L$; compactness
makes one of them all of $L$. This proves uniform convergence. Finally
\[
 \sum_j\frac{\tau(z_{m,j})}{q_{m,j}}
 \leq\tau(z_{m,N})+\frac1{N+1}.
\]
First choose $N$ large, and then $m$ large, to prove
\eqref{eq:diffuse-face-rounding-mass}.
\end{proof}

\begin{lem}[Regular finite-stage existence over a face]
\label{lem:diffuse-face-regular-existence}
Let $A,L$ be as in Lemma~\ref{lem:diffuse-face-small-blocks}, let
$F=\bigoplus_{i=1}^r\mathrm{M}_{d_i}$ have its standard matrix units, and let
$\lambda:L\to \T(F)$ be continuous and affine. Given $\delta>0$ and a
prescribed starting stage, there are a later stage $A_m$ and a regular,
possibly nonunital, homomorphism $\alpha:F\to A_m$ such that
\begin{equation}\label{eq:diffuse-face-existence-traces}
 \sup_{\tau\in L}\sup_{x\in(F)_1}
       |\tau(\alpha(x))-\lambda(\tau)(x)|<\delta.
\end{equation}
In particular, $\sup_{\tau\in L}\tau(1-\alpha(1))<\delta$.
\end{lem}

\begin{proof}
The trace space $\T(A)$ is a metrizable Choquet simplex and $L$ is
a closed face. Set $p_i(\tau)=\lambda(\tau)(1_{\mathrm{M}_{d_i}})$, so $p_i$ are nonnegative
continuous affine functions on $L$ with $\sum_i p_i=1$. The bounded
affine extension theorem of Edwards, in the form
\cite[Corollary~3.4]{Asimow}, extends them to a nonnegative continuous
affine partition $g_1,\ldots,g_r$ on $\T(A)$. To see that the scalar
extension theorem gives the whole partition, extend $p_1$ between zero
and one, then extend $p_2$ between zero and $1-g_1$, and continue;
take the last function to be the remaining affine difference. The
remaining upper bound is nonnegative at each step, and the required
inequalities hold on $L$.

The linear span of finite-stage trace profiles is uniformly dense in
$\operatorname{Aff}(\T(A))$. Here is a direct argument. If its closure
were proper, a separating functional on $\operatorname{Aff}(\T(A))$
could be extended to a signed measure on $\T(A)$. Since it vanishes on
the constant function, its positive and negative parts have equal mass.
Since it vanishes on all profiles $\tau\mapsto\tau(a)$ for
$a\in\bigcup_m(A_m)_{\mathrm{sa}}$, the barycenters of those two parts
are the same trace of $A$. Every continuous affine function integrates
to its value at the barycenter. The separating functional would
therefore vanish on all of $\operatorname{Aff}(\T(A))$, a contradiction.

Choose $b_1,\ldots,b_{r-1}$ in a fixed finite self-adjoint stage whose
trace profiles approximate $g_1,\ldots,g_{r-1}$ within $\eta$, and set
$b_r=1-\sum_{i<r}b_i$. All $r$ profiles then approximate the corresponding
$g_i$ within $(r-1)\eta$ (with the case $r=1$ immediate). At all sufficiently
late stages, every normalized block trace takes each $b_i$ to a number
at least $-r\eta$. Indeed, otherwise normalized block traces along a
sequence of later stages could be extended to states on $A$; a cluster
state is tracial on each fixed stage and hence is a trace of $A$,
contradicting the preceding profile bound and $g_i\geq0$.

At such a stage write
$a_{i,j}=\operatorname{tr}_{q_{m,j}}(b_i|_{\mathrm{M}_{q_{m,j}}})$.
Their sum over $i$ is one. Replace their negative entries by zero and
renormalize the resulting vector to have sum one, obtaining $w_{i,j}$.
If $s_j=\sum_i\max\{-a_{i,j},0\}$, then $s_j\leq r^2\eta$ and
\[
 \sum_i|w_{i,j}-a_{i,j}|=2s_j\leq2r^2\eta.
\]
Consequently, uniformly over $\T(A)$,
\[
 \left|\sum_j\tau(z_{m,j})w_{i,j}-g_i(\tau)\right|
 \leq(r-1)\eta+2r^2\eta.
\]
This estimate continues to hold if the stage is made later and the
clipping operation is repeated there.

In target block $j$ place
\[
 t_{i,j}=\left\lfloor\frac{q_{m,j}w_{i,j}}{d_i}\right\rfloor
\]
coordinate copies of $\mathrm{M}_{d_i}$. Their total rank is at most $q_{m,j}$,
so these choices define a regular homomorphism $\alpha:F\to A_m$;
unused coordinates give its complementary unit. For each minimal
projection $e^i_{kk}$,
\[
 \left|\tau(\alpha(e^i_{kk}))-\frac{p_i(\tau)}{d_i}\right|
 \leq\frac{(r-1)\eta+2r^2\eta}{d_i}
       +\sum_j\frac{\tau(z_{m,j})}{q_{m,j}}
 \qquad(\tau\in L).
\]
Off-diagonal matrix units have trace zero on both sides. For
$x\in(F)_1$ its diagonal coefficients have absolute value at most one;
thus summing these bounds over $\sum_i d_i$ minimal projections bounds
the error in \eqref{eq:diffuse-face-existence-traces}. Choose $\eta$
sufficiently small, and use
\eqref{eq:diffuse-face-rounding-mass} to make the rounding term sufficiently
small at a later stage. The desired error is then less than $\delta$.
The assertion about the complementary unit follows by taking $x=1$.
\end{proof}

\begin{lem}[Regular uniqueness over a trace face]
\label{lem:diffuse-face-regular-uniqueness}
Let $A$ be a unital AF algebra, let $L\subseteq \T(A)$ be a nonempty
closed face, and let $F=\bigoplus_i\mathrm{M}_{d_i}$. Suppose that
$\alpha,\beta:F\to A_k$ are regular, possibly nonunital, homomorphisms and
\[
 \sup_{\tau\in L}
 |\tau(\alpha(e^i_{11}))-\tau(\beta(e^i_{11}))|\leq\delta
 \qquad(i=1,\ldots,r).
\]
For every $\delta'>\delta$ there are a later stage $A_m$ and a
coordinate permutation unitary $u\in A_m$ such that
\begin{equation}\label{eq:diffuse-face-regular-uniqueness}
 \sup_{x\in(F)_1}\|u\alpha(x)u^*-\beta(x)\|_{2,L}
 \leq2\sqrt{\delta'\sum_i d_i}.
\end{equation}
Diffuseness is not required for this uniqueness assertion.
\end{lem}

\begin{proof}
For $m\geq k$ write the multiplicities of $\alpha$ and $\beta$ in
target block $j$ as $s_{i,m,j}$ and $t_{i,m,j}$. Define
\[
 a_{i,m}(\tau)=\sum_j\tau(z_{m,j})
                    \frac{|s_{i,m,j}-t_{i,m,j}|}{q_{m,j}}.
\]
These are nonnegative continuous affine functions, decreasing with
$m$: under a coordinate refinement, the new signed multiplicity
difference is the sum of old signed differences with nonnegative
integer coefficients, and the triangle inequality gives the claimed
monotonicity after taking traces.

Fix $\sigma\in\partial_eL$, put
$q_i=\alpha(e^i_{11})-\beta(e^i_{11})$, and work in its finite GNS factor
$\mathcal V_\sigma$. Averaging conjugations by the finite-stage unitary
group of $A_m$ is the trace-preserving conditional expectation onto
$A_m'\cap\mathcal V_\sigma$. Since $q_i\in A_m$, its average is
\[
 c_{i,m}=\sum_j
       \frac{s_{i,m,j}-t_{i,m,j}}{q_{m,j}}z_{m,j}.
\]
Thus $a_{i,m}(\sigma)=\sigma(|c_{i,m}|)$.
The commutants decrease to the center of $\mathcal V_\sigma$.
Equivalently, the corresponding orthogonal projections on its
$L^2$ space decrease to the projection onto the intersection of their
closed ranges. Hence $c_{i,m}\to\sigma(q_i)1$ in $L^2$, and
\[
 a_{i,m}(\sigma)\longrightarrow|\sigma(q_i)|\leq\delta.
\]
The last implication follows also directly from
$|\sigma(|c|)-|\sigma(q_i)||\leq\|c-\sigma(q_i)1\|_{2,\sigma}$,
by the scalar triangle inequality and its reverse.

The representing probability of every $\tau\in L$ gives full
measure to $\partial_eL$, which is contained in $\partial_e \T(A)$
because $L$ is a face. Integrating the
preceding limit gives $\lim_m a_{i,m}(\tau)\leq\delta$ for every
$\tau\in L$. Monotonicity and compactness, as in
Lemma~\ref{lem:diffuse-face-small-blocks}, therefore give a common later
stage at which $a_{i,m}(\tau)<\delta'$ for all $i$ and $\tau\in L$.

In each target block match $\min\{s_{i,m,j},t_{i,m,j}\}$ coordinate
copies of the $i$th source block. A coordinate permutation carries
the matched matrix units of $\alpha$ exactly to those of $\beta$.
Complete it arbitrarily on the remaining coordinates. To estimate
their rank, adjoin a one-dimensional source block to $F$, mapping its
unit to $1-\alpha(1)$ and $1-\beta(1)$ respectively. These extended
maps are unital. The absolute difference of their complementary
multiplicities is at most
$\sum_i d_i|s_{i,m,j}-t_{i,m,j}|$.
After matching all common copies, including the complementary block,
the remaining rank on either side is therefore at most
$\sum_i d_i|s_{i,m,j}-t_{i,m,j}|$.
Both maps agree on the matched reducing coordinates, and their
difference on the remaining coordinates has norm at most two on
$(F)_1$. Consequently
\[
 \|u\alpha(x)u^*-\beta(x)\|_{2,\tau}^2
 \leq4\sum_i d_i a_{i,m}(\tau)
 <4\delta'\sum_i d_i
 \qquad(x\in(F)_1,\ \tau\in L),
\]
which proves \eqref{eq:diffuse-face-regular-uniqueness}.
\end{proof}

\begin{thm}[Classification over diffuse closed trace faces]
\label{thm:diffuse-face-AF-classification}
Let $(E_C,C)$ and $(E_B,B)$ be standard unital AF diagonal pairs and
let $L_C\subseteq \T(C)$ and $L_B\subseteq \T(B)$ be nonempty closed faces
whose extremal GNS factors are diffuse. Every affine homeomorphism
$h:L_C\to L_B$ is implemented by an isomorphism
\[
 \Phi:C^{L_C}\longrightarrow B^{L_B},\qquad
 \Phi(E_C^{L_C})=E_B^{L_B},
\]
satisfying $h(\sigma)(\Phi(x))=\sigma(x)$ for all $\sigma\in L_C$.
It intertwines the completed diagonal expectations and is an isometry
for the designated uniform $2$-norms. Neither simplicity nor full trace
families are required, provided the designated families are faces.
\end{thm}

\begin{proof}
We give the summable intertwining, retaining its uniform estimates.
Choose a decreasing summable sequence $\varepsilon_n>0$.
Pass to increasing cofinal coordinate stages $F_n=C_{k_n}$ and
$G_n=B_{l_n}$, and construct regular maps
\[
 \alpha_n:F_n\to G_n,\qquad \beta_n:G_n\to F_{n+1}.
\]
Let $D_n$ and $Q_n$ be the sums of the matrix sizes in $F_n$ and
$G_n$ respectively. Require the trace errors on whole unit balls to
be at most $a_n$ for $\alpha_n$ and $b_n$ for $\beta_n$, with the markings
$h$ and $h^{-1}$ respectively, where
\[
 \begin{split}
 a_n&<\min\{\varepsilon_n^2/(32D_n),\ \varepsilon_n^2,\
                 \varepsilon_{n-1}^2/(32Q_{n-1})\},\\
 b_n&<\min\{\varepsilon_n^2/(32D_n),\
                 \varepsilon_n^2/(32Q_n),\ \varepsilon_n^2\}.
 \end{split}
\]
Omit the term involving $n-1$ when $n=1$.
The existence lemma constructs each map with any such error after its
source stage is fixed; its target stage can be taken beyond any prescribed
stage. These are therefore legitimate inductive requirements.

The composed maps $\beta_n\alpha_n$ and $\alpha_{n+1}\beta_n$ have trace
errors on their respective source unit balls at most $a_n+b_n$ and
$b_n+a_{n+1}$. For example, insert
$h(\sigma)(\alpha_n(x))$ between
$\sigma(\beta_n\alpha_n(x))$ and $\sigma(x)$ and use that
$\|\alpha_n(x)\|\leq\|x\|$.
Apply the uniqueness lemma to each composed map and the corresponding
coordinate inclusion. Its bound permits coordinate permutation
conjugacies with
\begin{equation}\label{eq:diffuse-face-intertwining-inverses}
 \begin{aligned}
 \|\beta_n\alpha_n(x)-x\|_{2,L_C}
       &\leq\varepsilon_n\|x\| &&(x\in F_n),\\
 \|\alpha_{n+1}\beta_n(y)-y\|_{2,L_B}
       &\leq\varepsilon_n\|y\| &&(y\in G_n).
 \end{aligned}
\end{equation}
First correct $\beta_n$ in its target; this does not change its trace
error. Then correct $\alpha_{n+1}$ in its target. All maps remain regular
and expectation preserving. Making these targets later preserves the
cofinality requirement and does not change previously constructed maps.

A trace error at most $a_n$ on $(F_n)_1$ gives the useful finite-stage
estimate
\begin{equation}\label{eq:diffuse-face-finite-continuity}
 \|\alpha_n(z)\|_{2,L_B}^2
 \leq\|z\|_{2,L_C}^2+a_n\|z\|^2
 \qquad(z\in F_n),
\end{equation}
by applying the trace estimate to $z^*z$; the analogous estimate holds
for $\beta_n$ with $b_n$. For $x\in F_n$, combine the two identities
in \eqref{eq:diffuse-face-intertwining-inverses} to get
\[
 \begin{split}
 \|\alpha_{n+1}(x)-\alpha_n(x)\|_{2,L_B}
 &\leq\|\alpha_{n+1}(x-\beta_n\alpha_n(x))\|_{2,L_B}
             +\varepsilon_n\|x\|\\
 &\leq\bigl(\sqrt{\varepsilon_n^2+4a_{n+1}}+\varepsilon_n\bigr)\|x\|
 <4\varepsilon_n\|x\|.
 \end{split}
\]
Here $\|x-\beta_n\alpha_n(x)\|\leq2\|x\|$.
Similarly, for $y\in G_n$,
\[
 \|\beta_{n+1}(y)-\beta_n(y)\|_{2,L_C}
 \leq\bigl(\sqrt{\varepsilon_n^2+4b_{n+1}}
                 +\varepsilon_{n+1}\bigr)\|y\|
 <4\varepsilon_n\|y\|.
\]
Thus the maps have bounded uniform $2$-limits on the finite-stage
unions. Since every map is a norm contraction, the limits extend
from these norm-dense unions to homomorphisms
\[
 \Phi_0:C\to B^{L_B},\qquad \Psi_0:B\to C^{L_C}.
\]
Multiplication passes to the limit on bounded sets.

Their exact trace markings follow by taking limits of $a_n,b_n\to0$:
\[
 h(\sigma)(\Phi_0(x))=\sigma(x),\qquad
 \sigma(\Psi_0(y))=h(\sigma)(y).
\]
Also $\|1-\alpha_n(1)\|_{2,L_B}^2\leq a_n$ and
$\|1-\beta_n(1)\|_{2,L_C}^2\leq b_n$, so both limiting maps are unital.
Applying the trace identities to $x^*x$ and $y^*y$ makes the limiting
maps uniform $2$-isometries. They consequently extend to the bounded
uniform completions, denoted $\Phi$ and $\Psi$.

Put $t_n=4\sum_{j\geq n}\varepsilon_j$. The Cauchy estimates give
\[
 \|\Phi(x)-\alpha_n(x)\|_{2,L_B}\leq t_n\|x\|
       \quad(x\in F_n),\qquad
 \|\Psi(y)-\beta_n(y)\|_{2,L_C}\leq t_n\|y\|
       \quad(y\in G_n).
\]
For an element $x$ of a fixed $F_k$, use arbitrarily late $n\geq k$,
the isometry of $\Psi$, and the first estimate of
\eqref{eq:diffuse-face-intertwining-inverses}, to obtain
\[
 \|\Psi\Phi(x)-x\|_{2,L_C}
 \leq(2t_n+\varepsilon_n)\|x\|\longrightarrow0.
\]
For $y\in G_k$, use the second estimate there and the isometry of
$\Phi$ to obtain
\[
 \|\Phi\Psi(y)-y\|_{2,L_B}
 \leq(t_n+t_{n+1}+\varepsilon_n)\|y\|\longrightarrow0.
\]
Density and continuity give both inverse identities on the completions.
Each finite-stage map carries the source diagonal into the target
diagonal and intertwines expectations. Their limits have the same
properties; the inverse identities give equality of the completed
diagonals and complete the proof.
\end{proof}

\subsection{Classification for full invariant-measure families}

\begin{lem}[The full invariant family gives a diffuse trace face]
\label{lem:full-family-AF-face}
Let $(X,R,K,D_K,\mathcal S)$ be a uniformly AF hyperfinite relation.
Assume that $K$ consists of all $R$-invariant Borel probability measures
on $X$, and that $R$ has infinite classes on a common invariant
$K$-universally conull set. In the reconstruction of
Theorem~\ref{thm:uniform-AF-reconstruction}, the designated family
$L=H(K)$ is a closed face of $\T(B)$, and every extreme trace in $L$
has a type $\mathrm{II}_1$ GNS factor.
\end{lem}

\begin{proof}
Let $\theta:X_0\to Y_0$ be the spatial tail-relation isomorphism of
Theorem~\ref{thm:uniform-AF-reconstruction}, after intersecting $X_0$
with the stated aperiodic conull set. The image $Y_0$ is Borel,
tail-invariant, and conull for every central measure in $L$.
If $\nu\in L$ and
\[
 \nu=t\nu_1+(1-t)\nu_2,
 \qquad 0<t<1,\quad \nu_1,\nu_2\in\T(B),
\]
positivity gives $\nu_1(Y_0)=\nu_2(Y_0)=1$. Pulling these central
measures back by $\theta$ gives $R$-invariant probabilities on $X$,
extended by zero off $X_0$. Fullness of $K$ puts both in $K$, so
$\nu_1,\nu_2\in L$. Thus $L$ is a face, and it is closed by compactness.

An extreme central measure in this face is extreme in $\T(B)$,
hence ergodic for tail equivalence. Its tail classes are infinite
almost everywhere. Centrality makes the masses of singleton points
in a tail class equal, so this measure is nonatomic. Its measured
tail-relation algebra is therefore an infinite-dimensional finite
factor, that is, a type $\mathrm{II}_1$ factor. This is the GNS
von Neumann algebra of $B$ for the corresponding trace.
\end{proof}

\begin{thm}[Classification of completed Cartan pairs]
\label{thm:uniform-Cartan-classification}
Let the regular uniform measure relations $(X_i,R_i,K_i,D_i,\mathcal S_i)$,
$i=1,2$, satisfy the following properties:
\begin{enumerate}
\item $K_i$ is the nonempty compact family of all $R_i$-invariant
Borel probability measures;
\item $(X_i,K_i)$ has relative SBP;
\item $(X_i,R_i,\mu)$ is amenable for every $\mu\in K_i$;
\item $R_i$ has infinite classes on a common invariant
$K_i$-universally conull set.
\end{enumerate}
Write $M_i=M_{R_i}$ for the untwisted completed relation algebra.
Every affine homeomorphism $H:K_1\to K_2$ is implemented by a unital
isomorphism of completed Cartan pairs
\[
 \Theta:(D_1,M_1)\longrightarrow(D_2,M_2),
 \qquad \tau_{H\mu}(\Theta(a))=\tau_\mu(a).
\]
It is an isometry for the uniform $2$-norms and intertwines the
canonical diagonal expectations. The same $H$ is implemented by a uniform measure relation
isomorphism, which conjugates the admissible full groups isometrically. Consequently the full
invariant-measure simplex is a complete invariant for marked
completed Cartan pairs and uniform measure relations in this class.

In particular, two Cartan subalgebras in one completed algebra,
arising in this class with the same designated traces, are conjugate
by a trace-preserving automorphism of that algebra.
\end{thm}

\begin{proof}
Theorem~\ref{thm:uniform-CFW-SBP} supplies compatible finite matrix
charts. Their AF reconstruction has the trace faces described in
Lemma~\ref{lem:full-family-AF-face}, and identifies the completed pairs with
\[
 ((E_i)^{L_i},B_i^{L_i}),
 \qquad H_i:K_i\longrightarrow L_i\subseteq\T(B_i),
\]
where $L_i$ are closed faces with type $\mathrm{II}_1$ extreme traces.
Apply Theorem~\ref{thm:diffuse-face-AF-classification} to the affine homeomorphism
$H_2HH_1^{-1}:L_1\to L_2$, and compose with these reconstructions.
Its trace identity gives the asserted marking and uniform $2$-norm
isometry, and its expectation identity gives the stated compatibility.
Theorem~\ref{thm:uniform-relation-pair-isomorphism} supplies the
spatial relation isomorphism. Conjugation preserves the displacement
sets in the full-group metric, so it is isometric.
Conversely, a marked pair isomorphism or a uniform measure relation
isomorphism gives an affine homeomorphism of the designated families.
For the final assertion, use the identity marking of the common trace
family and identify the two ambient completed algebras.
\end{proof}

\begin{cor}[Simple full-trace models for the full family]
\label{cor:full-family-simple-AF-model}
Every relation satisfying the hypotheses of
Theorem~\ref{thm:uniform-Cartan-classification} admits a marked model
with the standard diagonal of a simple unital infinite-dimensional AF
algebra and its full trace simplex. Thus it satisfies
\eqref{eq:uniform-AF-model-hypothesis} and is uniformly measure
isomorphic to a Cantor minimal $\mathbb Z$-system with the prescribed
full invariant-measure simplex.
\end{cor}

\begin{proof}
Lemma~\ref{lem:full-family-AF-face} identifies $K$ with a compact
face of an AF trace simplex, so $K$ is a metrizable Choquet simplex.
Choose an infinite Cantor minimal homeomorphism realizing this
simplex, using~\cite{Dow91}; for a singleton simplex one may use
an odometer. Use a properly ordered simple Bratteli diagram for this
Cantor minimal system~\cite[Corollary~1.3]{GPSA}. Its tail relation
is obtained by splitting one orbit into two infinite half-orbits.
Every tail-invariant probability is nonatomic, since its tail classes
are infinite. It is therefore also invariant under the homeomorphism,
which agrees off the split orbit with a countable union of tail maps.
Conversely, each prefix substitution is piecewise a power of the
homeomorphism. Thus the central probabilities of the diagram, or
equivalently the traces on its simple unital infinite-dimensional AF
algebra, give precisely the prescribed simplex $K$. Alternatively,
applying \cite[Theorem~4.1]{TC} to this Cantor minimal system gives
the same type of marked AF model. Apply
Theorem~\ref{thm:diffuse-face-AF-classification} to this pair and the reconstructed
pair over $L=H(K)$. This gives the required marked simple model.
Corollary~\ref{cor:uniform-CFW-Z-generator} gives the minimal
dynamical realization.
\end{proof}

\begin{cor}[Minimal principal amenable \'etale groupoids]
\label{cor:minimal-SBP-groupoid-classification}
Let $\mathcal G_i$, $i=1,2$, be second-countable locally compact
Hausdorff principal minimal \'etale groupoids with infinite compact
metrizable unit spaces $X_i$. Assume that their full invariant
probability families $K_i$ are nonempty, that $(X_i,K_i)$ has
relative SBP, and that every designated measured orbit relation is
amenable. In particular, topological amenability suffices.
Then every affine homeomorphism $K_1\to K_2$ is implemented by the
completed-pair isomorphism and uniform measure orbit equivalence of
Theorem~\ref{thm:uniform-Cartan-classification}. Each pair has the
simple full-trace AF model and Cantor minimal realization of
Corollary~\ref{cor:full-family-simple-AF-model}.
Here the completed groupoid pair is
\[
 \bigl(\mathrm{C}(X_i)^{K_i},(\mathrm{C}_r^*(\mathcal G_i))^{K_i}\bigr),
\]
with its canonical invariant-measure traces. The completed-pair
classification also holds for locally trivial circle twists of these
groupoids.
\end{cor}

\begin{proof}
The full invariant family is compact for these \'etale presentations.
Proposition~\ref{prop:principal-groupoid-SBP-atlas} supplies regular
countable atlases. Minimality on an infinite compact unit space
forces every orbit to be infinite: a finite orbit is closed and
cannot be dense. All hypotheses of
Theorem~\ref{thm:uniform-Cartan-classification} are therefore met.
For each atlas bisection, continuous cutoffs at its small-boundary
source approximate its graph operator in uniform $2$-norm.
These cutoffs define compactly supported continuous groupoid
coefficients. Conversely, a compactly supported continuous groupoid
function is a finite sum of continuous bisection coefficients, using
a partition of unity. The graph realization preserves the canonical
expectation and traces, and these two density statements identify
the groupoid completion with its completed relation pair.
For a locally trivial circle twist, Proposition~\ref{prop:SBP-regularizes-twist} gives
a regular scalar presentation, and
Theorem~\ref{thm:uniform-regular-untwisting} removes its cocycle while
preserving the diagonal and designated traces.
\end{proof}

\begin{rem}[The scope of the invariant]
\label{rem:uniform-classification-scope}
Minimality is not required in the relation classification theorem;
aperiodicity and the full compact invariant family suffice. More
generally, its proof applies whenever the reconstructed trace family
is a closed face with type $\mathrm{II}_1$ extreme traces.
An arbitrary designated compact convex family does not have this
property. For example, compare a uniquely ergodic Cantor minimal
homeomorphism with two disjoint copies, designating only the average
of their invariant measures in the second system. Both designated
families are singletons, and both relations are aperiodic, amenable,
and have relative SBP. Their completed algebras are respectively
the hyperfinite factor and its two-fold direct sum, with different
centers. Fullness, or the stated trace-face condition, is substantive.
Finite classes also require their matrix sizes as additional data.
The conjugacy assertion for Cartan subalgebras is by automorphisms;
no inner-unitary conjugacy assertion is made.

The small boundary hypothesis is substantive even for two inclusions
in the same ambient completion with the same designated traces.
Take a free minimal $\mathbb Z$-system of Giol--Kerr with crossed
product $A$ of positive radius of comparison~\cite[Theorem~2.2]{GiolKerr}.
Its canonical diagonal $D=\mathrm{C}(X)$ fails SBP, since SBP implies zero
mean dimension and hence $\mathcal Z$-stability of $A$
by the main theorem of~\cite{ENmdim}. By~\cite[Theorem~A]{Vaccaro},
$A$ has uniform property $\Gamma$. The ambient classification
in~\cite[Theorem~C]{CCEGSTW} therefore identifies $A^u$, with its trace
marking, with $B^u$ for a simple unital infinite-dimensional AF algebra
$B$ with the same trace simplex. Transport the completed standard
diagonal of $B$ into $A^u$. It and $D^u$ are uniform Cartan subalgebras
by~\cite[Propositions~2.5 and~2.6]{TC}. The transported diagonal has
real rank zero. Its failure for $D^u$, by~\cite[Theorem~8.4]{TC},
distinguishes the two diagonals. Thus they are nonconjugate
even by arbitrary automorphisms of $A^u$.
\end{rem}

\subsection{Uniform full groups and extreme amenability}
\label{subsec:uniform-full-group-amenability}

Throughout this subsection, full groups carry the topology of the
complete bi-invariant metric $d_K$. Topological amenability means the
fixed-point property for continuous affine actions on compact convex
subsets of locally convex spaces; extreme amenability means the
fixed-point property for all continuous actions on compact Hausdorff
spaces. A metrizable group is a \emph{L\'evy group} if it has increasing
compact subgroups with dense union whose normalized Haar measures
concentrate: for every $\eps>0$ and every sequence of Borel subsets
$A_n$ of the subgroups with Haar masses bounded away from zero, the
Haar masses of their $\eps$-neighborhoods tend to one. The subgroups
constructed below are finite.

\begin{lem}[Dense restriction to a designated measure]
\label{lem:uniform-full-group-measured-density}
Let $(X,R,K,D_K,\mathcal S)$ be a regular uniform measure relation
with relative $K$-SBP. For every $\mu\in K$, restriction gives a
continuous homomorphism
\[
 \rho_\mu:[R]^u_{K,\mathcal S}\longrightarrow[R]_\mu
\]
with dense image, where $[R]_\mu$ has its usual disagreement metric.
Fullness of $K$, aperiodicity, and measured amenability are not
required for this assertion.
\end{lem}

\begin{proof}
Continuity follows from $d_\mu\leq d_K$. By
Corollary~\ref{cor:uniform-relation-countable-group-presentation},
work on a common invariant $K$-conull Borel set carrying global
admissible maps $g_1,g_2,\ldots$ whose graphs cover $R$.
For $T\in[R]_\mu$, choose a Borel representative on an invariant
$\mu$-conull subset of this set and extend it by the identity on the
complement. This representative preserves $\mu$; no assertion that
it preserves the other designated measures is needed. Its disjoint
first-chart partition is
\[
 E_j=\{x:T(x)=g_j(x)\}
       \setminus\bigcup_{i<j}\{x:T(x)=g_i(x)\}.
\]
These sets need not be regular. Fix $\delta>0$ and choose $m$ such that
\[
 \mu\!\left(X\setminus\bigcup_{j\leq m}E_j\right)<\delta.
\]
The finite projection diagonals of
Proposition~\ref{prop:relative-SBP-finite-diagonals} are dense in
$L^2(X,\mu)$: their norm closure contains $\mathrm{C}(X)$, which is
$L^2(\mu)$-dense. Approximating $1_{E_j}$ by a finite-diagonal
positive contraction and taking its spectral projection for
$[1/2,1]$ therefore gives regular sets $P_j$ for which
\[
 \eta:=\sum_{j\leq m}\mu(P_j\triangle E_j)
\]
is arbitrarily small.

Disjointize the sources and then the ranges:
\[
 C_j=P_j\setminus\bigcup_{i<j}P_i,\qquad
 A_j=C_j\setminus
       g_j^{-1}\!\left(\bigcup_{i<j}g_i(A_i)\right).
\]
Finite Boolean operations and admissible pullbacks make all these
sets regular. The ranges $g_j(A_j)$ are regular and disjoint, so
\[
 \varphi=\bigsqcup_{j\leq m}g_j|_{A_j}
\]
is an admissible partial bijection.

Source disjointization loses at most $\eta$ of the correct points:
a lost point is either in its own $E_j\setminus P_j$ or in an earlier
$P_i\setminus E_i$. Moreover,
$\sum_j\mu(C_j\setminus E_j)\leq\eta$.
If a correct edge $x\mapsto T(x)=g_j(x)$ is removed by range
disjointization, its image equals $g_i(z)$ for an earlier retained
point $z\in A_i$. Necessarily $z\notin E_i$; otherwise injectivity
of $T$ gives $z=x$, contrary to source disjointness.
The images of the removed correct points thus lie in
$\bigcup_i g_i(A_i\setminus E_i)$, whose measure is at most $\eta$.
Since $T$ preserves $\mu$, the removed points have measure at most
$\eta$ as well. Lemma~\ref{lem:uniform-partial-map-extension} extends
$\varphi$ to $U\in[R]^u_{K,\mathcal S}$, with
\begin{equation}
\label{eq:uniform-full-group-measured-density}
 d_\mu(U,T)\leq\delta+2\eta.
\end{equation}
Letting $\delta$ and $\eta$ tend to zero proves density.
\end{proof}

\begin{lem}[Concentration for finite permutation products]
\label{lem:uniform-permutation-concentration}
Let $H=\prod_{j=1}^b\mathfrak S_{r_j}$ have normalized counting
measure $m_H$, and put
\[
 \rho(\sigma,\pi)=\max_j
 \frac{\#\{i:\sigma_j(i)\ne\pi_j(i)\}}{r_j},
 \qquad s=\sum_{j=1}^b r_j^{-1}.
\]
For $\eps>0$, its concentration function satisfies
\[
 \alpha_{H,\rho}(\eps)\leq
       \exp\!\left(-\frac{\eps^2}{32s}\right).
\]
The same bound holds for every bi-invariant metric $d\leq\rho$.
\end{lem}

\begin{proof}
Here
$\alpha_{H,\rho}(\eps)=1-\inf\{m_H(A_\eps^\rho):
m_H(A)\geq1/2\}$, where $A_\eps^\rho$ is the $\eps$-neighborhood
of $A$. The assertion is immediate if $H$ is trivial.
Otherwise grow each symmetric factor through
$\mathfrak S_1\subseteq\cdots\subseteq\mathfrak S_{r_j}$, completing
one factor before beginning the next. Each successive quotient has
diameter at most $2/r_j$ in the factor metric: its cosets have
representatives differing by a transposition. Thus the sum of
squared quotient diameters is at most
\[
 \sum_{j=1}^b(r_j-1)\frac4{r_j^2}\leq4s.
\]
The compact-group concentration inequality
of Giordano--Pestov~\cite[Theorem~2.9]{GP} gives the stated bound.
Since $A_\eps^\rho\subseteq A_\eps^d$ when $d\leq\rho$, the
concentration function for $d$ is no larger.
\end{proof}

\begin{thm}[Uniform Giordano--Pestov theorem]
\label{thm:uniform-full-group-extreme-amenability}
Let $(X,R,K,D_K,\mathcal S)$ be a regular uniform measure relation.
Assume relative $K$-SBP, that $K$ is the nonempty compact family of
all $R$-invariant Borel probability measures, and that $R$ has
infinite classes on a common invariant universally conull set.
For $G=[R]^u_{K,\mathcal S}$, the following are equivalent:
\begin{enumerate}
\item $(X,R,\mu)$ is amenable for every $\mu\in K$;
\item $(G,d_K)$ is topologically amenable;
\item $(G,d_K)$ is extremely amenable;
\item $(G,d_K)$ is a L\'evy group.
\end{enumerate}
In the fourth condition, the concentrating compact subgroups can
be chosen finite. No minimality or Bauer hypothesis is required.
\end{thm}

\begin{proof}
Assume the first condition. By
Corollary~\ref{cor:full-family-simple-AF-model}, the completed pair
has a marked standard diagonal model in a simple unital
infinite-dimensional AF algebra $B$, carrying its full trace
simplex. Choose coordinate stages
\[
 B_n=\bigoplus_{j=1}^{b_n}\mathrm{M}_{a_{n,j}}.
\]
Let $U=\mathrm{M}_{2^\infty}$ with its standard diagonal $E_U$ and unique
trace. The AF pair $(E_B\otimes E_U,B\otimes U)$ has the same full
trace simplex as $(E_B,B)$. Every extreme trace of either simple
infinite-dimensional AF algebra has type $\mathrm{II}_1$ GNS
closure, so Theorem~\ref{thm:diffuse-face-AF-classification}, applied
to the full trace faces, identifies their completed Cartan pairs
with the natural trace marking.

Choose dyadic integers $q_n$ with $q_n\mid q_{n+1}$ and
$q_n\geq2^n b_n$. Then
\[
 C_n=B_n\otimes \mathrm{M}_{q_n}
     =\bigoplus_{j=1}^{b_n}\mathrm{M}_{r_{n,j}},
 \qquad r_{n,j}=q_na_{n,j},
 \qquad s_n:=\sum_jr_{n,j}^{-1}\leq2^{-n},
\]
are coordinate stages exhausting $B\otimes U$.
Their permutation groups $H_n=\prod_j\mathfrak S_{r_{n,j}}$,
viewed in the completed normalizer quotient modulo diagonal
unitaries, are increasing. Lemma~\ref{lem:AF-normalizer-density}
makes their union dense in that quotient.
Transporting them to $G$ preserves nesting and density by
Theorem~\ref{thm:uniform-relation-pair-isomorphism} and
Proposition~\ref{prop:uniform-relation-isomorphism-transport}.
By \eqref{eq:uniform-relation-quotient-metric}, the transported
disagreement distance $d_K$ is one half of the square of the
quotient uniform $2$-distance. Allowing independent diagonal phases
on both representatives gives the same infimum, by unitary invariance
and normalization of the diagonal. We use $d_K$ in the concentration
estimate below.

For the unit $z_{n,j}$ of the $j$th block and the corresponding
normalized Hamming metric $d_{n,j}$, the metric on $H_n$ is
\[
 d_K(\sigma,\pi)
 =\sup_{\tau\in \T(B\otimes U)}
       \sum_j\tau(z_{n,j})d_{n,j}(\sigma_j,\pi_j)
 \leq\max_j d_{n,j}(\sigma_j,\pi_j).
\]
Lemma~\ref{lem:uniform-permutation-concentration} yields
\begin{equation}
\label{eq:uniform-full-group-Levy-bound}
 \alpha_{H_n,d_K}(\eps)
 \leq\exp\!\left(-\frac{\eps^2}{32s_n}\right)
 \leq\exp(-2^n\eps^2/32)\longrightarrow0.
\end{equation}
To see concentration also for sets of mass below $1/2$, fix
$r>0$ and let $m_n(A_n)\geq a>0$, where $m_n$ is normalized
counting measure. Once $\alpha_{H_n,d_K}(r)<a$, the
$r$-neighborhood of $A_n$ has mass at least $1/2$. Otherwise
its complement has mass at least $1/2$ and that complement's open
$r$-neighborhood is disjoint from $A_n$, contradicting the bound
on $\alpha_{H_n,d_K}(r)$. Applying the bound once more gives
$m_n((A_n)_{2r})\geq1-\alpha_{H_n,d_K}(r)\to1$.
This proves the fourth condition. A L\'evy group is extremely
amenable~\cite[Corollary~2.19]{GP}, and extreme amenability implies
topological amenability.

Finally assume the second condition and fix $\mu\in K$.
Lemma~\ref{lem:uniform-full-group-measured-density} gives a
continuous dense homomorphism $G\to[R]_\mu$. Restricting any
continuous affine action of $[R]_\mu$ on a compact convex set along
this homomorphism supplies a $G$-fixed point, which density makes
$[R]_\mu$-fixed. Hence $[R]_\mu$ is topologically amenable.
The measure $\mu$ is nonatomic: an atom in the aperiodic conull
set would have the same positive mass at every point of its
infinite $R$-class, contrary to finiteness of $\mu$.
Le Ma\^itre's nonergodic measured theorem~\cite[Theorem~8.4]{LM}
therefore makes $R$ $\mu$-hyperfinite, hence amenable by
Connes--Feldman--Weiss~\cite{CFW}. This holds for every $\mu\in K$.
\end{proof}

\begin{rem}
The full-family hypothesis is used in the forward proof to obtain
the simple full-trace model and transport its UHF amplification.
The dense restriction lemma and the implication from topological
group amenability to measured amenability do not require fullness.
For the principal minimal amenable groupoids of
Corollary~\ref{cor:minimal-SBP-groupoid-classification}, the admissible
uniform full group is consequently a L\'evy group and extremely
amenable. These assertions concern the admissible full group in
the uniform measure topology.
\end{rem}

\subsection{Reconstruction from the uniform full group}
\label{subsec:uniform-Dye-reconstruction}

The full group retains the regular Boolean algebra through its
involutions. Its topology then recovers the uniform diagonal topology.
For a full invariant family, the designated measures can also be
recovered intrinsically. We first separate these three steps.

Write $\mathcal P_K=\operatorname{Proj}(D_K)$ and
$c_K(p)=\sup_{\mu\in K}\mu(p)$. An element
$T\in G=[R]^u_{K,\mathcal S}$ acts on $D_K$ by
$\beta_T(f)=f\circ T^{-1}$, and its support is the regular projection
$p_T=1_{\{x:Tx\ne x\}}$. All equalities of projections and full-group
elements are understood modulo universally null sets.

\begin{lem}[Every regular projection is an involution support]
\label{lem:uniform-involution-supports}
Let $(X,R,K,D_K,\mathcal S)$ be a regular uniform measure relation
with relative $K$-SBP and infinite classes on a common invariant
universally conull set. For every $p\in\mathcal P_K$ there is an
involution $t_p\in G$ with $p_{t_p}=p$. Neither amenability nor
fullness of the designated family is needed.
\end{lem}

\begin{proof}
Choose a countable basis $(U_m)$ of open sets with universally null
boundaries, and a Borel regular representative $A_0$ of $p$ in the
aperiodic conull set. The projections $1_{U_m}$ belong to $D_K$.
Enumerate all pairs $(s_j,U_m)$, where the graphs of the atlas
$(s_j)$ cover $R$. Starting with the unmatched set $A_{n-1}$, for the
$n$th pair $(s,U)$ put
\[
 E_n=\{x\in A_{n-1}\cap\operatorname{dom}s\cap U:
                    sx\in A_{n-1}\setminus U\},\qquad
 F_n=s(E_n),\qquad
 A_n=A_{n-1}\setminus(E_n\cup F_n).
\]
These sets are regular, by the corner isomorphism of $s$.
The sets $E_n,F_n$ are disjoint. Flip them using $s$ and $s^{-1}$,
and fix their complement. This is an admissible involution.
The successive flip supports are disjoint, so their product through
stage $n$ is an involution $t_n$ supported on $A_0\setminus A_n$.
All representative identities may be imposed on one common invariant
universally conull set.

The residual Borel set $A_\infty=\bigcap_n A_n$ meets each class in at
most one point. Indeed, two distinct residual points would be joined
by an atlas map and separated by a basis element; the corresponding
pair would have removed them. A Borel transversal in an aperiodic
probability-measure-preserving relation is null for every invariant
probability. To see this directly, enumerate distinct points in each
infinite class above the transversal using the countable atlas.
For each index this gives an injective Borel partial $R$-map; its
image has the same measure as the transversal, and these images are
pairwise disjoint. Invariance for such maps follows by countably
partitioning them into atlas restrictions. Finite total measure
therefore forces the transversal to have measure zero.

For each $n$, the function $\mu\mapsto\mu(A_n)$ is continuous on $K$,
since $1_{A_n}\in D_K$. These functions decrease pointwise to zero.
Dini's theorem gives $c_K(A_n)\to0$. For $m\geq n$,
\[
 d_K(t_m,t_n)\leq c_K(A_n).
\]
Completeness of the full group
(Proposition~\ref{prop:uniform-relation-full-group-complete})
gives an involution $t=\lim_n t_n$. Moreover,
\[
 c_K(p_t\mathbin{\triangle}p)
 \leq d_K(t,t_n)+c_K(A_n)\longrightarrow0.
\]
Thus $p_t=p$. For $p=0$, take the identity.
\end{proof}

\begin{lem}[Boolean and diagonal reconstruction]
\label{lem:uniform-full-group-Boolean-reconstruction}
Let two regular uniform measure relations satisfy the hypotheses of
Lemma~\ref{lem:uniform-involution-supports}, and write $G_i$ for their
admissible full groups. Every abstract group isomorphism
$\Phi:G_1\to G_2$ induces a unique equivariant Boolean isomorphism
\[
 b:\mathcal P_{K_1}\longrightarrow\mathcal P_{K_2},
 \qquad b(p_T)=p_{\Phi(T)}.
\]
If $\Phi$ is a topological group isomorphism for the uniform
disagreement topologies, $b$ extends uniquely to a unital
$*$-isomorphism $\Theta:D_{K_1}\to D_{K_2}$ that is bicontinuous on
norm-bounded sets for the uniform $2$-topologies and satisfies
$\Theta\beta_T=\beta_{\Phi(T)}\Theta$.
\end{lem}

\begin{proof}
The action of $G_i$ on $\mathcal P_{K_i}$ is faithful, and $p_T$ is
its least Boolean support. Here is a useful verification of the
support assertion. A countable separating SBP basis expresses
$\{Tx\ne x\}$ as a countable union of regular sets on which $x$ and
$Tx$ lie on opposite sides of a basis element. Every nonzero regular
$r\leq p_T$ therefore contains a nonzero regular $e$ with
$e\wedge\beta_T(e)=0$. Consequently $T$ acts nontrivially on every
nonzero regular part of $p_T$.

The Boolean algebra $\mathcal P_{K_i}$ need not be Dedekind complete.
Pass to its Dedekind completion $\widehat{\mathcal P}_{K_i}$, in
which it is order dense. Boolean automorphisms extend uniquely.
The extended action is the identity below $1-p_T$, by order density.
Its support is still $p_T$: any nonzero
part of $p_T$ omitted from a proposed smaller support contains a
nonzero regular projection, and then a displaced regular subprojection
as above. Every nonzero element of the completion contains a nonzero
regular projection $p$ and hence, by
Lemma~\ref{lem:uniform-involution-supports}, the support of a
nontrivial involution. Thus these automorphism groups have many
involutions in the sense of Fremlin~\cite[\S384D]{Fremlin}.
That theorem gives a unique equivariant Boolean isomorphism between
the Dedekind completions implementing $\Phi$. Supports are preserved.
Since every regular projection is an involution support, this
isomorphism carries $\mathcal P_{K_1}$ onto $\mathcal P_{K_2}$,
giving $b$. Any equivariant Boolean isomorphism of the regular
algebras extends uniquely to their Dedekind completions, so the same
uniqueness applies before completion.

Suppose now that $\Phi$ is a topological isomorphism. For the
involutions supplied by the first lemma,
\[
 d_{K_1}(t_p,1)=c_{K_1}(p),\qquad
 d_{K_2}(\Phi(t_p),1)=c_{K_2}(b(p)).
\]
Continuity of $\Phi$ and $\Phi^{-1}$ at the identity therefore gives
uniform moduli in both directions for the capacities of projections.
Apply this also to symmetric differences to obtain bicontinuity of
$b$ for the Boolean capacity metrics.

Extend $b$ linearly to finite projection simple functions. It is a
norm-isometric unital $*$-isomorphism on these algebras. If
$\|h\|\leq M$ and $p_\delta=1_{\{|h|>\delta\}}$, then
\[
 c_{K_1}(p_\delta)\leq
       \delta^{-2}\|h\|_{2,K_1}^2,\qquad
 \|b(h)\|_{2,K_2}^2\leq
       \delta^2+M^2c_{K_2}(b(p_\delta)).
\]
The capacity moduli show uniform $2$-continuity on each norm-bounded
set, in both directions. Relative SBP gives bounded uniform
$2$-density of finite projection simple functions
(Proposition~\ref{prop:relative-SBP-finite-diagonals}). Bounded
completeness consequently gives the asserted extension $\Theta$,
with inverse. Multiplication and involution pass to bounded limits.
Boolean equivariance extends by the same density argument.
\end{proof}

\begin{lem}[Intrinsic recovery of the full measure family]
\label{lem:uniform-full-group-invariant-states}
Assume the hypotheses of
Theorem~\ref{thm:uniform-Cartan-classification}. The states on $D_K$
that are invariant under $G$ and continuous on norm-bounded sets for
the uniform $2$-topology are precisely the states
$f\mapsto\int f\,d\mu$, $\mu\in K$.
\end{lem}

\begin{proof}
Each designated integral has these properties. Conversely, transport
a state $\lambda$ with the stated properties to the standard
completed diagonal of the simple full-trace AF model supplied by
Corollary~\ref{cor:full-family-simple-AF-model}. Let
$B_n=\bigoplus_j\mathrm{M}_{r_{n,j}}$ be its unital coordinate stages.
Finite coordinate permutations belong to the admissible full group.
Invariance makes $\lambda$ equal on all diagonal minimal projections
within each matrix block. Its restriction to the finite diagonal
therefore extends uniquely to a trace $\tau_n$ on $B_n$.
These traces are compatible: a trace on a finite-dimensional algebra
is determined by its diagonal values, and the coordinate inclusions
preserve the diagonals. They define $\tau\in \T(B)$.

The model carries all of $\T(B)$ as its designated family, so $\tau$
has its canonical bounded uniform $2$-continuous extension.
This extension agrees with $\lambda$ on the finite diagonals, hence
on the original AF diagonal by norm density, and on its completion
by bounded uniform $2$-density. Transport back to the original model
identifies $\lambda$ with a measure in $K$.
\end{proof}

\begin{thm}[Uniform Dye reconstruction]
\label{thm:uniform-Dye-reconstruction}
For $i=1,2$, let $(X_i,R_i,K_i,D_{K_i},\mathcal S_i)$ be a regular
uniform measure relation with relative $K_i$-SBP. Assume that $K_i$
is the nonempty compact family of all $R_i$-invariant Borel
probability measures, that $R_i$ is amenable for every such measure,
and that its classes are infinite on a common invariant universally
conull set. Put $G_i=[R_i]^u_{K_i,\mathcal S_i}$. The following are
equivalent:
\begin{enumerate}
\item the two relations are uniformly measure isomorphic;
\item $G_1$ and $G_2$ are isomorphic as topological groups;
\item $(G_1,d_{K_1})$ and $(G_2,d_{K_2})$ are isomorphic as metric groups;
\item $K_1$ and $K_2$ are affinely homeomorphic.
\end{enumerate}
Every given topological group isomorphism $\Phi:G_1\to G_2$ is
implemented by a uniform measure relation isomorphism $\theta$:
\[
 \Phi(T)=\theta T\theta^{-1}\qquad(T\in G_1).
\]
The implementing map is unique modulo equality almost everywhere for
every designated measure. Its affine measure marking is unique, and
$\Phi$ is automatically an isometry.
\end{thm}

\begin{proof}
A uniform relation isomorphism induces an isometric full-group
isomorphism by
Proposition~\ref{prop:uniform-relation-isomorphism-transport}.
Theorem~\ref{thm:uniform-Cartan-classification} identifies uniform
relation isomorphism with affine homeomorphism of the full families.
It remains to reconstruct a uniform relation isomorphism from a
given topological group isomorphism $\Phi$.

Lemma~\ref{lem:uniform-full-group-Boolean-reconstruction} supplies the
equivariant boundedly bicontinuous diagonal isomorphism $\Theta$.
Lemma~\ref{lem:uniform-full-group-invariant-states} shows that
pullback by $\Theta$ carries the full designated state family onto
the other one. Thus
\[
 J:K_2\longrightarrow K_1,\qquad
 \int\Theta(f)\,d\nu=\int f\,d(J\nu)
\]
is an affine bijection. Integral profiles of completed diagonal
elements are continuous on the designated compact families.
Testing the displayed formula on $\mathrm{C}(X_1)$ proves continuity of $J$;
applying the inverse formula proves continuity of $J^{-1}$.
Put $H=J^{-1}$. In particular, $\Theta$ is an exact uniform
$2$-isometry with this marking.

The spatial diagonal theorem
(Theorem~\ref{thm:orbit-spatial-diagonal}) implements $\Theta$ by a
Borel bijection $\theta$ between common universally conull Borel
sets, with $\theta_*\mu=H(\mu)$. We explain the countable
bookkeeping needed to make these domains invariant and to identify
the orbits. Extend every source atlas chart and every target atlas
chart to an admissible full-group element, using
Lemma~\ref{lem:uniform-partial-map-extension}. Let $\Gamma$ be the
countable subgroup of $G_1$ generated by the source extensions and
the $\Phi^{-1}$-images of the target extensions. By
Corollary~\ref{cor:uniform-relation-countable-group-presentation},
these elements and their images have simultaneous Borel
representatives on common invariant conull sets. Their orbits are
respectively $R_1$ and $R_2$.

Equivariance of $\Theta$, tested on countable continuous separating
families, gives
$\theta(gx)=\Phi(g)\theta(x)$ for every $g\in\Gamma$ outside one
universally null Borel set. Include in this exceptional set the
complements of the spatial domains and the preimages needed for
these equations to be defined. Remove its $R_1$-saturation; this is
still universally null, by the countable atlas and invariance.
The resulting source set $X_1^0$ is $R_1$-invariant and conull.
Its image $X_2^0=\theta(X_1^0)$ is Borel, conull for every measure
in $K_2$, and invariant under $\Phi(\Gamma)$, hence under $R_2$.
The intertwining equations show that $\theta$ carries whole
$R_1$-classes onto whole $R_2$-classes.

A source partial chart is a restriction of its chosen full-group
extension to a regular projection. Its conjugate is therefore a
restriction of an admissible element of $G_2$ to the corresponding
projection under $\Theta$. It is admissible. The chosen target
extensions give the reverse assertion. Consequently $\theta$
satisfies Definition~\ref{defn:uniform-measured-relation-isomorphism}.

For each arbitrary $T\in G_1$, the same diagonal equivariance and a
countable separating family give
$\Phi(T)=\theta T\theta^{-1}$ as full-group classes. This assertion
does not require an intersection of exceptional sets over all $T$.
Uniqueness of the Boolean reconstruction gives uniqueness of
$\Theta$. A countable continuous separating family shows that two
spatial maps inducing $\Theta$ agree outside one universally null
set. The integral
identity determines $H$ uniquely. Finally spatial transport with the
full measure marking preserves $d_K$, so $\Phi$ is isometric.
\end{proof}

\begin{example}[Topology alone can lose designated weights]
\label{ex:full-group-topology-loses-designated-weights}
Fullness cannot be omitted from the preceding topological
reconstruction statement. Let $R$ be the disjoint union of the
orbit relations of two uniquely ergodic Cantor minimal
homeomorphisms, with invariant probabilities $\mu_1,\mu_2$.
For $0<a<1$, designate only
$K_a=\{a\mu_1+(1-a)\mu_2\}$. These relations are aperiodic,
amenable, and have relative SBP. Their regular diagonals and
admissible full groups agree for every $a$, since the weighted
single-measure $2$-norms are equivalent. The common group is
$[R_1]_{\mu_1}\times[R_2]_{\mu_2}$, with metric
\[
 d_a=a\,d_{\mu_1}+(1-a)\,d_{\mu_2}.
\]
All these metrics are bi-Lipschitz equivalent. However, a marked
probability orbit equivalence must preserve the two atoms of the
invariant measure algebra and their unordered weights
$\{a,1-a\}$. Thus $a=1/3$ and $a=1/4$ give topologically
isomorphic full groups but no uniform measure relation isomorphism
with the designated markings.
\end{example}

\begin{rem}\label{rem:uniform-Dye-abstract-scope}
The Boolean reconstruction in
Lemma~\ref{lem:uniform-full-group-Boolean-reconstruction} applies to
abstract group isomorphisms, without amenability or fullness.
The subsequent recovery of the diagonal topology uses continuity,
and the recovery of the designated measures above uses the full
invariant family and its amenable AF model.
Theorem~\ref{thm:uniform-Dye-reconstruction} does not assert automatic
continuity of arbitrary abstract full-group isomorphisms. For an
arbitrary designated subfamily, one must retain additional measure
marking; Example~\ref{ex:full-group-topology-loses-designated-weights}
already prevents reconstruction from the topological group alone.
\end{rem}

\subsection{The simple full-trace model interface}

Let $(X,R,K,D_K,\mathcal S)$ be a uniform measure relation and use its
untwisted completed relation algebra.  Put
\[
 D=D_K=\mathrm{C}(X)^K,\qquad M=M_R,\qquad G=[R]^u_{K,\mathcal S}.
\]
Assume that $(E_0,B)$ is the standard diagonal pair of a simple unital
infinite-dimensional AF algebra, and identify $L=\T(B)$ with the
central probability measures on its Bratteli path space $Y$.  Write
\[
 E=(E_0)^L,\qquad N=B^L.
\]
Our explicit hypothesis is a trace-preserving isomorphism of completed
pairs, marked by an affine homeomorphism $H:K\to L$:
\begin{equation}\tag{AFM}\label{eq:uniform-AF-model-hypothesis}
 \Psi:(D,M)\longrightarrow(E,N),\qquad
 \tau_{H\mu}(\Psi(a))=\tau_\mu(a)
 \quad(\mu\in K,\ a\in M).
\end{equation}
\stepcounter{equation}
Here the completions are bounded uniform tracial completions, and on $M$
we use the designated traces constructed from $K$.  Let
$B=\overline{\bigcup_n B_n}$,
$B_n=\bigoplus_{v\in V_n}\mathrm{M}_{k_{n,v}}$, be a coordinate Bratteli
presentation, and let $R_B$ denote tail equivalence on $Y$.

\begin{cor}[Classification under the simple full-trace model hypothesis]
\label{cor:AFM-uniform-Cartan-classification}
Suppose that two regular uniform measure relations admit the marked
completed-pair models \eqref{eq:uniform-AF-model-hypothesis}. Every affine
homeomorphism between their designated families is implemented by a
marked isomorphism of completed Cartan pairs and a uniform measure
relation isomorphism. The isomorphism intertwines the canonical
expectations and is isometric for the uniform $2$-norms.
\end{cor}

\begin{proof}
Compose the two model markings with the specified affine homeomorphism.
The AF-pair classification theorem~\cite[Theorem~5.3]{TC},
or Theorem~\ref{thm:diffuse-face-AF-classification} with full trace faces,
implements the resulting correspondence between the simple AF models.
Composition with the two model isomorphisms gives the marked pair
isomorphism. Each model isomorphism intertwines the expectations:
for $x$ in the ambient algebra, the diagonal element $\mathbb E(x)$
is uniquely determined by
$\tau(d^*\mathbb E(x))=\tau(d^*x)$ for every diagonal $d$ and every
designated trace $\tau$. Indeed, testing the difference of two candidates
against itself gives zero uniform $2$-norm.
Theorem~\ref{thm:uniform-relation-pair-isomorphism} gives
its spatial relation isomorphism.
\end{proof}

\begin{prop}\label{prop:AF-tail-spatialization}
Under \eqref{eq:uniform-AF-model-hypothesis}, the completed-pair
isomorphism induces a Borel isomorphism
\[
 \theta:X_0\longrightarrow Y_0
\]
between invariant conull Borel sets, carrying $K$ onto $L$ by the prescribed
affine homeomorphism, such that
\[
 (x,x')\in R\quad\Longleftrightarrow\quad
 (\theta x,\theta x')\in R_B.
\]
There are increasing finite subrelations $F_n\subset R|_{X_0}$,
with uniformly bounded class sizes at each fixed stage, for which
\[
 \sup_{\mu\in K}\mu\{x:(x,Tx)\notin F_n\}\longrightarrow0
 \qquad(T\in G).
\]
Their finite charts are admissible restrictions of elements of the uniform
measure full group.
\end{prop}

\begin{proof}
Theorem~\ref{thm:orbit-spatial-diagonal} applies to these two commutative
completions and gives an initial conull Borel bijection $\theta$, with
$\Psi(d)=d\circ\theta^{-1}$.  This use of that theorem involves only the
diagonals; no spatialization theorem for an AF normalizer is being assumed.

Let $P_n=\prod_{v\in V_n}\operatorname{Sym}(k_{n,v})$ be the coordinate
permutation group, regarded both as permutation unitaries in $B_n$ and as
homeomorphisms of $Y$ which replace a finite prefix and retain its tail.
The coordinate inclusions make $P_n\subset P_{n+1}$.  Set $P=\bigcup_nP_n$.
By Proposition~\ref{prop:uniform-full-group-normalizers}, each
$\Psi^{-1}(p)$ induces $T_p\in G$.  Covariance on the diagonal,
tested on a countable family of continuous functions separating points,
gives
\begin{equation}\label{eq:AF-permutation-spatial}
 \theta(T_px)=p\theta(x)
\end{equation}
outside a fixed universally null Borel set, for every $p\in P$.
The normalizer quotient also gives $T_pT_q=T_{pq}$ modulo universally null
sets.  There are only countably many such identities.

By Corollary~\ref{cor:uniform-relation-countable-group-presentation},
choose countably many $T_i\in G$ generating $R$ on a common invariant
conull set, each extending an original atlas chart.  For each $i$,
Lemma~\ref{lem:AF-normalizer-density}, applied to the image of its
canonical graph unitary $v_{T_i}$, gives $p_{i,r}\in P$ and diagonal
phases $a_{i,r}$ such that
\[
 \sum_r\|v_{T_i}-a_{i,r}\Psi^{-1}(p_{i,r})\|_{2,K}^2<\infty.
\]
The quotient metric in
Proposition~\ref{prop:uniform-full-group-normalizers} therefore implies
$\sum_r d_K(T_i,T_{p_{i,r}})<\infty$.
The Borel--Cantelli argument applies simultaneously to every $\mu\in K$:
the fixed Borel set on which $T_ix\ne T_{p_{i,r}}x$ infinitely often has
measure zero for every such $\mu$.  Thus, for every $i$, its graph is
contained almost everywhere in the relation generated by the $T_p$.

Here are the common-null-set choices explicitly.  Choose Borel
representatives of the countably many $T_p$, defined on common conull
$R$-invariant sets and extended by the identity elsewhere.  Collect
the complements of the initial spatialization sets, the exceptional sets
for \eqref{eq:AF-permutation-spatial}, the group identities, and the
preceding Borel--Cantelli conclusions.  Remove their saturation under the
countable group generated by the $T_i$ and these representatives.
Every generator preserves every $\mu\in K$, so the resulting exceptional
set is still universally null.  On the remaining invariant Borel set
$X_0$ all the identities hold pointwise.  Its image $Y_0=\theta(X_0)$ is
Borel, conull for $L$, and $P$-invariant by
\eqref{eq:AF-permutation-spatial}.  Since each $T_p$ lies in $R$ and every
$T_i$-edge lies in a $P$-orbit, $R|_{X_0}$ is
exactly the relation generated by the $T_p$.  The $P$-orbit relation on
$Y$ is exactly $R_B$: two paths with a common tail differ by a permutation
of their finite prefixes at some level.  This proves the asserted spatial
identification.

Let $R_{B,n}$ mean equality of all edges after level $n$, and pull its
restriction to $Y_0$ back to $F_n$.  Its class at a path ending at $v$ at
level $n$ has size $k_{n,v}$, so the stage has bounded finite classes.
Its prefix-substitution charts are restrictions of $p\in P_n$ to clopen
prefix cylinders.  Under $\theta$ their domain projections belong to
$D$, and their maps are restrictions of $T_p$; hence they are
admissible.  Finally, apply the same phase approximation argument to
$\Psi(v_T)$ for any $T\in G$.  This proves that the $T_p$ are dense in
$G$.  For $T$ and
$\varepsilon>0$, choose $p\in P_m$ with $d_K(T,T_p)<\varepsilon$.
For $n\geq m$, agreement with $T_p$ implies $(x,Tx)\in F_n$, giving the
claimed uniform estimate.
\end{proof}

\begin{cor}\label{cor:uniform-AF-charts-from-model}
Under \eqref{eq:uniform-AF-model-hypothesis}, the original uniform measure
relation is uniformly AF hyperfinite, in the sense of
Definition~\ref{defn:uniform-AF-hyperfinite-relation}.  Its finite matrix
charts can be chosen as the transports of the standard AF coordinates,
with all finite diagonals exhausting $D$ in the bounded uniform $2$-norm.
Moreover $\theta$ in Proposition~\ref{prop:AF-tail-spatialization} is a
uniform measure isomorphism between the original regular atlas
$\mathcal S$ of $R$ and the finite-prefix atlas of $R_B$.
\end{cor}

\begin{proof}
For each length-$n$ prefix $i$ ending at $v$, let $Q_{n,v,i}$ be its
clopen cylinder and put $Y_{n,v,i}=\theta^{-1}(Q_{n,v,i}\cap Y_0)$.
These give a finite partition of $X_0$, with
$1_{Y_{n,v,i}}=\Psi^{-1}(1_{Q_{n,v,i}})\in D$.  A fixed permutation in
$P_n$ sending prefix $1$ to prefix $i$ transports to an admissible map
$Y_{n,v,1}\to Y_{n,v,i}$.  Their matrix charts give precisely $F_n$.
At the next stage each prefix cylinder splits into longer cylinders;
each old substitution retains the additional suffix, and therefore
splits into the corresponding finite disjoint union of new matrix
charts.  This is the required coordinate refinement.

Write $E_n$ for the finite diagonal of $B_n$ and $D_n=\Psi^{-1}(E_n)$.
The union of the $E_n$ is norm dense in $E_0=\mathrm{C}(Y)$, and $E_0$ is bounded
uniformly $2$-norm dense in $E$.  Consequently, for each $b\in E$ there
are $b_n\in E_n$ with $\|b_n\|\leq\|b\|$ and
$\|b-b_n\|_{2,L}\to0$.  One way to retain the bound is to radially clip
bounded continuous approximants to the disk of radius $\|b\|$, then
approximate them by their values on the refining clopen cylinders.
Their transports prove the asserted bounded density of $D_n$ in $D$.
Uniform escape and the common conull exhaustion were proved in
Proposition~\ref{prop:AF-tail-spatialization}; hence all the uniform AF
hyperfiniteness conditions hold.

We also verify the atlas regularity of $\theta$.  For each $T_i$ in the
chosen global graphing, put $S_i=\theta T_i\theta^{-1}$.  The preceding
phase approximation gives permutations $p_r\in P$ such that
\[
 d_L(S_i,p_r)\longrightarrow0.
\]
For every fixed $p\in P$, their coincidence indicators satisfy
\[
 \bigl\|1_{\{S_i=p\}}-1_{\{p_r=p\}}\bigr\|_{2,L}^2
 \leq d_L(S_i,p_r)\longrightarrow0.
\]
The indicators on the right are clopen finite-prefix functions, so
$1_{\{S_i=p\}}\in E$.  Restricting $T_i$ to the domain of the original
chart it extends multiplies these projections by the transported domain
projection, which belongs to $E$.  The inverse transport sends $p$ to
$T_p\in G$, whose coincidence indicators with each original chart lie
in $D$.  Domain restrictions to prefix cylinders retain this property.
Together with the diagonal pullback isomorphism, these are the forward
and inverse regularity conditions for the atlases.  The same conclusion
for their full pseudogroups follows by decomposing an admissible map
into its atlas charts with diagonal domain projections: the uniform
tail estimate makes the finite decompositions converge on both sides.
\end{proof}

\begin{prop}\label{prop:AF-uniform-Vershik}
The simple AF model can be ordered so that its path space carries a
Cantor minimal homeomorphism $V$ with
\[
 \mathcal M_V(Y)=L.
\]
There are $c_n\in P_n$ such that, with $m_n=\min_v k_{n,v}$,
\begin{equation}\label{eq:uniform-Vershik-estimate}
 \sup_{\nu\in L}\nu\{c_n\ne V\}\leq m_n^{-1}\longrightarrow0.
\end{equation}
The relation $R_V$ agrees with $R_B$ outside a countable universally null
invariant set.  Moreover the canonical unitary for $V$ is the uniform
$2$-norm limit of the phase-free unitaries $c_n$ in $N$, and
\[
 \bigl(\mathrm{C}(Y)^L,(\mathrm{C}(Y)\rtimes_{r,V}\mathbb Z)^L\bigr)
 \cong\bigl((E_0)^L,B^L\bigr)
\]
as completed Cartan pairs, with the identity on the diagonal.
\end{prop}

\begin{proof}
First telescope to make all incidence matrices strictly positive, using
simplicity.  A further telescoping makes every multiplicity at least two:
passing through a level with at least two vertices doubles every path
count, and if there is eventually only one vertex, infinite dimensionality
provides factors at least two arbitrarily far out.  The dimension recurrence
then gives $m_{n+1}\geq2m_n$, so $m_n\to\infty$.
At level $n$, choose the minimal incoming edge at
every vertex to have one prescribed source $v_{n-1}$, and the maximal
incoming edge to have one prescribed source $w_{n-1}$.  When these sources
coincide, use distinct edges.  Order the remaining edges arbitrarily.
An infinite minimal path is forced, by the next incoming edge, to pass
through $v_n$ at every level; hence it is unique.  The analogous maximal
path is unique as well.  This is a proper order.  Its successor map,
completed by sending $y_{\max}$ to $y_{\min}$, is a homeomorphism:
near $y_{\max}$ the first nonmaximal edge tends to infinity, and the
successor's arbitrarily long minimal prefix converges to $y_{\min}$;
the inverse argument is identical.  Simplicity makes every tail class
dense, and the successor orbit contains its tail class.  Thus $V$ is
minimal.  This is also the standard proper-order construction recorded
in \cite[Proposition~1.1]{GPSA}.

Let $c_n$ cyclically permute the linearly ordered paths from the root to
each $v\in V_n$.  It agrees with $V$ unless the length-$n$ prefix is
maximal.  Every central measure assigns the same mass to all
$k_{n,v}$ prefix cylinders ending at $v$.  Consequently the union of the
maximal-prefix cylinders has mass
$\sum_v\nu(z_{n,v})/k_{n,v}\leq m_n^{-1}$, proving
\eqref{eq:uniform-Vershik-estimate}.  In particular every singleton has
mass at most $m_n^{-1}$ for every $n$, so all central measures are
nonatomic.  The set
\[
 Z=R_B[y_{\min}]\cup R_B[y_{\max}]
\]
is countable and universally null.  It is invariant under $R_B$ and $V$.
On its complement every successor and predecessor retains a tail.  If
two paths have the same tail beyond level $n$, their finite prefixes
have ranks $i,j$ in the same ordered set, and $V^{j-i}$ sends the first
path to the second: the intermediate ranks remain between $i$ and $j$,
so no cyclic wrap is needed.  This proves $R_V=R_B$ off $Z$.
Central measures preserve every $c_n$, and their pointwise limit is $V$
off $y_{\max}$; hence they preserve $V$.  Conversely every finite prefix
permutation is a Borel piecewise power of $V$, by the same rank argument,
so a $V$-invariant measure is central.  This proves
$\mathcal M_V(Y)=L$.

For finite coordinate permutations the diagonal expectation satisfies
$\mathbb{E}(c_n^*c_r)=1_{\{c_n=c_r\}}$ after passage to a common stage.
It follows from \eqref{eq:uniform-Vershik-estimate} that their unitaries
are uniform $2$-norm Cauchy.  Their limit $w\in N$ is unitary, and
$wfw^*=f\circ V^{-1}$ for $f\in \mathrm{C}(Y)$; this follows by the same estimate
and then extends to $E$ by continuity.  Thus $w$ normalizes
$E$ and is the phase-free graph unitary for $V$.

Covariance gives a homomorphism
$\iota:\mathrm{C}(Y)\rtimes_{r,V}\mathbb Z\to N$, sending its canonical unitary to
$w$.  For $k\ne0$, covariance and the diagonal expectation imply
\[
 f\mathbb{E}(w^k)=\mathbb{E}(w^k)(f\circ V^{-k})
 \qquad(f\in \mathrm{C}(Y)).
\]
Since the Cantor minimal map has no periodic points, a finite family of
clopen functions separates each $y$ from $V^ky$; thus
$\mathbb{E}(w^k)=0$.  The homomorphism therefore preserves every canonical
trace.  It is injective: an invariant measure has full support, and its
trace is faithful because the reduced crossed-product expectation is
faithful; full and reduced crossed products coincide for $\mathbb Z$.

If $e_{ij}^{(n,v)}$ are the matrix units indexed by the ordered finite
paths ending at $v$, then for every $r\geq n$,
\[
 c_r^{\,i-j}e_{jj}^{(n,v)}=e_{ij}^{(n,v)}.
\]
Indeed these iterates stay within the indicated length-$n$ prefixes.
Taking uniform $2$-norm limits gives
$w^{i-j}e_{jj}^{(n,v)}=e_{ij}^{(n,v)}$.  Hence
$B_n\subset\iota(\mathrm{C}(Y)\rtimes_{r,V}\mathbb Z)$ for every $n$.
The image is therefore uniformly $2$-norm dense in $N$, and the
trace-preserving embedding extends to the claimed completed-pair
isomorphism.
\end{proof}

\begin{cor}[Uniform CFW consequence]
\label{cor:uniform-CFW-Z-generator}
Under \eqref{eq:uniform-AF-model-hypothesis}, there are $T\in G$ and one
invariant conull Borel set $X_1$ such that
\[
 R|_{X_1}=R_T|_{X_1}.
\]
The original relation is uniform measure isomorphic to the orbit
relation of the Cantor minimal $\mathbb Z$-system $(Y,V)$ above, with
the prescribed affine identification $H:K\to\mathcal M_V(Y)$ of
designated measure simplexes.  Its completed untwisted Cartan pair is
also the completion of $(\mathrm{C}(Y),\mathrm{C}(Y)\rtimes_{r,V}\mathbb Z)$.
\end{cor}

\begin{proof}
Compose $\Psi$ in \eqref{eq:uniform-AF-model-hypothesis} with the
inverse of the isomorphism in
Proposition~\ref{prop:AF-uniform-Vershik}.  The diagonal spatialization
is the map $\theta$ of Proposition~\ref{prop:AF-tail-spatialization},
after restriction to common conull sets.  The image of $w$ in the
original completed algebra induces, by
Proposition~\ref{prop:uniform-full-group-normalizers},
$T\in G$, and diagonal covariance gives
$\theta T=V\theta$ almost everywhere simultaneously for all measures.
Remove the inverse image of $Z$, and then the saturation of the
countably many exceptional sets under the global graphing $T_i$ and $T$.
These maps preserve every measure; the resulting $X_1$ is $R$-invariant
and conull.  The two preceding propositions now give $R_T=R$ on $X_1$.

We verify uniform relation regularity as well.  Pullbacks identify the
two completed diagonals.  By
Lemma~\ref{lem:uniform-relation-algebra-action-case}, the completed
$\mathbb Z$-crossed product is the untwisted relation algebra for $R_V$.
The image of $v_{T_i}$ is a normalizer in this algebra, so
Proposition~\ref{prop:uniform-full-group-normalizers} puts all its
coincidence indicators with the powers of $V$ in $E$.  Restricting to the
domain of each original chart gives the forward atlas coefficients.
Conversely, the inverse image of every canonical unitary $w^n$ induces
an admissible element of $G$, so its coincidence indicators with every
original chart belong to $D$.  These give the inverse atlas
coefficients.  Thus $\theta$ is a uniform measure relation isomorphism.
Trace preservation gives
$\theta_*\mu=H\mu$ for every $\mu\in K$.  Finally the completed-pair
assertion is the composition already constructed.
\end{proof}

\begin{cor}[Free-action specialization]
\label{cor:AF-CFW-free-action}
Let a countable group $\Gamma$ act freely and continuously on a compact metrizable space
$X$, let $K=M_\Gamma(X)$ be nonempty, and suppose
\[
 \bigl(\mathrm{C}(X)^K,(\mathrm{C}(X)\rtimes_r\Gamma)^K\bigr)
 \cong\bigl((E_0)^{\T(B)},B^{\T(B)}\bigr)
\]
is a trace-preserving marked isomorphism with a standard simple unital
infinite-dimensional AF pair.  Then the system is uniformly measure
orbit equivalent to a Cantor minimal $\mathbb Z$-system with the prescribed
invariant-measure simplex.  It has an admissible single generator for
its orbit relation on one universally conull invariant set.
\end{cor}

\begin{proof}
Apply Lemma~\ref{lem:uniform-relation-algebra-action-case} and then
Corollary~\ref{cor:uniform-CFW-Z-generator}.  For the translation atlases,
uniform measure relation isomorphism is exactly uniform measure orbit
equivalence: the mixed graph coefficients are the forward and inverse
cocycle-fiber projections, and the diagonal pullbacks give the observable
conditions.
\end{proof}

The AF-model theorem~\cite[Theorem~4.1]{TC} supplies a simple model carrying the full trace
simplex in the free minimal action setting.  The following consequence
applies that theorem to the relation framework developed here.

\begin{cor}[Free minimal actions with small boundaries]\label{cor:SBP-uniform-CFW}
Let $\Gamma$ be countably infinite and amenable, and let its action on
$X$ be free and minimal with the small boundary property.  Then its
uniform measure relation is uniformly AF hyperfinite, and all the
conclusions of Corollary~\ref{cor:AF-CFW-free-action} hold.
\end{cor}

\begin{proof}
The completed AF-model theorem~\cite[Theorem~4.1]{TC} supplies
\eqref{eq:uniform-AF-model-hypothesis}, with a standard simple unital
infinite-dimensional AF pair and the prescribed trace correspondence.
Apply Corollaries~\ref{cor:uniform-AF-charts-from-model} and
\ref{cor:AF-CFW-free-action}, using the action specialization of the
untwisted relation algebra.
\end{proof}

\section{Local regular finite approximation}
\label{sec:local-regular-cuts}

The AF hypothesis can be expressed without requiring an initially nested
system of matrix charts.  Throughout this section the relation and atlas
are those of Definition~\ref{defn:regular-relation-atlas}.

\begin{defn}\label{defn:local-regular-cuts}
The relation has \emph{local regular finite approximation} if, for every
finite $\mathcal F\subset\mathscr P$ and $\varepsilon>0$, there is a
Borel subequivalence relation $S\subset R$ with uniformly bounded finite
classes and finitely many admissible partial generators such that
\[
 \sup_{\mu\in K}\sum_{s\in\mathcal F}
 \mu\{x\in\operatorname{dom}s:(x,sx)\notin S\}<\varepsilon.
 \tag{RF}\label{eq:local-regular-cuts}
\]
\stepcounter{equation}
Relations and generators may be restricted to one invariant universally
conull set.  No nesting or choice of transversals is part of this condition.
\end{defn}

\begin{lem}[Regular finite transversals]\label{lem:regular-finite-transversals}
Assume $K$-SBP.  Every bounded finite subrelation generated by finitely
many admissible partial maps has finite $D_K$-regular matrix charts.
Its charts can be refined to incorporate any prescribed finite partition
of regular projections.  For nested such relations, a new chart system
can be chosen to refine a prescribed old chart system in coordinates.
\end{lem}
\begin{proof}
Use the nested finite Borel partitions of mesh tending to zero from
Proposition~\ref{prop:relative-SBP-finite-diagonals}, and order points by
their lexicographic partition codes.  On one common conull set these
codes separate points.  If $\varphi,\psi$ are admissible partial maps,
let $a_m$ be the projection, on their common domain, where their first
$m$ codes agree.  Corner pullbacks make $a_m$ regular.  These projections
decrease to $c(\varphi,\psi)$, which is regular by
Corollary~\ref{cor:regular-relation-equalizers}.  The continuous functions
$\mu\mapsto\mu(a_m-c(\varphi,\psi))$ decrease to zero on compact $K$.
Dini's theorem gives uniform convergence.  Finite-code strict comparisons
are regular projections, and their difference from the exact lexicographic
comparison is supported in $a_m-c(\varphi,\psi)$.  Thus the exact comparison
is a projection of $D_K$.

Let $N$ bound the class sizes of $S$.  Words of length at most $N-1$ in
its finitely many generators and inverses, with the identity included,
enumerate each class: a shortest path to a class member visits no vertex
twice.  There are finitely many such partial words.  Their supports,
equalizers and order comparisons are regular.  Finite Boolean formulas
therefore select the least member of each class, detect its class size,
and rank the distinct word outputs.  On the regular transversal of
classes of size $k$, the map sending the least member to its $i$th ranked
member is a finite cut-and-paste of these words.  It is admissible, and
the rank-$i$ ranges have regular projections.  This gives finitely many
matrix charts with $k\le N$.

For refinement, label every member of a new class by its old coordinate
cell and by the prescribed projection partition.  Also record which
pairs of new coordinates are connected by each old matrix-chart map.
There are finitely many labels and tests.  Their projections on the new
transversal are regular, by corner pullback and coincidence regularity.
Subdivide that transversal into the resulting finite patterns and
transport each subdivision to every coordinate.  Each new coordinate
cell then lies in one old cell, and every old matrix map is a finite
disjoint union of new coordinate maps.  The same construction includes
the prescribed projection partition without changing $S$.
\end{proof}

\begin{lem}[Nesting with small loss]\label{lem:regular-cut-nesting}
Let $H,G$ be bounded finite regularly generated subrelations, and suppose
$H$ has class-size bound $M$.  Choose finitely many admissible arrow maps
$h_1,\ldots,h_q$ covering $H$, including its identity arrows.  Set
\[
 E=\bigcup_j\{x\in\operatorname{dom}h_j:(x,h_jx)\notin G\},
 \qquad B=[E]_H,\qquad C=X\setminus B.
\]
Then these sets are regular, $B$ is $H$-invariant,
$\mu(B)\le M\mu(E)$ for every $\mu\in K$, and
\[
 F=G|_C\ \sqcup\ H|_B
\]
is bounded finite, regularly generated, and contains $H$.  For $s\in\mathscr P$,
\[
 \mu\{(x,sx)\notin F\}
 \le\mu\{(x,sx)\notin G\}+2\mu(B),
\]
where escape sets are restricted to $\operatorname{dom}s$.
\end{lem}
\begin{proof}
Membership of an admissible arrow in $G$ is a finite union of coincidence
projections: finitely many words enumerate each bounded $G$-class.
Thus $E$ is regular.  The finite $H$-arrow cover shows that its saturation
$B$ is regular.  On an $H$-class of size $k\le M$, the indicator of its
meeting $E$ is at most the sum of the $E$-indicators of its members.
Invariance, or integration in finite matrix charts, gives
$\mu(B)\le M\mu(E)$.  Every $H$-class in $C$ lies in a single $G$-class;
therefore $H\subset F$.  Restricting finite arrow covers of $G,H$ to
$C,B$ gives finitely many admissible generators, and class sizes remain
bounded by the larger original bound.  A captured $G$-arrow can fail to
be an $F$-arrow only when its source or target lies in $B$.  Invariance
under $s$ bounds their total measure by $2\mu(B)$.
\end{proof}

\begin{thm}[A local criterion for AF charts]\label{thm:local-cuts-AF}
Assume $K$-SBP.  Condition \eqref{eq:local-regular-cuts} is equivalent to
uniform AF hyperfiniteness.  Consequently it supplies the completed AF
pair and spatial tail-relation model of
Theorem~\ref{thm:uniform-AF-reconstruction}, with its designated trace family.
\end{thm}
\begin{proof}
An AF exhaustion captures any finite collection of admissible arrows
uniformly, so it implies (RF).  Conversely enumerate the atlas as $(s_j)$
and fix positive summable $\eta_n\to0$.  Starting with the identity relation,
suppose $H_{n-1}$ has been constructed.  In (RF) test its finite arrow
cover together with $s_1,\ldots,s_n$.  Choose the approximation so accurate
that Lemma~\ref{lem:regular-cut-nesting} produces $H_n\supset H_{n-1}$
whose summed escape for these $n$ tests is below $\eta_n$.  This is possible
because the old class-size and arrow-count bounds are fixed before making
the new choice.

For every fixed $s_j$, the escape probabilities are summable uniformly
over $K$.  The Borel set of points escaping infinitely often is null for
every $\mu\in K$, by Borel--Cantelli.  Remove the saturation of the union
of these countably many exceptional sets.  On the remaining invariant
conull set every atlas arrow belongs to $\bigcup_n H_n$, and therefore
$R=\bigcup_n H_n$.  Escape for an arbitrary admissible partial map follows
by its first-chart decomposition and the uniform tail estimate
\eqref{eq:regular-relation-uniform-tails}.

Apply Lemma~\ref{lem:regular-finite-transversals} successively, refining the
old chart system and including the $n$th SBP partition at stage $n$.
The resulting diagonals contain the SBP step algebras, so they are bounded
uniformly dense in $D_K$.  All conditions of
Definition~\ref{defn:uniform-AF-hyperfinite-relation} now hold.
\end{proof}

\begin{cor}\label{cor:bounded-regular-exhaustion-AF}
If $K$-SBP holds and $R$ has an increasing common-conull exhaustion by
bounded finite subrelations with finitely many admissible generators at
each level, then $R$ is uniformly AF hyperfinite.  A separate uniform
escape assumption is unnecessary in this case.
\end{cor}
\begin{proof}
For each admissible $s$, its escape projection at a fixed level is regular
by the finite-word argument above.  These projections decrease to zero
universally almost everywhere.  Their continuous measure profiles decrease
to zero on compact $K$, so Dini gives uniform escape.  Thus (RF) holds,
and Theorem~\ref{thm:local-cuts-AF} applies.
\end{proof}

\begin{prop}\label{prop:Reiter-measurewise-equivalence}
Regular uniform Reiter amenability is equivalent to amenability for
every $\mu\in K$.
\end{prop}
\begin{proof}
The reverse implication is
Corollary~\ref{cor:measurewise-amenable-uniform-Reiter}.  For the forward
implication, pass to a subsequence whose uniform integrated errors on the
first $n$ atlas charts have sum below $2^{-n}$.  For each $\mu$, summability
gives pointwise asymptotic invariance for all charts almost everywhere.
After removing the saturation of the exceptional set, the kernels are
ordinary measurable Reiter witnesses on that relation.  The same Borel
exceptional set is null for every designated measure.  The usual measured
Reiter characterization gives amenability.
\end{proof}

\begin{rem}\label{rem:deterministic-cut-mechanism}
Proposition~\ref{prop:Reiter-measurewise-equivalence} requires no
amenability assumption stronger than amenability for each designated
measure.  Relative SBP and
Theorem~\ref{thm:regular-Reiter-rounding} convert its regular kernels
into deterministic regular cuts.  Convex averaging alone would not
do this: finite subrelations are not a convex family, and their
joins can have infinite classes.  The finite-kernel construction
and local derandomization provide the required step; the chart
and nesting lemmas then produce compatible matrix coordinates.
Aperiodicity is needed only for the diffuse dynamical realization.
\end{rem}

\section{Aperiodic AF models and dynamical realization}
\label{sec:general-AF-generator}

Simplicity in Section~\ref{subsec:uniform-CFW} supplies a minimal Cantor
system and its full invariant-measure simplex.  It is unnecessary for
the admissible single-generator conclusion.

\begin{prop}\label{prop:general-AF-successor}
Let $B$ be a standard unital AF algebra with path space $Y$ and compact
designated central-measure family $L$.  Suppose its tail relation has
infinite classes on one invariant $L$-conull set.  Then it has a single
uniformly admissible generator on an invariant $L$-conull set.
\end{prop}
\begin{proof}
Order the incoming edges arbitrarily.  Let $k_{n,v}$ count paths from the
root to $v$, and cyclically permute the ordered level-$n$ prefixes with
the same endpoint by $c_n$.  The maximal-prefix cylinder union $A_n$
decreases with $n$, and for a central measure
\[
 r_n(\nu)=\nu(A_n)=\sum_v\frac{\nu(z_{n,v})}{k_{n,v}}
 =\int_Y\frac{1}{k_{n,v_n(y)}}\,d\nu(y).
\]
The class sizes $k_{n,v_n(y)}$ increase to infinity on the stipulated
conull set.  Hence $r_n(\nu)\downarrow0$ for each $\nu\in L$.
The $r_n$ are continuous on compact $L$, so
$\sup_Lr_n\to0$.  Minimal-prefix cylinder unions have the same measures.
Thus the sets of all-maximal and all-minimal paths are universally null,
as are their tail saturations.

On the complement of these saturations the ordinary successor and
predecessor maps are inverse Borel bijections and retain tails.  Denote
the successor by $V$.  If two paths share a tail beyond level $n$, their
prefix ranks are $i,j$ in one finite ordered list.  Successor or predecessor
steps between these ranks do not wrap around that list, so $V^{j-i}$
connects the paths.  Its orbits are therefore exactly tail equivalence.
Moreover $c_n$ agrees with $V$ outside $A_n$ on this conull set, and
\[
 d_L(c_n,V)\le\sup_{\nu\in L}r_n(\nu)\longrightarrow0.
\]
The finite-prefix permutations are admissible.  Completeness of the
uniform full group proves admissibility of $V$.  Equivalently, their
phase-free graph unitaries converge in uniform $2$-norm to its graph
unitary.  A proper ordering and a continuous Vershik map on the original
path space are not asserted.
\end{proof}

\begin{thm}\label{thm:local-cuts-dynamical-model}
Let the regular relation have $K$-SBP and (RF), and suppose $R$ has infinite
classes on one invariant $K$-conull set.  There exist an admissible $T$
with $R_T=R$ on such a set, a Cantor space $Z$, a homeomorphism $W$ of $Z$,
and a compact designated invariant family $K_Z$ affine-homeomorphic to $K$
such that the original relation is uniform measure isomorphic to $R_W$.
The map $W$ is topologically free and aperiodic for every measure in $K_Z$,
and
\[
 (D_K,M_R)\cong
 \bigl(\mathrm{C}(Z)^{K_Z},(\mathrm{C}(Z)\rtimes_{r,W}\mathbb Z)^{K_Z}\bigr).
\]
Neither minimality, everywhere freeness, nor $K_Z=M_W(Z)$ is asserted.
\end{thm}
\begin{proof}
Theorem~\ref{thm:local-cuts-AF} and forward reconstruction give an AF tail
model with designated measures.  Transport the generator of
Proposition~\ref{prop:general-AF-successor}; its regular coincidence
projections with the original charts are supplied by the uniform relation
isomorphism.  This gives $T\in G$ with $R_T=R$ on a common invariant
conull set.

Take the countable SBP projection partitions and close their Boolean
algebra under pullback by every power of $T$.  Its norm closure is a
separable unital commutative AF algebra $\mathrm{C}(Z)\subset D_K$, contains the
continuous diagonal, and is bounded uniformly dense in $D_K$.  Its spectrum
$Z$ is compact metrizable and zero dimensional.  Pullback by $T$ induces
a homeomorphism $W$.  Restriction identifies $K$ homeomorphically with a
compact invariant family $K_Z$: the continuous diagonal separates the
original measures.  The designated measures are jointly faithful because
we use the common-null quotient.  Bounded completion and
Theorem~\ref{thm:orbit-spatial-diagonal} give a Borel conjugacy of $T,W$
on common conull sets.

All measures in $K$ are nonatomic: invariance makes singleton masses
constant on each infinite class.  Hence the transported measures are
nonatomic.  Joint faithfulness then rules out isolated points of $Z$,
so $Z$ is Cantor.  The sets $\operatorname{Fix}(W^m)$, $m\ne0$, are null
for every designated measure, since $T$ is aperiodic there.  They have
empty interior by joint faithfulness; this is topological freeness.

Covariance with the canonical graph unitary $v_T$ gives a homomorphism
$\mathrm{C}(Z)\rtimes_{r,W}\mathbb Z\to M_R$; full and reduced crossed products
coincide for $\mathbb Z$.  For $m\ne0$ the diagonal expectation of
$v_T^m$ is the fixed-point projection, which is zero.  Thus the map
preserves every designated canonical trace.  Joint faithfulness on the
diagonal and the faithful crossed-product expectation give injectivity.
For an original chart $s$, disjointize the regular projections
$1_{\{sx=T^m x\}}$, $m\in\mathbb Z$.  They exhaust its domain because
$R_T=R$.  Compactness and Dini give uniformly vanishing tails.  Consequently
the finite sums of $v_T^m$ times these source projections converge uniformly
in $2$-norm to $v_s$.  The image is bounded uniformly dense in $M_R$ and
extends to the asserted completed-pair isomorphism.  The same coefficients,
and those of powers of $T$ relative to the original atlas, verify forward
and inverse uniform relation regularity.
\end{proof}

\begin{rem}\label{rem:full-simplex-interface}
The simple AF hypothesis (AFM) remains useful: it supplies the minimal
Cantor realization and the entire invariant-measure simplex, as proved
in Section~\ref{subsec:uniform-CFW}.  Local reconstruction alone produces
a possibly nonsimple AF algebra with $L\subseteq \T(B)$. For an
aperiodic relation carrying its full compact invariant family,
Lemma~\ref{lem:full-family-AF-face} and
Theorem~\ref{thm:diffuse-face-AF-classification} give a simple full-trace
replacement and simplex classification. Thus (AFM) is automatic under
the hypotheses of Theorem~\ref{thm:uniform-Cartan-classification}.
For an arbitrary designated family, the face and diffuseness conditions
need not hold, and a compact convex family need not be a Choquet
simplex. Remark~\ref{rem:uniform-classification-scope} explains why the
designated family alone cannot classify all such models.
\end{rem}

\section{\'Etale groupoids and group-action presentations}
\label{sec:etale-presentations}

Let $\mathcal G$ be second countable, locally compact, Hausdorff and
\'etale, with compact metrizable unit space $X$.  Let $K$ be a nonempty
compact convex family of invariant probability measures.  Write
$R_{\mathcal G}=\{(r\gamma,s\gamma):\gamma\in\mathcal G\}$ for its
countable Borel orbit relation.  ``Principal'' means trivial isotropy
at every unit; ``ample'' means that compact-open bisections form a basis.
Example~\ref{ex:regular-atlas-no-original-etale-model} shows that a
regular relation need not admit an \'etale presentation on its original
unit space. The model changes below are therefore substantive.

\begin{thm}[A continuous group presentation]\label{thm:ample-group-presentation}
If $\mathcal G$ is principal and ample, there is a countable
$\Gamma\le\operatorname{Homeo}(X)$ with
\[
 R_\Gamma=R_{\mathcal G},\qquad M_\Gamma(X)=M_{\mathcal G}(X).
\]
The identity on $X$ is a uniform measure relation isomorphism for every
designated $K$, using the compact-open bisection atlas and the group-word
atlas.  If $\mathcal G$ is minimal, the action is minimal.  No amenability
assumption is needed.  The action is not asserted to be free.
\end{thm}
\begin{proof}
Every nonunit $\gamma$ has different source and range, by principality.
Choose a compact-open bisection $U$ containing it with disjoint clopen
source and range.  Its partial homeomorphism $\beta_U$ extends to the
involution acting by $\beta_U$ on $s(U)$, by $\beta_U^{-1}$ on $r(U)$,
and by the identity elsewhere.  The corresponding full bisection is
\[
 W_U=U\sqcup U^{-1}\sqcup\bigl(X\setminus(s(U)\cup r(U))\bigr).
\]
Second countability gives a countable such cover of the nonunits.
Let $\Gamma$ be generated by its involutions.  Every group map has graph
in $R_{\mathcal G}$, and each nonunit arrow is represented locally by
one generator, so the orbit relations agree everywhere.

Every word has a compact-open full bisection $W_g$.  For a compact-open
chart $V$, principality implies that its coincidence set with $g$ is
$s(V\cap W_g)$, a clopen subset of $X$.  This proves the regularity of
both atlases and all their transition coefficients.  Their pullbacks
preserve $\mathrm{C}(X)$ on the relevant clopen corners, hence also its bounded
uniform completion.  All group maps preserve every groupoid-invariant
measure.  Conversely the countable cover decomposes any bisection into
Borel restrictions of generators and the identity; group invariance
therefore implies bisection invariance.  This proves equality of the
full invariant-measure families.  Equality of orbits also preserves
minimality.
\end{proof}

\begin{rem}\label{rem:germs-not-crossed-product}
The presentation identifies $\mathcal G$ with the groupoid of germs of
the action, rather than generally with $\Gamma\ltimes X$.  Indeed,
equal germs give equal endpoints and hence the same arrow by principality;
equal arrows give equal germs by intersection of their open bisections.
The cover gives surjectivity and local homeomorphisms give the topology.
The involutions may fix open complements. Thus a nonunit stabilizer
arrow of $\Gamma\ltimes X$ can map to a unit germ.  Accordingly the action
presentation does not identify $\mathrm{C}_r^*(\mathcal G)$ with an ordinary
crossed product, or assert that its crossed-product diagonal is Cartan.
Nor is the presenting group automatically amenable.
\end{rem}

\begin{prop}\label{prop:principal-groupoid-SBP-atlas}
If $\mathcal G$ is principal and $(X,K)$ has relative SBP, its bisections
admit a regular countable atlas.  Hence its relation has an admissible
Borel group-action presentation on one common conull set.  It also has
a continuous countable-group presentation on a compact zero-dimensional
model, with the corresponding designated measure family.
\end{prop}
\begin{proof}
Take a countable open bisection cover.  In each source, take a countable
cover by small-boundary open sets $A$ with closure contained in that
source.  The image $\beta(A)$ has boundary $\beta(\partial A)$ because
the closure lies inside the bisection domain.  Invariance makes its
boundary null for every designated measure.  Thus both domain indicators
are regular.  Pullback by this chart preserves the completed corners:
for a continuous function on the range, cut its pullback off near the
source boundary, extend by zero, and use uniformly vanishing boundary
neighborhood measure to obtain bounded uniform continuous approximation.
Invariance then extends the corner isomorphism to the completions.

Consider all finite words in these restrictions and their inverses.
Each word has open domain $A$, universally null boundary, and an
extension as an open bisection to a neighborhood of $\overline A$.
To see this inductively, suppose $\beta:A\to B$ and
$\gamma:C\to D$ have these properties for the maps and their inverses.
The domain of $\gamma\beta$ is $A\cap\beta^{-1}(C)$. Its closure lies
in the domain of the composite of the extended bisections, since
$\overline A$ and $\overline C$ lie in their respective extended
domains. Its boundary is contained in $\partial A$ together with the
preimage of $\partial C$ under the extended $\beta$. Both are
universally null by invariance. Applying the same argument to the
inverse word gives the assertion for its range. Thus each word
satisfies the support and corner conditions by the preceding cutoff
argument.

For any two words, equality of endpoints means equality of arrows by
principality. This coincidence set is closed in the common domain by
Hausdorffness and open there as the source of the bisection intersection.
Its boundary lies in the boundaries of the two domains, so its indicator
belongs to $D_K$. The words, the identity and the empty map form the
required regular countable atlas.

Corollary~\ref{cor:uniform-relation-countable-group-presentation} gives
global admissible generators.  Use
Corollary~\ref{cor:relative-SBP-AF-diagonal-core}, then close its countable
projection algebra under all generator pullbacks.  Its spectrum is compact
metrizable and zero dimensional, with a continuous countable group action.
Its bounded diagonal completion is still $D_K$; spatialization transports
the original action on common conull sets, including its regular atlas.
Only the designated invariant measures are asserted for this model.
\end{proof}

\begin{prop}[Almost finiteness supplies regular cuts]
\label{prop:groupoid-almost-finite-cuts}
Let $\mathcal G$ be principal and ample.  Assume that for every compact
$Q\subset\mathcal G$ and $\varepsilon>0$ there is a compact-open elementary
subgroupoid $H$ with
\[
 |QH_x\setminus H_x|<\varepsilon|H_x|\quad(x\in H^{(0)}),
 \qquad\mu(X\setminus H^{(0)})<\varepsilon\quad(\mu\in K).
\]
Then its regular orbit relation satisfies (RF).  In particular an aperiodic
such relation has the dynamical Cartan realization of
Theorem~\ref{thm:local-cuts-dynamical-model}.
\end{prop}
\begin{proof}
These are the elementary-subgroupoid conditions of almost finiteness in
measure, relative here to $K$; compare~\cite[Definition~4.15]{Ma}.
Fix finitely many compact-open bisections and include them in $Q$.
Extend the finite relation of $H$ by singleton classes outside $H^{(0)}$.
Its finite matrix charts and all their projections are compact-open.
For a test bisection $U$, in an $H$-class based at $x$, each member $y$
whose $U$-image escapes that class gives a distinct arrow in
$QH_x\setminus H_x$.  Principality identifies the class size with
$|H_x|$, so the escaping proportion is at most $\varepsilon$.
Invariance integrated in the elementary matrix charts gives escaping
measure at most $\varepsilon$ inside $H^{(0)}$; the uncovered part costs
at most another $\varepsilon$.  The estimate is uniform over $K$.
For arbitrary finitely many admissible maps, their first-chart
decompositions reduce to finitely many compact-open charts with uniformly
small omitted domains.  Making $\varepsilon$ sufficiently small proves
(RF).  The compact-open basis gives relative SBP on the unit space,
so Theorem~\ref{thm:local-cuts-AF} applies.
\end{proof}

\begin{rem}\label{rem:groupoid-isotropy-and-amenability}
Orbit equivalence forgets isotropy.  For example the one-point groupoid
of $\mathbb Z/2$ has the equality orbit relation, while its group algebra
is $\mathbb C^2$ and the unit diagonal $\mathbb C1$ is not maximal abelian.
For a completed-Cartan conclusion in the nonprincipal case, a clean
additional assumption is that isotropy-free units form an invariant set
conull for every $\mu\in K$, together with a regular bisection atlas.
Topological principality alone does not imply this measure condition.
Renault's correspondence retains the groupoid and a possible twist,
not merely its orbits~\cite{RenaultC}.

Groupoid amenability gives ordinary measured amenability of the orbit
relation for invariant measures.  For principal groupoids, relative
SBP supplies a regular atlas, and
Theorem~\ref{thm:uniform-CFW-SBP} then supplies (RF).
The elementary approximations above are sufficient but are no longer
necessary for this completed-model conclusion.
Proposition~\ref{prop:SBP-regularizes-twist} and
Theorem~\ref{thm:uniform-regular-untwisting} also remove locally trivial
circle twists at the level of the uniform completion.
For uniformly Borel amenable F\o lner graphs, the simultaneous
Borel cuts of~\cite[Theorems~7 and~8]{ET} remain a related input;
their indicators need not belong to the prescribed $D_K$.
Related groupoid SBP and castle constructions
are developed in~\cite{MW}; its Theorem~6.12 and Remark~6.13 distinguish
the proved implication from almost elementariness in measure to SBP from
the unresolved converse in general.
\end{rem}

\section{Uniform Connes--Feldman--Weiss under relative SBP}
\label{sec:uniform-CFW-SBP}

The remaining finite-approximation step is a local rounding argument.
A finite regular probability row is replaced by a finite token set.
Shared priorities give a finite-fibre center map with small expected
escape.  Relative SBP separates the finitely many vertices in each
dependency window outside a uniformly small set.  Local conditional
expectations then choose the priorities deterministically.
The change of variables under admissible maps makes each choice
decrease the integrated cost for every invariant measure at once.

The token construction is related to the finite random tilings of
Elek--Tim\'ar~\cite[Proposition~1]{ETD}; their local replacement
method~\cite[Lemma~8]{ETD} also uses measure-monotone choices.
Here the input is the integrated regular Reiter condition for the
designated family, and the following conditional-expectation argument
provides the required regular deterministic cuts directly.

For the rounding argument through Remark~\ref{rem:uniform-CFW-scope},
the relation and pseudogroup are those of
Definition~\ref{defn:regular-relation-atlas}, and relative $K$-SBP is
assumed.  The finite projection-density characterization gives
finite regular step approximations to diagonal coefficients.
A countable small-boundary basis also gives refining finite regular
Borel partitions $\mathcal P_n$ separating points on a common
$K$-conull set.  Representatives and identities are chosen on one
invariant set conull for every designated measure.

\begin{lem}[Rational regular probability rows]\label{lem:regular-kernel-rational-rows}
Suppose $\phi_1,\ldots,\phi_J$ are admissible partial bijections and
\[
 m_x=\sum_{j=1}^J f_j(x)\delta_{\phi_jx},\qquad
 0\leq f_j\in1_{\operatorname{dom}\phi_j} D_K,\qquad \sum_jf_j=1.
\]
For every $q>0$ there are finite-valued rational functions
$g_j\in1_{\operatorname{dom}\phi_j} D_K$, with $g_j\geq0$ and $\sum_jg_j=1$,
and an integer $L\geq1$ such that $Lg_j$ is integer-valued and
\[
 \sup_{\mu\in K}\int\sum_j|g_j-f_j|\,d\mu<q.
\]
Thus, for $\widetilde m_x=\sum_jg_j(x)\delta_{\phi_jx}$,
$\sup_\mu\int\|\widetilde m_x-m_x\|_1\,d\mu<q$.
\end{lem}
\begin{proof}
Choose nonnegative rational step functions $h_j$ on a common finite
regular partition, supported on $\operatorname{dom}\phi_j$, such that
\[
 \sup_\mu\int\sum_j|h_j-f_j|\,d\mu<q/2.
\]
To preserve supports, enlarge the finite partition to contain the domain
projections before multiplying the step functions by them.  Clipping
and rational approximation preserve the required approximation.
Let $H=\sum_jh_j$.  Where $H>0$, put $g_j=h_j/H$.
Where $H=0$, give weight one to the first defined $\phi_j$.
At least one $\phi_j$ is defined because $\sum_jf_j=1$.
These choices are finite-valued rational regular functions.
Pointwise,
\[
 \sum_j|g_j-f_j|\leq 2\sum_j|h_j-f_j|.
\]
Indeed, where $H>0$, use
$\sum_j|g_j-h_j|=|1-H|\leq\sum_j|h_j-f_j|$; where $H=0$ the
right side is $2$ and the left side is at most $2$.
A common denominator of the finitely many values gives $L$.
\end{proof}

\begin{lem}[A finite regular coloring of local windows]
\label{lem:regular-window-colors}
Let $\psi_1,\ldots,\psi_Q$ be finitely many admissible partial
bijections and put
\[
 W_x=\{\psi_ax:x\in\operatorname{dom}\psi_a,\ 1\leq a\leq Q\}.
\]
For every $\delta>0$ there are a finite-valued regular coloring
$c:X\to\{1,\ldots,k\}$ and a regular Borel set $E$ such that
$\sup_{\mu\in K}\mu(E)<\delta$ and, for each $x\notin E$,
distinct vertices of $W_x$ have distinct colors.
\end{lem}
\begin{proof}
Let $c_n$ be the code of the finite partition $\mathcal P_n$.
For each pair $a,b$, consider
\[
 E_{n,a,b}=\{x\in\operatorname{dom}\psi_a\cap\operatorname{dom}\psi_b:
     \psi_ax\ne\psi_bx,\ c_n(\psi_ax)=c_n(\psi_bx)\}.
\]
Every indicator is in $ D_K$: use the corner pullbacks of the partition
projections and subtract the regular equalizer projection.
The sets decrease to the empty set on the common conull set because the
partition codes separate distinct points.  Set
$E_n=\bigcup_{a,b}E_{n,a,b}$.
For each $n$, $\mu\mapsto\mu(E_n)$ is continuous on $ K$.
This follows from bounded uniform $2$-approximation of its indicator by
continuous functions.  The functions decrease pointwise to zero, so
Dini's theorem gives $\sup_\mu\mu(E_n)\to0$.
Choose a suitable $n$ and take $c=c_n$ and $E=E_n$.
\end{proof}

\begin{lem}[Finite-window conditional-expectation derandomization]
\label{lem:regular-local-derandomization}
Use the windows and coloring of Lemma~\ref{lem:regular-window-colors}.
Let $A$ be a finite alphabet.  For each root $x$, let $h_x$ be a
nonnegative cost depending only on one alphabet label $a(v)\in A$
at each distinct $v\in W_x$.  Suppose the costs are uniformly bounded,
finite-valued, and regular in the following precise sense: after
specifying any finite pattern of chart coincidences, colors, and
previously assigned regular label functions, each finite local
conditional-expectation table is a finite-valued element of $ D_K$.
Write $\mathbb{E} h_x$ for the finite average obtained by assigning
independent uniform labels in $A$ to the distinct vertices of $W_x$.
There is a finite-valued regular label function $a:X\to A$ such that,
for every $\mu\in K$ simultaneously,
\[
 \int_{X\setminus E}h_x(a)\,d\mu(x)
 \leq \int_{X\setminus E}\mathbb{E} h_x\,d\mu(x).
\]
No random field on an uncountable space is used: all expectations are
finite local averages.
\end{lem}
\begin{proof}
Process the colors $1,\ldots,k$ in order.  At the start of a stage the
labels at earlier colors are fixed by regular functions.  Labels at
current and later colors remain formal independent uniform alphabet
variables in each local averaging table.
Every good root $x\notin E$ has at most one vertex of the current
color in $W_x$.  Therefore simultaneous assignment of this color
introduces no interaction between two currently assigned variables in
any good root cost.

Fix a current-color vertex $v$.  Its affected good roots are the
distinct points
\[
 Z_v=\{\psi_a^{-1}v:v\in\operatorname{ran}\psi_a\}\setminus E.
\]
There are at most $Q$ such roots.  Deduplicate this enumeration by
retaining the first inverse chart producing each root.  Explicitly,
with $x_a(v)=\psi_a^{-1}v$, use the indicator
\[
 I_a(v)=1_{\operatorname{ran}\psi_a}(v)1_{X\setminus E}(x_a(v))
 \prod_{b<a}\left(1-
 1_{\operatorname{ran}\psi_b}(v)1_{\{\psi_b^{-1}v=x_a(v)\}}(v)
 \right).
\]
Undefined corner terms are interpreted as zero.  These indicators are
regular by the inverse-chart equalizer hypothesis.

For $x\in Z_v$ and $u\in A$, let $H(x,v,u)$ be the finite average
of $h_x$ after fixing the current variable at $v$ to $u$, retaining
the earlier choices and averaging independently over the remaining
unassigned vertices of $W_x$.
Set
\[
 B(v,u)=\sum_{a=1}^Q I_a(v)H(x_a(v),v,u).
\]
Choose $a(v)$ to be the first minimizer of $B(v,u)$ in a fixed order on
$A$.  Each $B(\cdot,u)$ is finite-valued and regular.
Thus this argmin is regular: the indicator of any comparison between
two finite-valued regular functions is a polynomial of their
finite-valued difference.  In particular,
\begin{equation}
 B(v,a(v))-\frac1{|A|}\sum_{u\in A}B(v,u)\leq0.
\label{eq:regular-local-conditional-cost}
\end{equation}

Here is the exact integration argument.  If $x\notin E$ has no
current-color window vertex, its conditional expected cost is
unchanged.  Otherwise let $v$ be its unique current-color window
vertex.  The change is
\[
 H(x,v,a(v))-\frac1{|A|}\sum_{u\in A}H(x,v,u).
\]
For every bounded regular function $F$ on these root--window pairs,
invariance of $\mu$ under the $\psi_a$ gives the finite mass-transport
identity
\begin{equation}
 \int\sum_a I_a(v)F(x_a(v),v)\,d\mu(v)
 =\int_{X\setminus E}\sum_{v\in W_x}F(x,v)\,d\mu(x).
\label{eq:regular-conditional-mass-transport}
\end{equation}
The left side enumerates each distinct root exactly once for each
vertex.  To verify the identity, apply the change of variables
$v=\psi_ax$ to each summand with its first-occurrence indicator;
the transformed indicators retain exactly the first chart producing
each distinct vertex in $W_x$.
Apply \eqref{eq:regular-conditional-mass-transport} with the current-color indicator and the displayed
cost change.  The change in the integrated conditional expected cost
is
\[
 \int_{\{c(v)=\text{current color}\}}
 \left(B(v,a(v))-
       \frac1{|A|}\sum_{u\in A}B(v,u)\right)d\mu(v)\leq0
\]
by \eqref{eq:regular-local-conditional-cost}.  This calculation holds for every $\mu\in K$
with the same pointwise minimizing labels.  Iterate through the
finitely many colors.  At the end all labels are assigned and no
local averaging remains, proving the assertion.
\end{proof}

\begin{thm}[Quantitative regular rounding]\label{thm:regular-Reiter-rounding}
Let $m_x=\sum_{j=1}^J f_j(x)\delta_{\phi_jx}$ be a finite regular
probability kernel as in Lemma~\ref{lem:regular-kernel-rational-rows}.  Let $F$ be a finite
set of admissible partial bijections.  For every $\eta>0$ there is a
regular finite Borel subequivalence relation $S\subseteq R$ such
that its classes have size at most $J$, it is generated by finitely
many admissible partial bijections, and for every $\mu\in K$,
\begin{equation}
 \sum_{s\in F}\mu\{x\in\operatorname{dom} s:(sx,x)\notin S\}
 \leq
 \sum_{s\in F}\int_{\operatorname{dom} s}\|m_{sx}-m_x\|_1\,d\mu(x)+\eta.
\label{eq:regular-kernel-rounding-bound}
\end{equation}
\end{thm}
\begin{proof}
The case $F=\varnothing$ is immediate.  Put $d=|F|\geq1$.
Use Lemma~\ref{lem:regular-kernel-rational-rows} to obtain rational rows
$\widetilde m$ with
\[
 q:=\sup_\mu\int\|\widetilde m_x-m_x\|_1\,d\mu(x)
       <\eta/(6d).
\]
Since $s$ preserves $\mu$ on its domain,
\begin{equation}
 \sum_{s\in F}\int_{\operatorname{dom} s}
       \|\widetilde m_{sx}-\widetilde m_x\|_1\,d\mu
 \leq
 \sum_{s\in F}\int_{\operatorname{dom} s}\|m_{sx}-m_x\|_1\,d\mu+2dq.
\label{eq:regular-rational-variation}
\end{equation}

Let $L$ be a common denominator for the rational rows.  Aggregate
coincident chart outputs: for $z\in[x]_ R$ put
$n_x(z)=L\widetilde m_x(z)$, an integer in $\{0,\ldots,L\}$.
The finite token set
\[
 A_x=\{(z,r):n_x(z)>0,\ 1\leq r\leq n_x(z)\}
\]
has exactly $L$ elements.  Its vertices belong to
$\{\phi_jx:x\in\operatorname{dom}\phi_j\}$.
Coincidence patterns of the outputs are regular, so all the token
counts and first-chart enumerations of $A_x$ are finite-valued
regular functions.

Choose an integer $N\geq2L$ large enough that
\begin{equation}
 d\binom{2L}{2}/N<\eta/3.
\label{eq:regular-priority-ties}
\end{equation}
The label at a vertex will be a vector in the finite alphabet
$A=\{1,\ldots,N\}^L$.  Its $r$th coordinate is the priority of
the token $(z,r)$.  Given any label function, define $C(x)$ to be the
vertex of a minimum-priority token in $A_x$.  Resolve a tie by the
first support-chart index and then the copy index.  This is a
finite regular rule; no global ordering of $X$ is necessary.

First examine the finite local average in which the vectors at all
distinct vertices of $A_x\cup A_{sx}$ are independently uniform in
$A$.  Their coordinates give independent uniform priorities to all
distinct tokens.  The probability of any priority collision among
these at most $2L$ tokens is at most $\binom{2L}{2}/N$.
Conditioned on no collision, each token is equally likely to be
the minimum of $A_x\cup A_{sx}$.  If this minimum is in their
intersection then both rows choose its vertex.  Consequently
\begin{equation}
 \mathbb{E}1_{\{C(x)\ne C(sx)\}}
 \leq\frac{|A_x\triangle A_{sx}|}{|A_x\cup A_{sx}|}
           +\binom{2L}{2}/N
 \leq\|\widetilde m_x-\widetilde m_{sx}\|_1
           +\binom{2L}{2}/N.
\label{eq:regular-token-minimum-bound}
\end{equation}
The last inequality uses
$|A_x\triangle A_{sx}|=L\|\widetilde m_x-\widetilde m_{sx}\|_1$
and $|A_x\cup A_{sx}|\geq L$.
All expectations here are finite sums over the displayed token
variables; no sample of an uncountable random field is asserted.

For the derandomization use the common root window
\[
 W_x=\{\phi_jx:x\in\operatorname{dom}\phi_j\}
      \cup\bigcup_{s\in F}
        \{\phi_j(sx):x\in\operatorname{dom}(\phi_j s)\}.
\]
It is specified by at most $J(d+1)$ admissible partial bijections.
Lemma~\ref{lem:regular-window-colors} supplies a finite regular coloring and a
regular set $E$ such that
\begin{equation}
 d\sup_{\mu\in K}\mu(E)<\eta/3
\label{eq:regular-window-bad}
\end{equation}
and distinct window vertices have distinct colors for $x\notin E$.
Let the root cost be
\[
 h_x=\sum_{s\in F}1_{\operatorname{dom} s}(x)
                         1_{\{C(x)\ne C(sx)\}}.
\]
This cost is integer-valued between $0$ and $d$.
Its local conditional-expectation tables satisfy the regularity
requirement in Lemma~\ref{lem:regular-local-derandomization}: partition by the finite
chart coincidence pattern, token counts, vertex colors, and earlier
label values.  On each resulting regular piece, the table is the
average of finitely many explicitly specified finite decisions.
These tables have finitely many rational values.  Pulling them back
along inverse window charts uses only admissible compositions.
Apply that lemma to produce one finite-valued regular vector-label
function such that, for every $\mu\in K$,
\begin{equation}
 \int_{X\setminus E}h_x\,d\mu
 \leq \int_{X\setminus E}\mathbb{E} h_x\,d\mu
 \leq\sum_{s\in F}\int_{\operatorname{dom} s}
         \|\widetilde m_x-\widetilde m_{sx}\|_1\,d\mu
             +d\binom{2L}{2}/N.
\label{eq:regular-rounded-cost}
\end{equation}

For these deterministic labels, $C$ is a Borel map whose graph is
a finite union of regularly restricted $\phi_j$ graphs.  Define
\[
 (x,y)\in S\quad\Longleftrightarrow\quad C(x)=C(y).
\]
This is a Borel equivalence relation contained in $ R$.
For any $z$, each $x$ with $C(x)=z$ satisfies $\phi_jx=z$ for some
$j$.  Each $\phi_j$ is injective, so $|C^{-1}(z)|\leq J$.

For completeness, $S$ is regular and finitely admissibly generated.
Choose disjoint regular pieces $B_j$ on which $C=\phi_j$, assigning
each root to its first selected chart.  For each $j,k$, restrict
$\phi_k^{-1}\phi_j$ to
\[
 B_j\cap\operatorname{dom}(\phi_k^{-1}\phi_j)
          \cap(\phi_k^{-1}\phi_j)^{-1}(B_k).
\]
The restriction is admissible.  These finitely many maps enumerate
all arrows of $S$, because two selected chart outputs are equal
exactly when the indicated inverse composition sends one root to
the other.  They therefore generate $S$.

The escape cost for $S$ is exactly $h_x$.
On $E$ it is at most $d$.  Combine
\eqref{eq:regular-rational-variation}, \eqref{eq:regular-priority-ties}, \eqref{eq:regular-window-bad}, and
\eqref{eq:regular-rounded-cost} to obtain \eqref{eq:regular-kernel-rounding-bound}.
\end{proof}

\begin{cor}[Uniform regular Reiter implies regular finite cuts]\label{cor:uniform-Reiter-RF}
If for every finite set $F$ of admissible partial maps and every
$\varepsilon>0$ there is a finite regular probability kernel $m$
such that
\[
 \sup_{\mu\in K}\sum_{s\in F}\int_{\operatorname{dom} s}
                \|m_{sx}-m_x\|_1\,d\mu<\varepsilon,
\]
then $ R$ has the regular finite-cut property: each finite $F$ and
each $\varepsilon>0$ admit a bounded regular finite subrelation
$S\subseteq R$, generated by finitely many admissible maps, with
\[
 \sup_{\mu\in K}\sum_{s\in F}
          \mu\{x\in\operatorname{dom} s:(sx,x)\notin S\}<\varepsilon.
\]
\end{cor}
\begin{proof}
Choose the Reiter error below $\varepsilon/2$, and apply
Theorem~\ref{thm:regular-Reiter-rounding} with $\eta<\varepsilon/2$.
\end{proof}

\begin{thm}[Uniform Connes--Feldman--Weiss with relative SBP]
\label{thm:uniform-CFW-SBP}
Let $(X,R,K,D_K,\mathcal S)$ be a regular uniform measure relation,
with the nonempty compact convex designated family $K$ as above.
Assume relative $K$-SBP.  The following are equivalent:
\begin{enumerate}
\item $(X,R,\mu)$ is amenable for every $\mu\in K$;
\item $R$ has regular uniform Reiter amenability;
\item $R$ has local regular finite approximation \textup{(RF)};
\item $R$ is uniformly AF hyperfinite in the finite matrix-chart sense.
\end{enumerate}
The AF model retains the designated compact trace family.
No aperiodicity, Bauer hypothesis, simplicity, or identification
with the entire AF trace simplex is required for these equivalences.
\end{thm}
\begin{proof}
The first two conditions are equivalent by
Proposition~\ref{prop:Reiter-measurewise-equivalence}.
Atlas independence gives the Reiter estimates for each finite set
of admissible tests.  Corollary~\ref{cor:uniform-Reiter-RF} proves
the implication from the second condition to the third.
Theorem~\ref{thm:local-cuts-AF} identifies the third and fourth
under relative SBP.  Finally finite AF classes give regular Reiter
kernels by Proposition~\ref{prop:uniform-AF-implies-Reiter}.
\end{proof}

\begin{rem}\label{rem:uniform-CFW-scope}
The deterministic construction uses invariant measures in its
mass-transport identity; it is not asserted for arbitrary measure
families.  Compactness makes the regular window-collision errors
uniform through Dini's theorem, and convexity is used in the earlier
measurewise-to-Reiter minimax argument.
The resulting common exhaustion is regular in the prescribed
diagonal, rather than merely a separate measurable exhaustion
for each trace.  Aperiodicity subsequently gives the single
generator and Cantor dynamical model of
Theorem~\ref{thm:local-cuts-dynamical-model}.
For aperiodic relations carrying their full compact invariant family,
Corollary~\ref{cor:full-family-simple-AF-model} supplies a simple
full-trace replacement and the corresponding minimal realization.
For arbitrary designated families, the separate model hypothesis
remains necessary for that conclusion.
\end{rem}

\subsection{Property (C) for principal \'etale groupoids}
\label{subsec:SBP-groupoid-property-C}

There is a comparison consequence of small boundaries which is
independent of the amenability and finite-exhaustion arguments above.
It concerns the original standard Cartan pair, with witnesses in its
groupoid $\mathrm{C}^*$-algebra. Throughout this subsection, $\mathcal G$ is a
second-countable locally compact Hausdorff principal \'etale groupoid
with compact metrizable unit space $X$, and
$K=M(\mathcal G)\neq\varnothing$ is its \emph{full} family of invariant
probability measures. This family is compact and convex. Principality
identifies the traces of $A=\mathrm{C}_r^*(\mathcal G)$ with
$\tau_\mu=\mu\circ\mathbb E_D$, $\mu\in K$, where $D=\mathrm{C}(X)$ and
$\mathbb E_D:A\to D$ is the canonical expectation. In particular,
$\norm{\cdot}_{2,\T(A)}=\norm{\cdot}_{2,K}$.

We recall the relative comparison property introduced in
\cite[Definition~3.1]{EN-C}.

\begin{defn}[Property (C)]\label{defn:Cartan-property-C}
Let $D\subseteq A$ be a unital commutative sub-$\mathrm{C}^*$-algebra of a unital
$\mathrm{C}^*$-algebra with $\T(A)\neq\varnothing$. The pair has \emph{Property
\textup{(C)}} if, whenever $f,g,h\in D$ are positive contractions with
$f,g\in\overline{hDh}$ and
\[
 d_\tau(f)<d_\tau(g)\qquad(\tau\in\T(A)),
\]
then, for every $\varepsilon>0$, there is a contraction
$u\in\overline{hAh}+\mathbb C1_A$ such that
\begin{align*}
 \dist_{2,\T(A)}(ufu^*,\overline{gAg})&<\varepsilon,\\
 \sup_{d\in D_1}\max\bigl\{
 \dist_{2,\T(A)}(udu^*,D_1),
 \dist_{2,\T(A)}(u^*du,D_1)\bigr\}&<\varepsilon,\\
 \max\bigl\{\norm{uu^*-1}_{2,\T(A)},
                \norm{u^*u-1}_{2,\T(A)}\bigr\}&<\varepsilon.
\end{align*}
Here $d_\tau(a)=\lim_n\tau(a^{1/n})$ for $a\geq0$, and
$\dist_{2,\T(A)}(a,S)=\inf_{b\in S}\norm{a-b}_{2,\T(A)}$.
\end{defn}

The uniformity over the whole diagonal unit ball is part of the
definition. The following matching lemma is the ingredient that
allows us to retain this uniformity without assuming amenability.

\begin{lem}[Uniform finite matching along bisections]
\label{lem:SBP-uniform-bisection-matching}
Assume $(X,K)$ has relative SBP. Let $U,V\subseteq X$ be open sets
with $K$-null boundaries and
\[
                    \mu(U)<\mu(V)\qquad(\mu\in K).
\]
For every $\eta>0$, there are finitely many open bisections with
induced partial homeomorphisms $\theta_i$ and Borel source sets $P_i$
such that, writing $Q_i=\theta_i(P_i)$,
\begin{enumerate}
\item $P_i\subseteq U\setminus V$, $Q_i\subseteq V\setminus U$,
      and all the sets $P_i,Q_i$ are pairwise disjoint;
\item every $P_i$ and $Q_i$ has $K$-null boundary, and
      $\overline{P_i}$ is compactly contained in a larger bisection
      domain on which $\theta_i$ extends;
\item
\[
 \sup_{\mu\in K}\mu\left((U\setminus V)
                            \setminus\bigcup_iP_i\right)<\eta.
\]
\end{enumerate}
\end{lem}
\begin{proof}
Choose a countable open bisection atlas covering $\mathcal G$ whose
induced maps $\beta_n$ have domains $\Omega_n$ with $K$-null boundary
and compact closure inside a larger bisection domain. Such an atlas
exists by SBP and second countability: shrink the source of a
precompact bisection to a small-boundary neighborhood whose closure
lies inside its domain. Invariance and this compact containment show
that the range of every restricted bisection also has $K$-null
boundary. The same is true of finite Boolean combinations and of
their images and inverse images inside these bisection domains.

Set $S_0=U\setminus V$ and $T_0=V\setminus U$. At step $n$, match
every currently available source point whose $\beta_n$-image is an
available target point:
\[
 P_n=S_{n-1}\cap\Omega_n\cap\beta_n^{-1}(T_{n-1}),\qquad
 Q_n=\beta_n(P_n),
\]
and put $S_n=S_{n-1}\setminus P_n$,
$T_n=T_{n-1}\setminus Q_n$. These are Borel sets with $K$-null
boundaries. For example,
\[
 \partial P_n\subseteq\partial S_{n-1}\cup\partial\Omega_n
       \cup\beta_n^{-1}(\partial T_{n-1}),
\]
where the inverse image is taken in the larger bisection domain.
Invariance makes the last set null for every $\mu\in K$.
The matched source and target sets are pairwise disjoint, and each
pair has equal measure for every invariant measure.

Let $S=\bigcap_nS_n$ and $T=\bigcap_nT_n$. No orbit meets both $S$
and $T$. Indeed, an arrow from $x\in S$ to $y\in T$ lies in some
atlas bisection, and its step would have matched $x$ to $y$.
Let $Z$ be the orbit saturation of $S$. This is an invariant Borel
set, since the atlas is countable, and $Z\cap T=\varnothing$.
If $\mu(S)>0$ for some $\mu\in K$, then
\[
                \nu=\mu(Z)^{-1}\mu|_Z
\]
is another invariant probability measure, hence belongs to $K$.
Countable additivity and invariance cancel all matched pairs, giving
\[
 \nu(U)-\nu(V)=\nu(S)-\nu(T)=\nu(S)>0,
\]
contrary to the hypothesis. Thus $\mu(S)=0$ for every $\mu\in K$.

Since $S_n$ has $K$-null boundary, the function
$\mu\mapsto\mu(S_n)$ is continuous on $K$. These continuous
functions decrease pointwise to zero. Dini's theorem gives
$\sup_{\mu\in K}\mu(S_n)\to0$. Retain a sufficiently long finite
initial segment of the matching, omit empty pieces, and use the
corresponding atlas maps on the retained pieces.
\end{proof}

\begin{thm}[SBP implies Property (C) for the standard Cartan pair]
\label{thm:SBP-groupoid-property-C}
Let $\mathcal G$ be a second-countable locally compact Hausdorff
principal \'etale groupoid with compact metrizable unit space $X$
and nonempty full invariant-probability family $K=M(\mathcal G)$.
If $(X,K)$ has relative SBP, then
\[
                    \mathrm{C}(X)\subseteq \mathrm{C}_r^*(\mathcal G)
\]
has Property \textup{(C)}. Its witnesses may be chosen to be
unitaries $u\in1+\overline{h\mathrm{C}_r^*(\mathcal G)h}$. Neither
amenability nor minimality is required.
\end{thm}
\begin{proof}
Write $D=\mathrm{C}(X)$ and $A=\mathrm{C}_r^*(\mathcal G)$. Fix positive contractions
$f,g,h\in D$ as in Definition~\ref{defn:Cartan-property-C}, and fix
$\varepsilon>0$. For $\mu\in K$ the rank inequality reads
\[
             \mu(\{f>0\})<\mu(\{g>0\}).
\]
Choose $t>0$ and put $f_t=(f-t)_+$. By SBP and compactness of
$\{f\geq t\}$, choose a small-boundary open set $U$ with
\[
             \{f\geq t\}\subseteq U\Subset\{f>0\}.
\]
In particular $f_t$ vanishes outside $U$. Choose increasing
small-boundary open sets $V_m\Subset\{g>0\}$ whose union is
$\{g>0\}$. This can be done by taking finite unions from a countable
SBP basis with closures inside $\{g>0\}$. For each $\mu\in K$,
$\mu(U)<\mu(V_m)$ for some $m$. Both sides are continuous functions
of $\mu$, so compactness of $K$ and monotonicity give one $m$ for
which the inequality holds for every $\mu\in K$. Fix $V=V_m$.
The closures of $U$ and $V$ lie in $\{h>0\}$, since
$f,g\in\overline{hDh}$.

Apply Lemma~\ref{lem:SBP-uniform-bisection-matching} with a small
$\eta>0$. Let $P_i,Q_i=\theta_i(P_i)$ be the resulting finite
matching, and denote its source remainder by
\[
                S_N=(U\setminus V)\setminus\bigcup_iP_i.
\]
Choose compact cores $C_i\subseteq\operatorname{int}P_i$ and
$C_V\subseteq V$ such that, with $c_K(B)=\sup_{\mu\in K}\mu(B)$,
\begin{equation}\label{eq:property-C-core-losses}
 c_K(S_N)<\eta,\qquad
 \sum_i c_K(P_i\setminus C_i)<\eta,\qquad
 c_K(V\setminus C_V)<\eta.
\end{equation}
Here uniform compact extraction is elementary: the boundaries of the
finitely many sets are compact and $K$-null; by compactness of $K$,
their open neighborhoods have arbitrarily small $c_K$-measure.
Removing such neighborhoods from their closures leaves compact
subsets of their interiors.

Choose nested open sets
\[
 C_i\subseteq G_i\Subset W_i\Subset\operatorname{int}P_i
\]
and positive contractions $\varphi_i,\psi_i\in \mathrm{C}_c(W_i)$ such that
$\varphi_i=1$ on $C_i$, $\supp\varphi_i\subseteq G_i$,
$\psi_i=1$ on $\overline{G_i}$, and
$\supp\psi_i\subseteq W_i$. The open levels $W_i$ and
$\theta_i(W_i)$ have pairwise disjoint compact closures. The
bisection between these two levels gives the usual matrix-chart
inclusion $\mathrm{M}_2(\mathrm{C}_0(W_i))\subseteq A$. In its multiplier algebra put
\[
 P=\begin{pmatrix}0&1\\1&0\end{pmatrix},\qquad
 v_i=1+(e^{\mathrm i\pi\psi_i}-1)\frac{1-P}{2}.
\]
The matrix $v_i$ is unitary, is the swap $P$ on $G_i$, and is the
identity off $\supp\psi_i$. Thus $v_i-1\in \mathrm{M}_2(\mathrm{C}_c(W_i))$.
Extend these differences by zero into $A$ and define
\[
                         u=1+\sum_i(v_i-1).
\]
Disjointness of all the two-level footprints makes $u$ unitary.
Every matrix coefficient of $u-1$ has compact support with both
endpoints in $\{h>0\}$. Dividing the coefficients by the two
endpoint values of $h$ therefore gives
\begin{equation}\label{eq:property-C-unitary-localization}
                         u-1\in hAh.
\end{equation}

Choose $q\in \mathrm{C}_c(V)$ with $0\leq q\leq1$ and $q=1$ on $C_V$, and
set
\[
                         z=f_tq+\sum_i f_t\varphi_i\in D.
\]
Its summands have disjoint supports, so $0\leq z\leq f_t\leq1$.
The function $f_tq$ vanishes on both levels of every swap: the
source levels are outside $V$, and the target levels are outside
$U$, where $f_t=0$. Hence $u$ fixes this summand. On the support
of $\varphi_i$, $u$ is the exact swap, so it transports
$f_t\varphi_i$ compactly into $V$. Consequently
\begin{equation}\label{eq:property-C-exact-diagonal-comparison}
                         uzu^*\in\overline{gDg}.
\end{equation}
The function $f_t-z$ vanishes on $C_V\cup\bigcup_iC_i$ and
outside $U$. Its nonzero set is contained in
$S_N\cup(V\setminus C_V)\cup\bigcup_i(P_i\setminus C_i)$.
Equation~\eqref{eq:property-C-core-losses} therefore gives
\[
 \norm{f_t-z}_{2,K}\leq\sqrt{3\eta},\qquad
 \dist_{2,K}(ufu^*,\overline{gAg})
       \leq t+\sqrt{3\eta}.
\]

Finally $u$ is a diagonal permutation normalizer on the $G_i$
levels and the identity outside the $W_i$ levels. Its failure to
normalize $D$ is confined to
\[
 C=\bigsqcup_i\bigl((W_i\setminus G_i)
                   \sqcup\theta_i(W_i\setminus G_i)\bigr).
\]
Invariance and \eqref{eq:property-C-core-losses} imply
\[
 c_K(C)\leq2\sum_i c_K(P_i\setminus C_i)<2\eta.
\]
For every $d\in D_1$, the defect
$udu^*-\mathbb E_D(udu^*)$ has norm at most $2$ and is supported
on these transition levels. Since all traces factor through
$\mathbb E_D$, the matrix-chart calculation gives
\[
 \norm{udu^*-\mathbb E_D(udu^*)}_{2,K}^2\leq4c_K(C)<8\eta.
\]
The same estimate holds for $u^*du$, and both expectations belong
to $D_1$. Choose $t,\eta$ so that
$t+\sqrt{3\eta}<\varepsilon$ and $2\sqrt{2\eta}<\varepsilon$.
This proves both almost-normalizing estimates uniformly over
$D_1$, as well as the comparison estimate. Exact unitarity supplies
the remaining two conditions with zero error.
\end{proof}

\begin{rem}[The measure family and the converse]
\label{rem:groupoid-property-C-scope}
The proof uses the full invariant-probability family when restricting
a measure to the invariant saturation of an unmatched source set.
For a designated compact subfamily, the same matching argument
requires closure under normalized restrictions to invariant Borel
sets of positive measure. No such closure is imposed on the general
designated families elsewhere in this paper.

The converse fails in general. In
Example~\ref{ex:cover-AH-without-transversal}, the standard Cartan
pair of the double covering has Property \textup{(C)} but not SBP.
Indeed, at each two-point fibre the strict rank inequality forces
$g$ to be positive at both points whenever $f$ is positive at either
one. Thus $\{f>0\}\subseteq\{g>0\}$ and $u=1$ verifies
Property \textup{(C)} exactly. The invariant uniform fibre measures
detect every nonempty boundary, so SBP would force a clopen basis
on the connected circle. This example has finite orbits and does
not address a converse for minimal groupoids with infinite orbits.
\end{rem}

\subsection{Small subgroupoids from relative SBP}
\label{subsec:SBP-groupoid-castles}

The regular finite models can also be realized by relatively compact
open subgroupoids in the original topology.  The following construction
retains the original bisections; it does not change the unit space or
replace the groupoid by its abstract AF model.

Let $\mathcal G$ be a second-countable locally compact Hausdorff principal
\'etale groupoid with compact metrizable unit space $X$, and let $K$ be
a nonempty compact family of invariant probability measures.  Write
\[
 c_K(A)=\sup_{\mu\in K}\mu(A)
\]
for a Borel subset $A\subseteq X$.  By a \emph{small elementary
subgroupoid} below we mean an open subgroupoid with a finite multisection
presentation as follows.  For finitely many indices $a$, let $V_a$ be
an open parameter set in $X$ and let $\gamma_{a,i}:V_a\to\mathcal G$, $i\in F_a$,
be continuous bisection sections with
\[
 s(\gamma_{a,i}(x))=x,\qquad h_{a,i}(x)=r(\gamma_{a,i}(x)).
\]
The sections extend to a neighborhood of $\overline{V_a}$, and all the
compact level sets $h_{a,i}(\overline{V_a})$ are pairwise disjoint.
The subgroupoid is
\begin{equation}\label{eq:small-elementary-multisections}
 H=\bigsqcup_a\bigsqcup_{i,j\in F_a}
 \{\gamma_{a,i}(x)\gamma_{a,j}(x)^{-1}:x\in V_a\}.
\end{equation}
Its unit space, or footprint, is
$H^{(0)}=\bigsqcup_{a,i}h_{a,i}(V_a)$.  Every source fiber based in this
footprint has $|F_a|$ elements.  The closure $\overline H$, obtained by
using $\overline{V_a}$ in the same formula, is a compact elementary
subgroupoid.  In particular, $H$ is relatively compact.  Its footprint
need not be all of $X$.

For a subset $Q\subseteq\mathcal G$ and a finite nonempty set
$F\subseteq\mathcal G_x=s^{-1}(x)$, use the full left boundary
\begin{align*}
 \partial_Q^+F&=QF\setminus F,\\
 \partial_Q^-F&=\{f\in F:qf\notin F\text{ for some }q\in Q
                              \text{ composable with }f\},\\
 \partial_QF&=\partial_Q^+F\cup\partial_Q^-F.
\end{align*}
The outer and inner boundaries are disjoint.  Thus the full
F\o lner estimate is $|\partial_QF|<\eta|F|$, rather than only an
estimate for $QF\setminus F$; compare \cite[Definitions~2.7--2.9]{MW}.

\begin{lem}[Uniform compact bases]\label{lem:regular-uniform-compact-bases}
Let $A\subseteq X$ be Borel with $1_A\in D_K$, let $O\subseteq X$ be
open, and suppose $c_K(A\setminus O)=0$.  For every $\rho>0$ there
is a compact set $F\subseteq O$ such that
\[
 c_K(A\setminus F)<\rho,\qquad c_K(F\setminus A)=0.
\]
Also, every compact $K$-null set has open neighborhoods of arbitrarily
small $c_K$-measure.  No SBP assumption is needed for these assertions.
\end{lem}
\begin{proof}
Choose continuous functions $0\leq f_n\leq1$ such that
\[
 e_n=\norm{f_n-1_A}_{2,K}^{\,2},\qquad \sum_ne_n<\infty,
\]
and put $F_N=\bigcap_{n\geq N}\{f_n\geq1/2\}$.  These sets are
compact, and Chebyshev's inequality gives
\[
 c_K(A\setminus F_N)\leq4\sum_{n\geq N}e_n.
\]
For every $n\geq N$ and $\mu\in K$, one also has
$\mu(F_N\setminus A)\leq4e_n$.  Letting $n\to\infty$ shows that
$F_N\setminus A$ is one common $K$-null Borel set.

If $N_0\subseteq X$ is compact and $K$-null, choose decreasing
continuous functions $0\leq u_m\leq1$ with $u_m=1$ on $N_0$ and
$u_m\downarrow1_{N_0}$ pointwise.  The continuous functions
$\mu\mapsto\int u_m\,d\mu$ decrease to zero on compact $K$.
Dini's theorem makes the convergence uniform.  The open neighborhoods
$\{u_m>1/2\}$ therefore have arbitrarily small $c_K$-measure.

Choose $N$ with $4\sum_{n\geq N}e_n<\rho/2$.
The compact set $F_N\setminus O$ is $K$-null, since it is contained
in $(F_N\setminus A)\cup(A\setminus O)$.  Take an open neighborhood
$W$ of this compact set with $c_K(W)<\rho/2$, and set $F=F_N\setminus W$.
This proves the assertions.  Notice that $F$ need only be contained in
$A$ modulo a common null set; no literal compact subset of $A$ is asserted.
\end{proof}

\begin{lem}[Lifting regular finite charts to small subgroupoids]
\label{lem:regular-open-castle-lifting}
Assume relative $K$-SBP, and use the regular bisection atlas of
Proposition~\ref{prop:principal-groupoid-SBP-atlas}.  Suppose finitely
many regular Borel bases $A_a$ have admissible coordinate maps
$t_{a,i}:A_a\to X$, $i\in F_a$, whose levels $t_{a,i}(A_a)$ are
globally disjoint.  Choose nonempty subsets $I_a\subseteq F_a$, and
put
\[
 W_0=\bigsqcup_a\bigsqcup_{i\in I_a}t_{a,i}(A_a).
\]
Let $Q\subseteq\mathcal G$ be compact and symmetric.  Suppose, on
one common invariant $K$-conull set, that
\begin{equation}\label{eq:regular-chart-core-capture}
 r\bigl(Q\cap s^{-1}(t_{a,i}(x))\bigr)
 \subseteq\{t_{a,j}(x):j\in F_a\}
 \quad(x\in A_a,\ i\in I_a).
\end{equation}
For every $\rho>0$ and finite open cover $\mathcal U$ of $X$, there
are small elementary subgroupoids $H_0\subseteq H_1$ with the
following properties:
\begin{enumerate}
\item their presentations have common parameter sets, and on every
      refined piece the index sets are respectively $I_a$ and $F_a$;
\item every compact level closure of $H_1$ lies in a member of
      $\mathcal U$, and all level boundaries are $K$-null;
\item $QH_0\cup H_0Q\subseteq H_1$, and the same inclusions hold
      for $\overline{H_0}\subseteq\overline{H_1}$;
\item $c_K(X\setminus H_0^{(0)})
       \leq c_K(X\setminus W_0)+\rho$.
\end{enumerate}
Empty refined pieces may be discarded.  No index is removed from a
surviving piece.
\end{lem}
\begin{proof}
First realize the coordinate maps by finitely many original bisections.
The first-chart decomposition and the uniform tail estimate
\eqref{eq:regular-relation-uniform-tails} give a finite truncation
for each $t_{a,i}$.  Intersect their retained source sets on each
base.  By invariance, the loss of the old core footprint is bounded by
\[
 \sum_a|I_a|\sum_{i\in F_a}
 c_K(\text{omitted source set for }t_{a,i}).
\]
Make this less than $\rho/3$.  These accuracies are chosen after the
finite index sets are fixed.  Refine the retained bases by the finite
tuples of first-chart choices.  All pieces are regular, and on each
piece every coordinate agrees with the partial homeomorphism of an
original open bisection, including any compositions in the atlas.

Choose a finite small-boundary open cover $(U_b)_b$ whose closures
refine $\mathcal U$.  Its first-membership partition
$P_b=U_b\setminus\bigcup_{c<b}U_c$ has regular indicators.
Refine the bases further by the memberships of every coordinate
in these $P_b$.  Denote the resulting finitely many regular bases
by $A_\alpha$, and their open bisection sections by
$\gamma_{\alpha,i}$.  Write $h_{\alpha,i}=r\circ\gamma_{\alpha,i}$.
The original coordinate agreements hold on the fixed common conull
set.  Choose a common open domain $\Omega_\alpha$ for the sections
where the levels are pairwise distinct and
$h_{\alpha,i}(x)\in U_{b(\alpha,i)}$.  This domain contains
$A_\alpha$ modulo a common null set.

Restrict the sections to $\Omega_\alpha$, and define the open sets
\[
 L_{\alpha,i}=\bigcup_{j\in F_\alpha}
       \gamma_{\alpha,j}\gamma_{\alpha,i}^{-1},\qquad
 O_\alpha=\Omega_\alpha\cap
       \bigcap_{i\in I_\alpha}h_{\alpha,i}^{-1}
          \bigl(X\setminus s(Q\setminus L_{\alpha,i})\bigr).
\]
Since $Q\setminus L_{\alpha,i}$ is compact, $O_\alpha$ is open.
Condition~\eqref{eq:regular-chart-core-capture} and principality
give $c_K(A_\alpha\setminus O_\alpha)=0$: an arrow with endpoints
$h_{\alpha,i}(x),h_{\alpha,j}(x)$ must equal
$\gamma_{\alpha,j}(x)\gamma_{\alpha,i}(x)^{-1}$.
For every $x\in O_\alpha$, all $Q$-arrows sourced at a core level
belong to $L_{\alpha,i}$.  Injectivity of $h_{\alpha,i}$ ensures
that their representations in this set use the same parameter $x$.
Thus the exact extension condition holds throughout $O_\alpha$,
including points outside the original Borel base.

Apply Lemma~\ref{lem:regular-uniform-compact-bases} to obtain
compact $T_\alpha\subseteq O_\alpha$ with
$T_\alpha\setminus A_\alpha$ common null and with
\[
 \sum_\alpha |I_\alpha|c_K(A_\alpha\setminus T_\alpha)<\rho/3.
\]
The compact level sets $h_{\alpha,i}(T_\alpha)$ are disjoint
modulo common null sets.  Indeed, they differ from subsets of the
globally disjoint original levels only by images of common null sets.
Each intersection between two distinct compact levels is consequently
compact and $K$-null. Let $Z$ be the finite union of these intersections,
and set
\[
 N_\alpha=T_\alpha\cap\bigcup_{i\in F_\alpha}
                    h_{\alpha,i}^{-1}(Z).
\]
This is a compact $K$-null subset of $T_\alpha$.
Choose open neighborhoods $W_\alpha\supseteq N_\alpha$ with
\[
 \sum_\alpha|I_\alpha|c_K(W_\alpha)<\rho/3,
\]
using Lemma~\ref{lem:regular-uniform-compact-bases}, and put
$T'_\alpha=T_\alpha\setminus W_\alpha$.  All their compact level
images are now globally disjoint.

By normality and continuity, the nonempty $T'_\alpha$ have open
neighborhoods with closures inside $O_\alpha$, with globally
disjoint compact level images, and with each level closure contained
in its prescribed member of $\mathcal U$.  Cover each $T'_\alpha$
by finitely many SBP basic open sets with closures in that neighborhood.
Their union $V_\alpha$ has $K$-null boundary.  Define $H_1$ by
\eqref{eq:small-elementary-multisections} over these $V_\alpha$
using all indices $F_\alpha$, and define $H_0$ using $I_\alpha$ over
the same parameter sets.

The level closures are disjoint, so these are open elementary
subgroupoids and their closures are compact elementary subgroupoids.
For each matrix-coordinate map the closure of its image over
$V_\alpha$ is precisely its image over $\overline{V_\alpha}$:
the latter image is compact, and continuity gives both inclusions.
The exact extension established on $O_\alpha$ therefore applies
also to the closures.  Symmetry of $Q$ gives right extension by
inverting the left inclusion.  Moreover,
\[
 \partial h_{\alpha,i}(V_\alpha)
       =h_{\alpha,i}(\partial V_\alpha),
\]
since $\overline{V_\alpha}$ lies in the coordinate domain.
Invariance makes every level boundary $K$-null.

Finally, invariance and disjointness give
\[
 \mu(H_0^{(0)})=\sum_\alpha|I_\alpha|\mu(V_\alpha)
       \geq\sum_\alpha|I_\alpha|\mu(T'_\alpha).
\]
The three loss estimates prove the footprint assertion uniformly over
$K$.  All refinements have retained the original finite index sets.
\end{proof}

\begin{thm}[Small subgroupoids under relative SBP]
\label{thm:SBP-groupoid-castles}
Let $\mathcal G$ be a second-countable locally compact Hausdorff
principal \'etale groupoid with compact metrizable unit space $X$.
Let $K$ be a nonempty compact convex family of invariant probability
measures.  Assume relative $K$-SBP and measured amenability of its
orbit relation for every $\mu\in K$; topological amenability of
$\mathcal G$ is sufficient.  For every compact
$Q\subseteq\mathcal G$, finite open cover $\mathcal U$ of $X$,
and $\varepsilon,\eta>0$, there are small elementary subgroupoids
$H_0\subseteq H_1$ such that:
\begin{enumerate}
\item $c_K(X\setminus H_0^{(0)})<\varepsilon$;
\item every compact level closure of $H_1$ is contained in a member
      of $\mathcal U$, all level closures are disjoint, and all
      level boundaries are $K$-null;
\item $QH_0\cup H_0Q\subseteq H_1$, with the same inclusions for
      their compact closures;
\item for $j=0,1$ and every $x\in H_j^{(0)}$,
\[
             |\partial_Q(H_j)_x|<\eta|(H_j)_x|.
\]
\end{enumerate}
Minimality is not required.  The same boundary estimates hold for
the source fibers of the closure subgroupoids.
\end{thm}
\begin{proof}
Replace $Q$ by $Q\cup Q^{-1}\cup X$; the stronger conclusion for this
symmetric compact set implies the assertion for the original $Q$.
Choose a finite open bisection
cover of $Q$.  Its cardinality gives a bound $d$ for both the source
and range fiber cardinalities of $Q$.  By SBP and compactness, refine
this to finitely many open bisections $B_1,\ldots,B_m$ covering $Q$,
whose source domains have small boundary and closures inside the
original bisection domains.  Their induced maps $\beta_l$ are
admissible by Proposition~\ref{prop:principal-groupoid-SBP-atlas}.

Fix $0<\theta<1$ with
\begin{equation}\label{eq:small-subgroupoid-boundary-parameters}
       \frac{(d+1)\theta}{1-\theta}<\eta,
\end{equation}
and choose $\delta>0$ so small that
$\delta/\theta<\varepsilon/2$.  Theorem~\ref{thm:uniform-CFW-SBP}
supplies (RF), hence a bounded regular finite subrelation $S$ with
\[
 \sup_{\mu\in K}\sum_{l=1}^m
 \mu\{y\in\operatorname{dom}\beta_l:
                   (\beta_l(y),y)\notin S\}<\delta.
\]
Lemma~\ref{lem:regular-finite-transversals} supplies finite regular
matrix charts for $S$.  Call $y$ \emph{good} if every applicable
$\beta_l(y)$ lies in $[y]_S$.  Its complement is a regular set:
membership in $S$ is a finite union of coincidences with its
admissible matrix charts.  Its $c_K$-measure is less than $\delta$.
Every $Q$-neighbor of a good point stays in its $S$-class.

Split each regular base by the subset $I$ of good levels in its
class.  On such a piece let $k$ be the total number of levels and
$b=k-|I|$ the number of bad levels.  Invariance gives, for every
$\mu\in K$,
\[
       \sum_{\text{pieces }a}b_a\mu(A_a)<\delta.
\]
Discard the pieces with $b_a>\theta k_a$.  Their whole footprint
has measure bounded by
\[
 \sum_{b_a>\theta k_a} k_a\mu(A_a)
       \leq\frac1\theta\sum_a b_a\mu(A_a)<\frac\delta\theta.
\]
On the surviving pieces $|I_a|\geq(1-\theta)k_a>0$, and their
good-level footprint $W_0$ consequently satisfies
\[
                c_K(X\setminus W_0)<\delta/\theta.
\]
For this last estimate, the removed-class footprint and the remaining
bad levels have total measure
\[
 \sum_{b_a>\theta k_a}k_a\mu(A_a)
       +\sum_{b_a\leq\theta k_a}b_a\mu(A_a)
       \leq\frac1\theta\sum_a b_a\mu(A_a).
\]
In particular the estimate does not multiply the escape error by
the largest class size.

Apply Lemma~\ref{lem:regular-open-castle-lifting}, with
$\rho<\varepsilon/2$, to these surviving charts and their good
indices.  This gives the first three conclusions and preserves
every ratio $b_a/k_a$ on its refined pieces.  It remains to check
the full boundary estimates, which are pointwise and survive the
topological lifting.

Fix a parameter in one of the resulting towers and a source $x$ in
an enlarged level. Put $F=(H_1)_x$, and let $A\subseteq F$ consist of
the arrows whose ranges lie in the core levels of this tower.
Then $|F|=k$ and $|A|=k-b$. If $x$ is a core point, $A=(H_0)_x$.
The validity region in the lifting proof gives $QA\subseteq F$
for every enlarged source $x$. Only the $b$ omitted levels can have a
$Q$-neighbor outside $F$.  Consequently
\[
          |\partial_Q^+F|\leq db,\qquad
          |\partial_Q^-F|\leq b.
\]
For the core, exact extension gives
$\partial_Q^+A\subseteq F\setminus A$, hence
$|\partial_Q^+A|\leq b$.  If a core level belongs to
$\partial_Q^-A$, a $Q$-arrow carries it to an omitted level in
$F\setminus A$.  Reverse that arrow.  Since $Q$ is symmetric and
has source degree at most $d$, there are at most $db$ such core
levels.  Thus
\[
          |\partial_Q^+A|\leq b,\qquad
          |\partial_Q^-A|\leq db.
\]
The same argument works for source fibers based at any enlarged
level, and for any core source fiber.  It follows that
\[
 \frac{|\partial_QF|}{k}\leq(d+1)\theta<\eta,
 \qquad
 \frac{|\partial_QA|}{k-b}
      \leq\frac{(d+1)\theta}{1-\theta}<\eta.
\]
These bounds use only index sizes and exact extension, so they hold
also on the compact closures.  This proves the theorem.
\end{proof}

\Needspace{16\baselineskip}
\begin{cor}[SBP and COS imply $\mathcal Z$-absorption]
\label{cor:SBP-COS-groupoid-Z}
Let $\mathcal G$ be a second-countable locally compact Hausdorff
principal minimal \'etale groupoid with infinite compact metrizable
unit space $X$. Suppose that $\mathcal G$ is topologically amenable,
$K=M(\mathcal G)\neq\varnothing$, and $(X,K)$ has relative SBP.
Put $A=\mathrm{C}_r^*(\mathcal G)$. Assume Cuntz comparison of open sets
(COS) in the
following fixed-amplification form: there are $\lambda\in(0,1]$ and
$m\in\mathbb N$, $m\geq1$, such that, for all open $U,V\subseteq X$,
\begin{equation}\label{eq:groupoid-weak-COS}
 \mu(U)<\lambda\mu(V)\quad(\mu\in K)
 \quad\Longrightarrow\quad
 \varphi_U\precsim\varphi_V^{\oplus m}.
\end{equation}
Here $\varphi_W\in \mathrm{C}(X)_+$ has cozero set $W$, and comparison is
in $A\otimes\mathcal K$. Then the original untwisted algebra $A$
is tracially $\mathcal Z$-absorbing and
\[
                         A\cong A\otimes\mathcal Z.
\]
\end{cor}
\begin{proof}
We first obtain the norm-central cones from
Theorem~\ref{thm:SBP-groupoid-castles}. We give the relative nesting
and defect estimates explicitly, since the small remainder is to be
compared using \eqref{eq:groupoid-weak-COS}.

\emph{Relative nesting.}
A fortified elementary system $E_-\Subset E_+$ will mean two finite
multisection presentations with the same coordinate sections and
index sets, over bases $V_a^-\Subset V_a^+$, respectively.
Their level boundaries may be taken $K$-null. Given such a system,
a compact arrow control $Q$, an integer $R\geq1$, and $\eta>0$,
we claim that there is another fortified system
$F_-\Subset F_+$, faithfully nested in $E_-$ with multiplicity
at least $R$, such that $F_+$ has an exact $Q$-extension to a
small elementary system and
\begin{equation}\label{eq:SBP-COS-relative-nesting}
 c_K\bigl(\overline{E_-^{(0)}}\setminus F_-^{(0)}\bigr)<\eta.
\end{equation}
Nesting here is nesting of the finite multisection coordinates:
each $E_-$-coordinate acts on an entire selected level, and this
action is again a selected coordinate. Our construction will ensure
that every $F_+$-tower meets every old $E_-$-level; in particular
the nesting is faithful. The multiplicity counts the selected levels
in such an old level, as in \cite[Definitions~12.1 and~12.9]{MW}.

Replace $Q$ by $Q\cup Q^{-1}\cup X$. Minimality and principality,
with $X$ infinite, imply that every orbit is infinite and dense.
For each point of $X$ and each old $E_-$-level, choose $R$ arrows
from that point with distinct ranges in the level. Shrinking their
bisections and then using compactness gives finitely many source
patches covering $X$, each with these recurrence sections and
pairwise disjoint ranges in each old level. Their restricted
closures, together with $X$, form a compact set $K_0$. Set
$L_2=\overline{E_+}\cup X$ and choose a symmetric compact set $P$
containing $QL_2K_0$ and $X$.
Theorem~\ref{thm:SBP-groupoid-castles} gives $H_0\subseteq H_1$,
with $PH_0\subseteq H_1$ and
$c_K(X\setminus H_0^{(0)})<\delta$. Use their common parameter
presentations, with nonempty core index sets $I_a\subseteq J_a$.

Choose $i_0\in I_a$ in each tower. At each parameter, select
\emph{every} original core point that belongs to $E_-^{(0)}$,
the $R$ recurrence images of the $i_0$-point in each old level,
and the $E_-$-saturation of these points. These points all occur
among the $J_a$-levels at the same parameter. Indeed, a selected
point has the form $r(ew)$ with $s(w)$ a core point,
$w\in K_0$ and $e\in E_-$; identities include the original core
points. Both $ew$ and $qew$, for any composable $q\in Q$, belong
to $P$. Since their source is a core point, they belong to $H_1$.
Consequently
\[
                       q=(qew)(ew)^{-1}\in H_1
\]
at the same parameter. This gives the required exact extension.

To make the selection a whole-level selection, first remove the
pullbacks, through all $J_a$-sections, of the old level boundaries.
Their union is closed relative to each open base $W_a$ and
is common null by invariance. Removing it gives an open set
$O_a\subseteq W_a$ with $c_K(W_a\setminus O_a)=0$. Membership in the
old levels is then locally constant. Choose a recurrence source
patch for the $i_0$-point. The recurrence sections and their old
coordinate translates equal unique $H_1$-coordinate sections.
These equality conditions are open, since the sections are \'etale,
and closed on their common domains, since $\mathcal G$ is Hausdorff.
There are only finitely many possible target indices. Thus a finite
open parameter refinement fixes all selected indices and all old
coordinate actions on them.

The open bases $W_a$ have $K$-null boundaries, so $\mu\mapsto\mu(W_a)$
is continuous on $K$. Continuous compactly supported cutoffs
increasing to $1_{O_a}$ capture $W_a$ uniformly by Dini's theorem.
Cover a retained compact set by finitely many SBP opens whose closures
lie in the fixed-pattern cover. Disjointize them by subtracting
the closures of earlier sets; only their common-null boundaries
are omitted. Compact shrinking and two further SBP neighbourhood
choices give $V_t^-\Subset V_t^+$ with disjoint plus closures,
$K$-null boundaries, and the same fixed equations everywhere.
Choose the total lost core mass, weighted by $|I_a|$ after the
finite shapes are fixed, smaller than $\kappa$.
Use all pairs of selected indices over $V_t^\pm$ for $F_\pm$;
use all $J_a$-indices over $V_t^+$ for the extension system.
Every selected plus level lies in one old minus level. For an
$E_+$-coordinate whose source contains that old minus level, its
restriction to the selected level agrees with the corresponding $E_-$-coordinate. By saturation
and the fixed equations, it sends the entire selected level onto a
selected level of the same new tower, at the same new parameter.
Thus left multiplication of a new matrix coordinate by an old one
is either empty or one entire new coordinate. This verifies nesting
in both $E_-$ and $E_+$ in the sense of \cite[Definition~12.1]{MW}.
The recurrence sections supply at least $R$ selected levels in
every old level in each retained tower, proving faithfulness and
the asserted multiplicity.
Every retained original core point in $E_-^{(0)}$ was selected.
Invariance, the weighted loss estimate, and the null old boundaries
therefore give
\[
 c_K\bigl(\overline{E_-^{(0)}}\setminus F_-^{(0)}\bigr)
 \leq c_K(X\setminus H_0^{(0)})+\kappa<\delta+\kappa.
\]
Choose $\delta+\kappa<\eta$ to obtain
\eqref{eq:SBP-COS-relative-nesting}. In particular, the estimate
is relative to the old footprint.

\emph{Norm-central cones and their full defect.}
We show that for every finite $\mathcal F\subset A$, $\varepsilon>0$,
$n\geq1$, and $\rho>0$, there is a cpc order-zero map
$\psi:\mathrm{M}_n\to A$ satisfying
\begin{equation}\label{eq:SBP-COS-central-cone}
 \|[\psi(b),a]\|<\varepsilon
 \quad(a\in\mathcal F,\ \|b\|\leq1),\qquad
 q:=1-\psi(1)\in \mathrm{C}(X)_+,
\end{equation}
with $\{q>0\}\subset U$ for an open set $U$ having $c_K(U)<\rho$.

Approximate the finite tests by finite sums of functions
$f_i\in \mathrm{C}_c(\mathcal G)$ supported in precompact open bisections
$B_i$. Let $L$ be a symmetric compact set containing $X$, all
$\overline{B_i}$, and their inverses. Choose a large integer $k$.
Apply Theorem~\ref{thm:SBP-groupoid-castles} with control $L^k$.
Refine the resulting finite matrix-arrow coordinates by
$\{B_i,\mathcal G\setminus\operatorname{supp}f_i\}$ and by small
coefficient oscillation. This is a finite open parameter refinement.
After the index sets are fixed, compact extraction and SBP
thickening retain all but an arbitrarily small weighted core mass.
Take all full indices $J_a$ over outer and inner bases for
$E_+$ and $E_-$, respectively. Exact capture preserves the same
parameter, so every retained core level is $L^k$-deep in $E_+$.
Let $W_k$ be the union of these levels over the inner bases.
We may require
\begin{equation}\label{eq:SBP-COS-deep-footprint}
                c_K(X\setminus W_k)<\eta_0.
\end{equation}
The coordinate refinement also ensures that any matrix-arrow
coordinate meeting $\operatorname{supp}f_i$ is wholly in $B_i$.

Apply the relative nesting construction with control $Q=L$,
multiplicity $nN$, and loss $\eta_1$, for a large integer $N$, obtaining $F_\pm$.
The finite matrix-coordinate algebra of a presentation is the
direct sum of its constant matrix algebras. Faithful nesting
gives a reduced copy of the $E_+$-coordinate algebra in the
$F_+$-coordinate algebra; its commutant has the finite
multiplicity matrix blocks described in \cite[Section~12]{MW}.
Group the multiplicity indices into $n$ equal parts in each block.
This defines a matrix embedding
\[
 \phi:\mathrm{M}_n\longrightarrow
 \mathrm{C}^*\bigl(\operatorname{Comm}(F_+,E_+)\bigr),
\]
leaving fewer than $n$ indices per multiplicity block unused.
Write $\iota$ for the canonical inclusion of this commutant
algebra in the $F_+$-coordinate algebra. Since every multiplicity
is at least $nN$, invariance and the finite counting estimate in
\cite[Lemma~13.10]{MW} give
\begin{equation}\label{eq:SBP-COS-unused-levels}
                     c_K(Q_2)\leq 1/N,
\end{equation}
where $Q_2$ is the union of unused $F_+$-levels.

We make the finite staircase explicit. Let $K'$ be the symmetric
collection of $E_+$ matrix-arrow coordinates meeting some
$\operatorname{supp}f_i$, together with the unit coordinates.
Each tested coordinate is contained in some $B_i$ or $B_i^{-1}$,
and the unit coordinates lie in $X$. Hence every coordinate in $K'$
is contained in $L$. On each finite index set $J_a$, call an index safe
if every original $L$-arrow from that level is captured by another
$E_+$ coordinate at the same parameter, throughout the plus base.
Let $D_a$ be the remaining indices and join indices by the $K'$ arrows.
Exact $L^k$ capture implies that every $K'$ path of length less than
$k$ from a retained core index ends at a safe index. Thus its graph
distance to $D_a$ is at least $k$. Define
\[
 e_a(j)=\min\{1,\operatorname{dist}_{K'}(j,D_a)/k\},
\]
with value one when the distance is infinite. The resulting diagonal
contraction $e$ is one on the retained core, zero on the unsafe levels,
and changes by at most $1/k$ across a tested coordinate. In particular,
original $L$ capture holds on its support. This is the finite-coordinate
staircase of \cite[Lemma~13.9]{MW}. Choose invariant
coordinate cutoffs $f$ and $g$, supported on $F_+$ and $E_+$,
respectively, with $f=1$ on $\overline{F_-^{(0)}}$ and
$g=1$ on $\overline{E_-^{(0)}}$. For a constant matrix-coordinate
algebra, $\Psi_f$ denotes its cpc order-zero realization in $A$,
obtained by multiplying each matrix-arrow section by its
invariant base cutoff; define $\Psi_g$ in the same way.
Then set
\begin{equation}\label{eq:SBP-COS-cone-formula}
                  \psi(b)=\Psi_f(\iota\phi(b))\Psi_g(e).
\end{equation}
The two images commute by \cite[Lemma~12.7]{MW}, since
$F_+$ is nested in $E_-$ and $g=1$ there. Multiplication by
the commuting positive contraction $\Psi_g(e)$ consequently
preserves complete positivity, contractivity, and order zero.
For $f_i$, let $\varphi(f_i)$ be its sampled constant coordinate
kernel and let $\omega_i$ bound its coefficient oscillation.
For clarity, the two estimates in \cite[Lemma~13.6]{MW} here are
\[
 \|\Psi_g(e)f_i-\Psi_g(e\varphi(f_i))\|\leq\omega_i,
 \qquad
 \|f_i\Psi_g(e)-\Psi_g(\varphi(f_i)e)\|\leq\omega_i .
\]
Indeed, every arrow in $\operatorname{supp}f_i$ with source or
range in a level on which $e$ is nonzero is captured in $E_+$,
by safety and symmetry of $L$. On a captured coordinate the actual
and sampled values differ by at most $\omega_i$. Each difference
is supported in the single bisection $B_i$, on which the reduced
norm is the supremum norm. Since $\Psi_f(\iota\phi(b))$ commutes
with all of $\Psi_g(\mathrm{C}^*(E_+))$ by \cite[Lemma~12.7]{MW}, these
two bounds give, for $\|b\|\leq1$,
\[
 \|[\psi(b),f_i]\|
 \leq2\omega_i+\|[e,\varphi(f_i)]\|
 \leq2\omega_i+\frac{\|f_i\|_\infty}{k}.
\]
The last estimate is \cite[Lemma~13.8]{MW}: the coordinates belonging
to the one bisection $B_i$ form a partial permutation, so no matrix-rank
factor enters the bound. Sum these estimates for the fixed Fourier
tests. Taking the oscillations small and $k$ large proves the norm estimate in
\eqref{eq:SBP-COS-central-cone}. This part of the construction
uses only finite coordinates, exact extension, and nesting.

The unit of the matrix cone is diagonal, and the formula shows
\begin{equation}\label{eq:SBP-COS-whole-defect-support}
 \{1-\psi(1)>0\}\subseteq Q_0\cup Q_1\cup Q_2,
 \qquad
 \begin{cases}
 Q_0=X\setminus W_k,\\
 Q_1=\overline{E_-^{(0)}}\setminus F_-^{(0)}.
 \end{cases}
\end{equation}
Indeed, outside this union all the coordinate cutoffs, the
staircase, and $\phi(1)$ equal one. The first two sets are compact,
with capacities less than $\eta_0$ and $\eta_1$, by
\eqref{eq:SBP-COS-deep-footprint} and
\eqref{eq:SBP-COS-relative-nesting}; the third is open and
satisfies \eqref{eq:SBP-COS-unused-levels}. A compact set of
capacity below $b$ has an open neighborhood of capacity below
any prescribed larger bound. To see this, take continuous
$u_j\downarrow1_Z$ and use compactness of $K$ to choose one
index uniformly, then take a superlevel set of $u_j$.
Choose $\eta_0,\eta_1$ small and $N$ large. Open neighborhoods
of $Q_0,Q_1$, together with $Q_2$, give the required open $U$
with $c_K(U)<\rho$. This proves the small-support assertion
for the entire diagonal defect in
\eqref{eq:SBP-COS-central-cone}.

\emph{Comparison of the defect.}
Fix $0\neq a\in A_+$. Principality implies effectiveness, so
\cite[Lemma~13.5]{MW} gives $0\neq b\in \mathrm{C}(X)_+$ with
$b\precsim a$. Choose a nonempty open set
$O\subseteq\{b>c\}$ for some $c>0$.
By minimality and infinite principal orbits, there are $m$
bisections with common open source $V$, extending over
$\overline V$, and with pairwise disjoint ranges contained
in $O$. Choose $\varphi_V\in \mathrm{C}(X)_+$ with cozero set $V$.
The corresponding cutoff normalizers $z_j\in A$ satisfy
\[
 z_j^*z_j=\varphi_V,\qquad
 z_jz_j^*\perp z_lz_l^*\ (j\neq l),\qquad
 \sum_{j=1}^m z_jz_j^*\in\operatorname{Her}(b).
\]
Their row matrix therefore gives
\begin{equation}\label{eq:SBP-COS-reserved-amplification}
 \varphi_V^{\oplus m}\sim\sum_{j=1}^m z_jz_j^*
                         \precsim b\precsim a.
\end{equation}
Every invariant probability measure has full support. Since
$\mu\mapsto\mu(V)$ is lower semicontinuous on compact $K$,
\[
                         \alpha:=\inf_{\mu\in K}\mu(V)>0.
\]
Apply the cone construction with $\rho<\lambda\alpha$.
Its full defect belongs to $\mathrm{C}_0(U)_+$ and hence is Cuntz-below
the strictly positive function $\varphi_U$ of $\mathrm{C}_0(U)$.
Equations \eqref{eq:groupoid-weak-COS} and
\eqref{eq:SBP-COS-reserved-amplification} now give
\[
 1-\psi(1)\precsim\varphi_U
       \precsim\varphi_V^{\oplus m}\precsim a.
\]
Together with \eqref{eq:SBP-COS-central-cone}, this is tracial
$\mathcal Z$-absorption. Finally, $A$ is simple, separable,
unital and nuclear. Theorem~4.1 of Hirshberg--Orovitz
\cite{HO} gives $A\cong A\otimes\mathcal Z$.
\end{proof}

\begin{cor}[Almost elementariness in measure in the principal amenable case]
\label{cor:SBP-groupoid-AE-measure}
Let $\mathcal G$ be second countable, locally compact, Hausdorff,
principal and \'etale, with compact metrizable unit space.  Suppose
$\mathcal G$ is minimal, $M(\mathcal G)\neq\varnothing$, and its
orbit relation is amenable for every invariant probability measure;
topological amenability is sufficient.  Then relative SBP for
$M(\mathcal G)$ is equivalent to almost elementariness in measure
in the sense of \cite[Definition~5.7]{MW}.
\end{cor}
\begin{proof}
The invariant-measure family $M(\mathcal G)$ is compact and convex.
Apply Theorem~\ref{thm:SBP-groupoid-castles} to the symmetric compact
set $Q\cup Q^{-1}\cup X$ for any prescribed compact $Q$.
The finite multisection presentations of $H_0$ and $H_1$ are the
open castles used in \cite[Definition~5.7]{MW}.  Their common
parameters and exact extension make the core presentation
$Q$-extendable to the enlarged one, including the left and right
conditions of \cite[Definitions~5.1 and~5.3]{MW}.  It refines the prescribed
cover and has uniformly small remainder.  If $\varepsilon<1$,
its footprint is nonempty and open, hence full by minimality.
These are all the conditions of almost elementariness in measure.
The converse is \cite[Theorem~6.12]{MW}.
\end{proof}

\begin{rem}[Small subgroupoids and relation URP]
\label{rem:SBP-small-subgroupoids-URP}
The small subgroupoids in
Theorem~\ref{thm:SBP-groupoid-castles} are open principal elementary
subgroupoids of the original $G$: their finite multisections are
closed under multiplication and inverse.  Their unit spaces are the
corresponding footprints, and their closures are compact subgroupoids
because the closed bases remain inside the coordinate domains and
the closed levels are disjoint.

Under the same SBP and amenability hypotheses, uniform CFW already
gives strong uniform AF approximation and hence split relation URP
in Definition~\ref{defn:relation-URP}; see
Corollary~\ref{cor:local-AH-Reiter-and-SBP}.  By
Theorem~\ref{thm:relation-URP-local-AH}, this is equivalent to split
uniform local AH Cartan approximation.  The small-subgroupoid theorem
adds a realization by open multisections inside the given topological
groupoid.  The local approximants are not asserted to be nested.

Without SBP, relation URP and local AH approximation retain their
parameter spaces.  One must also relax the requirement that the levels
refine every prescribed open cover: disjoint open levels of arbitrarily
small uniform remainder refining every cover already imply relative
SBP by \cite[Lemma~6.11]{MW}; in particular, almost elementariness in
measure implies SBP by \cite[Theorem~6.12]{MW}.  Thus the arbitrary-cover
refinement in Corollary~\ref{cor:SBP-groupoid-AE-measure} belongs to the
SBP conclusion, rather than to a general formulation of relation URP.
\end{rem}

\begin{rem}[The group-action case]
\label{rem:SBP-small-subgroupoids-actions}
For a transformation groupoid $G=\Gamma\ltimes X$, the group
coordinates of bisections are locally constant.  The compact bases
in the construction admit finite subdivision on which these
coordinates are constant, giving ordinary group-labelled towers.
Freeness identifies the source-fibre coordinates with finite subsets
of $\Gamma$, so the corresponding F\o lner estimates take the usual
finite-shape form.

For free actions of countable amenable groups on compact metrizable
spaces, the action counterpart is the equivalence of SBP and almost
finiteness in measure in \cite[Theorem~5.6]{KS}.  Its proof uses an
Ornstein--Weiss tower construction.  This implies action URP, as stated
in \cite[Theorem~3.3]{NiuURP}.  The latter property, defined in
\cite[Definition~3.1 and Lemma~3.2]{NiuURP}, does not require arbitrarily
fine levels.  Minimality is unnecessary for this
Kerr--Szab\'o implication.  Relation URP above is the intrinsic
completed-Cartan formulation, and its local AH interpretation retains
the parameter spaces when SBP is not assumed.
\end{rem}

\section{Regular twists and uniform untwisting}
\label{sec:regular-twists}

The uniform finite charts also remove the scalar cohomology obstruction.
Feldman--Moore~\cite[Theorem~6]{FMI} prove measured cohomology vanishing
in degrees at least two for hyperfinite relations.  Here the issue is
to obtain one gauge whose graph coefficients lie in the designated
completed diagonal, simultaneously for all measures in $K$.

\begin{defn}\label{defn:regular-two-cocycle}
Write an arrow as $(y,x)$, with source $x$ and range $y$.
A \emph{regular scalar twist} is a normalized Borel function
$\sigma(y,z,x)\in\mathbb T$ on triples in one $R$-class, satisfying
\[
 \sigma(y,z,x)\sigma(y,x,w)
   =\sigma(y,z,w)\sigma(z,x,w),
 \qquad
 \sigma(y,y,x)=\sigma(y,x,x)=1,
\]
such that, for all admissible composable charts $s,t$,
\[
 \sigma_{s,t}(x)=\sigma(stx,tx,x)
       \in 1_{\operatorname{dom}(st)}D_K.
\]
A \emph{regular gauge} is a normalized Borel function
$b:R\to\mathbb T$ with
$b(sx,x)\in1_{\operatorname{dom}s}D_K$ for every admissible $s$.
Its coboundary is
\[
 (\delta b)(y,z,x)=\frac{b(y,z)b(z,x)}{b(y,x)}.
\]
All these identities are understood on one invariant set conull for
every measure in $K$.
\end{defn}

The regularity condition makes twisted multiplication of finite
regular chart kernels well defined:
\[
 (f*_\sigma g)(y,x)
       =\sum_{z\in[x]_R}f(y,z)g(z,x)\sigma(y,z,x).
\]
The inverse formula uses the regular coefficient
$\overline{\sigma(y,x,y)}$.
The diagonal expectation is restriction to $(x,x)$, and
$\tau_\mu(f)=\int f(x,x)\,d\mu(x)$ are the designated traces.
Denote the bounded uniform tracial completion by $M_K^\sigma(R)$.
It can equivalently be constructed in the twisted regular
representation as in Section~\ref{sec:uniform-relation-algebra}.

\begin{thm}[Uniformly regular untwisting]
\label{thm:uniform-regular-untwisting}
Suppose that $R$ has nested finite regular matrix-chart subrelations
$R_n$ which exhaust $R$ on one invariant $K$-conull set and have
uniformly vanishing escape along admissible charts.
Every regular scalar twist is a regular coboundary:
\[
 H^2_{\mathrm{reg},K}(R,\mathbb T)=0.
\]
Consequently there is a diagonal-fixing isomorphism
\[
 (D_K,M_K^\sigma(R))\cong(D_K,M_K^1(R))
\]
preserving the expectation and every designated trace.
Coordinate refinement, diagonal density, simplicity, and aperiodicity
are unnecessary for this assertion.
\end{thm}

\begin{proof}
Let $r_n(x)$ be the base representative of the $R_n$-class of $x$.
This map is given, on finitely many regular level projections, by
inverses of the matrix-coordinate maps.  For $(x,y)\in R_n$ put
\[
 \gamma_n(x,y)=\sigma(x,y,r_n(x)).
\]
The cocycle identity, and $r_n(x)=r_n(y)$, give
$\delta\gamma_n=\sigma|_{R_n}$.  Evaluating along the finitely many
coordinate pieces shows that its chart coefficients are regular.

We choose these trivializations compatibly.  Start with $b_1=\gamma_1$.
If $b_n$ is chosen, then
\[
 \eta_n=\frac{b_n}{\gamma_{n+1}|_{R_n}}
\]
is a regular $1$-cocycle on $R_n$.  Put
$a_n(x)=\eta_n(x,r_n(x))$.  Finite chart regularity gives
$a_n\in\mathcal U(D_K)$, and the $1$-cocycle identity yields
$\eta_n(x,y)=a_n(x)\overline{a_n(y)}$ on $R_n$.  Define, on $R_{n+1}$,
\[
 b_{n+1}(x,y)=
       \gamma_{n+1}(x,y)a_n(x)\overline{a_n(y)}.
\]
Then $\delta b_{n+1}=\sigma|_{R_{n+1}}$ and
$b_{n+1}|_{R_n}=b_n$.  The vertex phase $a_n$ is already defined on
all units, so no extension across the larger finite classes is needed.

The compatible union $b=\bigcup_n b_n$ is a normalized Borel gauge on
the common conull relation.  For an admissible partial map $s$, set
\[
 q_n(s)=1_{\{x\in\operatorname{dom}s:(sx,x)\in R_n\}}.
\]
This is a regular projection.  Indeed, the finite coordinate maps
cover $R_n$, and membership of $(sx,x)$ is a finite union of their
regular coincidence tests with $s$.  The coefficient
$q_n(s)b_n(sx,x)$ is regular for the same reason.  Extending it by
$1$ on the rest of $\operatorname{dom}s$ gives a regular phase
$b_{n,s}$, and
\[
 \norm{b_{n,s}-b(s\,\cdot,\cdot)}_{2,K}
 \leq
 2\sqrt{\sup_{\mu\in K}
       \mu\{x\in\operatorname{dom}s:(sx,x)\notin R_n\}}
 \longrightarrow0.
\]
Thus $b$ is regular along every admissible chart, and
$\delta b=\sigma$ on the exhausted relation.

Finally define $J(f)(y,x)=b(y,x)f(y,x)$.  The identity
$\delta b=\sigma$ gives
$J(f*_\sigma g)=J(f)*J(g)$ and preserves the involution.
Multiplication by $b(y,x)$ on each regular-representation Hilbert
space implements $J$, so reduced norms agree.  Since $b(x,x)=1$,
$J$ fixes the diagonal and the expectation and preserves each
$\tau_\mu$.  Regularity of $b$ and $\overline b$ makes $J$ and its
inverse maps of the regular chart algebras.  They extend to the
asserted bounded uniform completions.
\end{proof}

\begin{prop}[Regular scalar presentation of a groupoid twist]
\label{prop:SBP-regularizes-twist}
Let $G$ be a principal, Hausdorff, second-countable \'etale groupoid with
compact metrizable unit space $X$, and let $K$ be a nonempty compact family of
invariant probability measures.  If $(X,K)$ has relative SBP, every
locally trivial circle twist $\Sigma\to G$ has a normalized Borel
section whose scalar cocycle is regular for the regular atlas of
Proposition~\ref{prop:principal-groupoid-SBP-atlas}.
The resulting regular twisted relation completion agrees with the
designated bounded uniform completion of $\mathrm{C}_r^*(G,\Sigma)$.
\end{prop}

\begin{proof}
Cover $G$ by countably many open bisections carrying continuous local
twist sections.  Refine their sources by small-boundary open sets $U$
with compact closure inside the source of the ambient bisection.
Their arrow maps and local sections extend continuously over
$\overline U$.  Domain and range indicators are regular; invariance
makes the range boundaries universally null.

The coincidence sets of two such arrow maps are relatively clopen.
They are closed by Hausdorffness and open because, in a principal
groupoid, equality of endpoints means equality of arrows, so the
coincidence set is the source of a bisection intersection.
Their boundaries are therefore contained in the domain boundaries.
Include inverses and finite compositions to obtain a regular atlas.

Put the canonical section on the unit arrows first, and assign every
other arrow to the first local-section chart containing it.  This
defines a normalized Borel section $q$.  On a fixed chart its
assignment pieces are regular coincidence projections.
On a multiplication piece where the selected left, right, and product
sections are $i,j,k$, their product-comparison phase is continuous.
The compact-closure choice keeps all three local sections defined on
the closure of that piece.  A bounded complex-valued extension of the
phase into the closed unit disk, multiplied by the piece projection,
belongs to $D_K$.  A circle-valued continuous extension is not needed.

The countable assignment pieces exhaust each chart domain.  Their
finite-union projections have continuous measure profiles on $K$,
and compactness and Dini give uniformly vanishing tails.  Applying
this to the left, right, and product chart assignments shows that
$\sigma_{s,t}$ is a bounded uniform $2$-norm limit of the regular
finite-piece phases.  Thus it belongs to $D_K$.
The same decomposition and corner pullbacks treat arbitrary
admissible $s,t$.

Local graph lifts on small-boundary supports belong to the
completion of $\mathrm{C}_r^*(G,\Sigma)$ by continuous cutoffs inside those
supports; compactness and Dini again give the uniform tails.
The graph lift from the Borel section $q$ differs from the chosen
local continuous lift by a source-corner phase in
$\mathcal U(1_{\operatorname{dom}s}D_K)$.  Multiplying by this
regular phase therefore gives the same completed generators.
Conversely compactly supported continuous twisted sections decompose
into finitely many local bisection sections with diagonal
coefficients.  These two density statements identify the completions.
\end{proof}

\begin{cor}[Amenable twisted dynamical realization]
\label{cor:amenable-twisted-dynamical-model}
For a regular uniform measure relation with relative SBP and
amenability for every $\mu\in K$, every regular scalar twist has a
diagonal-fixing, trace-preserving untwisting.
For a principal groupoid as in
Proposition~\ref{prop:SBP-regularizes-twist}, this conclusion applies
to every locally trivial circle twist if $K$ is convex and all designated measured
orbit relations are amenable.
If the relation is also aperiodic on one common $K$-conull set, its
twisted completed Cartan pair has the untwisted Cantor
$\mathbb Z$-model of Theorem~\ref{thm:local-cuts-dynamical-model},
with a compact designated measure family affinely identified with $K$.
\end{cor}
\begin{proof}
Theorem~\ref{thm:uniform-CFW-SBP} supplies the finite regular charts.
Apply Theorem~\ref{thm:uniform-regular-untwisting}, and, in the
aperiodic case, Theorem~\ref{thm:local-cuts-dynamical-model}.
\end{proof}

\begin{rem}\label{rem:untwisting-boundaries}
The gauge is uniformly regular and Borel; it need not be continuous.
Thus the theorem does not assert continuous triviality of the
original topological twist or an isomorphism of the original
$\mathrm{C}^*$-algebras.  Circle bundles can already be nontrivial over a
principal two-point-class groupoid; see the example discussed
in~\cite[Section~2]{ANSZ}.
Genuine isotropy is also outside the relation statement: the
amenable one-unit groupoid $\mathbb Z^2$ has nontrivial scalar
multipliers.  For example
$\omega((m,n),(p,q))=\exp(2\pi i\theta np)$ is not a coboundary
when $\theta\notin\mathbb Z$, since a scalar coboundary on an
abelian group is symmetric whereas this multiplier is not.
Finally, untwisting is not canonical: two regular gauges differ
by a regular $1$-cocycle, so the residual phases of
Remark~\ref{rem:uniform-relation-pair-residual-phases} remain relevant.
\end{rem}

\section{Controlled quasi-invariant families and uniform measure classes}
\label{sec:controlled-RN}

This section extends the finite-cut, AF, and twist constructions to
quasi-invariant probabilities with a common controlled Radon--Nikodym
cocycle. The completion is a state completion rather than a tracial
completion. Some assertions using invariance require replacement:
the full-group metric is no longer bi-invariant, and arbitrary partial
maps need not extend to full-group elements.

\subsection{The uniform structure on a family of measure classes}

For a family $K$ of probabilities set
\[
 c_K(A)=\sup_{\mu\in K}\mu(A)\qquad(A\text{ Borel}).
\]
Two representing families $K,L$ define the same uniform measure class
if their capacities are uniformly absolutely continuous with respect
to each other: for every $\varepsilon>0$ there is $\delta>0$ such that
$c_K(A)<\delta$ implies $c_L(A)<\varepsilon$, and conversely.
The common null ideal is then the same.

\begin{prop}[Representative independence on bounded sets]
\label{prop:RN-capacity-topology}
Uniformly equivalent capacities give the same uniform $2$-topology
on bounded Borel functions, the same regular diagonal $D_K=D_L$,
and the same relative SBP. For diagonal states $\mu\circ \mathbb{E}$, they
also give the same symmetric uniform state topology on bounded
regular relation kernels. Equality of the particular states or of
the displayed seminorms is not asserted.
\end{prop}
\begin{proof}
For $0\leq f\leq M$ and $t>0$,
\[
 c_K\{f>t\}\leq t^{-1}\sup_K\int f\,d\mu,
 \qquad
 \sup_L\int f\,d\nu\leq t+M c_L\{f>t\}.
\]
These estimates in both directions prove equivalence of bounded
uniform integral convergence. Apply them to $|b|^2$, and then to
$\mathbb{E}(a^*a)$ and $\mathbb{E}(aa^*)$. Boundary nullity depends only on the common
null ideal, proving the SBP assertion.
\end{proof}

Bare measure classes do not determine this uniform structure. For
one probability $\mu$ and a bounded Borel function $f$,
\[
 \sup_{\nu\sim\mu,\ \nu(X)=1}\|f\|_{2,\nu}
       =\|f\|_{\infty,\mu}.
\]
The lower bound follows by concentrating almost all the mass on a
positive-measure set where $|f|$ is close to its essential supremum,
while keeping a positive background density. Thus taking all
representatives imposes a different topology.

\subsection{Controlled charts and the nonsingular pseudogroup}
\leavevmode\par

Let $X$ be compact metrizable, $R$ countable Borel, and $K$ a nonempty
compact convex family of $R$-quasi-invariant probabilities. Retain
the diagonal, support, coincidence, and corner-pullback regularity
conditions of Definition~\ref{defn:regular-relation-atlas}, replacing
measure preservation by the following condition. For every chart $s$
in the countable atlas, closed under inverse and composition, there
is a common positive regular function $j_s$ on its range such that
\[
 s_*(\mu|_{\operatorname{dom}s})
   =j_s\,\mu|_{\operatorname{ran}s}\quad(\mu\in K),
 \qquad C_s=\|j_s\|_\infty<\infty,
 \quad C_{s^{-1}}<\infty.
\]
The cocycle identities hold on one invariant common conull set.
There is no bound required uniformly over the entire countable atlas.
In particular, for a bounded function $f$,
\begin{equation}\label{eq:rn-change-variable}
 \int_{\operatorname{dom}s}f(sx)\,d\mu(x)
  =\int_{\operatorname{ran}s}f(v)j_s(v)\,d\mu(v),
\end{equation}
and
\[
 \mu(sA)\leq C_{s^{-1}}\mu(A),\qquad
 \mu(s^{-1}B)\leq C_s\mu(B).
\]
For a regular source cut $P\subseteq\operatorname{dom}s$ the
Jacobian of $s|_P$, extended by zero, is $1_{sP}j_s$. The range
projection is regular by the chart corner-pullback hypothesis.
For a finite injective pasting with disjoint sources and ranges,
these densities add on disjoint ranges. For a composite $st$ the
density at $v$ is $j_s(v)j_t(s^{-1}v)$ on its range. These formulas
give common bounded regular Jacobians, also for the inverse maps.

For an arbitrary partial map $\varphi$, require regular domain and range
projections and regular graph coefficients for the map and its inverse.
Let $\varphi_N$ be its first $N$ chart pieces, with omitted source
and range sets $P_N,Q_N$. These are regular projections decreasing
to zero, so Dini gives $c_K(P_N)+c_K(Q_N)\to0$. For Borel $A$,
\[
 c_K(\varphi A)\leq c_K(\varphi_N A)+c_K(Q_N),\qquad
 c_K(\varphi^{-1} A)\leq c_K(\varphi_N^{-1}A)+c_K(P_N).
\]
The finite-pasting estimates imply uniform nonsingularity in both
directions. For $f\in D_K$, the zero-extended pullbacks of $f$ along
$\varphi$ and $\varphi_N$ differ in uniform $2$-norm by at most
$\|f\|c_K(P_N)^{1/2}$. The analogous inverse bound uses $Q_N$.
Thus both corner pullbacks preserve $D_K$. These maps form the
nonsingular regular pseudogroup.
Bounded Jacobians on arbitrary countable pastings are not required.
The finite approximations constructed below have controlled generators.

\begin{prop}[Full groups and group presentations]
\label{prop:RN-full-group}
On the nonsingular regular full group the metric
\[
 d_K^{\pm}(T,S)=c_K\{T\ne S\}+c_K\{T^{-1}\ne S^{-1}\}
\]
is complete and separable, and gives a topological group. The source
disagreement metric remains left invariant, but bi-invariance is not
asserted. With relative SBP, the relation has a countable group
presentation by controlled regular involutions, without amenability.
\end{prop}
\begin{proof}
For a Cauchy sequence select a subsequence with summable consecutive
source and inverse disagreements. Borel--Cantelli and quasi-invariant
saturation give eventual pointwise limits that are inverse Borel
maps on a common invariant conull set. Their graph coefficients are
bounded uniform $2$-limits of regular coefficients. Approximation by
a fixed late map in both disagreement coordinates proves uniform
nonsingularity of the limit. For instance
\[
 c_K(TE)\leq c_K(T_nE)+c_K\{T^{-1}\ne T_n^{-1}\}.
\]
The analogous inverse inequality proves the other transport estimate.
Thus the limit is regular. Continuity of products follows from
uniform nonsingularity of each fixed translating map; inversion
exchanges the two coordinates. Separability follows by embedding the
forward and inverse first-chart coefficient sequences in the two
separable uniform $L^2$ sequence spaces used in
Proposition~\ref{prop:uniform-relation-full-group-complete}.

For the presentation, separate the source and target of each nonfixed
chart arrow by the nested finite regular SBP partitions. Disjointize
according to the first separating partition and its two different
cells. This gives countably many regular restrictions whose domains
and ranges are disjoint. Extend each restriction to an involution by
its inverse on its range and the identity elsewhere. These are finite
controlled pastings. The group they generate is countable and its
orbits are exactly $R$ on a common invariant conull set.
\end{proof}

\begin{example}[Partial maps need not extend]
\label{ex:RN-partial-extension-fails}
Take $X=\mathbb N\cup\{\infty\}$ with its one-point compactification,
$\mathbb N=\{0,1,\ldots\}$, the relation with classes $\mathbb N$
and $\{\infty\}$, and $\mu(n)=2^{-n-1}$, $\mu(\infty)=0$.
The shift $s(n)=n+1$, $s(\infty)=\infty$ maps $X$ onto
$X\setminus\{0\}$. Its Jacobian on its range is $2$, and its inverse
Jacobian is $1/2$, modulo the null point. It is a controlled regular
partial map with total domain. It cannot agree almost everywhere
with a full-group element: every $n$ has positive mass, and the
prescribed images omit $0$. The countable atlas of words in the shift and its inverse is
regular, and SBP holds. Thus Lemma~\ref{lem:uniform-partial-map-extension}
does not extend to the nonsingular setting, even with controlled
Jacobians. In the completed algebra its graph is a proper isometry.
\end{example}

\subsection{Measurewise amenability gives regular uniform Reiter kernels}

For one $\mu\in K$, measured amenability gives hyperfiniteness by the
ordinary nonsingular Connes--Feldman--Weiss theorem~\cite{CFW}.  Uniform probabilities
on the finite measured classes therefore give measurable probability
kernels with arbitrarily small integrated Reiter errors on any finite
test set $F$.  The fixed-measure regularization can be performed in
\[
 \lambda=\mu+\sum_{s\in F}s_*(\mu|_{\operatorname{dom}s}).
\]
Quasi-invariance gives $\lambda\ll\mu$.  Truncate to the first finitely
many disjoint atlas arrow coordinates, place the tail on the identity,
and approximate the finite coefficients in $L^1(\lambda)$ by regular
coefficients obtained from continuous functions.  Normalize using the
continuous formula in the proof of
Theorem~\ref{thm:measurewise-hyperfinite-uniform-Reiter}.  Its pointwise row error is
at most twice the sum of the coefficient errors.  If $e(x)$ denotes
the resulting row $\ell^1$ error, the change in a test error is at most
\[
 \int e(x)\,d\mu(x)
 +\int e(v)\,d\big(s_*(\mu|_{\operatorname{dom}s})\big)(v)
 \leq\int e\,d\lambda.
\]
Thus the infimum over finite regular kernels is zero for each $\mu$.
Their integrated error profiles remain continuous affine functions on
$K$, because the pointwise error functions lie in $D_K$.  The convex
upward-set separation and barycenter proof of
Theorem~\ref{thm:measurewise-hyperfinite-uniform-Reiter} applies without
change and gives regular uniform Reiter kernels.  This part does not
require a common cocycle or bounded Jacobians.

\subsection{Weighted local derandomization}
\leavevmode\par

Consider finitely many controlled window charts $\psi_1,\ldots,\psi_Q$
and $W_x=\{\psi_ax:x\in\operatorname{dom}\psi_a\}$.  Relative SBP
supplies a finite regular coloring $c:X\to\{1,\ldots,k\}$ and a regular
exceptional root set $E$ of arbitrarily small uniform measure, so that
distinct vertices of $W_x$ have distinct colors for $x\notin E$.
Let $A$ be a finite alphabet, and let bounded nonnegative regular local
costs $h_x$ satisfy the finite conditional-table hypotheses of
Lemma~\ref{lem:regular-local-derandomization}.
Then, for every $\alpha>0$, one can choose a finite-valued regular
label function $a:X\to A$ such that, simultaneously for every $\mu\in K$,
\begin{equation}\label{eq:weighted-derand}
 \int_{X\setminus E}h_x(a)\,d\mu(x)
 \leq\int_{X\setminus E}\mathbb{E} h_x\,d\mu(x)+\alpha.
\end{equation}

To prove this, process colors in order as in the proof of
Lemma~\ref{lem:regular-local-derandomization}.  For a current
vertex $v$, deduplicate the affected good roots
$x_a(v)=\psi_a^{-1}v$ using the regular first-occurrence indicators
$I_a(v)$.  If $H(x,v,u)$ is the current conditional average when the
label at $v$ is fixed to $u$, define
\[
 B(v,u)=\sum_a j_{\psi_a}(v)I_a(v)H(x_a(v),v,u).
\]
Every $B(\cdot,u)$ is a bounded element of $D_K$; the common cocycle
makes it independent of $\mu$.  It need not be finite-valued, so an
exact regular argmin must not be assumed.  Instead choose finite-valued
regular step functions $\widetilde B(\cdot,u)$ with
\[
 \sup_{\mu\in K}\int|B(v,u)-\widetilde B(v,u)|\,d\mu(v)
        <\frac{\alpha}{2k|A|}
 \quad(u\in A).
\]
Choose the first minimizer of $\widetilde B(v,u)$ on the current color.
This is a regular finite-valued rule.  Pointwise,
\[
 B(v,a(v))-\frac1{|A|}\sum_uB(v,u)
 \leq2\sum_u|B(v,u)-\widetilde B(v,u)|.
\]

The required mass transport is now weighted:
\begin{equation}\label{eq:weighted-mass-transport}
 \int\sum_a j_{\psi_a}(v)I_a(v)F(x_a(v),v)\,d\mu(v)
 =\int_{X\setminus E}\sum_{v\in W_x}F(x,v)\,d\mu(x).
\end{equation}
Apply the chart change of variables to each term; the first-occurrence
indicators retain exactly one copy of each distinct root--vertex pair.
Restricting this identity to the current color shows that its total
integrated conditional-cost increase is below $\alpha/k$ for every
$\mu$ simultaneously.  After the $k$ stages one obtains
\eqref{eq:weighted-derand}.  No global random field is used.

\subsection{Quantitative regular rounding}

Let $m_x=\sum_{j=1}^J f_j(x)\delta_{\phi_jx}$ be a finite regular
probability kernel supported on controlled charts.  Let $F$ be a finite
set of controlled tests.  For every $\eta>0$, there is a regular finite
subrelation $S\subset R$ with classes of size at most $J$, generated by
finitely many controlled regular partial maps, such that
\begin{equation}\label{eq:rn-rounding}
 \sum_{s\in F}\mu\{x\in\operatorname{dom}s:(sx,x)\notin S\}
 \leq\sum_{s\in F}\int_{\operatorname{dom}s}
                  \|m_{sx}-m_x\|_1\,d\mu(x)+\eta
 \qquad(\mu\in K).
\end{equation}

Indeed, choose rational regular rows $\widetilde m$ with uniform
integrated row error below $q$.  By \eqref{eq:rn-change-variable},
\[
 \sum_{s\in F}\int_{\operatorname{dom}s}
       \|\widetilde m_{sx}-\widetilde m_x\|_1\,d\mu
 \leq\sum_{s\in F}\int_{\operatorname{dom}s}\|m_{sx}-m_x\|_1\,d\mu
       +q\sum_{s\in F}(1+C_s).
\]
Choose $q$ so that the last summand is below $\eta/4$.  The token sets,
shared finite priorities, and minimum-token center map $C(x)$ from the
proof of Theorem~\ref{thm:regular-Reiter-rounding} give the same
pointwise finite-average bound, independently of
invariance:
\[
 \mathbb{E}1_{\{C(x)\ne C(sx)\}}
 \leq\|\widetilde m_x-\widetilde m_{sx}\|_1
       +\binom{2L}{2}/N.
\]
Choose $N$ so that the summed collision term is below $\eta/4$.
Use the window consisting of the controlled charts $\phi_j$ and
$\phi_js$, and choose its bad-root set $E$ with
$|F|\sup_\mu\mu(E)<\eta/4$.  Apply
\eqref{eq:weighted-derand} with $\alpha<\eta/4$ to the root cost
$\sum_s1_{\operatorname{dom}s}1_{\{C(x)\ne C(sx)\}}$.  The four
errors give \eqref{eq:rn-rounding}.  Define $S$ by $C(x)=C(y)$.
Injectivity of the $J$ support charts bounds each fibre by $J$.
The finite restrictions of $\phi_k^{-1}\phi_j$ between selected-chart
pieces enumerate all $S$-arrows and are controlled and regular.

For arbitrary admissible tests, keep only finitely many of their
regular first-chart pieces, so that the omitted source sets have
arbitrarily small total uniform measure.  These finite restrictions
are controlled.  Apply \eqref{eq:rn-rounding} there and add the source
tail error.

If a finite Reiter kernel itself uses admissible support maps
$\phi_1,\ldots,\phi_J$ with unbounded Jacobians, first fix the
controlled restrictions of the test maps. Truncate each $\phi_j$ to
finitely many of its regular first-chart pieces and transfer its
omitted weight to the identity. For the resulting controlled kernel
$m'$ the row error $e(x)=\|m'_x-m_x\|_1$ is bounded by twice the
sum of the indicators of the omitted source sets. Thus
$\sup_{\mu\in K}\int e\,d\mu$ can be made arbitrarily small.
For each fixed controlled test $s$, the change in its integrated
variation is at most $(1+C_s)\sup_{\mu\in K}\int e\,d\mu$ by
\eqref{eq:rn-change-variable}. Applying \eqref{eq:rn-rounding} to
$m'$ now proves (RF) for the full admissible pseudogroup, with
controlled finite generators.

\subsection{Nesting and the uniform AF conclusion}

Let $H,G$ be bounded regular finite subrelations whose finite arrow
covers are controlled.  Let $h_1,\ldots,h_q$ cover $H$, including the
identity, and let $E$ be the sources of $H$-arrows missed by $G$.
Set $B=[E]_H$ and $C=X\setminus B$.  All these sets are regular and
\[
 \mu(B)\leq A_H\mu(E),\qquad
 A_H=\sum_{i=1}^q C_{h_i^{-1}}.
\]
The repaired relation $F=G|_C\sqcup H|_B$ contains $H$, remains bounded
finite and controlled regularly generated, and satisfies
\[
 \mu\{(sx,x)\notin F\}
 \leq\mu\{(sx,x)\notin G\}+(1+C_s)\mu(B).
\]
The image estimate follows by covering $B$ by the $h_i$-images of $E$;
the escape estimate follows because only a source or target in $B$
can lose a captured arrow.  Constants are fixed before choosing $G$.
If (RF) gives total error below $\delta$ on the old arrow cover and
the new tests $F_0$, the repaired summed error on $F_0$ is at most
\[
 \left(1+A_H\sum_{s\in F_0}(1+C_s)\right)\delta.
\]
Thus choose $\delta$ to achieve any prescribed summable new tolerance.
Take the new tests at stage $n$ to be the first $n$ controlled
atlas charts. Summable errors and Borel--Cantelli give an increasing
bounded regular finite exhaustion on one invariant set conull for
every designated measure, since quasi-invariance preserves nullity
under saturation. Every admissible test has uniformly vanishing
first-chart source tails. For each fixed finite truncation its escape
is bounded by the sum of the corresponding atlas-chart escapes, so
the exhaustion has uniformly vanishing escape for the full
pseudogroup. In Lemma~\ref{lem:regular-finite-transversals} the
matrix coordinates are finite regular pastings of bounded-length
words in the stage generators, hence remain controlled under all
the coordinate refinements.

\begin{thm}[Weighted uniform Connes--Feldman--Weiss]
\label{thm:RN-uniform-CFW}
\leavevmode\par
Under the controlled common-co\-cycle hypotheses and relative SBP, the following are equivalent:
\begin{enumerate}
\item measured amenability for every designated measure;
\item regular uniform Reiter amenability;
\item regular finite cuts (RF), with controlled finite generators;
\item a regular AF exhaustion with coordinate refinement and dense
      finite diagonals.
\end{enumerate}
\end{thm}
\begin{proof}
The constructions above prove the forward implications. A regular AF exhaustion supplies regular uniform class probabilities and uniformly small escape, hence regular Reiter kernels. A summable subsequence is a measured Reiter witness for every designated measure. No invariant measure is needed.\end{proof}
\begin{rem}
The finite-cut condition in this theorem requires generators that are finite controlled pastings of atlas charts. This is the condition supplied by the rounding construction; it implies the original (RF) condition. No assertion here identifies arbitrary countable pastings with bounded-Jacobian charts.\end{rem}

\subsection{Scalar twists}

The compatible finite-class cochain argument of
Theorem~\ref{thm:uniform-regular-untwisting} does not use invariance.
On a regular AF exhaustion, every regular circle-valued
$2$-cocycle is a regular coboundary.  Its chartwise gauge approximation
still has error at most twice the square root of the uniform source
escape.  The resulting algebraic untwisting fixes the diagonal and the
expectation, hence preserves each designated state $\mu\circ \mathbb{E}$.
It retains the positive Radon--Nikodym $1$-cocycle.  Extension to the
chosen nonsingular operator-algebra completion requires that completion
to have been defined with a suitable uniform strong-star topology.

\subsection{State completions, AF reconstruction, and scalar twists}
\label{subsec:controlled-RN-Cartan}

We use the controlled quasi-invariant setting of this section.  In
particular, every designated measure is quasi-invariant, and the
jacobians of every controlled atlas chart, its inverse, and every
finite controlled regular pasting are bounded uniformly over the
designated family. A controlled finite-chart kernel is a finite sum
of regular diagonal coefficients times these graph operators.
Arbitrary admissible countable pastings need not have bounded
Jacobians. The common-cocycle assumption
is needed for the preceding simultaneous rounding argument.  The
completion results below only need these uniform chart bounds.

Let $A_R$ denote the norm closure of the controlled finite-chart kernel
algebra, represented by convolution on the relation counting spaces.
Its diagonal expectation is $\mathbb{E}:A_R\to D_K$, and its designated states
are $\phi_\mu=\mu\circ \mathbb{E}$.  These states are generally not traces.
Put
\[
 p_K(a)^2=\sup_{\mu\in K}
       \bigl(\phi_\mu(a^*a)+\phi_\mu(aa^*)\bigr).
\]
On the diagonal this is $\sqrt{2}$ times the uniform $2$-norm.
The faithful common-kernel quotient is understood throughout.

For an arbitrary admissible partial map $\varphi$, let $\varphi_N$
retain its first $N$ controlled atlas pieces. The omitted source and
range projections have uniformly vanishing measure. Computed in the
measured relation representations,
\[
 p_K(v_\varphi-v_{\varphi_N})^2
 =\sup_{\mu\in K}\bigl(
   \mu(\operatorname{dom}\varphi\setminus\operatorname{dom}\varphi_N)
  +\mu(\operatorname{ran}\varphi\setminus\operatorname{ran}\varphi_N)
                     \bigr)\longrightarrow0.
\]
These truncations are contractions. Thus the bounded state completion
below contains every admissible graph operator, although such an
operator need not belong to the norm core $A_R$.

\begin{lem}[Fixed-chart multiplication]
\label{lem:RN-fixed-chart-multiplication}
For every controlled finite-chart element $b$ there is a finite constant
$C_b$ such that
\[
 p_K(ab)\leq C_b p_K(a),\qquad
 p_K(ba)\leq C_b p_K(a)\qquad(a\in A_R).
\]
Multiplication by a fixed element of $A_R$ is therefore continuous
for $p_K$ on norm-bounded sets.
\end{lem}

\begin{proof}
Write $j_s$ for the density of $s_*(\mu|_{\operatorname{dom}s})$
relative to $\mu|_{\operatorname{ran}s}$.  The formulas below also
apply when this density depends on $\mu$, provided its bounds are
uniform.  For a positive $t\in A_R$, covariance of the diagonal
expectation and change of variables give
\[
 \phi_\mu(v_s^*tv_s)\leq\|j_s\|_\infty\phi_\mu(t),
 \qquad
 \phi_\mu(v_stv_s^*)\leq\|j_s^{-1}\|_\infty\phi_\mu(t).
\]
The integrals on the right include the appropriate corner support;
discarding that support only enlarges them.  If
$b=\sum_{i=1}^m f_i v_{s_i}$, the operator inequality
$(\sum_i x_i)(\sum_i x_i)^*\leq m\sum_i x_ix_i^*$, applied to
$x_i=f_i v_{s_i}t^{1/2}$, yields constants $c_b,c_{b^*}<\infty$
such that
\[
 \phi_\mu(btb^*)\leq c_b\phi_\mu(t),\qquad
 \phi_\mu(b^*tb)\leq c_{b^*}\phi_\mu(t)
\]
for every $\mu$.  Together with
$b^*b,bb^*\leq\|b\|^2 1$, these prove the displayed estimates.
For a fixed $b\in A_R$, approximate $b$ in operator norm by such
controlled finite-chart elements. On a bounded set of $a$'s, the resulting
error is uniformly small because $p_K(x)\leq\sqrt{2}\|x\|$.
\end{proof}

\begin{prop}[Bounded uniform state completion]
\label{prop:RN-state-completion}
The operator-norm bounded $p_K$-Cauchy sequences in $A_R$ form a
$\mathrm{C}^*$-subalgebra of $\ell^\infty(A_R)$.  Its $p_K$-null sequences
form a closed two-sided ideal.  The quotient, denoted $M_K(R)$,
has a faithful realization in the product of the designated
measured-relation von Neumann algebras.  Each $\phi_\mu$ extends,
and its closed norm balls are complete for $p_K$.

The unit ball is the $p_K$ completion of the unit ball of $A_R$.
The diagonal completion is $D_K$, and the diagonal expectation
extends to a faithful conditional expectation $M_K(R)\to D_K$.
For a singleton designated family this construction is the usual
nonsingular measured-relation von Neumann algebra.
\end{prop}

\begin{proof}
We give the product argument since the tracial estimate
\eqref{eq:two-mult} is unavailable.  Let $(a_n),(b_n)$ be bounded
$p_K$-Cauchy sequences, and let $M_a,M_b$ bound their norms.  In
each designated GNS representation, with cyclic vector $\Omega_\mu$,
fix a late index $N$ and write
\[
 \begin{split}
 (a_nb_n-a_mb_m)\Omega_\mu
  ={}&(a_n-a_m)b_N\Omega_\mu\\
     &+(a_n-a_m)(b_n-b_N)\Omega_\mu
       +a_m(b_n-b_m)\Omega_\mu.
 \end{split}
\]
The first term tends uniformly to zero for fixed $N$ by
Lemma~\ref{lem:RN-fixed-chart-multiplication}.  The remaining terms
are bounded by $2M_a p_K(b_n-b_N)$ and
$M_a p_K(b_n-b_m)$.  For the adjoint use the separate expansion
\[
 \begin{split}
 (b_n^*a_n^*-b_m^*a_m^*)\Omega_\mu
  ={}&(b_n^*-b_m^*)a_N^*\Omega_\mu\\
     &+(b_n^*-b_m^*)(a_n^*-a_N^*)\Omega_\mu
       +b_m^*(a_n^*-a_m^*)\Omega_\mu.
 \end{split}
\]
First choose $N$ to make the tails small, and then choose $n,m$
large.  Both vector norms tend uniformly to zero, proving that
$(a_nb_n)$ is $p_K$-Cauchy.  The same argument with one null
sequence proves the ideal assertion.  Norm closure follows from
$p_K(x)\leq\sqrt{2}\|x\|$.  Involution is $p_K$-isometric.

For any fixed $c\in A_R$, the vectors
$a_nc\Omega_\mu$ and $a_n^*c\Omega_\mu$ are Cauchy.  Such vectors
are dense in each GNS space, so $(a_n)$ has a bounded strong-star
limit in each designated measured algebra.  The resulting product
representation is faithful: the functions
$\phi_\mu(a_n^*a_n)$ and $\phi_\mu(a_na_n^*)$ are uniformly Cauchy
in $\mu$, by Cauchy--Schwarz and the common norm bound.  If every
fiber limit is zero, their pointwise limits are zero, hence their
uniform limits are zero and the sequence is $p_K$-null.
The limits of $\phi_\mu(a_n)$ give the extended states.  A standard
diagonal choice of bounded representatives proves bounded
$p_K$-completeness.

The quotient is a $\mathrm{C}^*$-algebra.  Every contraction in it has a
contractive lift in the sequence algebra: apply radial truncation
to any lift, using continuous functional calculus in that algebra.
This proves the unit-ball assertion without a tracial truncation
inequality.  Schwarz's inequality for $\mathbb{E}$ gives
$p_K(\mathbb{E}a)\leq p_K(a)$, so $\mathbb{E}$ extends and has range $D_K$.
Faithfulness follows in each designated measured algebra from its
faithful diagonal expectation, and then from the faithful product
realization.  For one measure, bounded $p_K$ convergence is the
ordinary strong-star convergence for its faithful diagonal state;
Kaplansky density gives the last assertion.
\end{proof}

\begin{rem}[Meaning of uniform strong-star convergence]
The seminorm $p_K$ is equivalent, on bounded sets, to uniform
convergence on the common test vectors $c\Omega_\mu$ and their
adjoint tests, for each fixed $c$ in the regular core.  It implies pointwise strong-star convergence in the product.
The converse can fail even for the identity relation: take
$X=[0,1]$, $K=\operatorname{Prob}(X)$, and continuous positive
contractions of norm one with pairwise disjoint supports.
Dominated convergence gives strong-star convergence to zero in
every designated representation, whereas $p_K$ stays equal to
$\sqrt{2}$ by the point masses.
\end{rem}

\begin{thm}[Nonsingular forward AF reconstruction]
\label{thm:RN-AF-reconstruction}
Suppose the controlled quasi-invariant relation has nested finite
regular matrix charts whose coordinates are finite controlled
pastings, with coordinate refinement, uniform chart
escape, and uniform diagonal density, as in
Definition~\ref{defn:uniform-AF-hyperfinite-relation}.  There are a
standard unital AF pair $E_B\subset B$ and a compact convex family
$L\subset S(B)$ of diagonal states such that
\[
 (D_K,M_K(R))\cong(E_B^L,B^L).
\]
The state correspondence $H:K\to L$ is an affine homeomorphism.
Quasi-invariant AF path measures with prescribed cocycles are discussed
in Renault~\cite{RenaultAF}.
The original relation is regularly isomorphic, on one universally
conull invariant set, to the AF tail relation with the designated
quasi-invariant path measures.  Its Radon--Nikodym cocycle is
transported by this isomorphism.  Centrality of the path measures,
traciality of the states, simplicity of $B$, and $L=\T(B)$ are not
asserted.
\end{thm}

\begin{proof}
Use the same phase-free graph matrix units and coordinate embeddings
as in Theorem~\ref{thm:uniform-AF-reconstruction}.  They generate
$B_n\subset B_{n+1}$ and the standard AF pair $E_B\subset B$.
Put $H(\mu)=\phi_\mu|_B$.  Finite diagonal integrals and diagonal
density give continuity, affinity, and injectivity exactly as in
that theorem.  Thus $L=H(K)$ has the asserted marking.  These are
diagonal states, not necessarily traces.

We check bounded density, the analytic step affected by losing
invariance.  If $s_n$ is the restriction of a chart $s$ to arrows
in $R_n$, and $F_n\subset\operatorname{dom}s$ is its omitted
source set, then
\[
 p_K(v_s-v_{s_n})^2
 =\sup_{\mu\in K}\bigl(\mu(F_n)+\mu(sF_n)\bigr)
 \leq(1+\|j_s^{-1}\|_\infty)\sup_{\mu\in K}\mu(F_n)
 \longrightarrow0.
\]
For fixed $n$, a contraction $z$ in
$\mathrm{C}^*(B_n,D_K)$ is a finite matrix of regular diagonal entries.
Entrywise approximation by later finite diagonals yields
$b_m\in B_m$ with $p_K(b_m-z)\to0$ and an operator-norm bound
depending only on this fixed finite matrix size.  By the bounded
algebra continuity just proved, continuous functional calculus is
$p_K$-continuous along this sequence.  Therefore its radial
truncations at norm one are contractions converging to $z$.
Applying this to $z=v_{s_n}$ and choosing a sufficiently accurate
approximation at each $n$ gives contractive AF approximations to
$v_s$.  Diagonal density is already assumed.  The bounded complete
AF closure consequently contains all regular graph generators and
equals $M_K(R)$.

Conjugation by a diagonal unitary is $p_K$-isometric because
$\phi_\mu=\mu\circ \mathbb{E}$.  Thus the finite diagonal-phase averages
used in the AF proof remain contractions for $p_K$.  They give the
faithful diagonal expectation and the maximal abelian property by
bounded density.  The spatial coding and its simultaneous
Borel--Cantelli injectivity argument use diagonal integrals only.
Removing relation saturations of null sets is valid by
quasi-invariance.  The rest of the spatial proof of
Theorem~\ref{thm:uniform-AF-reconstruction} therefore gives the
asserted regular tail-relation isomorphism.  Change of variables
under that isomorphism transports the designated cocycle.
\end{proof}

\begin{prop}[Scalar untwisting and spatial transport]
\label{prop:RN-untwisting-transport}
Under the finite-chart hypotheses of
Theorem~\ref{thm:uniform-regular-untwisting}, every regular scalar
$2$-cocycle is a regular coboundary in the controlled quasi-invariant
setting.  The induced diagonal-fixing isomorphism of bounded state
completions preserves the expectation and every designated state.
Also, a regular spatial relation isomorphism transporting the
designated probabilities induces an isomorphism of their untwisted
completed Cartan pairs preserving the corresponding states.
\end{prop}

\begin{proof}
The compatible finite-relation gauges in
Theorem~\ref{thm:uniform-regular-untwisting} use no change of
variables and no invariance.  Their regular convergence estimates
are diagonal estimates on uniformly small omitted source sets.
The resulting gauge multiplication is implemented by a unitary on
every counting-space representation, and fixes the diagonal
expectation.  It therefore preserves $p_K$ and extends to the
bounded state completions.  For spatial transport, the unitary
$(U_\theta\xi)(\theta y,\theta x)=\xi(y,x)$ is still unitary
when $\theta$ transports the designated probabilities; it carries
the graph operators and their expectation to the corresponding
ones.  The same state-preserving completion argument applies.
This argument concerns the prescribed controlled nonsingular relation
models; it does not supply an intrinsic nonsingular inverse theorem.
\end{proof}

\begin{prop}[Common modular flow]
\label{prop:RN-common-modular-flow}
If the designated measures have the same regular chart Jacobians,
their modular automorphism groups on the measured relation algebras
come from one strongly continuous flow on the regular core and one
bounded $p_K$-continuous flow on its completion.  With the convention
$j_s=d s_*(\mu|_{\operatorname{dom}s})/d\mu|_{\operatorname{ran}s}$,
its formulas are
\[
 \sigma_t(f)=f\quad(f\in D_K),\qquad
 \sigma_t(v_s)=j_s^{-it}v_s.
\]
Every controlled finite-chart normalizer is analytic for this flow. The
diagonal-enlarged finite-relation algebra $\mathrm{C}^*(B_n,D_K)$ is invariant;
invariance of the finite-dimensional stage $B_n$ itself is not
asserted.
\end{prop}

\begin{proof}
The positive controlled-chart Jacobians satisfy the multiplicative
Radon--Nikodym cocycle identities.  Their boundedness and bounded
inverses give $\log j_s\in1_{\operatorname{ran}s}D_K$ by continuous
functional calculus.  Consequently the displayed formulas preserve
the graph multiplication, involution, and covariance identities.
On each designated counting-space representation they are the usual
modular automorphisms of the faithful diagonal state: an arrow from
$x$ to $sx$ has modular weight equal to the target/source measure
ratio, namely $j_s(sx)^{-1}$.  These fiber automorphisms have the same
formula because the cocycle is common.  They preserve operator norms,
the expectation, and every designated state, and are therefore
$p_K$-isometries.  Finite-chart operator-norm continuity gives a
strongly continuous flow on the core.  Bounded density then gives
$p_K$ continuity on the completion.  Analyticity follows from
$\exp(-iz\log j_s)$ for $z\in\mathbb C$.  Multiplication by this
regular diagonal coefficient keeps a level-$n$ graph operator in
$\mathrm{C}^*(B_n,D_K)$, proving the final assertion.
\end{proof}

\begin{rem}[Changing representatives]
If $\nu=h\mu$ is an equivalent probability, then on the range of
$s$ its new Jacobian is
\[
 j_{s,\nu}(v)=j_{s,\mu}(v)
                 \frac{h(s^{-1}v)}{h(v)}.
\]
Thus changing representatives changes the positive cocycle by a
vertex coboundary.  For uniformly equivalent representing families,
the bounded completion topology is unchanged, but its particular
states and modular-flow formulas can change.  There is no reason for
an arbitrary representative correspondence to be affine or to be a
$p_K$ isometry.  The exact state marking above refers to the chosen
representatives; its spatial isomorphisms also implement uniform
measure-class orbit equivalences.
\end{rem}

\begin{rem}[The marking must retain modular information]
Controlled Jacobians do not restore classification by the abstract
designated simplex alone.  On the binary tail relation of
$X=\{0,1\}^{\mathbb N}$, let $\mu_p$ be the Bernoulli product
measure, $0<p<1$.  Its Jacobians on every fixed prefix-substitution
chart are positive constants with bounded inverses.  Relative SBP
holds since the space has a clopen basis.  Every singleton
$K_p=\{\mu_p\}$ is the same abstract compact convex space.

Nevertheless, if $p\ne1/2$, its measure class contains no invariant
probability.  Invariance under all equal-length prefix substitutions
forces every cylinder of length $n$ to have mass $2^{-n}$, so the
only invariant probability is $\mu_{1/2}$.  The strong law of large
numbers makes $\mu_p$ and $\mu_{1/2}$ mutually singular.
Thus the fair and biased product measures cannot be related by a
measure-class orbit equivalence: transport of the fair invariant
probability would supply an equivalent invariant probability on
the biased side.  The positive Radon--Nikodym $1$-cocycle must be
retained; removal of a scalar circle-valued $2$-twist does not remove
this modular information.
\end{rem}

\subsection{Aperiodic generators and representative changes}
\begin{prop}[Common-cocycle endpoint ordering]\label{prop:RN-endpoint-ordering}
Let $(X,R,K,D)$ have a regular nested AF exhaustion $(H_m)$ on a
common invariant $K$-conull set, with bounded class sizes at each stage.
Assume $K$ is compact, the regular projection profiles are continuous on
$K$, and the designated probabilities are quasi-invariant.  Assume also
that finite regular projection diagonals are bounded uniformly
$2$-norm dense in $D$, equivalently that relative $K$-SBP holds.
Suppose
their Radon--Nikodym cocycles agree, and every fixed controlled finite chart and
its inverse have bounded regular Jacobians.  If $R$ has infinite classes
on that conull set, then $R$ has a uniformly regular nonsingular
single generator.  Conservativity is not required.
\end{prop}

\begin{proof}
Let $F$ be a previously ordered stage whose classes have size at
most $M$, with its order given by finitely many regular pieces.
For a later stage $H\supseteq F$, write its finite matrix charts as
$Y_{v,i}=t_{v,i}(Z_v)$, with $t_{v,1}=\mathrm{id}_{Z_v}$.
Fix an error allowance $\delta>0$. A given $t_{v,i}$ is
admissible and can be a countable pasting. On each base retain the
intersection of finite first-chart truncations for its finitely many
coordinates, forming a regular set $Z'_v\subseteq Z_v$. Uniform
source-tail convergence and uniform nonsingularity of each fixed
$t_{v,i}$ allow these truncations to be chosen so that
\[
 c_K\left(\bigcup_{v,i}t_{v,i}(Z_v\setminus Z'_v)\right)<\delta/2 .
\]
Subdivide each tower over $Z'_v$ and $Z_v\setminus Z'_v$.
On the omitted towers concatenate the old ordered blocks in order
of their least coordinate index. Their entire contribution to the
new endpoint mass is less than $\delta/2$. On the retained bases all
coordinates are finite controlled pastings. We continue to write
$Z_v$ for such a base. For each of its charts there is a common
positive bounded regular function
\[
 h_{v,i}(z)
 =\frac{d\,(t_{v,i}^{-1})_*(\mu|_{Y_{v,i}})}
        {d\,(\mu|_{Z_v})}(z),\qquad
 W_v=\sum_i h_{v,i}.
\]
The function is independent of $\mu\in K$.  Thus for every bounded
Borel $f$,
\begin{equation}\label{eq:rn-generator-disintegration}
 \int_{\bigcup_iY_{v,i}} f\,d\mu
 =\int_{Z_v}\sum_i h_{v,i}(z)f(t_{v,i}z)\,d\mu(z).
\end{equation}
In particular $h_{v,i}/W_v$ are common conditional point weights on
the finite $H$-classes.

Refine $Z_v$ into finitely many regular pieces
according to which of the points $t_{v,i}z$ belong to the same $F$-class.
These tests are regular by the finite-chart coincidence conditions.
On such a piece the index sets $I_1,\ldots,I_q$ of the old blocks are
fixed.  Their conditional masses are the regular functions
\[
 \beta_b(z)=\frac{\sum_{i\in I_b}h_{v,i}(z)}{W_v(z)},
 \qquad \sum_{b=1}^q\beta_b=1.
\]
For $q\geq2$, the sum of the two smallest $\beta_b$ is at most $2/q$.
Choose these two old blocks as the first and last blocks of the new
order, retain their internal old orders, and place all other old blocks
between them in a fixed order.  The union of the two new endpoint sets
is contained in the two selected blocks, so its conditional mass is at
most $2/q$.  For $q=1$ preserve the old order; the same bound, though
coarse, is valid.

An exact comparison of regular real functions need not have regular
level sets.  To make the choice regular, approximate all the
$\beta_b$ by finite-valued functions $\widetilde\beta_b$ on a common
finite regular projection partition.  Select the distinct pair that
minimizes $\widetilde\beta_b+\widetilde\beta_c$, with a fixed tie rule.
The true mass of the selected pair is at most
\[
 \min_{b\ne c}(\beta_b+\beta_c)
       +2\sum_b|\beta_b-\widetilde\beta_b|.
\]
Finite projection density makes the last term arbitrarily small in
uniform $L^1(K)$ on every base.  Multiplying by $W_v$ in
\eqref{eq:rn-generator-disintegration} costs only its fixed sup-norm,
and there are finitely many bases and combinatorial pieces.  We may
therefore bound this integrated excess by $\delta/2$. Together
with the omitted footprint this costs at most $\delta$.  The resulting order uses finitely many regular pieces.
If $A_H,B_H$ denote its maximal and minimal endpoint sets, then
\begin{equation}\label{eq:rn-generator-endpoint-estimate}
 \sup_{\mu\in K}\mu(A_H\cup B_H)
 \leq 2M\sup_{\mu\in K}\int_X\frac{1}{|H[x]|}\,d\mu(x)+\delta,
\end{equation}
because $q\geq |H[x]|/M$.

For the original exhaustion, the finite-valued regular functions
$|H_m[x]|^{-1}$ decrease to zero on the common aperiodic conull set.
Their integrals are continuous on compact $K$ and decrease pointwise
to zero there.  Dini gives
\[
 \sup_{\mu\in K}\int_X |H_m[x]|^{-1}\,d\mu(x)\longrightarrow0.
\]
Starting with the singleton relation, apply
\eqref{eq:rn-generator-endpoint-estimate} recursively, choosing a
sufficiently late new stage after fixing the old size bound $M$, and
then choosing $\delta$ sufficiently small.  This telescopes the
exhaustion to ordered relations $F_n$ with
\[
 \sup_{\mu\in K}\mu(A_n\cup B_n)\longrightarrow0.
\]
Each new order concatenates the previously ordered blocks.  Hence
$A_{n+1}\subset A_n$ and $B_{n+1}\subset B_n$.

Let $c_n$ cyclically permute the finite ordered $F_n$-classes.  For
$m>n$, block compatibility gives
\[
 \{c_m\ne c_n\}\subset A_n,\qquad
 \{c_m^{-1}\ne c_n^{-1}\}\subset B_n.
\]
The intersections $A_\infty,B_\infty$ are common $K$-null.  Remove their
$R$-saturations, which remain common null by quasi-invariance.  On the
remaining invariant set the $c_n$ and $c_n^{-1}$ eventually stabilize
at each point.  Their limits are inverse Borel maps $V,V^{-1}$.
The finite-prefix rank argument shows that every $F_n$-class lies in
a $V$-orbit, so $R_V=R$ there.  Moreover
\[
 c_K\{c_n\ne V\}\leq c_K(A_n),\qquad
 c_K\{c_n^{-1}\ne V^{-1}\}\leq c_K(B_n).
\]

The generator has uniform nonsingularity as well.  For every Borel $E$,
\[
 c_K(V^{-1}E)\leq c_K(c_n^{-1}E)+c_K(A_n),
 \qquad
 c_K(VE)\leq c_K(c_nE)+c_K(B_n).
\]
For fixed $n$, $c_n$ and $c_n^{-1}$ are finite pastings of
admissible coordinate maps. Both therefore transport uniformly
small sets to uniformly small sets, by the first-chart tail
estimates for the nonsingular regular pseudogroup.
First choosing $n$ large and then $c_K(E)$ small proves uniform absolute
continuity for $V$ and $V^{-1}$.

Regular coincidence projections of $V$ with any fixed chart are uniform
$L^2(K)$ limits of those for $c_n$.  Pullbacks of bounded regular
functions are obtained in the same way; uniform nonsingularity permits
passing from continuous or finite-step approximants to all of $D$.
Thus $V$ and $V^{-1}$ preserve the regular diagonal and are uniformly
admissible.  For the phase-free graph unitaries one also has
\[
 \phi_\mu((v_{c_n}-v_V)^*(v_{c_n}-v_V))
    =2\mu\{c_n\ne V\},\qquad
 \phi_\mu((v_{c_n}-v_V)(v_{c_n}-v_V)^*)
    =2\mu\{c_n^{-1}\ne V^{-1}\},
\]
so they converge in the designated uniform strong-star seminorm.
\end{proof}

\begin{rem}[Bounded Jacobian for the new generator]
The conclusion is uniform nonsingularity, not a uniform essential bound
on the Jacobian of $V$.  This distinction is necessary.  Let
$X=\mathbb N\cup\{\infty\}$ be the one-point compactification, let $R$
have the infinite class $\mathbb N$ and the singleton class $\{\infty\}$,
and put $\mu(\{n\})=c\,2^{-2^n}$ and $\mu(\{\infty\})=0$.
Use finite coordinate permutations as the regular atlas.  Every fixed
chart has bounded continuous Jacobian and bounded inverse Jacobian;
SBP and regular finite exhaustion hold, and $R$ is $\mu$-aperiodic.
But for every fixed $C$, sufficiently late point weights are separated
from all other point weights by a factor exceeding $C$.  Any bijection
whose Jacobian and inverse Jacobian are bounded by $C$ must therefore
fix all sufficiently late points.  It cannot generate the infinite
class.  Thus a demand for bounded Jacobian on the resulting single
$\mathbb Z$-generator itself would be stronger than the chartwise
controlled-cocycle hypotheses.

\end{rem}

\begin{prop}[Bounded generator cocycle within the uniform measure class]
\label{prop:RN-generator-reweighting}
Let $V$ be the uniformly nonsingular generator just constructed.
There is an affine homeomorphism from $K$ onto a compact convex family
$K'$ of equivalent representatives with the same uniform measure
class, for which the common regular generator Jacobian and its
inverse are bounded by $2$.
\end{prop}
\begin{proof}
Extend $V$ by the identity off its common invariant conull domain and
put $P=V_*$. Define
\begin{equation}\label{eq:RN-geometric-reweighting}
 \nu_\mu=\frac13\sum_{k\in\mathbb Z}2^{-|k|}P^k\mu.
\end{equation}
The coefficients sum to one. Quasi-invariance gives
$\nu_\mu\sim\mu$, and comparison of consecutive coefficients gives
\[
 \tfrac12\nu_\mu\leq P\nu_\mu\leq2\nu_\mu,
 \qquad \nu_\mu\geq\tfrac13\mu.
\]
For the capacities, $c_K(A)\leq3c_{K'}(A)$. Conversely, split
\eqref{eq:RN-geometric-reweighting} into finitely many powers and a
geometric tail. Each fixed power of $V$ is uniformly nonsingular,
and the tail has arbitrarily small total mass independently of
$\mu$. This proves uniform absolute continuity in the other direction.

The parameter map is affine and weakly continuous, since
$f\circ V^k\in D_K$ for every continuous $f$, and the geometric
series converges uniformly in its integrals. It is injective, with
the exact inverse on finite signed measures
\[
 \mu=5\nu_\mu-2(P\nu_\mu+P^{-1}\nu_\mu).
\]
This follows by multiplying the absolutely convergent operator series
by $5I-2(P+P^{-1})$. Compactness therefore makes the parameter map
an affine homeomorphism onto compact convex $K'$.

We check regularity, not just boundedness, of the new Jacobian.
The original common cocycle gives common Borel densities
$d_k=d(P^k\mu)/d\mu$. First-chart truncations of $V^k$ express $d_k$
as a sum of positive bounded regular density pieces. The omitted
uniform $L^1(K)$ mass is exactly the omitted source mass of that
chart decomposition, and tends to zero. Thus each $d_k$ is uniformly
$L^1(K)$ approximable by bounded positive elements of $D_K$.
The functions
\[
 h=\tfrac13\sum_k2^{-|k|}d_k,\qquad
 h_1=\tfrac13\sum_k2^{-|k|}d_{k+1}
\]
are common densities of $\nu_\mu$ and $P\nu_\mu$ relative to $\mu$.
They have the same uniform $L^1$ approximation property, $h\geq1/3$,
and $h/2\leq h_1\leq2h$. Choose bounded positive regular
approximants $h_n,h_{1,n}$ with $h_n\geq1/3$. Then
\[
 \left|\frac{h_{1,n}}{h_n}-\frac{h_1}{h}\right|
 \leq3|h_{1,n}-h_1|+6|h_n-h|.
\]
Clip the first ratio to $[1/2,2]$. These regular bounded functions
converge uniformly in $L^1(K)$, hence in $L^2(K)$, to $h_1/h$.
Proposition~\ref{prop:RN-capacity-topology} transfers this regularity
to $D_{K'}=D_K$. The new generator Jacobian is $h_1/h$;
its inverse is its reciprocal composed with $V$ and is regular and
bounded by $2$ as well. Every fixed power of $V$ consequently has a
bounded common regular Jacobian in the new representing family.
\end{proof}

\begin{thm}[Uniform measure-class dynamical realization]
\label{thm:RN-dynamical-model}
Under the controlled common-cocycle hypotheses, relative SBP and
measurewise amenability give an AF Cartan model with designated
diagonal states and a regular trivialization of every scalar twist.
If $R$ has infinite classes on a common invariant conull set, it is
uniformly measure-class orbit equivalent to a Cantor $\mathbb Z$-system.
After the uniformly equivalent representative change
\eqref{eq:RN-geometric-reweighting}, its generator and inverse have
common regular Jacobians in $[1/2,2]$. Conservativity is unnecessary.
The model is not asserted minimal or to carry the full quasi-invariant
probability family.
\end{thm}
\begin{proof}
Use the weighted finite cuts, the forward AF reconstruction for states,
the construction of ordered endpoints, and
Proposition~\ref{prop:RN-generator-reweighting}. Close a countable
uniformly dense family of SBP projection partitions under all powers
of $V$. Their separable commutative AF core has compact metrizable
zero-dimensional spectrum $Z$; its induced homeomorphism $W$ is
spatially conjugate to $V$ on common conull sets by diagonal
spatialization. The partitions separate the original measures, so
the designated family is transported affinely and homeomorphically.
Regular graph coincidence coefficients give the required atlas
transport in both directions.

Covariance sends the crossed-product generator to the graph unitary
of $V$ in the state completion. The fixed-point projections of its
nonzero powers vanish for every designated measure. Hence the map
preserves diagonal expectations and the designated states. Every
original chart is a disjoint regular pasting of powers of $V$;
compactness gives uniformly small source and range tails. Its graph
operator is therefore the bounded strong-star limit of finite such
pastings. The crossed product has the same completed Cartan pair.

The designated measures on $Z$ are jointly faithful. If they are
nonatomic, this rules out isolated points and $Z$ is Cantor.
In general take a Cantor space $C$, $p\in C$, the product
$Z\times C$ with $W\times\mathrm{id}$, and measures
$\nu\times\delta_p$. The product is Cantor and its invariant
common conull section $Z\times\{p\}$ carries the required relation.
Evaluation at $p$ identifies the completed crossed products after
the common-null quotient. This product realization loses joint
faithfulness on the uncompleted topological diagonal. The fixed sets
of $W^m$, $m\ne0$, have empty interior by aperiodicity and joint
faithfulness on $Z$, so the product homeomorphism remains
topologically free.
\end{proof}

\begin{cor}[Principal \'etale groupoids and circle twists]
\label{cor:RN-groupoid-twists}
Let a principal second countable Hausdorff \'etale groupoid have compact
metrizable units and a nonempty compact convex designated family of
quasi-invariant probabilities. Suppose relative SBP holds and a
countable bisection atlas has common bounded regular Jacobians and
inverse Jacobians. If its measured orbit relation is amenable for
every designated measure, every locally trivial circle twist has the
same untwisted bounded state Cartan model. In the aperiodic case the
orbit relation has the uniform measure-class dynamical realization
of Theorem~\ref{thm:RN-dynamical-model}.
\end{cor}
\begin{proof}
Restrict a trivializing atlas bisection $s$ to an SBP open set
$V$ with $\overline V\subseteq\operatorname{dom}s$. Then
$\partial(sV)=s(\partial V)$ is null for every designated measure
by quasi-invariance, so both support projections are regular. For
$f\in \mathrm{C}(X)$ the pullback on $V$ is the product of $1_V$ and a
continuous extension of $f\circ s|_{\overline V}$. The chart bound
extends this pullback to $D_K$, and the inverse bound gives the
reverse corner map. The restricted Jacobian is $1_{sV}j_s$.
Using a countable cover by such restrictions, the local sections
and finite first-chart approximations of
Proposition~\ref{prop:SBP-regularizes-twist} give the regular
scalar cocycle. Apply
Theorem~\ref{thm:RN-uniform-CFW} and
Proposition~\ref{prop:RN-untwisting-transport}, and then
Theorem~\ref{thm:RN-dynamical-model} when applicable. The gauge is
regular in the completion and need not be continuous on the original
twist. Nontrivial isotropy is not covered by this conclusion.
\end{proof}

\section{Intrinsic reconstruction and uniform CFW consequences}
\label{sec:intrinsic-reconstruction}

The preceding relation constructions start with a regular atlas. We now
recover such a relation, together with its twist, from the marked
completed Cartan inclusion. The fibrewise Cartan condition supplies
exactly the missing principality hypothesis under countable
partial-normalizer regularity. With SBP, measurewise injectivity then
produces one compatible AF Cartan model.

\subsection{Intrinsic hypotheses and the fibrewise Cartan condition}
\label{subsec:intrinsic-reconstruction-hypotheses}

We first formulate the hypotheses on the completed pair itself.  The
designated traces remain part of the data; no assertion below replaces
them by all norm-continuous traces on the completion.

\begin{defn}[Marked uniformly separable completed Cartan pair]
\label{defn:intrinsic-completed-Cartan}
Extending the completed-pair convention of Section~6, a
\emph{marked uniformly separable completed Cartan pair} $(D,M,\mathbb{E},K)$ consists of
a unital $\mathrm{C}^*$-algebra $M$, a nonempty compact convex jointly faithful
family $K\subseteq \T(M)$, a unital maximal abelian $\mathrm{C}^*$-subalgebra
$D\subseteq M$, and a conditional expectation $\mathbb{E}:M\to D$ such that
$\tau\circ \mathbb{E}=\tau$ for every $\tau\in K$.  We require the unit ball
$M_1$ to be complete and separable for $\|\cdot\|_{2,K}$, and require
\emph{partial-normalizer regularity}:
\[
 \overline{\mathrm{C}^*(D,\mathcal N_{\mathrm{pi}}(D,M))_1}
          ^{\|\cdot\|_{2,K}}=M_1,
\]
where $\mathcal N_{\mathrm{pi}}(D,M)$ denotes the partial isometries
$v\in M$ with $vDv^*,v^*Dv\subseteq D$.
\end{defn}

Joint faithfulness makes $\|\cdot\|_{2,K}$ a norm and makes $\mathbb{E}$
faithful.  The expectation is $2$-norm contracting, so $D_1$ is also
complete and separable.  Partial-normalizer regularity and separability
allow us to choose a sequence $(v_j)$ in
$\mathcal N_{\mathrm{pi}}(D,M)$ with
\begin{equation}\label{eq:intrinsic-countable-normalizers}
 \overline{\mathrm{C}^*(D,v_1,v_2,\ldots)_1}^{\|\cdot\|_{2,K}}=M_1.
\end{equation}
Indeed, approximate a countable dense subset of $M_1$ to successively
smaller errors by bounded elements of the algebra generated by all
partial normalizers, then approximate these elements in norm by
polynomials.  Collect the finitely many normalizers occurring in each
polynomial.  Radial truncation, as in Section~2, keeps the approximants
in the unit ball of the algebra generated by the collected normalizers.

For $\tau\in K$, write
$M_\tau=\pi_\tau(M)''$ and $D_\tau=\pi_\tau(D)''$.  The expectation
extends to the trace-preserving normal expectation
$\mathbb{E}_\tau:M_\tau\to D_\tau$: on $L^2(M_\tau,\tau)$ its extension is the
orthogonal projection onto $L^2(D_\tau,\tau)$.  We shall use the separate
\emph{fibrewise Cartan condition}
\begin{equation*}\label{eq:intrinsic-fibre-Cartan}
 \tag{F}
 D_\tau'\cap M_\tau=D_\tau\qquad(\tau\in K).
\end{equation*}
Regularity of $D_\tau\subseteq M_\tau$ follows from
\eqref{eq:intrinsic-countable-normalizers}; thus \eqref{eq:intrinsic-fibre-Cartan}
makes every fibre a von Neumann Cartan pair.  Joint faithfulness shows
that (F) implies the global maximal abelian condition, although we
retain that condition in Definition~\ref{defn:intrinsic-completed-Cartan}.

\begin{lem}[Uniqueness of the designated-trace-preserving expectation]
\label{lem:intrinsic-expectation-unique}
The expectation $\mathbb{E}$ is determined by $(D,M,K)$.
\end{lem}
\begin{proof}
If $\mathbb{E}_1,\mathbb{E}_2$ are two such expectations, then for $a\in M$, $d\in D$,
and $\tau\in K$, bimodularity and trace preservation give
\[
 \tau(d^*\mathbb{E}_1(a))=\tau(d^*a)=\tau(d^*\mathbb{E}_2(a)).
\]
Take $d=\mathbb{E}_1(a)-\mathbb{E}_2(a)$ and use joint faithfulness.
\end{proof}

\begin{lem}[Diagonal corners of partial normalizers]
\label{lem:intrinsic-normalizer-corners}
For $v\in\mathcal N_{\mathrm{pi}}(D,M)$,
\[
 e_v:=|\mathbb{E}(v)|^2=v^*\mathbb{E}(v)
\]
is a projection in $D$, $ve_v=\mathbb{E}(v)$, and $e_v$ is the largest projection
$e\leq v^*v$ such that $ve\in D$.  Equivalently, it is the largest
source corner on which the partial automorphism implemented by $v$
is the identity.  For two partial normalizers $v,w$, the projection
\[
 e(v,w):=|\mathbb{E}(v^*w)|^2
\]
is the largest projection $e\leq(v^*v)(w^*w)$ on which their corner
isomorphisms agree.  On this corner,
\[
 we(v,w)=v\mathbb{E}(v^*w),
\]
and $\mathbb{E}(v^*w)$ is a diagonal partial unitary with support $e(v,w)$.
\end{lem}
\begin{proof}
This is the algebraic diagonal-corner argument of
\cite[Theorem~5.12, Claim~1]{DFP}; we include it because no normality of
$\mathbb{E}$ is needed.  Put $p=v^*v$, $q=vv^*$, and $b=\mathbb{E}(v)$.  Bimodularity
gives $b=qb=bp$ and
$bd=v(pd)v^*b$ for $d\in D$.  Hence $r=v^*b$ commutes with $D$ and
belongs to $D$ by maximal abelianness.  Since $vr=b$,
\[
 br=\mathbb{E}(vr)=\mathbb{E}(b)=b,
 \qquad r^2=v^*br=v^*b=r.
\]
An idempotent in the abelian $\mathrm{C}^*$-algebra $D$ is a projection.  Also
$r\leq p$ and $b^*b=rv^*vr=r$, proving the stated identities.
If $e\leq p$ and $ve\in D$, then $ve=\mathbb{E}(ve)=be=vre$, whence
$e=re$.  Conversely, $ve\in D$ is exactly the assertion that $ve$
commutes with $D$, which is exactly the identity-corner assertion.

Products of partial normalizers are again partial normalizers.  Apply
the preceding conclusion to $v^*w$.  Its diagonal corner
$e=e(v,w)$ has both source and range $e$, lies below $(v^*v)(w^*w)$,
and satisfies $v^*we=\mathbb{E}(v^*w)$.  Its source identity gives
$w^*vv^*we=e$, so $vv^*we=we$ and $we=v\mathbb{E}(v^*w)$.
The diagonal coefficient is a unitary in $eD$, and therefore the two
corner isomorphisms agree on $e$.  Conversely, if they agree on a
projection $f\leq(v^*v)(w^*w)$, their range projections agree there,
and $v^*wf$ commutes with $D$.  Thus $v^*wf\in D$ and the maximal
diagonal-corner assertion gives $f\leq e(v,w)$.
\end{proof}

\begin{defn}[Intrinsic small boundary property]
\label{defn:intrinsic-Cartan-SBP}
The marked diagonal $(D,K)$ has the \emph{intrinsic small boundary
property} if the finite-spectrum self-adjoint contractions in $D$ are
uniformly $2$-norm dense in $(D_{\mathrm{sa}})_1$.
\end{defn}

By \cite[Lemma~8.3]{TC}, bounded
uniform completeness makes this property equivalent to norm real rank zero
of $D$.  The converse follows from norm density of finite-spectrum
self-adjoint elements and clipping to the contraction interval.
Since $D_1$ is uniformly $2$-separable, it is also equivalent to the existence
of increasing finite-dimensional unital subalgebras $F_n\subseteq D$
whose unit balls have uniformly $2$-dense union in $D_1$.  To see this,
choose countably many finite-spectrum approximants to a dense sequence
in $(D_{\mathrm{sa}})_1$, and let $F_n$ be generated by the first $n$
finite spectral partitions.  Their projections commute.  Approximating
real and imaginary parts and then clipping to the unit disk gives the
required bounded approximation of $D_1$.  The converse is immediate.

\begin{lem}[Small boundaries are independent of the diagonal core]
\label{lem:intrinsic-SBP-core-independent}
Let $B\subseteq D$ be a separable unital $\mathrm{C}^*$-subalgebra with
$\overline{B_1}^{\|\cdot\|_{2,K}}=D_1$, put $X=\operatorname{Spec}(B)$,
and let $K_X$ be the restrictions of $K$ to $B=\mathrm{C}(X)$.  Then
Definition~\ref{defn:intrinsic-Cartan-SBP} holds if and only if
$(X,K_X)$ has relative SBP.  In particular, this property does not
depend on the chosen uniformly dense diagonal core or on (F).
\end{lem}
\begin{proof}
Restriction identifies $K$ affinely and homeomorphically with $K_X$.
Indeed, agreement on $B$ implies agreement on $D$ by bounded uniform
$2$-norm density, and then agreement on $M$ by trace preservation of
$\mathbb{E}$.  Compactness gives the homeomorphism.  The bounded uniform
completion of $B$ identifies with $D$.

Suppose first that the intrinsic property holds.  Given $f\in B_{\mathrm{sa}}$
and $\eps>0$, choose a finite-spectrum $h\in D_{\mathrm{sa}}$ with
$\|f-h\|_{2,K}<\eps/3$.  Perturb $h$ by a real scalar of absolute
value less than $\eps/3$ to obtain an invertible finite-spectrum
$h'$.  Put $\gamma=\min\{|\lambda|:\lambda\in\operatorname{Spec}(h')\}>0$.
Choose $g\in B_{\mathrm{sa}}$ such that
\[
 \|g-h'\|_{2,K}<\min\{\eps/3,\gamma\sqrt{\eps}\}.
\]
In the commutative completion, on $\{g=0\}$ one has
$|g-h'|\geq\gamma$.  Consequently
\[
 \sup_{\mu\in K_X}\mu(\{g=0\})
 \leq \gamma^{-2}\|g-h'\|_{2,K}^2<\eps,
 \qquad \|f-g\|_{2,K}<\eps.
\]
The characterization in \cite[Theorem~2.9]{EN} gives relative
SBP.  Conversely, relative SBP gives the finite projection diagonals
of Proposition~\ref{prop:relative-SBP-finite-diagonals} in
$\mathrm{C}(X)^{K_X}=D$, and these give the intrinsic property.
\end{proof}

\begin{prop}[Intrinsic tests for the fibrewise Cartan condition]
\label{prop:intrinsic-fibre-Cartan-tests}
For a pair as in Definition~\ref{defn:intrinsic-completed-Cartan}, (F) is equivalent to
the following uniform diagonal averaging condition: for every $a\in M$
and $\eps>0$, there are $u_1,\ldots,u_r\in\mathcal U(D)$ and
$t_1,\ldots,t_r\geq0$, with $\sum_jt_j=1$, such that
\begin{equation}\label{eq:intrinsic-uniform-diagonal-averaging}
 \left\|\sum_{j=1}^r t_j u_jau_j^*-\mathbb{E}(a)\right\|_{2,K}<\eps.
\end{equation}
If the intrinsic small boundary property holds, choose any increasing
finite diagonals $F_n\subseteq D$ with bounded uniformly dense union,
write their atoms as $p_{n,1},\ldots,p_{n,r_n}$, and define
\[
 \mathbb{E}_n^{\mathrm{pin}}(a)=\sum_i p_{n,i}ap_{n,i}.
\]
Then (F) is also equivalent to
\begin{equation}\label{eq:intrinsic-uniform-diagonal-pinching}
 \|\mathbb{E}_n^{\mathrm{pin}}(a)-\mathbb{E}(a)\|_{2,K}\longrightarrow0\qquad(a\in M).
\end{equation}
For the words $w$ in the generators of
\eqref{eq:intrinsic-countable-normalizers} and their adjoints, it is
enough to check
\begin{equation}\label{eq:intrinsic-normalizer-coincidence-test}
 \sup_{\tau\in K}\tau\bigl(f_n(w)-e_w\bigr)\longrightarrow0,
 \qquad f_n(w)=\sum_i p_{n,i}w^*p_{n,i}w.
\end{equation}
The assertions concern each fixed element or word, rather than
convergence uniformly over $M_1$.
\end{prop}
\begin{proof}
Let $\mathbb{E}_\tau^{\mathrm{comm}}:M_\tau\to D_\tau'\cap M_\tau$ be the trace-preserving
normal expectation.  The $L^2$-closed convex hull of the
$\pi_\tau(\mathcal U(D))$-conjugates of $\pi_\tau(a)$ contains
$\mathbb{E}_\tau^{\mathrm{comm}}(\pi_\tau(a))$.  Indeed, its unique minimum-norm element is
fixed by these unitaries.  The convex hull lies in the operator ball
of radius $\|a\|$, which is $L^2$-closed, so that element lies in
$D_\tau'\cap M_\tau$.  Moreover $\mathbb{E}_\tau^{\mathrm{comm}}$ is constant on the orbit,
identifying the minimum-norm element with $\mathbb{E}_\tau^{\mathrm{comm}}(\pi_\tau(a))$.

Under (F), this element is $\mathbb{E}_\tau(\pi_\tau(a))$.  For each $\tau$
choose a finite diagonal-unitary averaging map $T_\tau$ with error
less than $\eps$ at $\tau$.  The squared error is continuous on $K$,
so this inequality holds in a neighbourhood of $\tau$.  Choose a
finite subcover and its averaging maps $T_1,\ldots,T_m$.  They commute
because $D$ is abelian, contract every tracial $2$-norm, and fix $\mathbb{E}(a)$.
Their composition $T$ is another finite diagonal-unitary average, and
on the neighbourhood corresponding to $T_i$ it satisfies
\[
 \|T(a)-\mathbb{E}(a)\|_{2,\tau}
       \leq\|T_i(a)-\mathbb{E}(a)\|_{2,\tau}<\eps.
\]
This proves \eqref{eq:intrinsic-uniform-diagonal-averaging}.
Conversely, $\mathbb{E}_\tau^{\mathrm{comm}}(uau^*)=\mathbb{E}_\tau^{\mathrm{comm}}(a)$ for $u\in\mathcal U(D)$.
Applying $\mathbb{E}_\tau^{\mathrm{comm}}$ to the approximants in that equation gives
$\mathbb{E}_\tau^{\mathrm{comm}}(\pi_\tau(a))=\mathbb{E}_\tau(\pi_\tau(a))$ for every $a\in M$.
Their $2$-norm density in $M_\tau$ gives (F).

For finite diagonals, uniform density gives
$(\bigcup_n\pi_\tau(F_n))''=D_\tau$.  The maps $\mathbb{E}_n^{\mathrm{pin}}$ are the
$L^2$-orthogonal projections onto
$L^2(\pi_\tau(F_n)'\cap M_\tau,\tau)$;
these subspaces decrease to $L^2(D_\tau'\cap M_\tau,\tau)$. Thus their
$L^2$-limit is $\mathbb{E}_\tau^{\mathrm{comm}}$.  Under (F), the continuous functions
\[
 \tau\longmapsto\|\mathbb{E}_n^{\mathrm{pin}}(a)-\mathbb{E}(a)\|_{2,\tau}^{\,2}
\]
decrease pointwise to zero on compact $K$.  Dini's theorem proves
\eqref{eq:intrinsic-uniform-diagonal-pinching}.  Conversely, its
fibrewise limit identifies $\mathbb{E}_\tau^{\mathrm{comm}}$ with $\mathbb{E}_\tau$ on $\pi_\tau(M)$,
and hence on $M_\tau$.

For a partial normalizer $w$, the elements $f_n(w)$ are decreasing
projections in $D$, $e_w\leq f_n(w)$, and
\[
 \mathbb{E}_n^{\mathrm{pin}}(w)=wf_n(w),\qquad
 \|\mathbb{E}_n^{\mathrm{pin}}(w)-\mathbb{E}(w)\|_{2,K}^{\,2}
        =\sup_{\tau\in K}\tau(f_n(w)-e_w).
\]
These follow by writing $p_{n,i}wp_{n,i}=w(w^*p_{n,i}w)p_{n,i}$
and using Lemma~\ref{lem:intrinsic-normalizer-corners}.  Thus
\eqref{eq:intrinsic-uniform-diagonal-pinching} implies
\eqref{eq:intrinsic-normalizer-coincidence-test}.  Conversely, the
latter gives pinching convergence on each word and on its diagonal
multiples, since $\mathbb{E}_n^{\mathrm{pin}}(dw)=d\mathbb{E}_n^{\mathrm{pin}}(w)$ and $\mathbb{E}(dw)=d\mathbb{E}(w)$.
Finite sums of these multiples are dense in $\mathrm{C}^*(D,v_1,v_2,\ldots)$;
bounded uniform density and $2$-norm contractivity of $\mathbb{E}_n^{\mathrm{pin}}$ and $\mathbb{E}$
give the convergence on every $a\in M$.
\end{proof}

\begin{rem}\label{rem:intrinsic-pinching-pointwise}
Convergence in Proposition~\ref{prop:intrinsic-fibre-Cartan-tests} is
pointwise in $a$.  Uniform convergence over $M_1$ fails even for the
standard diagonal Cartan subalgebra in the hyperfinite factor: for the
finite tensor diagonals $F_n$, an off-diagonal unitary in a later tensor
factor commutes with $F_n$, has expectation zero, and has $2$-norm one.
The intrinsic averaging condition, and hence (F), does not require
small boundaries.
\end{rem}

\subsection{Separable cores and inverse reconstruction}
\label{subsec:intrinsic-inverse-reconstruction}

The countable normalizer presentation permits an inverse construction
without a small boundary assumption.  The additional condition is
maximal abelianness in every designated tracial fibre, rather than
maximal abelianness only in the completed algebra.

\begin{lem}[A compatible separable Cartan core]
\label{lem:intrinsic-separable-Cartan-core}
Let $(D\subseteq M,\mathbb{E},K)$ be a completed Cartan pair as in
Definition~\ref{defn:intrinsic-completed-Cartan}.  Suppose that the
unit ball of $D$ is separable for $\norm{\cdot}_{2,K}$ and that there
are countably many partial-isometry normalizers $(v_j)$ such that
\[
 \overline{\mathrm{C}^*(D,v_j:j\geq1)_1}^{\,\norm{\cdot}_{2,K}}=M_1.
\]
There are separable unital $\mathrm{C}^*$-subalgebras $D_0\subseteq D$ and
$A_0=\mathrm{C}^*(D_0,v_j:j\geq1)\subseteq M$ such that $(D_0,A_0)$ is a
$\mathrm{C}^*$-Cartan pair, $\mathbb{E}(A_0)\subseteq D_0$, and
\[
 D_0^K=D,\qquad A_0^K=M.
\]
Restriction identifies $K$ affinely and homeomorphically with the
designated traces on $A_0$, and also with their restrictions to
$D_0$.  The algebra $D_0$ can be chosen to contain any prescribed
separable $\mathrm{C}^*$-subalgebra of $D$.
\end{lem}

\begin{proof}
Choose a countable set $S\subseteq D_1$ that is uniformly
$2$-norm dense in $D_1$.  Let $W$ be the countable semigroup of words
in the $v_j$ and their adjoints, including $1$.  Products of partial
normalizers are partial normalizers, since their initial and final
projections lie in the abelian algebra $D$.  Start with the unital
diagonal algebra generated by $S$, the prescribed separable
subalgebra, and
\[
 w^*w,\quad ww^*,\quad \mathbb{E}(w)\qquad(w\in W).
\]
At each step adjoin the images of a countable norm-dense subset under
$d\mapsto v_jdv_j^*$ and $d\mapsto v_j^*dv_j$ for every $j$.
Take the norm closure of the increasing union.  The resulting $D_0$
is separable and invariant under these corner maps, by their norm
continuity.  Its bounded uniform completion is $D$.

Moving diagonal coefficients through words shows that a dense
*-subalgebra of $A_0$ consists of finite sums $wd$, with $w\in W$
and $d\in D_0$.  Diagonal bimodularity gives
$\mathbb{E}(wd)=\mathbb{E}(w)d\in D_0$, so $\mathbb{E}(A_0)\subseteq D_0$.  Its restriction
is faithful.  If $a\in A_0$ commutes with $D_0$, bounded uniform
$2$-norm approximation and
\[
 \norm{[a,d-d_0]}_{2,K}\leq2\norm a\norm{d-d_0}_{2,K}
\]
show that it commutes with $D$.  Maximal abelianness in $M$ gives
$a\in D$, and hence $a=\mathbb{E}(a)\in D_0$.  Thus $D_0$ is maximal
abelian in $A_0$.  The normalizers $v_j$ and the diagonal generate
$A_0$, so $(D_0,A_0)$ is a $\mathrm{C}^*$-Cartan pair.

Replacing the finitely many diagonal coefficients in a polynomial
by bounded approximants from $D_0$ proves uniform $2$-norm density
of $A_0$ in $\mathrm{C}^*(D,v_j:j\geq1)$, using
\eqref{eq:two-mult}.  Radial truncation, as in Section~2, makes the
approximants contractions.  The assumed bounded density therefore
gives $A_0^K=M$.

Restriction of traces is continuous and affine.  If two members of
$K$ agree on $D_0$, uniform $2$-norm density makes them agree on
$D$.  Since both satisfy $\tau=\tau\circ \mathbb{E}$, they agree on $M$.
Compactness of $K$ now gives the asserted affine homeomorphisms.
No assertion that the designated traces exhaust $\T(A_0)$ is needed.
\end{proof}

\begin{thm}[Inverse uniform Feldman--Moore reconstruction]
\label{thm:intrinsic-uniform-FM}
Let $(D\subseteq M,\mathbb{E},K)$ satisfy the hypotheses of
Lemma~\ref{lem:intrinsic-separable-Cartan-core}.  The following
conditions are equivalent:
\begin{enumerate}
\item The fibre Cartan condition
\eqref{eq:intrinsic-fibre-Cartan} holds.
\item There are a compact metrizable space $X$, a nonempty compact
convex family $K_X$ of invariant probability measures, a countable
Borel equivalence relation $R$ with a regular countable atlas, and
a regular scalar twist $\sigma$, such that
\[
 (D\subseteq M,\mathbb{E},K)
 \cong
 (D_{K_X}\subseteq M_{K_X}^{\sigma}(R),\mathbb{E}_R,K_X).
\]
The isomorphism preserves the expectation and implements an affine
homeomorphism from $K$ to the designated trace family indexed by
$K_X$.
\end{enumerate}
No amenability, aperiodicity, or small boundary assumption is needed.
If the completed diagonal has the intrinsic small boundary property,
the model $(X,K_X)$ has relative SBP.
\end{thm}

\begin{proof}
Choose the core $(D_0,A_0)$ from the preceding lemma and put
$X=\operatorname{Spec}(D_0)$.  Write
$K_X=\{\mu_\tau:\tau\in K\}$, where
$\tau(d)=\int_Xd\,d\mu_\tau$ for $d\in D_0=\mathrm{C}(X)$.
The preceding lemma makes this a compact convex family with the
specified trace marking.  Renault's reconstruction
\cite[Theorem~5.9]{RenaultC} gives a second-countable locally compact
Hausdorff topologically principal \'etale groupoid $G$ and its Weyl
twist $\Sigma$, with
\[
 (\mathrm{C}(X)\subseteq A_0,\mathbb{E}|_{A_0})
 \cong(\mathrm{C}(X)\subseteq \mathrm{C}_r^*(G,\Sigma),\mathbb{E}_G).
\]
Every designated trace is $\mu_\tau\circ \mathbb{E}_G$.  In particular,
each $\mu_\tau$ is invariant: applying the trace identity to a
normalizer and a diagonal coefficient gives invariance under its
bisection map.

Enumerate the word semigroup as $W=(u_i)_{i\geq1}$, with $u_1=1$.
Let $B_i\subseteq G$ be the Weyl bisection of $u_i$.  Its source
and range are clopen, with indicators $u_i^*u_i$ and $u_iu_i^*$.
The coefficient of $u_i$ has modulus one on $B_i$, so it supplies
a continuous section of the twist there.  These bisections cover
$G$.  Indeed, given $\gamma\in G$, choose a compactly supported
local section $a$ with $a(\gamma)\ne0$.  Approximate $a$ in norm
by a finite diagonal-word polynomial $P$.  Coefficient evaluation
in the regular representation gives
\[
 |a(\gamma)-P(\gamma)|\leq\norm{a-P}.
\]
For a sufficiently close approximation $P(\gamma)\ne0$, forcing
$\gamma$ to lie in one of the word bisections occurring in $P$.

For two words put
\[
 e_{ij}=|\mathbb{E}(u_i^*u_j)|^2.
\]
Lemma~\ref{lem:intrinsic-normalizer-corners} makes this a projection
in $D_0$.  In the Weyl model its support is
$s(B_i\cap B_j)$: the expectation coefficient of $u_i^*u_j$
is a unit-modulus phase exactly where the two word arrows agree,
and is zero elsewhere.  Thus intersections of word bisections have
clopen source projections.

Assume condition (1).  We verify principality on one set conull
for all designated measures.  For a word $u_i$, let $F_i$ be the
Borel fixed-point set of its partial homeomorphism and let
$U_i=s(B_i\cap G^{(0)})$, whose indicator is $e_{1i}$.  For a
fixed $\tau\in K$, the projection
$p=1_{F_i\setminus U_i}$ belongs to
$\pi_\tau(D_0)''=\pi_\tau(D)''$.  The normalizer
$\pi_\tau(u_i)p$ commutes with this diagonal because its underlying
partial map fixes every point of that support.  The
trace-preserving normal expectation satisfies
\[
 \mathbb{E}_\tau(\pi_\tau(u_i)p)=\pi_\tau(\mathbb{E}(u_i))p=0.
\]
Condition (1) consequently makes $\pi_\tau(u_i)p$ zero.  Since
$p\leq\pi_\tau(u_i^*u_i)$, this gives $p=0$, or
$\mu_\tau(F_i\setminus U_i)=0$.  This is also the measured
principality criterion in \cite[Proposition~3.11]{Renault}.

Because the word bisections cover $G$, the set of units with
nontrivial isotropy is
\[
 N=\bigcup_{i\geq1}(F_i\setminus U_i).
\]
It is Borel, invariant, and null for every member of $K_X$.
Put $X_0=X\setminus N$.  The restriction $G|_{X_0}$ is principal;
its arrows are therefore identified with the countable Borel orbit
relation on $X_0$.  Extend this relation by the equality relation
on $N$, obtaining a relation $R$ on the compact space $X$.
Restrict the word maps to $X_0$, and adjoin the identity on $X$
and the empty map.  Their graphs cover $R$ and they are closed
under composition and inversion.  Their domain and range
projections agree, modulo $K_X$-null sets, with the clopen word
projections.  Word corner maps preserve $\mathrm{C}(X)$ and extend
isometrically to $D_{K_X}$.  On $X_0$ endpoint coincidence is
arrow coincidence, so the coincidence projections are the $e_{ij}$.
This gives a regular countable atlas.  Invariance under its charts,
and countable Borel disjointization, give invariance under every
Borel partial bijection with graph in $R$.

To obtain a scalar presentation of the twist, select on each arrow
the section of its first word bisection.  Its disjoint source
pieces are the clopen projections
\[
 r_j=u_j^*u_j\prod_{i<j}(1-e_{ij}).
\]
The section is normalized because the first word is $1$.
On a fixed word chart its first-word pieces have index at most
that word's index.  For a product of two fixed word charts, each
resulting word has a finite index, and its first matching word has
no larger index.  Hence the scalar multiplication coefficient on
that chart pair is a finite patch of the diagonal phases
$\mathbb{E}(u_k^*u_i u_j)$ and the corresponding word-section phases.
It belongs to $D_0$ on the composition domain.  The scalar
cocycle is Borel and normalized and satisfies the cocycle identity
by associativity in $\Sigma$.  Extend it trivially on the equality
relation on $N$.

This cocycle is regular also along arbitrary admissible chart
pairs.  Disjointize each partial map into its first-word pieces
as in Proposition~\ref{prop:regular-relation-pseudogroup-closure}.
On finitely many pieces the cocycle coefficient is a regular
diagonal function.  The omitted composition-domain mass is
bounded by the two omitted source masses, which tend uniformly
to zero.  Since all phases have modulus one, the resulting finite
coefficients converge uniformly in $2$-norm to the required
coefficient.  Thus Definition~\ref{defn:regular-two-cocycle} holds.

Finally, the designated regular groupoid representations become
the twisted relation representations after restriction to $X_0$.
The trace obtained from a dense barycentre of $K$ is faithful on
$A_0$, by joint faithfulness and continuity of the trace profiles.
Its regular representation therefore identifies the core faithfully
with its relation image.
The word operators correspond to regular word graph operators
with the diagonal phases just constructed.  Conversely those
graph operators are obtained from the word operators by such
phases.  The original core has diagonal $\mathrm{C}(X)$, whose bounded
completion is $D_{K_X}=D$.  These two descriptions consequently
generate the same bounded uniform completion.  The expectation
is restriction to the units on both sides, and its composition
with $\mu_\tau$ is the original trace.  They therefore give the
claimed marked expectation-preserving isomorphism.

Conversely, suppose (2) holds.  For each designated measure the
tracial GNS von Neumann algebra of the twisted relation model is
the measured twisted relation algebra: finite graph kernels are
dense in its regular $L^2$ space.  The strong closure of the
diagonal is $L^\infty(X,\mu)$.  It is maximal abelian, since a
kernel commuting with a countable family of continuous functions
separating points is supported on the diagonal.  Equivalently,
this is the principal-relation direction of the Feldman--Moore
construction \cite{FMII}, or \cite[Proposition~3.11]{Renault}.
Transporting this assertion through the marked pair isomorphism
gives condition (1).  If intrinsic SBP holds, relative SBP of
$(X,K_X)$ follows from
Lemma~\ref{lem:intrinsic-SBP-core-independent}.
\end{proof}

\begin{rem}[Scope of the inverse theorem]
\label{rem:intrinsic-inverse-scope}
Partial-isometry normalizer generation supplies the clopen domain
and germ-coincidence projections needed in this construction; it
does not assert norm density of projections in $D_0$.  The fibre
condition removes isotropy simultaneously for the designated
measures.  Amenability is a separate property: when it holds and
intrinsic SBP is present, the uniform CFW and untwisting theorems
apply to this reconstructed model.
\end{rem}

\subsection{Independence and uniqueness of the twisted model}
\label{subsec:intrinsic-model-uniqueness}

The invariant reconstructed from a completed pair is a marked regular
relation together with a regular scalar cohomology class.  The compact
unit space and the topology of the auxiliary Cartan core need not be
unique.  We first extend the normalizer correspondence to the twisted
models of Section~\ref{sec:regular-twists}.

\begin{prop}[Twisted uniform normalizer correspondence]
\label{prop:twisted-uniform-normalizers}
Let $(X,R,K,D_K,\mathcal S)$ be a regular uniform measure relation,
and let $\sigma$ be a regular scalar twist.  For an admissible partial
map $\varphi:A\to B$, the phase-one graph kernel defines a partial
isometry $u_\varphi^\sigma\in M_K^\sigma(R)$, with initial and final
projections $1_A$ and $1_B$, implementing the corner pullback of
$\varphi$.  Every partial-isometry normalizer $v$ of $D_K$ in
$M_K^\sigma(R)$ has a unique factorization
\[
 v=u_\varphi^\sigma a,
 \qquad a\in1_AD_K,\qquad a^*a=aa^*=1_A,
\]
where $\varphi$ is admissible and is unique modulo $K$-null sets.
Thus partial normalizers modulo diagonal source phases recover the
full admissible pseudogroup.  The products of the graph lifts are
\[
 u_s^\sigma u_t^\sigma
   =u_{st}^\sigma\sigma_{s,t},\qquad
 \sigma_{s,t}(x)=\sigma(stx,tx,x),
\]
on the domain of the composite.
\end{prop}

\begin{proof}
First-chart disjointization expresses an admissible graph as a
countable union of regular chart restrictions with disjoint sources
and disjoint ranges.  The associated finite sums of twisted graph
operators have norm at most one.  Their uniformly vanishing source
tails show that their bounded uniform limit is the phase-one graph
kernel $u_\varphi^\sigma$.  The corner covariance and the displayed
product formula follow directly from twisted convolution.

For the converse, the disjoint graph coordinates $w_j$ of
Section~\ref{sec:uniform-relation-algebra}, now interpreted as twisted
graph operators, still satisfy
\[
 \mathbb{E}(w_i^*w_j)=\delta_{ij}r_j.
\]
Twisting changes the product coefficients but not the graph supports
or this orthogonality identity.  Hence the proof of
Lemma~\ref{lem:uniform-relation-Fourier-approximation} gives
\[
 b_j(x)=\mathbb{E}(w_j^*x),\qquad
 \norm{x-\sum_{j\leq n}w_jb_j(x)}_{2,K}\longrightarrow0
\]
and its Fourier-tail identity, with exactly the same proof.
The diagonal-kernel argument also proves that $D_K$ is maximal
abelian in the twisted completion.

Let $p=v^*v$ and $q=vv^*$.  Choose the faithful dense-barycenter
representation as in Section~\ref{sec:uniform-relation-algebra}.
The corner automorphism induced by $v$ extends normally to the
corresponding corners of $L^\infty(X,\mu_0)$: it is trace preserving,
and $D_K$ is strongly dense in this diagonal.  A countable
point-separating family from $\mathrm{C}(X)$ spatializes that normal corner
automorphism as a Borel partial bijection.  Covariance forces the
kernel $v\delta$ to be supported on its graph, and its squared column
norm is $p$.  Consequently
\[
 e_j=|b_j(v)|^2\in D_K
\]
are orthogonal projections describing the first matching chart of
this graph.  These identities initially hold in the faithful
$\mu_0$-representation and therefore hold in $D_K$.
The Fourier-tail identity gives the essential uniform conclusion
\[
 \sup_{\mu\in K}\mu\left(p-\sum_{j\leq n}e_j\right)
       \longrightarrow0.
\]
The chart images of these source projections are orthogonal regular
range projections.  Since every designated trace satisfies
$\tau_\mu(p)=\tau_\mu(q)$, their range tails vanish uniformly as well.
Thus these chart restrictions define a Borel partial bijection on
one invariant $K$-conull set.  Coincidence with any original chart
is the bounded uniform limit of the corresponding finite regular
coincidence sums, so the partial bijection is admissible.

The proof is the same support-and-tail argument as in
Proposition~\ref{prop:uniform-full-group-normalizers}; the scalar
twist only changes the diagonal phases of the coefficients.
The graph lift and $v$ implement the same corner automorphism, so
$(u_\varphi^\sigma)^*v$ is a unitary in $pM_K^\sigma(R)p$ commuting
with $pD_K$.  Maximal abelianness puts it in $pD_K$.
Faithfulness and the point-separating family give uniqueness of
both the partial map and its source phase.  Equivalently, the
coincidence projections are the intrinsic projections of
Lemma~\ref{lem:intrinsic-normalizer-corners}.
\end{proof}

\begin{thm}[Twisted model isomorphisms]
\label{thm:twisted-uniform-pair-isomorphism}
Let $(X_i,R_i,K_i,D_i,\mathcal S_i)$, $i=1,2$, be regular uniform
measure relations with regular scalar twists $\sigma_i$, and put
$M_i=M_{K_i}^{\sigma_i}(R_i)$.  Fix an affine homeomorphism
$H:K_1\to K_2$.  The following conditions are equivalent:
\begin{enumerate}
\item There is a unital isomorphism $\Phi:M_1\to M_2$ carrying
      $D_1$ onto $D_2$ and satisfying
      \[
       \tau_{H(\mu)}(\Phi(x))=\tau_\mu(x)
       \qquad(x\in M_1,\ \mu\in K_1).
      \]
\item There are a uniform measure relation isomorphism $\theta$
      implementing $H$ and a regular gauge $b:R_1\to\mathbb T$ such that
      \begin{equation}
      \label{eq:twisted-model-cohomology-transport}
       \sigma_1(y,z,x)
        =\sigma_2(\theta y,\theta z,\theta x)
          \frac{b(y,z)b(z,x)}{b(y,x)}.
      \end{equation}
      The spatial map and the identity are defined on one invariant
      set conull for every measure in $K_1$ and its image.
\end{enumerate}
For a given $\Phi$, the spatial map $\theta$ and gauge $b$ are unique
modulo universally null sets, relative to the chosen scalar
presentations.  Every such isomorphism preserves the expectations.
\end{thm}

\begin{proof}
Suppose first that $\Phi$ satisfies (i).  Its diagonal restriction
and Theorem~\ref{thm:orbit-spatial-diagonal} give a Borel bijection
$\theta$ between universally conull sets, implementing $H$ and the
diagonal pullback.  For every source chart $s$, the preceding
proposition supplies the admissible target partial map of
$\Phi(u_s^{\sigma_1})$.  Corner covariance, tested on a countable
continuous separating family, identifies this map with
$\theta s\theta^{-1}$ outside a universally null set.
Apply the same argument to $\Phi^{-1}$ and all target charts.
There are only countably many chart and separating-function
identities.  Removing their exceptional sets and taking the
invariant cores on both sides, exactly as in the proof of
Theorem~\ref{thm:uniform-relation-pair-isomorphism}, makes $\theta$
a uniform measure relation isomorphism on common invariant conull
sets.  In particular both directions transport admissible maps.

Write the unique source phases as
\[
 \Phi(u_s^{\sigma_1})
    =u_{\theta s\theta^{-1}}^{\sigma_2}\Phi(a_s),
 \qquad a_s\in p_sD_1,\qquad |a_s|=p_s.
\]
On a regular coincidence projection, two source graph kernels
agree.  Their transported graph kernels also agree there, so the
phases $a_s$ agree on graph overlaps.  Choose simultaneous Borel
representatives for these countably many phases and patch them by
first-chart disjointization to a normalized Borel function
\[
 b(sx,x)=a_s(x)
\]
on $R_1$.  The overlap identities make this independent of the
selected chart.  Along an arbitrary admissible partial map,
finite first-chart phase sums belong to $D_1$ and approximate its
coefficient uniformly in $2$-norm, by the uniformly vanishing
source tails.  Thus $b$ is a regular gauge.

Multiplicativity of $\Phi$ and the product formula in the preceding
proposition give, on the domain of $st$,
\[
 a_{st}(x)\sigma_1(stx,tx,x)
  =\sigma_2(\theta stx,\theta tx,\theta x)
     a_s(tx)a_t(x).
\]
This is~\eqref{eq:twisted-model-cohomology-transport}.  Since chart
pairs are countable and cover all composable arrows, it holds on
one common invariant conull set.  Expectation preservation follows
from Lemma~\ref{lem:intrinsic-expectation-unique}: both transported
expectations have the same designated trace pairing with the
diagonal.

Conversely assume (ii).  On regular graph kernels define
\[
 (Jf)(\theta y,\theta x)=b(y,x)f(y,x).
\]
Equation~\eqref{eq:twisted-model-cohomology-transport} shows that
$J(f*_{\sigma_1}g)=J(f)*_{\sigma_2}J(g)$ and that $J$ preserves
the involution.  Counting-space transport by $\theta$, followed by
multiplication by $b$, implements $J$ unitarily in compatible
dense-barycenter regular representations.  Consequently it
preserves reduced norms.  Admissibility of transported charts and
regularity of $b$ show that the images of regular graph operators
and diagonal coefficients belong to the target completed algebra;
the same holds in the reverse direction using $\overline b$.
Since $b(x,x)=1$, $J$
preserves the diagonal, the expectations, and all designated traces.
It therefore extends isometrically for the uniform $2$-norm to
the required isomorphism of bounded completions.

Finally, the diagonal action of a fixed $\Phi$ determines $\theta$
modulo universal null sets by a countable separating family.  Its
action on the countable graph lifts then determines every $a_s$,
and hence $b$, on the same conull relation.  This proves uniqueness.
\end{proof}

\begin{cor}[Independence of the Cartan core]
\label{cor:intrinsic-core-independence}
The marked regular twisted relation model of a completed pair
satisfying the equivalent conditions of
Theorem~\ref{thm:intrinsic-uniform-FM} is independent of
the chosen separable Cartan core, up to uniform relation isomorphism
and regular coboundary.  More generally, any marked isomorphism
between such completed pairs induces a unique spatial transport
and regular gauge, modulo universally null sets and relative to
the selected scalar presentations.  These transports are functorial.
\end{cor}

\begin{proof}
The models reconstructed from two cores have canonical marked
identifications with the original completed pair.  Compose one
identification with the inverse of the other and apply the preceding
theorem.  The same argument applies to any fixed marked pair
isomorphism.  If successive isomorphisms induce $(\theta_{12},b_{12})$
and $(\theta_{23},b_{23})$, their composite induces
\[
 \theta_{13}=\theta_{23}\theta_{12},\qquad
 b_{13}(y,x)=b_{12}(y,x)
              b_{23}(\theta_{12}y,\theta_{12}x).
\]
The identities follow on countably many graph lifts, hence on the
common conull relation.  Intersecting the transported invariant
conull sets makes the composition well defined simultaneously for
all designated measures.
\end{proof}

\begin{prop}[The intrinsic normalizer extension and uniform metric]
\label{prop:intrinsic-normalizer-extension}
For a pair satisfying the equivalent conditions of
Theorem~\ref{thm:intrinsic-uniform-FM}, let
$\mathcal S(D,M)$ be the inverse semigroup of corner isomorphisms
$\Ad(v)$ induced by all $v\in\mathcal N_{\mathrm{pi}}(D,M)$.
The map from partial normalizers to these corner maps has kernel
$\operatorname{PI}(D)$, the diagonal partial isometries. Thus its
inverse-semigroup extension, including the diagonal phases, is determined
by the marked inclusion. In every reconstruction it is the extension
of the admissible pseudogroup by its diagonal source phases.

For unitary normalizers, the quotient
$\mathcal G(D,M)=\mathcal N_{\mathcal U(M)}(D)/\mathcal U(D)$
has the intrinsic displacement metric
\begin{equation}\label{eq:intrinsic-normalizer-displacement}
 d_K([u],[v])
 =\sup_{\tau\in K}\tau\bigl(1-|\mathbb{E}(u^*v)|^2\bigr)
 =\frac12\inf_{a\in\mathcal U(D)}\|u-va\|_{2,K}^2.
\end{equation}
Reconstruction identifies it isometrically with the admissible uniform
full group equipped with its disagreement metric. The extension and
metric are preserved by marked pair isomorphisms.
\end{prop}

\begin{proof}
A partial normalizer inducing the identity on its source corner
commutes with $D$ and hence belongs to $D$. Conversely every diagonal
partial isometry induces the identity on its support. The coincidence
corners and their transition phases are those of
Lemma~\ref{lem:intrinsic-normalizer-corners}; the twisted normalizer
correspondence identifies the quotient corner maps with admissible
partial bijections.

For unitary normalizers put $b=\mathbb{E}(u^*v)$ and $e=|b|^2$.
The same lemma makes $b$ a partial unitary supported on $e$, and
identifies $e$ with the coincidence projection of the induced maps.
For any $a\in\mathcal U(D)$,
\[
 \|u-va\|_{2,\tau}^2
 =2-2\operatorname{Re}\tau(ba)
 \geq2\tau(1-e).
\]
The single diagonal unitary $a=b^*+(1-e)$ gives equality for every
$\tau\in K$, proving the infimum formula after taking suprema.
Expectation uniqueness and preservation of the trace marking give
invariance under marked isomorphisms.
\end{proof}

\subsection{Scope of reconstruction}
\label{subsec:intrinsic-reconstruction-scope}

The fibrewise maximal abelianness condition cannot be replaced by
maximal abelianness in the uniform completion, even in the presence
of SBP, nuclearity, and the full trace family of the original
$\mathrm{C}^*$-algebra.

\begin{example}[A uniform MASA with a non-Cartan tracial fibre]
\label{ex:uniform-MASA-not-fibre-MASA}
Let $X_n=(\mathbb Z_2)^2$ for $n\geq1$, where $\mathbb Z_2$ is the
compact group of $2$-adic integers, and let
\[
 X=\Bigl(\coprod_{n\geq1}X_n\Bigr)\cup\{\infty\}
\]
be the one-point compactification. The group $\mathbb Z^2$ acts on
each $X_n$ by coordinatewise translation and fixes $\infty$.
Each component action is free, minimal, and uniquely ergodic, with
Haar probability $\mu_n$. The space $X$ is compact, metrizable, and
zero-dimensional. Its action is topologically free, since the fixed
set of every nonzero group element is $\{\infty\}$ and has empty
interior.

Fix an irrational $\theta$, choose the constant group cocycle
$\sigma_\theta$ whose canonical unitaries satisfy
$UV=e^{2\pi i\theta}VU$, and put
\[
 A=\mathrm{C}(X)\rtimes_{\sigma_\theta}\mathbb Z^2,
 \qquad D_0=\mathrm{C}(X),\qquad K=\T(A).
\]
The pair $D_0\subseteq A$ is a separable unital $\mathrm{C}^*$-Cartan pair, and
$A$ is nuclear. Write $\mathbb{E}_0:A\to D_0$ for its canonical expectation.
The corner over $X_n$ has the unique trace
$\tau_n=\mu_n\circ \mathbb{E}_0$, while the quotient over $\infty$ is the
irrational rotation algebra $A_\theta$, with its unique trace
$\tau_\infty=\delta_\infty\circ \mathbb{E}_0$.
For completeness, freeness on $X_n$ forces every trace there to
vanish on nonidentity Fourier coefficients, and unique ergodicity
determines its diagonal restriction. In $A_\theta$, conjugation by
$U$ and $V$ forces the trace of $U^pV^q$ to vanish unless
$(p,q)=(0,0)$, which gives uniqueness of its trace.

The component projections $p_n=1_{X_n}$ are central in $A$.
Decomposing a trace along them shows that every trace of $A$ has
the form
\[
 \tau=\sum_{n\geq1}t_n\tau_n+t_\infty\tau_\infty,
 \qquad t_n,t_\infty\geq0,\quad
 \sum_{n\geq1}t_n+t_\infty=1.
\]
Indeed, the residual positive functional after removing the
component traces vanishes on the ideal generated by the $p_n$, so
it factors through $A_\theta$. On continuous Fourier polynomials,
$\tau_n\to\tau_\infty$, because $\mu_n\to\delta_\infty$.
Norm density gives the same convergence on $A$. The family $K$ is
jointly faithful: a convex combination giving positive weight to
every component and to $\infty$ has a full-support diagonal measure,
and $\mathbb{E}_0$ is faithful.

Let $M=A^K$ and $D=D_0^K$, retaining the designated extensions of
the traces in $K$, and let $\mathbb{E}:M\to D$ be the extended expectation.
The bounded unit balls are uniformly $2$-separable. Increasing
finite clopen partitions of $X$ give finite-dimensional diagonals
whose bounded union is uniformly $2$-dense in $D$. Thus the
intrinsic SBP condition holds. The two unitaries $U,V$ normalize
$D$, and $\mathrm{C}^*(D,U,V)_1$ is bounded uniformly $2$-dense in $M_1$.

Nevertheless, $D$ is a MASA in $M$. To see this, let
$a\in D'\cap M$ and put $b=a-\mathbb{E}(a)$. For every $n$, the diagonal in
the $\tau_n$-GNS closure is the Cartan algebra
$L^\infty(X_n,\mu_n)$, since the component action is free. Hence
$\|b\|_{2,\tau_n}=0$. The function
\[
 K\longrightarrow[0,\infty),\qquad
 \tau\longmapsto\tau(b^*b)
\]
is continuous: bounded uniform $2$-approximation of $b$ by elements
of $A$ approximates these functions uniformly, using
\[
 \sup_{\tau\in K}|\tau(b^*b-c^*c)|
 \leq(\|b\|+\|c\|)\|b-c\|_{2,K}.
\]
Therefore it vanishes also at $\tau_\infty$ and, by the displayed
trace decomposition, on all of $K$. Joint faithfulness gives
$a=\mathbb{E}(a)\in D$.

At $\tau_\infty$, however, the completed diagonal has scalar GNS
closure, whereas the ambient GNS closure is the hyperfinite
$\mathrm{II}_1$ factor generated by $A_\theta$. Thus
\[
 D_{\tau_\infty}=\mathbb C
 \quad\text{is not a MASA in}\quad M_{\tau_\infty}.
\]
All tracial fibres are injective, by nuclearity of $A$~\cite{Connes76}, but the
fibrewise Cartan condition fails. Here $K$ is the full trace family
of $A$, extended to $M$; no assertion that $K=\T(M)$ is needed.
\end{example}

\begin{example}[The topology of a Cartan core is not canonical]
\label{ex:noncanonical-Cartan-core-topology}
Take $D=M=L^\infty([0,1],m)$ with its Lebesgue trace. The interval
core $\mathrm{C}([0,1])$ and the binary Cantor core
$\mathrm{C}(\{0,1\}^{\mathbb N})$, equipped with fair Bernoulli measure,
both have this bounded $2$-norm completion. To identify the latter
core with a subalgebra of $D$, use the binary expansion map
\[
 (\varepsilon_k)_{k\geq1}\longmapsto
 \sum_{k\geq1}\varepsilon_k2^{-k}.
\]
It pushes Bernoulli measure to Lebesgue measure and is one-to-one
after removing the countable ambiguity set of dyadic expansions.
Both cores have relative SBP: interval boundaries can be chosen
finite, and the Cantor space has a clopen basis. The equality
relations and trivial twists give the same completed Cartan pair,
although the two compact spectra are not homeomorphic. Thus
uniqueness concerns the regular uniform measure relation and the
regular cohomology class of its twist, rather than a compact
topological presentation.
\end{example}

\subsection{Intrinsic uniform CFW and Cartan classification}
\label{subsec:intrinsic-CFW}

Reconstruction turns measurewise amenability into a property of the
marked inclusion itself.  The small boundary property supplies the
diagonal approximation needed to make the resulting finite models
compatible.  Throughout this subsection the reconstructed relation
and its cocycle are those of Theorem~\ref{thm:intrinsic-uniform-FM}.

\begin{thm}[Intrinsic uniform Connes--Feldman--Weiss theorem]
\label{thm:intrinsic-uniform-CFW}
Let $(D\subseteq M,K)$ satisfy
Definition~\ref{defn:intrinsic-completed-Cartan}, intrinsic SBP, and
\textup{(F)}.  Write $M_\tau=\pi_\tau(M)''$.  The following conditions
are equivalent:
\begin{enumerate}
\item $M_\tau$ is injective, equivalently hyperfinite, for every
      $\tau\in K$;
\item the reconstructed measured relation is amenable for every
      designated measure;
\item the reconstructed relation has regular uniform Reiter
      amenability, equivalently regular finite-cut approximation
      \textup{(RF)}, equivalently strong uniform AF hyperfiniteness
      in the finite matrix-chart sense;
\item there are a standard unital AF diagonal pair $(E_B,B)$, a
      compact convex designated family $L\subseteq \T(B)$, and a
      marked, expectation-preserving isomorphism
      \[
       (D,M,K)\cong(E_B^L,B^L,L).
      \]
\end{enumerate}
Equivalently, there are increasing unital finite-dimensional
subalgebras $B_n\subseteq M$, with compatible coordinate matrix
units normalizing $D$, such that
\[
 \mathbb{E}(B_n)=B_n\cap D,
 \qquad
 \overline{\bigcup_n(B_n)_1}^{\,\|\cdot\|_{2,K}}=M_1,
 \qquad
 \overline{\bigcup_n(B_n\cap D)_1}^{\,\|\cdot\|_{2,K}}=D_1.
\]
Here $B_n\cap D$ is the coordinate diagonal of $B_n$, and the
coordinate embeddings preserve these diagonals.  No additional
uniform amenability, aperiodicity, Bauer, minimality, or simplicity
assumption is required.
\end{thm}

\begin{proof}
Theorem~\ref{thm:intrinsic-uniform-FM} identifies the marked pair
with a regular twisted principal relation model.  In each tracial
fibre, \textup{(F)} and normalizer generation give a measured Cartan
pair. If the fibre is injective, a conditional expectation
from all bounded operators in its measured regular representation
onto the fibre sends multiplication operators on the arrow space
into the Cartan diagonal, by bimodularity and maximal abelianness.
Equivariance under the graph normalizers gives an invariant mean
as in \cite[Proposition~7]{CFW}: the scalar phases cancel in their
conjugation of multiplication operators. Conversely, measured
amenability gives a hyperfinite exhaustion by \cite{CFW}. The
twisted algebras of its finite subrelations are type I and generate
the measured fibre, so the fibre is hyperfinite. The tracial GNS
preduals are separable by bounded uniform separability. Thus
injectivity and hyperfiniteness agree here \cite{Connes76}, proving
the equivalence of the first two conditions.

Intrinsic SBP gives relative SBP on a separable reconstruction
core.  Theorem~\ref{thm:uniform-CFW-SBP} therefore identifies the
second and third conditions, including the compatible matrix-chart
and finite-diagonal conclusions.  Theorem~\ref{thm:uniform-regular-untwisting}
removes the regular scalar twist, and
Theorem~\ref{thm:uniform-AF-reconstruction} gives the fourth
condition.  Transporting the coordinate AF stages, after quotienting
by their common trace-null ideal if necessary, gives the displayed
finite-dimensional approximation of the marked inclusion.  Their
matrix units normalize the completed diagonal, and the transported
expectation restricts to their coordinate diagonal expectations.

Conversely, the displayed finite-dimensional system gives a
standard AF pair by taking its norm closure and then its designated
uniform completion.  Every designated von Neumann fibre of that
pair is hyperfinite.  This proves the converse implications and the
algebraic formulation.
\end{proof}

\begin{cor}[Untwisting, generators, and the uniform full group]
\label{cor:intrinsic-CFW-consequences}
Under the equivalent conditions of
Theorem~\ref{thm:intrinsic-uniform-CFW}, every regular scalar twist
on the reconstructed relation is a regular coboundary:
\[
 H^2_{\mathrm{reg},K}(R,\mathbb T)=0.
\]
The twisted and untwisted completions are isomorphic by a
diagonal-fixing map preserving the expectation and all designated
traces.  The intrinsic uniform full group
\[
 \mathcal G(D,M)=
 \mathcal N_{\mathcal U(M)}(D)/\mathcal U(D)
\]
has increasing finite subgroups with dense union, and is
topologically amenable for its uniform displacement topology.
If $R$ has infinite classes on one invariant universally conull
set, it has one admissible $\mathbb Z$-generator and a Cantor
homeomorphism model with the designated measure family.
\end{cor}

\begin{proof}
Uniform untwisting is Theorem~\ref{thm:uniform-regular-untwisting}.
Lemma~\ref{lem:AF-normalizer-density}, transported through the
marked AF model, gives the dense increasing finite permutation
groups.  For a continuous affine action of $\mathcal G(D,M)$ on a
compact convex set, the fixed sets of these finite groups are
nonempty compact sets and decrease.  Their intersection is fixed
by their dense union and hence by the full group.  The last
assertion is Theorem~\ref{thm:local-cuts-dynamical-model}.
\end{proof}

The stronger classification statement has an intrinsic condition
on the designated traces.  It does not require that their spatial
presentation initially contain every invariant probability.

\begin{lem}[Trace faces and factorial extreme traces]
\label{lem:intrinsic-trace-face}
Let $(E_B^L,B^L,L)$ be a marked AF model of $(D,M,K)$, and identify
$B$ with its image in $M$.  Restriction identifies $K$ with $L$.
Then $K$ is a closed face of $\T(M)$ if and only if $L$ is a closed
face of $\T(B)$.  Moreover, if $M_\tau$ is a factor for every extreme
point $\tau$ of $K$, then $L$ is automatically a closed face of
$\T(B)$.  If these factors are of type $\mathrm{II}_1$, this is a
diffuse trace face.
\end{lem}

\begin{proof}
Suppose first that $K$ is a face.  For $\tau\in K$, put
$\nu=\tau|_B$, and let
$\nu=t\nu_1+(1-t)\nu_2$, where $0<t<1$ and $\nu_i\in \T(B)$.
Positivity gives $\nu_i\leq c_i\nu$ on $B_+$, with
$c_1=t^{-1}$ and $c_2=(1-t)^{-1}$.  Cauchy--Schwarz gives
\[
 |\nu_i(b)|\leq\sqrt{c_i}\,\|b\|_{2,\tau}
              \leq\sqrt{c_i}\,\|b\|_{2,K}.
\]
Thus $\nu_i$ extends to $M$ by bounded uniform $2$-approximation.
Bounded approximation also preserves products, so the extension
$\rho_i$ is a positive norm-bounded tracial state.  Density gives
$\tau=t\rho_1+(1-t)\rho_2$.  The face property puts $\rho_i$ in
$K$, and hence $\nu_i$ in $L$.  Compactness makes $L$ closed.

Conversely, any convex decomposition of $\tau\in K$ by tracial
states on $M$ has components dominated by a scalar multiple of
$\tau$.  They are therefore uniformly $2$-continuous by the same
inequality.  If $L$ is a face, their restrictions belong to $L$;
bounded uniform density identifies them with the corresponding
traces of $K$.  Compactness again gives closedness.

For the final assertion, bounded density gives
$\pi_\tau(B)''=M_\tau$.  Hence factoriality makes every extreme
point of $L$ an extreme point of $\T(B)$.  The latter is a
metrizable Choquet simplex.  A Choquet representing probability
for $\nu\in L$, supported on the extreme points of $L$, is also
its unique extreme-point representing probability in $\T(B)$.
In an ambient convex decomposition of $\nu$, uniqueness forces
both representing probabilities to be supported on those same
extreme points.  Their barycentres belong to the closed convex set
$L$.  Thus $L$ is a face.  The type $\mathrm{II}_1$ assumption is
exactly the stated diffuseness of its extreme fibres.
\end{proof}

The same dominated-extension argument shows that a compact closed
trace face of a separable $\mathrm{C}^*$-core remains a closed face after its
designated bounded uniform completion.  The designated traces need
not exhaust all norm tracial states of that completion.

\begin{thm}[Intrinsic uniform Cartan classification]
\label{thm:intrinsic-Cartan-classification}
For $i=1,2$, let $(D_i\subseteq M_i,K_i)$ satisfy the hypotheses
and equivalent conditions of
Theorem~\ref{thm:intrinsic-uniform-CFW}.  Assume that
$(M_i)_\tau$ is a type $\mathrm{II}_1$ factor for every extreme
point $\tau$ of $K_i$.  Every affine homeomorphism $K_1\to K_2$
is implemented by an expectation-preserving isomorphism of the
marked completed Cartan pairs, and by a uniform measure isomorphism
of their reconstructed relations.  The designated trace simplex
is therefore a complete invariant in this class.

Each pair admits a simple full-trace AF model and a Cantor minimal
$\mathbb Z$-model.  Its intrinsic uniform full group is a L\'evy
group, with concentrating finite subgroups, and is extremely
amenable.  In particular, two such Cartan diagonals in one ambient
completed algebra, with the same designated trace marking, are
conjugate by a marking-preserving automorphism.
\end{thm}

\begin{proof}
Take the marked standard AF models of
Theorem~\ref{thm:intrinsic-uniform-CFW}.
Lemma~\ref{lem:intrinsic-trace-face} makes their designated families
diffuse closed faces.  Apply
Theorem~\ref{thm:diffuse-face-AF-classification} to the transported
affine homeomorphism and compose with the marked reconstructions.
Theorem~\ref{thm:twisted-uniform-pair-isomorphism} supplies the spatial
isomorphism; uniform untwisting removes the cocycle obstruction.

The designated face is a metrizable Choquet simplex.  The argument
of Corollary~\ref{cor:full-family-simple-AF-model}, using a Cantor
minimal realization of that simplex \cite{Dow91} and diffuse-face
AF classification, supplies a simple full-trace AF replacement.
Corollary~\ref{cor:uniform-CFW-Z-generator} then supplies the
minimal dynamical model.  For this replacement the proof of
Theorem~\ref{thm:uniform-full-group-extreme-amenability} applies:
UHF amplification preserves its full trace simplex, diffuse-face
classification identifies the marked completed pairs, and the
large finite matrix stages have dense permutation groups with
concentrating normalized counting measures.  Thus the intrinsic
full group is L\'evy and extremely amenable.  For the final
assertion use the identity marking of the common trace family.
\end{proof}

\begin{rem}[Scope of the intrinsic conclusions]
\label{rem:intrinsic-CFW-scope}
Uniform AF approximation is approximation in the bounded uniform
$2$-norm; it does not assert norm AF structure of $M$.
The original-topology small subgroupoids of
Theorem~\ref{thm:SBP-groupoid-castles} additionally use a principal
\'etale presentation and the regular-to-open lifting lemma.
The intrinsic completed-pair theorem alone does not specify such a
topological presentation. Scalar untwisting
alone is not an amenability criterion: a nonamenable relation can
already have the trivial cocycle.  An arbitrary designated compact
convex family need not classify its completed pair; the factorial
extreme-trace condition is substantive.  Finally, \textup{(F)}
remains necessary even when every fibre is injective.  Nuclearity
of a separable core implies injectivity of its tracial fibres, but
does not by itself make their designated diagonals maximal abelian.
\end{rem}

\section{Uniform Rokhlin properties and approximation without small boundaries}
\label{sec:AH-without-SBP}

Relative SBP supplies regular transversals and makes finite projection
diagonals uniformly dense in $D_K$. Combined with amenability for each
designated measure, it also gives deterministic regular cuts. The
projection-density conclusion is specific to an AF model. If coefficients are retained as continuous functions on
parameter spaces, the natural building blocks are homogeneous Cartan
pairs, with recursive subhomogeneous blocks allowing classes to split
at boundaries. Regular bounded finite cuts alone also give an ASH
Cartan core, without any choice of global regular transversals. This
section gives precise versions of these replacements. We begin with
explicit homogeneous Cartan models. In
Subsection~\ref{subsec:uniform-local-AH-definition} we define relation
URP and its split version, identify their local AH approximation
formulations, and then compare the two properties on these models.

For actions, we use the uniform Rokhlin property (URP) of
\cite[Definition~3.1]{NiuURP}, recalled in
Subsection~\ref{subsec:URP-local-RSH}.
For a free action of a countably infinite amenable group, with its entire
invariant-probability family, one has~\cite{KS,LW,NiuZ}
\[
 \mathrm{SBP}\quad\Longleftrightarrow\quad
 \mathrm{URP}\ \text{and}\ \operatorname{mdim}(X,\Gamma)=0.
\]
Consequently, URP alone should be distinguished from the strong AF
condition in Definition~\ref{defn:uniform-AF-hyperfinite-relation}.
The latter includes finite projection-diagonal density and therefore
forces relative SBP. A regular atlas by itself provides neither property;
for a free continuous action the canonical atlas exists without SBP or
URP. The action results below retain URP after dropping SBP. Their
estimates also hold for any compact designated subfamily of invariant
measures.

In this section, an \emph{AH algebra} means an inductive limit of finite
direct sums of homogeneous corners $p\mathrm{M}_m(\mathrm{C}(Z))p$, including
$\mathrm{M}_k(\mathrm{C}(Z))$, where $Z$ is compact metrizable and has finite
covering dimension. We do not impose finite CW spectra. An \emph{ASH
algebra} is an inductive limit of subhomogeneous $\mathrm{C}^*$-algebras; the compact
base spaces of the recursive subhomogeneous (RSH) blocks used below may
have infinite dimension. We use the pullback convention for RSH
algebras of \cite{PhillipsRSH}. All assertions concern the Cartan pair together
with its designated traces, rather than an unmarked algebra alone.

\subsection{Elementary homogeneous relations without small boundaries}
\label{subsec:AH-elementary-examples}

The absence of relative SBP need not prevent finite orbit coordinates.
What can fail is approximation of the transverse observables by finite
projection diagonals.  Retaining those observables as commutative
coefficient algebras gives homogeneous blocks instead.

\begin{example}[A split homogeneous block]
\label{ex:split-AH-without-SBP}
Fix $k\geq2$, put $Z=[0,1]$ and
\[
 X=Z\times\{1,\ldots,k\},\qquad
 (t,i)\mathrel R(t',j)\ \Longleftrightarrow\ t=t',
\]
and designate the compact convex family
\[
 K=\left\{\mu_\nu=\frac1k\sum_{i=1}^k
                 \nu\otimes\delta_i:\nu\in\operatorname{Prob}(Z)\right\}.
\]
These measures are invariant under $R$.  For $f\in \mathrm{C}(X)$,
\[
 \|f\|_{2,K}^2
   =\sup_{t\in Z}\frac1k\sum_{i=1}^k|f(t,i)|^2,
 \qquad
 \frac1{\sqrt k}\|f\|_\infty\leq\|f\|_{2,K}
                       \leq\|f\|_\infty.
\]
Consequently $D_K=\mathrm{C}(X)\cong \mathrm{C}(Z)^k$.  The maps
$s_i(t,1)=(t,i)$ give a finite regular matrix chart system for all of
$R$.  Taking this same relation at every level gives an exact finite
exhaustion with zero escape.  The completed untwisted Cartan pair is
already the homogeneous pair
\[
 \bigl(\mathrm{C}(Z)^k,\mathrm{M}_k(\mathrm{C}(Z))\bigr),
 \qquad
 \tau_\nu(a)=\int_Z\operatorname{tr}_k(a(t))\,d\nu(t).
\]
Indeed, on this matrix algebra the designated uniform $2$-norm is
equivalent to the $\mathrm{C}^*$-norm by the same factor $\sqrt k$.

Every projection of $D_K$ is constant on each connected copy of $Z$.
The algebra generated by all these projections is therefore only
$\mathbb C^k$, and cannot approximate, for example, the function
$(t,i)\mapsto t$ in uniform $2$-norm.  Relative SBP fails, as does the
finite-diagonal density clause of uniform AF hyperfiniteness.
Thus an exact finite matrix relation can have a homogeneous Cartan
model while admitting no standard AF Cartan model with the prescribed completed
diagonal and designated traces.  By
Corollary~\ref{cor:relative-SBP-uniform-invariance}, it cannot be
uniform measure isomorphic to a relation on a zero-dimensional compact
space with its corresponding completed continuous diagonal.
\end{example}

\begin{example}[A finite relation without a regular transversal]
\label{ex:cover-AH-without-transversal}
Let $X=\mathbb T$, let $q(z)=z^2$, and let $R$ be equality of the
$q$-images.  For $\nu\in\operatorname{Prob}(\mathbb T)$ designate the
measure obtained by integrating the uniform measure on
$q^{-1}(w)=\{z,-z\}$ against $\nu(w)$.  These measures form a compact
convex invariant family $K$, and
\[
 \|f\|_{2,K}^2
    =\sup_{w\in\mathbb T}\frac12\sum_{z^2=w}|f(z)|^2.
\]
Thus $D_K=\mathrm{C}(\mathbb T)$, again by equivalence with the supremum norm.
The identity and the fixed-point-free involution $z\mapsto-z$ form a
regular finite atlas; $R$ itself is a bounded finite relation.
Nevertheless it has no regular Borel transversal.  Its characteristic
function would be a nontrivial projection of $\mathrm{C}(\mathbb T)$, whereas
this connected space has only the projections $0$ and $1$.  In
particular, the exact finite relation cannot be represented by the
regular split matrix charts used in
Definition~\ref{defn:uniform-AF-hyperfinite-relation}.

Its natural homogeneous Cartan model is explicit. On the fiber
$q^{-1}(w)$ use normalized counting measure. The resulting rank-two
bundle has the global orthonormal frame given by the functions $1$ and $z$.
For the unique expression
\[
 f(z)=a(z^2)+z\,b(z^2),\qquad a,b\in \mathrm{C}(\mathbb T),
\]
multiplication by $f$ is represented in this frame by
\[
 \begin{pmatrix}a(w)&w\,b(w)\\b(w)&a(w)\end{pmatrix}.
\]
The deck involution is represented by
$U=\operatorname{diag}(1,-1)$.  These multiplication operators and $U$
generate $\mathrm{M}_2(\mathrm{C}(\mathbb T))$, and the coefficient-of-identity
expectation for the finite crossed product is faithful.  Fiberwise
commutation with the multiplication operators shows that their
commutant is precisely the displayed copy of $\mathrm{C}(\mathbb T)$; the deck
involution normalizes it.  Hence this is a Cartan pair
\[
 \bigl(\mathrm{C}(\mathbb T),\mathrm{M}_2(\mathrm{C}(\mathbb T))\bigr).
\]
Its diagonal is not isomorphic to the standard split diagonal
$\mathrm{C}(\mathbb T)\oplus \mathrm{C}(\mathbb T)$: the former has connected spectrum
and the latter does not.  The covering monodromy is part of the Cartan
data. The connected completed diagonal also prevents relative SBP.
Thus this bounded finite relation has an exact homogeneous Cartan
model, despite having neither a regular transversal nor relative SBP.
This finite-covering model retains the transverse parameter space
and provides an exact instance of homogeneous Cartan approximation.
\end{example}

\begin{rem}[Finite covering blocks and twists]
\label{rem:finite-covering-AH-blocks}
More generally, let $q:P\to Z$ be a finite covering of constant degree
$k$ between compact metrizable spaces.  The elementary relation
$P\times_ZP$ has the natural homogeneous Cartan pair
\[
 \mathrm{C}(P)\subseteq\Gamma\bigl(\operatorname{End}(q_*\mathbb C)\bigr).
\]
Here $q_*\mathbb C$ is the rank-$k$ bundle with fiber
$\mathbb C^{q^{-1}(z)}$.  Locally the pair is
$\mathbb C^k\subseteq \mathrm{M}_k(\mathbb C)$ over the parameter space; its
coordinate permutations on overlaps retain the finite covering.
Finite kernels act as endomorphisms of this bundle, giving the stated
algebra, and the diagonal expectation and local bisection normalizers
give the Cartan structure.  With all mixtures of the uniform fiber
measures designated, the uniform completions agree with these $\mathrm{C}^*$-norm
algebras.  A globally split covering recovers
Example~\ref{ex:split-AH-without-SBP}.

If $L\to P$ is a Hermitian line bundle, replace the block by
$\Gamma(\operatorname{End}(q_*L))$.  Multiplication still embeds $\mathrm{C}(P)$
as its Cartan diagonal.  The unitary maps from $L_y$ to $L_x$, for
$(x,y)\in P\times_ZP$, form a circle twist whose convolution algebra
is this block.  These endomorphism algebras are homogeneous corners
$p\mathrm{M}_N(\mathrm{C}(Z))p$, since the underlying bundles are finite projective.
They exhibit the additional local coordinates and bundle information
that a nonsplit AH formulation should retain; global regular
transversals and scalar untwisting cannot simply be assumed from
finite class sizes.
\end{rem}

\subsection{Subhomogeneous reconstruction without small boundaries}
\label{subsec:RF-ASH-without-SBP}

Here an \emph{ASH algebra} means an inductive limit of separable
subhomogeneous $\mathrm{C}^*$-algebras; the bounds on their irreducible
representation dimensions may increase with the stage.  This
conclusion does not require the stages to have globally split matrix
coordinates or homogeneous coefficient spaces.

\begin{thm}[Regular finite cuts give an ASH Cartan core]
\label{thm:RF-ASH-without-SBP}
Let $(X,R,K,D_K,\mathscr P)$ be a regular uniform measure relation in
the probability-measure-preserving setting.  Assume the regular
finite-cut condition \eqref{eq:local-regular-cuts}.  There are a
separable unital ASH algebra $B$, a Cartan subalgebra $F=\mathrm{C}(Y)\subseteq B$,
and a compact convex designated family $L\subseteq \T(B)$ such that
restriction gives an affine homeomorphism $H:K\to L$ and
\[
 (F^L,B^L)\cong(D_K,M_R)
\]
as completed Cartan pairs, preserving the diagonal expectation and
each corresponding designated trace.  More precisely, $B$ has an
increasing presentation $B=\overline{\bigcup_nB_n}^{\|\cdot\|}$ by
unital subhomogeneous algebras, with $F\cap B_n$ Cartan in $B_n$.
No relative SBP assumption is made, and no assertion that $L=\T(B)$
is made.
\end{thm}

\begin{proof}
First construct increasing bounded finite regularly generated
subrelations $S_n$ with uniform escape tending to zero.  This is
the nesting part of Theorem~\ref{thm:local-cuts-AF}, before its use
of regular transversals.  For clarity, enumerate the regular atlas
as $(s_j)$, choose positive summable $\eta_n\to0$, and start with the
identity relation.  At stage $n$, apply (RF) to a finite arrow cover
of $S_{n-1}$ together with $s_1,\ldots,s_n$.  Make the approximation
sufficiently accurate that Lemma~\ref{lem:regular-cut-nesting} gives
$S_n\supseteq S_{n-1}$ with total escape for these $n$ atlas tests
below $\eta_n$.  The old class-size bound is fixed before this
choice, so the nesting error can be made arbitrarily small.
Neither this step nor the finite-class integration estimate in that
lemma requires regular transversals or SBP.  Borel--Cantelli,
followed by removal of the saturation of the exceptional sets,
gives $R=\bigcup_nS_n$ on one invariant $K$-conull set.  The
first-chart decomposition and
\eqref{eq:regular-relation-uniform-tails} extend the uniform escape
conclusion to every admissible partial map.

Let $N_n$ bound the class sizes of $S_n$.  Choose a finite inverse-closed
generating collection $\mathcal T_n$ of admissible maps for $S_n$,
including the identity and every previously chosen generator.  Put
\[
 A_n=\mathrm{C}^*(D_K,v_t:t\in\mathcal T_n)\subseteq M_R.
\]
This algebra is $N_n$-subhomogeneous.  Indeed, for each $\mu\in K$,
choose an ordinary Borel transversal for the finite relation $S_n$.
The finite-relation kernel representation decomposes over that
transversal into matrices of sizes $|[x]_{S_n}|\leq N_n$.  Its
integrated normalized matrix trace is $\tau_\mu|_{A_n}$.  These
representations are jointly faithful because the designated traces
are jointly faithful on $M_R$.  Consequently $A_n$ satisfies the
standard polynomial identity
\[
 \sum_{\pi\in\operatorname{Sym}(2N_n)}
       \operatorname{sgn}(\pi)\,
       a_{\pi(1)}\cdots a_{\pi(2N_n)}=0
       \qquad(a_1,\ldots,a_{2N_n}\in A_n),
\]
which holds in every matrix algebra of size at most $N_n$.  The
standard polynomial-identity criterion for $\mathrm{C}^*$-algebras bounds the
dimension of every irreducible representation by $N_n$.  Equivalently,
the converse bound follows by applying irreducible transitivity to
the failure of this identity in $\mathrm{M}_{N_n+1}$.  A regular transversal
has not been used.

We next choose nested separable diagonals $F_n\subseteq D_K$, starting
with $F_0$ equal to the image of $\mathrm{C}(X)$.
At stage $n$, start with the image of $\mathrm{C}(X)$, the preceding $F_{n-1}$,
and the support and coincidence projections of all partial words in
$\mathcal T_n$.  Close under every corner pullback associated to
$\mathcal T_n$ and their inverses, iterating the countable dense-set
construction of
Proposition~\ref{prop:uniform-relation-separable-Cartan-core}.
The norm closure is a separable unital algebra $F_n$ invariant under
these corner maps.  Define
\[
 B_n=\mathrm{C}^*(F_n,v_t:t\in\mathcal T_n).
\]
The choices imply $B_{n-1}\subseteq B_n$, and the preceding
polynomial identity makes $B_n$ subhomogeneous.  Moving diagonal
factors through graph words, and then taking their identity
coefficients, shows that $\mathbb{E}(B_n)\subseteq F_n$.  Thus
$B_n\cap D_K=F_n$.  The restricted expectation is faithful, and
the generating graph operators normalize $F_n$.  If $b\in B_n$
commutes with $F_n$, it commutes with $\mathrm{C}(X)$. Bounded uniform density
of $\mathrm{C}(X)$ in $D_K$ and continuity of multiplication then show that it
commutes with $D_K$. The MASA assertion of
Lemma~\ref{lem:uniform-relation-Fourier-approximation} gives
$b=\mathbb{E}(b)\in F_n$.  These facts show that $F_n$ is Cartan in $B_n$.

Set
\[
 B=\overline{\bigcup_nB_n}^{\|\cdot\|},\qquad
 F=\overline{\bigcup_nF_n}^{\|\cdot\|}.
\]
Then $B$ is separable, unital and ASH, and $\mathbb{E}(B)\subseteq F$.
The same completed-diagonal argument proves maximal abelianness of $F$.
Every graph generator from stage $n$ normalizes all $F_m$ for $m\geq n$
by construction, hence normalizes $F$.  These normalizers and $F$
generate $B$, so $F\subseteq B$ is a Cartan pair.  In particular
$F\cap B_n=F_n$.

To check the completed pair, first note that $\mathrm{C}(X)\subseteq F\subseteq
D_K$, so the bounded uniform completion of $F$ is $D_K$.  Fix an
admissible partial map $\varphi$, and restrict it to the regular set
on which $(x,\varphi x)\in S_n$, obtaining $\varphi_n$.  Words of
length at most $N_n-1$ in $\mathcal T_n$ enumerate each $S_n$-class.
Successively assigning each arrow of $\varphi_n$ to the first
matching word expresses its graph operator as
\[
 v_{\varphi_n}=\sum_{\ell=1}^{q_n}v_{w_\ell}p_\ell\in A_n,
\]
where the $p_\ell\in D_K$ are regular coincidence projections with
disjoint source supports.  Since $\mathrm{C}(X)\subseteq F_n$ is bounded
uniformly dense in $D_K$, approximate these finitely many coefficients
by elements of $F_n$.  This approximates $v_{\varphi_n}$ by elements
of $B_n$ in uniform $2$-norm.  Radial truncation at norm two retains
a bound independent of $n$, by the estimate in the proof of
Theorem~\ref{thm:uniform-AF-reconstruction}.  Uniform escape gives
\[
 \|v_\varphi-v_{\varphi_n}\|_{2,K}^2
 =\sup_{\mu\in K}\mu\{x\in\operatorname{dom}\varphi:
                           (x,\varphi x)\notin S_n\}
 \longrightarrow0.
\]
Thus the bounded uniform completion of $B$ contains $D_K$ and every
original atlas graph operator, and is exactly $M_R$.  The restricted
expectation completes to its original faithful diagonal expectation.

Finally, put $H(\mu)=\tau_\mu|_B$ and $L=H(K)$.  For $b\in B$,
the function $\mathbb{E}(b)\in D_K$ has a continuous integral profile on $K$,
by bounded uniform approximation from $\mathrm{C}(X)$.  Hence $H$ is continuous
and affine.  It is injective because $F$ contains $\mathrm{C}(X)$, and compactness
of $K$ makes it a homeomorphism onto the compact convex family $L$.
The preceding inclusions induce the claimed marked completed-pair
isomorphism.
\end{proof}

\begin{rem}[Retaining a regular scalar twist]
\label{rem:RF-ASH-retains-twist}
The ASH conclusion also applies with a regular scalar $2$-cocycle
retained in the completed pair. In the construction of $F_n$, adjoin
the countably many regular phase coefficients in products and adjoints
of the chosen graph words, and close under their pullbacks as before.
On a bounded finite measured relation an ordinary Borel transversal
untwists each class, so its twisted measured representations are still
matrices of sizes at most $N_n$. The polynomial identity, expectation,
normalization and bounded-density arguments therefore remain valid.
This gives a twisted ASH Cartan core, without claiming a regular
global coboundary. Scalar untwisting from compatible regular matrix
charts is the stronger conclusion recorded below.
\end{rem}

\subsection{AH reconstruction from regular matrix charts}
\label{subsec:AH-without-SBP}

The pro\-jec\-tion-di\-ag\-o\-nal-den\-sity condition in
Definition~\ref{defn:uniform-AF-hyperfinite-relation} is the part of the
AF construction that forces relative SBP.  The finite orbit coordinates
and the approximation of diagonal observables can be separated.  If the
coordinates are retained, a commutative coefficient algebra can replace
the scalar coefficient in each finite matrix block.

The theorem below gives full matrix blocks $\mathrm{M}_k(\mathrm{C}(Z))$ with finite
covering-dimensional compact metrizable coefficient spectra. No finite
CW structure, uniform bound on $\dim Z$, or dimension-growth estimate
is asserted.

\begin{thm}[Conditional AH reconstruction without SBP]
\label{thm:AH-reconstruction-without-SBP}
Let $(X,R,K,D_K,\mathscr P)$ be a uniform measure relation.  Suppose
that there are nested finite $D_K$-regular matrix-chart subrelations
$(R_n)$ with coordinate refinement, exhausting $R$ on one invariant
$K$-conull set, and that for every $\varphi\in\mathscr P$,
\[
 \lim_n\sup_{\mu\in K}
 \mu\{x\in\operatorname{dom}\varphi:
             (x,\varphi x)\notin R_n\}=0.
\]
No density condition on the finite projection diagonals is assumed.
Then there is a unital separable AH algebra $B$ with a $\mathrm{C}^*$-Cartan
subalgebra $F=\mathrm{C}(Y)$ and a compact convex family $L\subseteq \T(B)$ such
that restriction gives an affine homeomorphism $K\to L$ and
\[
                  (D_K,M_R)\cong(F^L,B^L).
\]
The expectation and every designated trace are preserved.  The diagonal
$F$ can be chosen to contain the image of $\mathrm{C}(X)$ in $D_K$.
\end{thm}

\begin{proof}
Write $e^{(n,a)}_{ij}$ for the level $n$ matrix units and
$p_{n,a}=e^{(n,a)}_{11}$.  Discard common-null blocks.  Coordinate
refinement makes the matrix-unit algebras nested, and
\[
 C_n=\mathrm{C}^*(D_K,\{e^{(n,a)}_{ij}\}_{a,i,j})
       \cong\bigoplus_a \mathrm{M}_{k_{n,a}}(p_{n,a}D_K)
\]
are nested as well.  The coefficient $d\in p_{n,a}D_K$ occupies the
$ij$ entry as $e^{(n,a)}_{i1}d e^{(n,a)}_{1j}$.

Choose a sequence $(f_j)$ of continuous functions whose linear span is
norm dense in $\mathrm{C}(X)$.  We recursively choose finitely generated unital
commutative subalgebras
$A_{n,a}\subseteq p_{n,a}D_K$ and put
\[
 B_n=\bigoplus_a \mathrm{M}_{k_{n,a}}(A_{n,a})\subseteq C_n,
 \qquad F_n=B_n\cap D_K.
\]
At stage $n$, include in the coefficient algebras every matrix entry,
in level $n$ coordinates, of a finite generating set of $B_{n-1}$ and
of $f_1,\ldots,f_n$.  There are finitely many such entries, all lying
in the indicated diagonal corners.  This ensures
$B_{n-1}\subseteq B_n$ and $f_1,\ldots,f_n\in F_n$.
Each $A_{n,a}$ is generated by finitely many self-adjoint elements, so
its spectrum $Z_{n,a}$ is a compact subset of some $\mathbb R^{m_{n,a}}$.
Consequently
\[
               B_n\cong\bigoplus_a \mathrm{M}_{k_{n,a}}(\mathrm{C}(Z_{n,a})),
 \qquad \dim Z_{n,a}\leq m_{n,a}<\infty.
\]
Set $B=\overline{\bigcup_n B_n}^{\|\cdot\|}$ and
$F=\overline{\bigcup_n F_n}^{\|\cdot\|}$.  Thus $B$ is AH and
$\mathrm{C}(X)\subseteq F\subseteq D_K$.  Since $D_K$ is the bounded uniform
$2$-norm completion of $\mathrm{C}(X)$, the bounded completion of $F$ is $D_K$.

The diagonal expectation of $M_R$ maps each $B_n$ onto $F_n$ and hence
restricts to a faithful expectation $B\to F$.  To prove the MASA
property in norm, let $b\in B$ commute with $F$ and choose
$b_n\in B_n$ with $\|b-b_n\|\to0$.  If
$P_n(z)=\sum_{a,i}e^{(n,a)}_{ii}z e^{(n,a)}_{ii}$, then
$P_n(b)=b$ and $P_n(b_n)$ belongs to $F_n$, so
\[
                    \|b-P_n(b_n)\|\leq\|b-b_n\|\longrightarrow0.
\]
It follows that $b\in F$.  Every level $n$ matrix unit normalizes
every $F_m$ for $m\geq n$: its coordinate refinement is a sum of
matching level $m$ matrix units with pairwise disjoint sources and
ranges.  Therefore these matrix units normalize $F$.  Together with
$F$ they generate $B$, proving norm regularity and the $\mathrm{C}^*$-Cartan
assertion.

For $\mu\in K$ let $H(\mu)=\tau_\mu|_B$.  Integration of a regular
diagonal coefficient is continuous on $K$, since that coefficient is
a bounded uniform $2$-norm limit of continuous functions.  The finite
matrix-entry formula and norm approximation therefore show that $H$
is continuous and affine.  It is injective because $F$ contains
$\mathrm{C}(X)$: equality of the restricted traces implies equality of the
integrals of every continuous function.  Taking $L=H(K)$ gives the
asserted affine homeomorphism and the equality of the diagonal
completions.

For the ambient completion, restrict an admissible map $\varphi$ to
$\varphi_n$ on the set where its graph belongs to $R_n$.  The domain
of this restriction is regular, being a finite union of coincidence
sets with level $n$ coordinate maps.  Its graph operator is a
contraction in $C_n$, and
\[
 \|v_\varphi-v_{\varphi_n}\|_{2,K}^{\,2}
 =\sup_{\mu\in K}
   \mu\{x\in\operatorname{dom}\varphi:(x,\varphi x)\notin R_n\}
 \longrightarrow0.
\]
Every element of $C_n$ is a finite matrix with entries in the corners
$p_{n,a}D_K$.  Approximate each entry by later $F_m$, $m\geq n$,
and compress the approximation by $p_{n,a}\in F_m$.  The compressed
entry lies in $p_{n,a}F_m$ and has no greater $2$-norm error.  Inserting
it between the old matrix units, which belong to $B_m$, gives
uniform $2$-norm approximations in $B_m$.  Bounded approximations can
be retained as
in the proof of Theorem~\ref{thm:uniform-AF-reconstruction}: if a
contraction $z\in C_n$ is approximated by $b\in B_m$, its radial
truncation $b^{[2]}$ satisfies
\[
 \|b^{[2]}\|\leq2,
 \qquad
 \|b^{[2]}-z\|_{2,K}\leq3\|b-z\|_{2,K}.
\]
Thus the completed image of $B$ contains $D_K$ and every original
atlas graph operator.  It is all of $M_R$.  The inherited expectation
and traces give the required isomorphism of completed Cartan pairs.
\end{proof}

\begin{rem}[The spatial AH relation]
\label{rem:AH-spatial-relation}
The level $n$ matrix units induce partial homeomorphisms of the clopen
subsets of $Y$ determined by their diagonal projections.  They form
a compact finite equivalence relation $S_n$ on $Y$, and coordinate
refinement gives $S_n\subseteq S_{n+1}$ as an open inclusion of \'etale
relations.  Their union, with the inductive-limit groupoid topology,
is the corresponding spatial analogue of an AF tail relation.
The bases of its finite matrix blocks are compact subspaces of $Y$;
they need not be zero-dimensional or finite-dimensional.  The
finite-dimensional spaces $Z_{n,a}$ above are the coefficient
spectra of the chosen AH core, rather than these full spatial bases.

The original relation is carried onto this union on invariant sets
conull for the designated measures.  Indeed, evaluating a countable
norm-dense algebra in $F$ gives a Borel coding map $X\to Y$ outside
a common null set.  The included continuous functions separate points,
so the map is injective there.  The matrix units intertwine the finite
coordinate maps, exactly as in the spatial proof of
Theorem~\ref{thm:uniform-AF-reconstruction}; saturation of the common
exceptional set gives invariant domains.  Thus $Y$ can have positive
topological dimension without changing the designated completed
relation.  This conclusion does not assert a Cantor model.
\end{rem}

\begin{rem}[What the hypothesis and the twist still require]
\label{rem:AH-without-SBP-scope}
The theorem retains regular finite transversals and compatible matrix
coordinates.  Without SBP, one-stage regular finite cuts alone do not
supply these data by Lemma~\ref{lem:regular-finite-transversals}; the
lemma uses SBP to make its ordering tests regular.  The theorem is
therefore a reconstruction result once the chart structure is given,
not an implication from amenability or URP alone.

Scalar untwisting already follows from nested regular finite matrix charts
with common conull exhaustion and uniformly vanishing escape.
Theorem~\ref{thm:uniform-regular-untwisting} does not require coordinate
refinement, SBP, or density of the finite projection diagonals.  Accordingly a regular scalar
$2$-cocycle yields the same completed AH Cartan pair after a regular
diagonal-fixing gauge.  This does not trivialize an arbitrary
topological circle bundle without first establishing a regular
presentation of its twist.
\end{rem}

\subsection{Local RSH approximation from URP}
\label{subsec:URP-local-RSH}

Here $X$ is compact metrizable, $\Gamma$ is countable amenable, the
action is free, and $K$ is a nonempty compact family of invariant
probability measures.  Write $A_0=\mathrm{C}(X)\rtimes_r\Gamma$ and
$\norm{a}_{2,K}^2=\sup_{\mu\in K}\tau_\mu(a^*a)$, where
$\tau_\mu=\mu\circ \mathbb{E}$. The weak-star closed convex hull
$\overline{\operatorname{co}}K$ gives the same uniform seminorm, since
$\mu\mapsto\tau_\mu(a^*a)$ is continuous and affine for $a\in A_0$.
Thus the completion for a compact subfamily is the same algebra, with the
traces in the original $K$ designated. We use right-action notation $xg$.
The uniform Rokhlin property (URP) supplies finitely many closed
tower bases $B_s$ and finite shapes $F_s$ such that their translated
levels are disjoint, the shapes are as invariant as prescribed, and
the complement has arbitrarily small orbit capacity
\cite[Definition~3.1]{NiuURP}. Orbit capacity bounds the measure of
the complement for every invariant probability measure.

\begin{prop}\label{prop:URP-local-RSH}
Suppose $(X,\Gamma)$ has URP.  For every finite set $L\subseteq\Gamma$
and every $\eps>0$, there is an open relatively compact principal
subgroupoid $G\subseteq X\rtimes\Gamma$, containing all units, such
that its finite equivalence relation $S_G$ satisfies
\[
 \sup_{\mu\in K}\mu\{x:(x,xg)\notin S_G\}<\eps,
 \qquad g\in L.
\]
The inclusion $\mathrm{C}(X)\subseteq \mathrm{C}^*_r(G)\subseteq A_0$ is a Cartan
inclusion, and $\mathrm{C}^*_r(G)$ is recursively subhomogeneous with diagonal
gluing maps.  Moreover, for every finite $\mathcal F\subseteq A_0$
and every $\eps>0$, such a subgroupoid can be chosen so that
\[
 \dist_{2,K}(a,\mathrm{C}^*_r(G))<\eps,\qquad a\in\mathcal F.
\]
No small boundary property is required.
\end{prop}

\begin{proof}
Fix $0<\eta<1$ and choose URP towers whose uncovered measure is
less than $\eta$ for every invariant measure, with
\[
 |\{h\in F_s:hg\notin F_s\}|<\eta|F_s|,
 \qquad g\in L.
\]
Thicken the closed bases to open sets $U_s\supseteq B_s$ so that
all the sets $\overline{U_s}h$, $h\in F_s$, remain mutually
disjoint.  There are only finitely many disjoint compact levels,
so these neighborhoods exist.  Enlarging the bases does not increase
the uncovered measure.  Let $G$ contain all units and all arrows
between the points of each finite tower fibre
$\{bh:h\in F_s\}$, $b\in U_s$.  In transformation-groupoid
coordinates these arrows are
\[
 (bh,h^{-1}k),\qquad b\in U_s,\quad h,k\in F_s.
\]
This is an open subgroupoid with finite group-coordinate support;
hence it is relatively compact.  Its classes have cardinality at
most $\max_s|F_s|$, and points outside the towers are singleton
classes.  Invariance and disjointness of the levels give
\begin{align*}
 \mu\{x:(x,xg)\notin S_G\}
 &\leq \mu\bigl(X\setminus\textstyle\bigcup_{s,h}U_sh\bigr)
       +\sum_s|\{h\in F_s:hg\notin F_s\}|\mu(U_s)\\
 &<\eta+\eta\sum_s|F_s|\mu(U_s)\leq2\eta.
\end{align*}

The recursively subhomogeneous structure is explicit.  Put
\[
 Z_s=\overline{U_s},\qquad
 Y=X\setminus\bigcup_{s,h\in F_s}U_sh.
\]
Then $\mathrm{C}^*_r(G)$ is the pullback algebra
\begin{equation}\label{eq:tower-RSH-pullback}
 \left\{\bigl(f,(a_s)_s\bigr):
 \begin{array}{l}
 f\in \mathrm{C}(Y),\quad a_s\in \mathrm{M}_{|F_s|}(\mathrm{C}(Z_s)),\\
 a_s(b)=\operatorname{diag}\bigl(f(bh)\bigr)_{h\in F_s},
 \quad b\in\partial U_s
 \end{array}\right\}.
\end{equation}
Indeed, the off-diagonal matrix entries are functions in $\mathrm{C}_0(U_s)$,
while diagonal entries glue to the function on $Y$ at the tower
boundaries.  Conversely, these entries extend to continuous
groupoid functions and generate the stated pullback.  This is a
recursive subhomogeneous decomposition with diagonal attaching
maps; it is a special case of \cite[Theorem~6.15]{NiuURP}.
The standard matrix diagonal identifies exactly with $\mathrm{C}(X)$.
The groupoid is principal and \'etale, so $\mathrm{C}(X)$ is maximal abelian;
the continuous functions on open bisections normalize it and span
a dense subalgebra.  The restriction expectation is faithful.
These facts prove the Cartan assertion.

For the approximation assertion, first take a finite collection of
Fourier polynomials and include their group coordinates in $L$.
For $g\in L$, let
\[
 D_g=\{x:(x,g)\in G\}.
\]
This set is open and $\inf_{\mu\in K}\mu(D_g)>1-2\eta$.
Choose continuous functions $0\leq\chi_n\leq1$, compactly
supported in $D_g$, increasing pointwise to $1_{D_g}$.
The open sets
\[
 \left\{\mu\in K:\int\chi_n\,d\mu>1-3\eta\right\}
\]
increase and cover the compact space $K$.  Consequently one
$\chi_g=\chi_n$ satisfies this inequality for every $\mu\in K$.
Write $\chi_gu_g$ for the corresponding bisection function, transporting
the cutoff to the coefficient side if required by convention. It
belongs to $\mathrm{C}^*_r(G)$ and
\[
 \norm{u_g-\chi_gu_g}_{2,K}^2
 \leq\sup_{\mu\in K}\int(1-\chi_g)^2\,d\mu<3\eta.
\]
The choice of source or range for the coefficient does not change
this estimate, by invariance.  Thus, for a Fourier polynomial,
\[
 \norm{\sum_{g\in L}f_gu_g-
             \sum_{g\in L}f_g\chi_gu_g}_{2,K}
 \leq\sqrt{3\eta}\sum_{g\in L}\norm{f_g}.
\]
Take $\eta$ sufficiently small and then use norm density of Fourier
polynomials to obtain the assertion for $\mathcal F$.
\end{proof}

\begin{rem}\label{rem:URP-RSH-boundaries}
The finite relation in Proposition~\ref{prop:URP-local-RSH} has
open arrow domains, but their characteristic functions need not
belong to $D_K$.  Thus these stages are not asserted to be finite
regular matrix charts in Definition~\ref{defn:uniform-AF-hyperfinite-relation}.
Their boundary behavior is retained by the pullback
\eqref{eq:tower-RSH-pullback}: a matrix block splits into singleton
diagonal blocks at a tower boundary.  SBP permits regular projection
cuts instead.  Without SBP, positive cutoff functions supply the
uniform approximation while the subhomogeneous stage retains the
topology of its coefficient spaces.
\end{rem}

There is also a local homogeneous formulation with dimension
control. Given a finite collection $a_1,\ldots,a_m\in A_0$ and
$\delta>0$, apply \cite[Theorem~3.9]{NiuURP} with tolerance $\delta/2$.
Under URP this gives a subalgebra
\[
 H\cong\bigoplus_s \mathrm{M}_{r_s}(\mathrm{C}_0(W_s)),
\]
whose standard diagonal lies in $\mathrm{C}(X)$, a positive contraction
$p\in H\cap \mathrm{C}(X)$, and elements $a_i'\in A_0$ such that
\[
 \norm{a_i-a_i'}<\delta,\qquad pa_i'p\in H,\qquad
 \sup_{\mu\in K}\mu\{p\ne1\}<\delta.
\]
It follows that
\begin{equation}\label{eq:URP-homogeneous-uniform-two}
 \dist_{2,K}(a_i,H)
 \leq\delta+2(\norm{a_i}+\delta)\sqrt\delta.
\end{equation}
The same theorem gives compact cores $Q_s\subseteq W_s$ with
\[
 \frac{\dim Q_s}{r_s}<\operatorname{mdim}(X,\Gamma)+\delta.
\]
This is a dimension bound on the compact cores, not a bound on all
of $W_s$.  Neither this theorem nor
Proposition~\ref{prop:URP-local-RSH} supplies nested stages or
compatible connecting maps.  A general AH or ASH Cartan
inductive-limit model therefore remains an additional reconstruction
problem. The following local definition retains homogeneous
coefficient algebras without imposing those connecting maps.

\subsection{Relation URP and local AH Cartan approximation}
\label{subsec:uniform-local-AH-definition}

Local approximation does not require nested finite relations or
connecting maps.  It can be formulated directly for the completed,
marked Cartan pair.  Homogeneous blocks with nonsplit Cartan diagonals
retain the finite-covering and bundle examples above; a split version
records the stronger form supplied by Rokhlin towers.

We first describe the finite homogeneous relation models used for
relation URP.  Let $q:Y\to Z$ be a finite covering of constant degree
$r$, where $Z$ is finite-dimensional, locally compact and metrizable,
and let $L\to Y$ be a Hermitian line bundle.  Put
\[
 Q=Y\times_ZY,\qquad
 (\delta L)_{(y,y')}=\operatorname{Hom}(L_{y'},L_y),\qquad
 \mathcal V=q_*L.
\]
Composition and adjoint of the fiber maps define the convolution and
involution of this line bundle over $Q$.  Require $\mathcal V$ to be
isomorphic to the range of a continuous constant-rank projection
$p\in \mathrm{M}_m(\mathrm{C}_b(Z))$ for some finite $m$.  Finite matrix convolution gives
\begin{equation}\label{eq:finite-cover-homogeneous-kernels}
 A(q,L):=\mathrm{C}^*_r(Q,\delta L)
 \cong\Gamma_0(Z,\operatorname{End}\mathcal V)
 \cong p\mathrm{M}_m(\mathrm{C}_0(Z))p.
\end{equation}
Multiplication on the individual lines embeds $\mathrm{C}_0(Y)$ as its Cartan
diagonal; write $\mathbb{E}_Q$ for the diagonal expectation.  These models retain
covering monodromy and line-bundle data.  The specified coefficient bundle $\delta L$ is part of the model. Its
compactly supported bisection sections can serve as positive-cutoff
normalizers even when the covering relation with its own fibre
measures admits no regular atlas. For instance, for connected $Y$
and compact $Z$, designation of all mixtures of uniform fibre measures
gives $D_K=\mathrm{C}(Y)$ by norm equivalence. Every nonempty regular chart
then has full source and range, since their indicators are continuous
projections. Its corner pullback on $\mathrm{C}(Y)$ makes it a homeomorphism
of $Y$ over $Z$. There are no nonempty common-null sets, so such
charts cover the relation only if the deck transformations act
transitively on each fibre. The embedded-block definition below uses
local continuous bisection sections of the covering itself.

\begin{defn}[Uniform Rokhlin property for relations]
\label{defn:relation-URP}
A completed regular prob\-abil\-ity-mea\-sure-pre\-serv\-ing relation
$(D_K\subseteq M_R,\mathbb{E}_R,K)$ has the \emph{uniform Rokhlin property
for relations}, or \emph{relation URP}, if the following holds.
For every finite $\mathcal D\subseteq(D_K)_1$, every finite
$F\subseteq\mathscr P$, and every $\varepsilon>0$, there are finitely
many models $A(q_s,L_s)$ as above and an injective *-homomorphism
\[
 \Phi:A=\bigoplus_s A(q_s,L_s)\longrightarrow M_R,
 \qquad H=\Phi(A),\qquad B=\bigoplus_s \mathrm{C}_0(Y_s),
\]
with these properties:
\begin{enumerate}
\item $\Phi(B)=H\cap D_K$ and $\mathbb{E}_R\circ\Phi=\Phi\circ \mathbb{E}_Q$;
\item $A$ has a generating family of compactly supported bisection
sections whose images under $\Phi$ normalize $D_K$;
\item there are contractions $h_f,h_t\in H$ such that
\[
 \norm{f-h_f}_{2,K}<\varepsilon\quad(f\in\mathcal D),
 \qquad
 \norm{v_t-h_t}_{2,K}<\varepsilon\quad(t\in F).
\]
Here $v_t$ is the graph partial isometry of the admissible map $t$.
\end{enumerate}
The property is \emph{split} if the models can be chosen with
$Y_s=Z_s\times\{1,\ldots,r_s\}$ and trivial $L_s$.
Neither nesting of models nor unitality of $\Phi$ is required.
\end{defn}

This definition permits positive cutoff normalizers; their support
projections need not belong to $D_K$.  It therefore differs from
regular finite-cut approximation \textup{(RF)}.  The model relations
are finite coverings over parameter spaces, represented in the
completed pair; they need not be admissible finite Borel subrelations
of the original $R$.  This is an intrinsic formulation, rather than
the group-labelled tower formulation of action URP.  The equivalence
below makes it independent of the chosen equivalent regular atlas.

\begin{defn}\label{defn:uniform-local-AH-Cartan}
Let $(D_K\subseteq M_R,\mathbb{E}_R,K)$ be a completed regular uniform measure
relation in the probability-measure-preserving setting.  It has
\emph{uniform local AH Cartan approximation} if, for every finite
$\mathcal F\subseteq(M_R)_1$ and every $\eps>0$, there is a
$\mathrm{C}^*$-subalgebra $H\subseteq M_R$ with the following properties.
\begin{enumerate}
\item There are finitely many finite-dimensional locally compact
metrizable spaces $Z_s$, integers $m_s,r_s\geq1$, and continuous
projections $p_s\in \mathrm{M}_{m_s}(\mathrm{C}_b(Z_s))$ of constant rank $r_s$ such
that
\[
 H\cong\bigoplus_s p_s\mathrm{M}_{m_s}(\mathrm{C}_0(Z_s))p_s.
\]
\item $D_H=H\cap D_K$ is Cartan in $H$, and
$\mathbb{E}_R(H)=D_H$.  There is a family of $\mathrm{C}^*$-normalizers generating $H$
whose members normalize both $D_H$ inside $H$ and $D_K$ inside
$M_R$: for each such $n$,
$nD_Kn^*\subseteq D_K$ and $n^*D_Kn\subseteq D_K$.
\item For every $a\in\mathcal F$ there is $b_a\in H_1$ with
$\norm{a-b_a}_{2,K}<\eps$.
\end{enumerate}
The property is \emph{split} if the approximants can be chosen with
\[
 (D_H\subseteq H)\cong
 \left(\bigoplus_s \mathrm{C}_0(Z_s)^{r_s} \subseteq\
             \bigoplus_s \mathrm{M}_{r_s}(\mathrm{C}_0(Z_s))\right).
\]
The blocks may be nonunital.  They are embedded subalgebras of the
original marked completion, and do not replace its designated trace
family by all traces on the blocks.
\end{defn}

The expectation condition supplies the diagonal approximation
automatically.  For $f\in D_K$ and $b\in H$,
\[
 \norm{f-\mathbb{E}_R(b)}_{2,K}\leq\norm{f-b}_{2,K},
 \qquad \mathbb{E}_R(b)\in D_H.
\]
The norm bound in the definition passes to the expectation.  Normalizer
compatibility retains the relation's Cartan structure, rather than
only a local approximation of its ambient algebra.  It does not
assert that the support projections of these normalizers belong to
$D_K$; positive cutoff normalizers are allowed.  They are not
required to be the graph partial isometries of the admissible
pseudogroup.  The definition asserts an embedded norm Cartan pair
$D_H\subseteq H$, not the strong partial-isometry Cartan property
for its designated completion as a subpair of $D_K\subseteq M_R$.

For the split version, \emph{locally bounded dimension ratio} means
that, for every $\mathcal F\subseteq(M_R)_1$ finite and $\eps>0$,
there is $N=N(\mathcal F,\eps)<\infty$ such that for every integer
$k\geq1$ the approximants in Definition~\ref{defn:uniform-local-AH-Cartan}
can be chosen with
\begin{equation}\label{eq:local-AH-bounded-ratio-quantifiers}
 \min_s r_s\geq k,
 \qquad\max_s\frac{\dim Z_s}{r_s}<N.
\end{equation}
The constant is independent of $k$, but may depend on the finite
test set and its approximation tolerance.  This is a local bound,
not a slow-dimension-growth assertion for an inductive system.
For split relation URP, locally bounded dimension ratio means the
same bounds in \eqref{eq:local-AH-bounded-ratio-quantifiers}, with
$N=N(\mathcal D,F,\varepsilon)$ fixed for the finite test request in
Definition~\ref{defn:relation-URP} and models available for every $k$.

Local dimension/rank control appears in the study of mean dimension for
AH algebras \cite[Theorem~4.3]{NiuAH} and of minimal homeomorphisms
of zero mean dimension
\cite[Theorem~5.1]{ENmdim}.

\begin{thm}[Relation URP is equivalent to local AH approximation]
\label{thm:relation-URP-local-AH}
A regular uniform measure relation has relation URP if and only if it
has uniform local AH Cartan approximation.  Split relation URP is
likewise equivalent to split uniform local AH Cartan approximation.
No relative SBP assumption is needed.  The correspondence retains
block ranks and parameter-space dimensions, so the locally bounded
dimension-ratio condition is the same in the split formulations.
\end{thm}

\begin{proof}
Suppose first that relation URP holds.  By
Lemma~\ref{lem:uniform-relation-Fourier-approximation}, every finite
set $\mathcal F\subseteq(M_R)_1$ is uniformly approximated by
\[
 Q_na=\sum_{j\leq n}w_jb_j(a),\qquad a\in\mathcal F,
\]
for one $n$.  Each $w_j$ is the graph operator of an admissible
restriction of an atlas map, and $b_j(a)\in(D_K)_1$.
Apply relation URP to these finitely many maps and diagonal
coefficients, with tolerance $\delta$.  If $u_j,d_{j,a}\in H_1$ are
the resulting approximants, then
\[
 \left\|Q_na-\sum_{j\leq n}u_jd_{j,a}\right\|_{2,K}
 \leq2n\delta.
\]
Choose $n$ and then $\delta$ for the desired error.  Radial truncation
of the sum into $H_1$ does not increase its distance from the
contraction $a$, by the tracial truncation estimate in Section~2.
Equation~\eqref{eq:finite-cover-homogeneous-kernels}, diagonal
expectation compatibility, and the generating bisection normalizers
now give all conditions of
Definition~\ref{defn:uniform-local-AH-Cartan}.

Conversely, take a local AH approximant for
$\mathcal D\cup\{v_t:t\in F\}$. First decompose the Cartan diagonal
along the summands. Projection onto a summand is multiplication by
a central multiplier $z_s$ of $H$. For $d\in D_H$, the element
$z_s d\in H$ commutes with $D_H$, hence belongs to $D_H$. It follows
that $D_H=\bigoplus_s(D_H\cap H_s)$. Bimodularity and an approximate
unit in each diagonal summand restrict the expectation accordingly,
and the components of the generating normalizers still normalize the
summand diagonals. They also normalize $D_K$: for a normalizer $n$,
$z_s n$ is the norm limit of $d_\lambda n$ for an approximate unit
$d_\lambda$ of $D_H\cap H_s$, which is an approximate unit of $H_s$.
Thus each $D_H\cap H_s$ is Cartan in $H_s$.

For a homogeneous summand $H_s=p_s\mathrm{M}_{m_s}(\mathrm{C}_0(Z_s))p_s$, put $\mathcal V_s=\operatorname{ran}p_s$.
Its Cartan pair is locally the standard matrix-diagonal pair by
\cite[Lemma~2.1]{LiRenault}.  Hence the spectrum $Y_s$ of its diagonal
is an $r_s$-sheeted covering of $Z_s$.  The rank-one diagonal
projections define the line bundle
$L_y\subseteq(\mathcal V_s)_{q_s(y)}$, and
$\mathcal V_s\cong(q_s)_*L_s$.  Thus this summand has precisely the
finite-cover kernel presentation above.  In these coordinates a
Cartan normalizer has at most one nonzero matrix entry in each row
and column, so its nonzero coefficients form an open bisection.
Multiplication by diagonal cutoffs
with compact support inside its source domain approximates it in
norm by compactly supported bisection sections.  These cutoffs
preserve normalization of $D_K$.  Applying this to the generating
normalizers of the local AH block supplies the family required in
the definition.  The canonical diagonal expectation agrees with
$\mathbb{E}_R|_{H_s}$ by uniqueness of the Cartan expectation.  This gives
the relation URP witness.  For split pairs choose the trivial cover
and trivial line bundle.  These constructions keep the same $Z_s$
and $r_s$, proving the dimension-ratio assertion as well.
\end{proof}

\begin{prop}\label{prop:local-AH-uniform-OE-invariant}
Relation URP and uniform local AH Cartan approximation, their split
versions, and locally bounded dimension ratio for the split versions
are invariant under
isomorphisms of completed Cartan pairs that preserve the expectation
and the designated traces.  In particular they are invariant under
the uniform measure orbit equivalences considered here.
\end{prop}

\begin{proof}
Use Theorem~\ref{thm:relation-URP-local-AH} and transport each embedded
block, its diagonal and its generating
normalizers through the marked isomorphism.  The expectation and
the uniform $2$-norm are preserved.  The same homogeneous
presentation, rank and coefficient-space dimension are retained,
so all quantifiers and estimates are unchanged.
\end{proof}

\begin{prop}\label{prop:local-AH-measurewise-hyperfinite}
If a regular uniform measure relation has uniform local AH Cartan
approximation, then $(R,\mu)$ is amenable and hyperfinite for every
$\mu\in K$.  Its finite von Neumann relation algebra is hyperfinite;
in the ergodic aperiodic case it is the hyperfinite
$\mathrm{II}_1$ factor.
\end{prop}

\begin{proof}
Fix $\mu$.  The image of $M_R$ in its $\tau_\mu$-GNS von Neumann
algebra $M_\mu$ is strongly dense.  This algebra has separable
predual: the separable Cartan core constructed earlier has the same
GNS closure.  Given a finite set of contractions from this image,
choose a local block $H$ approximating it in $\norm{\cdot}_{2,\mu}$.
Let $N\subseteq M_\mu$ be the unital von Neumann algebra generated
by the image of $H$.  A representation of a homogeneous corner
has a type I von Neumann closure with bounded finite matrix sizes.
Adjoining the ambient unit adds at most a scalar summand on the
complement of its support.  Thus $N$ is a finite type I von Neumann
algebra and is injective.

Let $P_N:M_\mu\to N$ be the trace-preserving conditional
expectation.  Orthogonal projection in $L^2(M_\mu,\tau_\mu)$ gives
\[
 \norm{a-P_N(a)}_{2,\mu}\leq\norm{a-b_a}_{2,\mu}<\eps.
\]
Finite collections in $N$ can be approximated in $2$-norm by finite
dimensional subalgebras: partition their central measurable matrix
coefficients into finitely many ranges.  Compose the corresponding
trace-preserving finite-dimensional approximations with $P_N$.
This gives finite-dimensional completely positive approximations
to each finite test set in $M_\mu$.  Strong density extends the
criterion to every finite set in its unit ball.  The finite von
Neumann approximation theorem therefore gives injectivity and
hyperfiniteness \cite{Connes76}.  The measured Cartan relation is
then amenable and hyperfinite by the classical Cartan equivalences
and Connes--Feldman--Weiss \cite{CFW}.  Ergodicity and aperiodicity
give the last assertion.
\end{proof}

\begin{cor}[SBP and AF for relations with URP]
\label{cor:local-AH-Reiter-and-SBP}
Let $(X,R,K,D_K,\mathcal S)$ be a regular uniform measure relation
with a nonempty compact convex designated family $K$ of invariant
probability measures.  Relation URP implies regular uniform Reiter
amenability.  For a relation with URP, one has
\[
 K\text{-SBP}\quad\Longleftrightarrow\quad
 \text{strong uniform AF hyperfiniteness of }R.
\]
The AF condition is that of
Definition~\ref{defn:uniform-AF-hyperfinite-relation}.
Moreover, under relative SBP, measurewise amenability, relation URP,
and regular finite-cut approximation \textup{(RF)} are equivalent.
In this case split relation URP holds, with finite-dimensional local
AH approximants and one-point coefficient spaces.  No additional
uniform amenability assumption is required.
\end{cor}

\begin{proof}
Theorem~\ref{thm:relation-URP-local-AH},
Proposition~\ref{prop:local-AH-measurewise-hyperfinite}, and
Proposition~\ref{prop:Reiter-measurewise-equivalence} give the first
assertion.  If relative SBP holds, apply
Theorem~\ref{thm:uniform-CFW-SBP} to obtain strong uniform AF
hyperfiniteness and regular finite cuts.  Conversely, the finite
projection diagonals of a strong AF exhaustion imply relative SBP
by Proposition~\ref{prop:relative-SBP-finite-diagonals}.
For the final assertion, the AF Cartan reconstruction in
Theorem~\ref{thm:uniform-AF-reconstruction} supplies embedded
finite-dimensional stages
$H_n\cong\bigoplus_s \mathrm{M}_{r_s}(\mathbb C)$ with bounded uniform
$2$-norm density.  Their matrix units normalize $D_K$, and $\mathbb{E}_R$
restricts to their diagonal expectations.  They therefore give split
local AH approximants and, by
Theorem~\ref{thm:relation-URP-local-AH}, split relation URP.
\end{proof}

\begin{rem}[Relation URP and its split version]
\label{rem:relation-URP-nonsplit-example}
The explicit homogeneous models in
Examples~\ref{ex:split-AH-without-SBP}
and~\ref{ex:cover-AH-without-transversal} now illustrate
Definition~\ref{defn:relation-URP}. The first has split relation URP.
The second has relation URP, witnessed by its single finite-cover
model, but fails split relation URP. Both examples fail relative SBP;
in particular, relation URP does not require a regular global
transversal or relative SBP.

To see the failure of split relation URP in
Example~\ref{ex:cover-AH-without-transversal}, test the unit and the
deck involution $U$ with tolerance $1/2$. Since $K$ contains the
uniform measure on each fibre, $\|a\|\leq\sqrt{2}\,\|a\|_{2,K}$ for
every $a\in \mathrm{M}_2(\mathrm{C}(\mathbb T))$. An approximating split subalgebra
$H$ would therefore contain an element at norm distance less than one
from the unit, so $1\in H$ by functional calculus. The abstract split
model is consequently unital. The characteristic functions of its
individual sheets give nonzero, pairwise orthogonal projections in
$D_K$ summing to one. Connectedness of $\mathbb T$ permits only one
sheet, so $H\subseteq D_K$. But $\mathbb E_R(U)=0$ and $U$ is
unitary, whence $\|U-d\|_{2,K}\geq1$ for every $d\in D_K$.
\end{rem}

\begin{question*}
The orbits in Example~\ref{ex:cover-AH-without-transversal} are the
pairs $\{z,-z\}$. Does there exist a minimal, second-countable,
locally compact Hausdorff principal \'etale groupoid $\mathcal G$
with infinite compact metrizable unit space $X$ and nonempty full
family $K$ of invariant probability measures whose completed orbit
relation has relation URP but fails split relation URP? Can such
a groupoid be topologically amenable?
\end{question*}

Any such example must have infinite-dimensional unit space. Indeed,
minimality makes every orbit dense and hence infinite. Relation URP
gives measurewise amenability by
Proposition~\ref{prop:local-AH-measurewise-hyperfinite} and
Theorem~\ref{thm:relation-URP-local-AH}. If $X$ were finite-dimensional,
Theorem~\ref{thm:finite-dimensional-relation-SBP} would give relative
SBP, and Corollary~\ref{cor:local-AH-Reiter-and-SBP} would then give
split relation URP.

\begin{example}[Small boundaries without amenability]
\label{ex:SBP-without-relation-Rokhlin}
Let $F_2=\langle a,b\rangle$ act by left translation on its metrizable
profinite completion $X=\widehat{F_2}$, and let $R$ be the orbit relation.
Residual finiteness makes this action free.  Its only invariant
probability measure is Haar measure $\mu$: invariance under the dense
subgroup $F_2$ extends, by continuity on $\mathrm{C}(X)$, to invariance under
all translations of $X$.  Take $K=\{\mu\}$.  The clopen basis gives
relative SBP, and the global action maps give a regular atlas.

Nevertheless, every finite Borel subrelation $S\subseteq R$ satisfies
\[
 \mu\{x:(x,ax)\notin S\}
 +\mu\{x:(x,bx)\notin S\}\geq1.
\]
Indeed, set
\[
 f_S(x)=\mathbf1_{\{ax\in[x]_S\}}
          +\mathbf1_{\{bx\in[x]_S\}}.
\]
Since $S$ preserves $\mu$, finite-class averaging gives
\[
 \int_X f_S\,d\mu
 =\int_X\frac1{|[x]_S|}\sum_{y\in[x]_S}f_S(y)\,d\mu(x).
\]
Each orbit has the Cayley graph of $F_2$, which is a tree.  For a finite
class $C$, the sum $\sum_{y\in C}f_S(y)$ counts its internal
undirected generator edges once and is therefore at most $|C|-1$.
Thus $\int_X f_S\,d\mu\leq1$, and the displayed escape sum is
$2-\int_X f_S\,d\mu\geq1$.  Condition \textup{(RF)} fails already
for the two tests $a,b$ and $\varepsilon<1$.  Since relative SBP holds,
Corollary~\ref{cor:local-AH-Reiter-and-SBP}
shows that this relation fails both relation URP and uniform local
AH Cartan approximation.  Thus relative SBP
alone implies neither relation URP nor local AH approximation;
measured amenability is an essential hypothesis.
\end{example}

\begin{rem}\label{rem:local-AH-intrinsic-vs-continuous}
There is a stronger \emph{continuous} version, relative to a chosen
Cartan core $D_0\subseteq A_0$, requiring all the blocks to lie in
$A_0$, their diagonals to lie in $D_0$, and their generating
normalizers to normalize $D_0$.  For an action this means
$D_0=\mathrm{C}(X)$ and $A_0=\mathrm{C}(X)\rtimes_r\Gamma$.  It implies the intrinsic
completed property by bounded density.  It retains more information:
an isomorphism of completed pairs need not carry the chosen
continuous core onto another chosen continuous core.  Thus the
intrinsic property is an invariant of uniform measure OE, whereas
the continuous version is a statement about the specified
topological presentation.

The broad definition permits the connected finite-covering diagonal
of Example~\ref{ex:cover-AH-without-transversal}.  The split version
requires more: its individual homogeneous blocks have global split
diagonals.  Neither version imposes nested blocks, and neither is
claimed to reconstruct an AH or ASH inductive-limit Cartan model.
\end{rem}

\begin{prop}\label{prop:URP-continuous-local-AH}
Let $\Gamma$ be countably infinite and amenable, and let its action
on compact metrizable $X$ be free and have URP.  Then its continuous
Cartan pair has split uniform local AH Cartan approximation with
locally bounded dimension ratio, uniformly over all invariant
probability measures and hence over every compact designated
subfamily.  The spaces $Z_s$ in these approximants have finite
covering dimension.  The completed pair consequently has the same
intrinsic properties, after discarding any common-null blocks.
In particular, its orbit relation has split relation URP.
\end{prop}

\begin{proof}
Fix a finite set of contractions and a tolerance. Approximate them
in norm by finite Fourier polynomials. Write the coefficients on the
source side, so that a term supported on the arrows $x\mapsto xg$ has
value $f(x)$ at source $x$. Let $f_1,\ldots,f_m$ list all coefficient
functions and let $L$ be the finite set of group labels. These are
fixed before the rank lower bound $k$ is prescribed.

Given $k$ and $0<\eta<1$, take the closed URP towers of
Proposition~\ref{prop:URP-local-RSH}, with uncovered measure less than
$\eta$ and $|\{h\in F_s:hg\notin F_s\}|<\eta|F_s|$ for $g\in L$.
We may require $|F_s|\geq k$: include an auxiliary finite set $P$
of cardinality greater than $2k$ in the invariance test and require
$|F_sP|<2|F_s|$. Choose open thickenings $U_s$ with pairwise
disjoint translated closures and cutoffs $\chi_s\in \mathrm{C}_c(U_s)$,
$0\leq\chi_s\leq1$, equal to one on the respective closed bases.
The tower bisections identify
$\bigoplus_s \mathrm{M}_{|F_s|}(\mathrm{C}_0(U_s))$ with a subalgebra of the crossed
product: its $k',h$ entry acts on the arrows $bh\mapsto bk'$,
$b\in U_s$. In particular, individual entries normalize $\mathrm{C}(X)$.

In $\mathrm{C}_0(U_s)$ let $C_s$ be the $\mathrm{C}^*$-algebra generated, without adjoining
a unit, by $\chi_s$ and the functions
\[
 b\longmapsto \chi_s(b) f_j(bh)
 \qquad(1\le j\le m,\ h\in F_s).
\]
Thus $C_s=\mathrm{C}_0(Z_s)$, where $Z_s$ is the joint range on the one-point
compactification with the zero character removed. There are at most
$1+2m|F_s|$ real coordinates. The covering dimension of this
metrizable subspace of Euclidean space is at most $1+2m|F_s|$.
Set $H=\bigoplus_s \mathrm{M}_{|F_s|}(C_s)$. Its matrix diagonal belongs
to $\mathrm{C}(X)$, and the crossed-product expectation is the coordinate
diagonal expectation on $H$. Hence $H\cap \mathrm{C}(X)$ is its coordinate
diagonal. Its individual matrix entries are generating normalizers
of both diagonals.

For a term with coefficient $f_j$ and label $g$, put
$\chi_s(b) f_j(bh)$ in the $(hg,h)$ entry when $h,hg\in F_s$ and
use zero elsewhere. This gives an element of $H$ agreeing with the
term on these arrows over the closed tower bases. Its squared
$2$-norm error is at most $2\eta\|f_j\|^2$: the uncovered sources
have measure less than $\eta$, and the omitted boundary levels above
the closed bases have total measure at most $\eta$, by invariance.
Sum over the finitely many terms and choose $\eta$ to give the desired
error. Radial truncation into $H_1$ preserves the error bound relative
to the original contractions. Since
$\dim Z_s/|F_s|\leq 2m+1$ for every $k$, this proves
\eqref{eq:local-AH-bounded-ratio-quantifiers} with $N=2m+2$.
The compact-core estimate in \cite[Theorem~3.9]{NiuURP} gives the
separate bound in terms of mean dimension stated above.

For a smaller designated family, passage to the completed pair may
factor a block through a quotient.  An ideal of a full matrix
algebra over $\mathrm{C}_0(Z_s)$ gives another full matrix block over a
closed subspace of $Z_s$, or the zero block.  Removing zero blocks
preserves the rank lower bound and cannot increase the dimensions.
Bounded uniform density extends the finite-set property from $A_0$
to its marked completion.  Apply
Theorem~\ref{thm:relation-URP-local-AH} for the relation URP assertion.
\end{proof}

\subsection{Diagonal-fixing completely positive local approximation}
\label{subsec:URP-Cartan-cp}
\leavevmode\par

The local subgroupoids in Proposition~\ref{prop:URP-local-RSH} give
more than distance estimates.  Their matrix coordinates support
completely positive approximation maps which retain the whole
diagonal and the designated measures exactly.

\begin{thm}[URP gives marked Cartan completely positive approximation]
\label{thm:URP-Cartan-cp}
In the setting of Proposition~\ref{prop:URP-local-RSH}, let
$A_0=\mathrm{C}(X)\rtimes_r\Gamma$ and $M=A_0^K$.  For every finite
$\mathcal F\subseteq A_0$ and $\eps>0$ there are an open relatively
compact principal subgroupoid $G\subseteq X\rtimes\Gamma$
containing all units, its recursively subhomogeneous Cartan algebra
$B_G=\mathrm{C}_r^*(G)$, and a unital completely positive map
\[
                    \Phi_G:A_0\longrightarrow B_G
\]
such that
\[
 \Phi_G|_{\mathrm{C}(X)}=\operatorname{id}_{\mathrm{C}(X)},\qquad
 \mathbb{E}\Phi_G=\mathbb{E},\qquad
 \tau_\mu\Phi_G=\tau_\mu\quad(\mu\in K),
\]
and
\[
              \|a-\Phi_G(a)\|_{2,K}<\eps
              \qquad(a\in\mathcal F).
\]
The map is $\mathrm{C}(X)$-bimodular and contractive for every designated
$2$-norm.  It extends to a unital completely positive map
\[
 \overline\Phi_G:M\longrightarrow B_G^{L_G}\subseteq M,
 \qquad L_G=\{\tau_\mu|_{B_G}:\mu\in K\},
\]
which fixes $D_K$, is $D_K$-bimodular, and preserves the expectation
and each designated trace. Restriction gives a homeomorphism
$K\to L_G$ that extends to an affine homeomorphism of their weak-star
closed convex hulls. The approximation conclusion also
holds for any prescribed finite subset of $M$ after choosing $G$
accordingly.  These conclusions require neither SBP nor nested
subgroupoids.
\end{thm}

\begin{proof}
First fix a finite set $L\subseteq\Gamma$ and $0<\eta<1$.  Choose
closed URP bases $B_s$ and shapes $F_s$ whose uncovered measure is
less than $\eta$ for every invariant probability measure, and such
that
\[
       |\{h\in F_s:hg\notin F_s\}|<\eta|F_s|
       \qquad(g\in L).
\]
Choose open thickenings $U_s\supseteq B_s$ with all the compact
levels $\overline{U_s}h$ mutually disjoint, and continuous functions
$0\leq\chi_s\leq1$, compactly supported in $U_s$, equal to $1$ on
$B_s$.  Let $G$ be the wide finite subgroupoid from
Proposition~\ref{prop:URP-local-RSH}.  Its nontrivial finite classes
are the tower fibres $\{bh:h\in F_s\}$, $b\in U_s$.

Define a scalar kernel on the transformation groupoid by
\[
 \kappa(x,x)=1,\qquad
 \kappa(bh,bk)=\chi_s(b)\quad(h\ne k),
\]
and set it equal to zero on all other nonidentity arrows.  Here
the notation records the range and source pair; freeness identifies
each such pair with a unique transformation-groupoid arrow.
The kernel is continuous and compactly supported in $G$.  Indeed,
the unit component is clopen, and every nonidentity group-coordinate
coefficient is a finite sum of translates of functions
$\chi_s\in \mathrm{C}_c(U_s)$, extended continuously by zero.  There are
only finitely many group coordinates $h^{-1}k$.

On each full orbit, partition the points into their finite
$G$-classes.  On a tower class the matrix of $\kappa$ is
\[
                   \chi_s(b)J+(1-\chi_s(b))I,
\]
where $J$ is the all-ones matrix.  This matrix is positive, and
the singleton classes have matrix $[1]$.  Cross-class coefficients
are zero.  Thus $\kappa$ is positive definite on every full orbit,
not merely on the arrows of $G$.

Multiply regular kernels by $\kappa$.  In each reduced regular
representation this is the unital positive Schur multiplier
associated to the preceding positive-definite matrix.  It is
completely positive and contractive.  The faithful family of
reduced regular representations therefore gives a unital completely
positive map $\Phi_G$ on $A_0$.  For a finite Fourier polynomial,
its image is a continuous compactly supported groupoid function
on $G$, so belongs to $B_G$.  Norm continuity gives the same range
inclusion for all of $A_0$.  The algebra $B_G$ has the diagonal
recursive subhomogeneous description
\eqref{eq:tower-RSH-pullback}.

The diagonal value of $\kappa$ is $1$, giving
$\Phi_G|_{\mathrm{C}(X)}=\operatorname{id}$ and $\mathbb{E}\Phi_G=\mathbb{E}$.  Alternatively,
the diagonal lies in the multiplicative domain of the unital
completely positive map; in particular the map is $\mathrm{C}(X)$-bimodular.
Since $\tau_\mu=\mu\circ \mathbb{E}$, it preserves every designated trace.
Schwarz's inequality now gives
\[
 \|\Phi_G(a)\|_{2,\mu}^{\,2}
 \leq\tau_\mu\bigl(\Phi_G(a^*a)\bigr)
 =\tau_\mu(a^*a),\qquad \mu\in K.
\]

For $g\in L$, the coefficient of $\Phi_G(u_g)$ is $1$ on every
source $x=bh$ with $b\in B_s$ and $h,hg\in F_s$.  Elsewhere this
coefficient still lies in $[0,1]$.  Invariance therefore gives
\begin{align*}
 \|u_g-\Phi_G(u_g)\|_{2,K}^{\,2}
 &\leq\sup_{\mu\in K}\left[
       \mu\bigl(X\setminus\textstyle\bigcup_{s,h}B_sh\bigr)
       +\sum_s|\{h\in F_s:hg\notin F_s\}|\mu(B_s)\right]\\
 &\leq 2\eta.
\end{align*}
Source versus range coefficients have the same integral by
invariance.  Bimodularity consequently gives, for a Fourier polynomial,
\[
 \left\|\sum_{g\in L}f_gu_g-
          \Phi_G\left(\sum_{g\in L}f_gu_g\right)\right\|_{2,K}
 \leq\sqrt{2\eta}\sum_{g\in L}\|f_g\|.
\]
Choosing $\eta$ small and using norm density of finite Fourier
polynomials proves the asserted finite-set approximation.

On the compact convex set $\overline{\operatorname{co}}K$, the map
$\mu\mapsto\tau_\mu|_{B_G}$ is continuous and affine. Since $B_G$
contains $\mathrm{C}(X)$, it is injective. Its image is
$\overline{\operatorname{co}}L_G$ by compactness and affinity, so it
is a homeomorphism onto that hull and restricts to the stated
homeomorphism $K\to L_G$.  The faithful-completion construction
embeds $B_G^{L_G}$ into $M$, with its diagonal completion equal
to $D_K$.  A bounded uniform $2$-norm Cauchy sequence in $A_0$
is carried, by the displayed contraction estimate, to a bounded
uniform $2$-norm Cauchy sequence in $B_G$.  This defines
$\overline\Phi_G$ with the claimed range.  Positivity at every
matrix level passes to the bounded limits in the faithful finite
von Neumann representations, so this extension is unital and
completely positive.  Diagonal density, bounded continuity of
multiplication, and trace continuity give its $D_K$-bimodularity,
diagonal-fixing property, and preservation of expectation and traces.

Finally, for a finite $\mathcal F\subseteq M$, choose bounded
$a_f\in A_0$ close to each $f\in\mathcal F$ in uniform $2$-norm.
Apply the core approximation to these finitely many $a_f$.  The
$2$-norm contraction of the extension bounds
\[
 \|f-\overline\Phi_G(f)\|_{2,K}
 \leq2\|f-a_f\|_{2,K}
       +\|a_f-\Phi_G(a_f)\|_{2,K},
\]
which proves the completed finite-set assertion.
\end{proof}

\begin{rem}\label{rem:URP-Cartan-cp-scope}
The maps in Theorem~\ref{thm:URP-Cartan-cp} need not be idempotent.
They provide diagonal-fixing, trace-preserving completely positive
local approximation whose norm cores are subhomogeneous Cartan pairs.
An exact inductive-limit presentation requires additional
compatibility of the stages.  The approximation topology here is
the bounded uniform $2$-norm topology; norm approximation of the
identity map is not asserted.
\end{rem}

\subsection{Uniform Borel finite exhaustion}
\label{subsec:local-finite-Borel-nesting}

Let $\mathcal G$ be a second-countable, locally compact, Hausdorff,
principal \'etale groupoid with compact unit space $X$, and let $K$ be a
nonempty family of invariant probability measures on $X$.  In this
subsection compactness or convexity of $K$ is not required.  For a Borel
set $A\subseteq\mathcal G$, put
\[
 \nu_\mu(A)=\int_X\#\bigl(A\cap s^{-1}(x)\bigr)\,d\mu(x).
\]
Invariance gives $\nu_\mu(A)=\nu_\mu(A^{-1})$.  Every compact
$Q\subseteq\mathcal G$ has finite source and range multiplicity bounds
\[
 d_s(Q)=\sup_{x\in X}\#\bigl(Q\cap s^{-1}(x)\bigr),\qquad
 d_r(Q)=\sup_{x\in X}\#\bigl(Q\cap r^{-1}(x)\bigr),
\]
because $Q$ is covered by finitely many open bisections.
By principality, a full-unit Borel subgroupoid is identified with a
Borel subequivalence relation of $R_{\mathcal G}$.

\begin{lem}[Borel repair of a finite stage]
\label{lem:Borel-finite-stage-repair}
Let $S,H\subseteq\mathcal G$ be full-unit Borel subgroupoids whose
classes have cardinality at most $M,N$, respectively.  Define
\[
 E=s(S\setminus H),\qquad B=[E]_S,
 \qquad F=(H|_{X\setminus B})\sqcup(S|_B).
\]
The subgroupoid $F$ is Borel, contains every unit, and satisfies
$S\subseteq F$. Its classes have size at most $\max\{M,N\}$.
For every $\mu\in K$,
\[
 \mu(B)\leq M\mu(E)\leq M\nu_\mu(S\setminus H).
\]
For every compact $Q\subseteq\mathcal G$,
\begin{equation}
\label{eq:Borel-finite-stage-counting-repair}
 \nu_\mu(Q\setminus F)
 \leq \nu_\mu(Q\setminus H)
       +\bigl(d_s(Q)+d_r(Q)\bigr)M\mu(E).
\end{equation}
For a Borel partial $R_{\mathcal G}$-bijection $t$ one also has
\begin{equation}
\label{eq:Borel-finite-stage-chart-repair}
 \mu\{x\in\operatorname{dom}t:(tx,x)\notin F\}
 \leq \mu\{x\in\operatorname{dom}t:(tx,x)\notin H\}
       +2M\mu(E).
\end{equation}
If $S$ and $H$ are relatively compact in $\mathcal G$, then so is $F$.
\end{lem}

\begin{proof}
The source and range maps are countable-to-one Borel maps, so $E$ and
$B$ are Borel.  Since $B$ is $S$-invariant, an $S$-class is either
contained in $B$ or disjoint from $B$.  In the latter case no arrow of
that class lies in $S\setminus H$, so the whole class is contained in
$H|_{X\setminus B}$.  This proves $S\subseteq F$; the remaining
subgroupoid and class-size assertions follow from its disjoint-union
description.

Counting mass transport for the finite Borel relation $S$ gives
\[
 \begin{aligned}
 \mu(B)
 &\leq\int_X\sum_{y\in[x]_S}1_E(y)\,d\mu(x)\\
 &=\int_E|[x]_S|\,d\mu(x)\leq M\mu(E).
 \end{aligned}
\]
This identity follows from invariance by decomposing $S$ into Borel
partial bijections; no regularity of their domain projections is used.
Also $\mu(E)\leq\nu_\mu(S\setminus H)$ by the definition of $E$.
An arrow in $H$ is retained in $F$ whenever both endpoints lie outside
$B$.  Consequently
\[
 Q\setminus F\subseteq(Q\setminus H)
       \cup(Q\cap s^{-1}(B))\cup(Q\cap r^{-1}(B)).
\]
The source contribution is at most $d_s(Q)\mu(B)$.  Inversion and
invariance bound the range contribution by $d_r(Q)\mu(B)$, proving
\eqref{eq:Borel-finite-stage-counting-repair}.
For a partial bijection $t$, the additional exceptional source set
is contained in $B\cup t^{-1}(B)$.  Invariance gives
$\mu(t^{-1}(B))\leq\mu(B)$, proving
\eqref{eq:Borel-finite-stage-chart-repair}.  Finally $F\subseteq S\cup H$,
which proves the assertion concerning relative compactness.
\end{proof}

\begin{thm}[Uniform Borel exhaustion from local topological stages]
\label{thm:local-topological-uniform-Borel-exhaustion}
Assume that for every compact $Q\subseteq\mathcal G$ and every
$\eps>0$ there is an open relatively compact subgroupoid
$H\subseteq\mathcal G$, containing all units, with
\begin{equation}
\label{eq:local-topological-counting-capture}
 \sup_{\mu\in K}\nu_\mu(Q\setminus H)<\eps.
\end{equation}
Then there are increasing full-unit finite Borel subgroupoids
$S_1\subseteq S_2\subseteq\cdots\subseteq\mathcal G$, each relatively
compact and with a uniform bound on its class sizes, such that
\[
 \sup_{\mu\in K}\nu_\mu(Q\setminus S_n)\longrightarrow0
 \qquad\text{for every compact }Q\subseteq\mathcal G.
\]
There is one $R_{\mathcal G}$-invariant $K$-conull Borel set $X_0$ on
which
\[
 R_{\mathcal G}|_{X_0}=\bigcup_n S_n|_{X_0}.
\]
In particular, every designated measured relation is hyperfinite.
The construction makes no assertion that the $S_n$ are open or that
their graph, domain, or coincidence projections belong to $D_K$.
\end{thm}

\begin{proof}
An open relatively compact full-unit subgroupoid $H$ has uniformly
finite classes: its source fibres have cardinality at most
$d_s(\overline H)$, and principality identifies a source fibre with
the corresponding equivalence class.  Choose increasing compact sets
$Q_n$ whose interiors cover $\mathcal G$, and include $X$ in each of
them.  Set $S_0=X$.

Suppose $S_{n-1}$ has been constructed, with class bound $M_{n-1}$.
Its closure in the ambient groupoid is compact; the closure need not
itself be a subgroupoid.  Choose $\eta_n>0$ so small that
\[
 \eta_n\bigl[1+
    \bigl(d_s(Q_n)+d_r(Q_n)\bigr)M_{n-1}\bigr]<2^{-n}.
\]
Apply \eqref{eq:local-topological-counting-capture} to the compact set
$Q_n\cup\overline{S_{n-1}}$ with tolerance $\eta_n$, and call the
resulting open subgroupoid $H_n$.  Repair $S_{n-1}$ and $H_n$ by
Lemma~\ref{lem:Borel-finite-stage-repair}, obtaining $S_n$.
Since
\[
 \sup_{\mu\in K}\mu\bigl(s(S_{n-1}\setminus H_n)\bigr)<\eta_n,
\]
the lemma gives $S_{n-1}\subseteq S_n$ and
\[
 \sup_{\mu\in K}\nu_\mu(Q_n\setminus S_n)<2^{-n}.
\]
It also supplies a finite class bound for $S_n$ and preserves relative
compactness.  For fixed $j$ and $n\geq j$ the same estimate holds
with $Q_j$ in place of $Q_n$.  Every compact subset of $\mathcal G$
is contained in some $Q_j$, proving the uniform capture assertion.

Put $Z=\mathcal G\setminus\bigcup_nS_n$.  For every $\mu\in K$ and
every $j$, continuity from above on the finite-measure set $Q_j$ gives
$\nu_\mu(Q_j\cap Z)=0$.  Hence $\nu_\mu(Z)=0$.  The set
$s(Z)\cup r(Z)$ is Borel, since both maps have countable fibres, and
it is null for every designated measure.  Its $R_{\mathcal G}$-saturation
is again Borel and $K$-null, by countability and invariance.  Let $X_0$
be the complement of this saturation.  Every arrow with endpoints
in $X_0$ belongs to some $S_n$, which proves the stated common
conull exhaustion.
\end{proof}

\begin{cor}[The Borel consequence of URP]
\label{cor:URP-uniform-Borel-exhaustion}
For a free action with URP as in
Proposition~\ref{prop:URP-local-RSH}, its orbit relation has the
uniform Borel exhaustion of
Theorem~\ref{thm:local-topological-uniform-Borel-exhaustion}, for every
designated family of invariant probability measures.  Neither SBP nor
a nested family of topological Rokhlin stages is required.
\end{cor}

\begin{proof}
A compact subset $Q$ of $X\rtimes\Gamma$ meets only finitely many
group coordinates, say $L$.  For a subgroupoid $H$ one has
\[
 \nu_\mu(Q\setminus H)
 \leq\sum_{g\in L}
       \mu\{x:(x,g)\notin H\}.
\]
Apply Proposition~\ref{prop:URP-local-RSH} to the compact family of
all invariant probability measures.  With a sufficiently small
coordinatewise tolerance it gives
\eqref{eq:local-topological-counting-capture}.
\end{proof}

\begin{rem}[Borel nesting and regular finite approximation]
\label{rem:Borel-nesting-not-regular-AF}
The theorem supplies actual nesting after passing from the local
topological stages to Borel cuts.  It does not supply nested
topological stages.  The saturation and restriction operations in
the repair can create nonregular characteristic functions, so the
result is weaker than regular finite approximation \textup{(RF)} and
does not provide a regular AF or AH Cartan core.  The uniform capture
assertion also gives uniform escape estimates for partial maps whose
graphs are relatively compact in the ambient groupoid.  Extending
those estimates to every admissible partial map requires a regular
atlas and the corresponding uniform chart-tail argument.
\end{rem}

\subsection{Orbit-breaking ASH models}
\label{subsec:orbit-breaking-ASH}

For a minimal homeomorphism the nesting problem has a direct
solution through orbit breaking.  Here ASH means an inductive limit
of subhomogeneous $\mathrm{C}^*$-algebras; the compact coefficient spaces are
allowed to have arbitrary covering dimension.

\begin{thm}\label{thm:orbit-breaking-ASH-uniform}
Let $T$ be a minimal homeomorphism of an infinite compact metrizable
space $X$, and let $K$ be a nonempty compact family of
$T$-invariant probability measures.  For any $y\in X$, the
orbit-breaking algebra
\[
 A_y=\mathrm{C}^*\bigl(\mathrm{C}(X),\,u g:g\in \mathrm{C}(X),\ g(y)=0\bigr)
 \subseteq \mathrm{C}(X)\rtimes_T\mathbb Z
\]
is ASH, $\mathrm{C}(X)$ is a Cartan subalgebra of $A_y$, and the inclusion
induces an isomorphism of the designated uniform tracial Cartan
completions.  Its principal \'etale groupoid agrees with the orbit
relation of $T$ away from the common $K$-null orbit of $y$.
No SBP assumption is needed.
\end{thm}

\begin{proof}
Choose closed neighborhoods $Y_n$ of $y$ with
$Y_{n+1}\subseteq\operatorname{int}Y_n$ and
$\bigcap_nY_n=\{y\}$. Define
\[
 A_{Y_n}=\mathrm{C}^*\bigl(\mathrm{C}(X),\,u g:g\in \mathrm{C}(X),\ g|_{Y_n}=0\bigr).
\]
These algebras increase.  The first-return decomposition shows that
each is recursively subhomogeneous with diagonal gluing maps
\cite[Theorem~5.5]{NiuURP}.  Functions vanishing on a neighborhood
of $y$ are norm dense among continuous functions vanishing at $y$.
It follows that $A_y=\overline{\bigcup_nA_{Y_n}}^{\norm{\cdot}}$,
so $A_y$ is ASH.

Every invariant probability measure is nonatomic: minimality on
an infinite space makes every orbit infinite, whereas an atom
would have the same positive mass at every point of its orbit.
Hence $\mu(Y_n)\downarrow0$ for every $\mu\in K$.  The functions
$\mu\mapsto\mu(Y_n)$ are upper semicontinuous and decrease on
compact $K$, so
\[
 \sup_{\mu\in K}\mu(Y_n)\longrightarrow0.
\]
For $n\geq2$, choose a continuous $0\leq g_n\leq1$ vanishing on
$Y_n$ and equal to $1$ off $\operatorname{int}Y_{n-1}$. The error
is supported in $Y_{n-1}$. Then $u g_n\in A_{Y_n}$ and
\[
 \norm{u-u g_n}_{2,K}^2
 =\sup_{\mu\in K}\int(1-g_n)^2\,d\mu\longrightarrow0.
\]
Together with $\mathrm{C}(X)\subseteq A_y$, this gives uniform $2$-norm
density in the crossed product.  Functional-calculus truncation
gives bounded approximants, so the designated bounded completions
coincide.

The groupoid of $A_{Y_n}$ is generated by the units and the
generator edges based outside $Y_n$. Their increasing union is
generated by the units and the generator edges based outside $\{y\}$.
On the orbit of $y$, the omitted edge separates the orbit into two
half-orbits, with no arrows between them. The union is an open
principal \'etale subgroupoid of the transformation groupoid,
whose $\mathrm{C}^*$-algebra is $A_y$.  The diagonal $\mathrm{C}(X)$ is therefore
Cartan.  A finite orbit segment avoiding $y$ eventually avoids
$Y_n$, so every arrow off the orbit of $y$ belongs to this union.
The countable orbit of $y$ is null for every invariant measure,
which proves the final assertion.
\end{proof}

Theorem~\ref{thm:orbit-breaking-ASH-uniform} supplies an actual
nested subhomogeneous Cartan model in this case, rather than only
local approximation.  It does not assert an AH decomposition with
homogeneous stages, nor an AF diagonal approximation.  In particular,
it remains applicable when the designated relative SBP fails.

\subsection{Compatibility and dimension}
\label{subsec:AH-compatibility-question}

The local and nested conclusions have different hypotheses.
Definition~\ref{defn:uniform-local-AH-Cartan} specifies approximation
of finite observable sets by embedded homogeneous Cartan pairs,
without connecting maps. Proposition~\ref{prop:URP-continuous-local-AH}
proves its split continuous version for free URP actions.
Theorem~\ref{thm:URP-Cartan-cp} additionally gives diagonal-fixing,
trace-preserving completely positive approximation through completed
RSH tower algebras. These statements retain positive cutoffs instead
of requiring their support characteristic functions to be regular.

Theorem~\ref{thm:local-topological-uniform-Borel-exhaustion} and
Corollary~\ref{cor:URP-uniform-Borel-exhaustion} give an increasing
finite Borel exhaustion from the topological local-capture hypothesis.
Their repair need not preserve the completed diagonal. They therefore
do not imply regular finite approximation \textup{(RF)}.
In contrast, regular finite cuts give an actual ASH Cartan core by
Theorem~\ref{thm:RF-ASH-without-SBP}, and compatible regular matrix
coordinates give an AH core by
Theorem~\ref{thm:AH-reconstruction-without-SBP}.
For minimal $\mathbb Z$-systems,
Theorem~\ref{thm:orbit-breaking-ASH-uniform} gives an actual nested
ASH model by orbit breaking.

\begin{rem}[The remaining reconstruction problem]
\label{rem:local-AH-reconstruction-question}
For a general free URP action, does its designated uniform completed
Cartan pair admit a separable ASH Cartan core whose designated
bounded completion is the given pair? Neither local AH approximation
nor the Borel exhaustion proves this assertion. A construction would
need compatible diagonal-preserving embeddings or regular nested
finite stages with appropriate boundary algebras. Uniform $2$-norm
approximate containment alone supplies no such exact inclusions.
The intrinsic local AH property remains useful independently of that
problem: it is preserved by uniform measure OE, implies measurewise
amenability, and, with relative SBP, gives the strong AF conclusion by
Corollary~\ref{cor:local-AH-Reiter-and-SBP}.
\end{rem}

The dimension estimates also have distinct scopes. The local bound on
dimension ratio in \eqref{eq:local-AH-bounded-ratio-quantifiers} is finite for each fixed
approximation request and independent of the required minimum rank.
It is not one global dimension/rank bound and does not assert zero
dimension growth. The compact-core bound of
\cite[Theorem~3.9]{NiuURP} retains the finer dependence on mean
dimension. An inductive-limit reconstruction preserving that bound
is a further dimension question. The conditional AH reconstruction
above supplies no such bound on its connecting system.

The broad local AH definition retains homogeneous corner bundles and
nonsplit Cartan diagonals, including as single blocks the covering
models of Remark~\ref{rem:finite-covering-AH-blocks} with
finite-dimensional base. Independent local
trivializations of a scalar or bundle twist do not produce one
uniformly regular global gauge. The regular untwisting theorem
still requires its stated compatible finite-chart hypotheses; the
Borel repair alone does not establish them.

All the local AH, completely positive, and Borel exhaustion results
above are stated for invariant probability measures. The controlled
Radon--Nikodym hypotheses of Section~\ref{sec:controlled-RN} govern
transport of observables and source/target errors. They do not supply
small boundaries or compatible topological blocks. A nonsingular
extension of the Borel repair requires uniform transport control on
the finite stages and the tested arrows, while an analytic extension
must retain the states and common modular cocycle. No such extension
is inferred merely from the chartwise Jacobian bounds.

\section{Further questions and related work}

A regular atlas specifies the uniform observable and graph-coefficient
topology.  Relative SBP gives the diagonal projection-density condition.
For an aperiodic regular relation on a finite-dimensional unit space,
Theorem~\ref{thm:finite-dimensional-relation-SBP} supplies this condition
automatically.
By Proposition~\ref{prop:Reiter-measurewise-equivalence} and
Theorem~\ref{thm:uniform-CFW-SBP}, amenability for every designated
measure then supplies deterministic regular finite cuts.
Theorem~\ref{thm:local-cuts-AF} supplies compatible matrix charts
and an AF model.  Aperiodicity gives a single admissible generator
and the designated Cantor dynamical realization of
Theorem~\ref{thm:local-cuts-dynamical-model}.
Example~\ref{ex:Reiter-without-AF-coordinates} shows that relative
SBP cannot be dropped from the AF conclusion.
For principal \'etale presentations, Theorem~\ref{thm:SBP-groupoid-castles}
also gives small subgroupoids in the original topology, with uniformly
large unit spaces, F\o lner fibres, and exact compact-arrow extension.
For minimal topologically amenable principal groupoids with infinite
compact units and relative SBP for the full nonempty family of
invariant probability measures, Corollary~\ref{cor:SBP-COS-groupoid-Z}
also gives $\mathcal Z$-absorption
of the original untwisted groupoid algebra under Cuntz comparison of open
sets. Removing SBP from the implication from amenability to intrinsic relation
URP remains open in this generality. Such an extension must retain
homogeneous parameters: demanding uniformly almost full disjoint levels
subordinate to every open cover would already imply SBP.
For an aperiodic regular relation with relative SBP, amenability for
every invariant probability measure, and the full nonempty compact
invariant family, the reconstructed AF traces form a closed face with
diffuse extreme traces. The classification theorem
of Section~\ref{subsec:uniform-CFW} then gives a simple full-trace
replacement and a minimal Cantor realization. Minimality of the
original relation is unnecessary; aperiodicity suffices. Arbitrary
designated subfamilies can fail the face condition and are not
classified by their abstract measure families alone.

Theorem~\ref{thm:uniform-Dye-reconstruction} also reconstructs
these amenable full-family relations from their topological uniform
full groups; each topological isomorphism is uniquely spatial and
isometric. Whether every abstract group isomorphism is automatically
continuous is a separate question. The Boolean reconstruction already
holds algebraically, but transferring the uniform topology requires
an additional argument. For designated subfamilies,
Example~\ref{ex:full-group-topology-loses-designated-weights} shows
that group topology alone cannot recover the measure marking.

Theorem~\ref{thm:uniform-regular-untwisting} removes every regular
scalar $2$-cocycle on a relation with the finite-chart hypotheses,
using one gauge regular for the entire designated family.
Relative SBP also supplies regular scalar presentations of locally trivial
circle twists on principal \'etale groupoids.
These results concern uniform completed pairs; they do not assert
continuous triviality of the original twist.
Section~\ref{sec:AH-without-SBP} records what can replace the AF
conclusion when relative SBP is dropped: regular finite cuts give an
ASH Cartan core, and exact matrix coordinates give a sharper AH core.
Relation URP is equivalent to local AH approximation without SBP;
within this class relative SBP is equivalent to strong uniform AF
hyperfiniteness. Action URP supplies split local AH blocks and
diagonal-fixing completely positive RSH approximations. Its local
finite topological stages also
admit a common uniform Borel exhaustion, with no assertion of regular
projections. Compatibility for a general ASH Cartan core requires
further work beyond the minimal $\mathbb Z$ case.
Section~\ref{sec:intrinsic-reconstruction} completes the inverse
uniform Feldman--Moore construction for uniformly separable tracial
completed pairs with partial-isometry normalizer regularity and
fibrewise maximal abelianness. The relation and regular twist class are
independent of the core, and intrinsic SBP with injective fibres gives
a compatible AF Cartan model. The untwisted graph section is a consequence
of the amenable SBP case, rather than part of the definition of a general
completed pair. Example~\ref{ex:uniform-MASA-not-fibre-MASA} shows that
SBP and global maximal abelianness cannot replace the fibre condition.
Reconstruction outside these hypotheses, including ordinary normalizer
regularity without a verified partial-isometry generation condition,
and a nonsingular inverse theorem retaining the modular cocycle,
remain separate questions.

There are also type restrictions on any prospective classification by a
measure simplex alone.  Finite single-class relations of different
cardinalities all have a singleton invariant-probability simplex.  Their
class sizes are additional data.  The aperiodicity hypothesis in the general single-generator theorem excludes those finite classes; no classification of arbitrary finite or mixed-type
relations is asserted.

Existing uses of ``uniform hyperfiniteness'' should be distinguished
from the definition here.  Elek~\cite{Elek} studies measured graphs of
bounded degree with quantitative finite-com\-po\-nent bounds valid on all
positive-measure subgraphs.  Elek and Tim\'ar~\cite{ET} relate uniform
Borel amenability to randomized hyperfiniteness and prove, for uniformly
Borel amenable F\o lner graphs, hyperfiniteness outside an invariant
Borel set null for every invariant probability measure.  Their
Theorem~8 gives the latter conclusion, and Theorem~7 gives almost
finiteness in measure uniformly over all invariant measures.  Our
definitions additionally specify the completed diagonal and require
regularity of graph coefficients, transversals, and finite matrix
coordinates.  Any comparison between the two frameworks must retain
those extra regularity requirements.

\end{document}